\documentclass[a4paper,10pt]{article}
\usepackage[margin=0.80in,bottom=0.90in,top=0.85in]{geometry}

\usepackage[english]{babel}
\usepackage{chngcntr}

\usepackage[utf8]{inputenc}
\usepackage{amsmath,mathrsfs}
\usepackage{indentfirst}
\usepackage{fancyhdr}
\usepackage{amssymb}
\usepackage{amsthm}
\usepackage{graphicx}
\usepackage{color}

\usepackage{cancel}
\usepackage{subcaption}
\usepackage[dvipsnames]{xcolor}
\usepackage{eucal}
\usepackage{latexsym}
\usepackage{chngcntr}

\usepackage[expansion=false]{microtype}

\usepackage{enumitem}

\usepackage{geometry}

\usepackage{tikz}

\usepackage{todonotes}

\usepackage{amsmath, amsfonts, amsthm,amssymb}

\usepackage{tikz-cd}
\usepackage{multicol}
\usepackage{cancel}
\usepackage{mathtools}
\usepackage{graphicx}

\usepackage{nicematrix}

\usepackage{comment}

\usepackage[colorlinks=true,
            linkcolor=blue,
            citecolor=blue,
            urlcolor=blue]{hyperref}

\newcommand{\e}{\epsilon}
\newcommand{\ve}{\varepsilon}
\newcommand{\ck}{\mathtt{c}_{\kappa}}

\newcommand{\pa}{\partial}

\newcommand{\de}{\textrm{d}}

\newcommand{\tts}{\mathtt{s}}
\newcommand{\cA}{{\mathcal A}}
\newcommand{\cB}{{\mathcal B}}

\newcommand{\cJ}{{\mathcal J}}

\newcommand{\cL}{{\mathcal L}}

\newcommand{\cO}{{\mathcal O}}

\newcommand{\cR}{{\mathcal R}}

\newcommand{\cV}{{\mathcal V}}
\newcommand{\cW}{{\mathcal W}}

\newcommand{\N}{{\mathbb{N}}}
\newcommand{\R}{{\mathbb{R}}}
\newcommand{\Z}{{\mathbb{Z}}}
\newcommand{\T}{{\mathbb{T}}}

\newcommand{\vet}[2]{\begin{bmatrix}#1 \\ #2 \end{bmatrix}}

\def\s{{\rm s}}

\newcommand{\fa}{\mathfrak{a}}
\newcommand{\fb}{\mathfrak{b}}
\newcommand{\fc}{\mathfrak{c}}
\newcommand{\fd}{\mathfrak{d}}

\newcommand{\im}{\textup{i}}

\newcommand{\ad}{{\rm ad}}
\newcommand{\di}{{\rm d}}

\newcommand{\vertiii}[1]{{\left\vert\kern-0.25ex\left\vert\kern-0.25ex\left\vert #1
    \right\vert\kern-0.25ex\right\vert\kern-0.25ex\right\vert}}
\newcommand{\opnorm}[1]{{\vert\kern-0.25ex\vert\kern-0.25ex\vert #1
    \vert\kern-0.25ex\vert\kern-0.25ex\vert}}

\newcommand{\sL}{\mathscr{L}}

\theoremstyle{plain}
\newtheorem{lemma}{Lemma}[section]
\newtheorem{theorem}[lemma]{Theorem}

\newtheorem{proposition}[lemma]{Proposition}
\newtheorem{remark}[lemma]{Remark}
\newtheorem{definition}[lemma]{Definition}

\theoremstyle{definition}

\numberwithin{equation}{section}
\counterwithin{lemma}{section}

\allowdisplaybreaks

\title{Modulational instability of Wilton ripples in deep water}

\begin{document}

\author{Ting-Yang Hsiao\footnote{International School for Advanced Studies (SISSA), Via Bonomea 265, 34136, Trieste, Italy. \newline
	\textit{Emails:}  \texttt{thsiao@sissa.it}, \texttt{amaspero@sissa.it}, \texttt{stefano.pasquali@sissa.it}},  Alberto Maspero$^*$,  Stefano Pasquali$^*$ and Changzhen Sun\footnote{Université Marie et Louis-Pasteur, CNRS, Laboratoire de Mathématiques de Besancon (UMR CNRS 6623), F-25000 Besancon, France \newline
	\textit{Emails:}  \texttt{changzhen.sun@cnrs.fr}}}

\maketitle

\begin{abstract}
We investigate the modulational instability of
small-amplitude Wilton ripples for the two-dimensional
gravity--capillary water-wave equations in infinite depth.
We consider the two branches of $2\pi$-periodic
traveling waves arising from the second-order harmonic
resonance, corresponding to the surface tension
$\kappa=1/2$. On the positive branch, we provide a complete description of the six
eigenvalues of the linearized water-wave operator near
the origin, for Floquet exponents in a parabolic
neighborhood of $\mu=\frac98\epsilon^2$, where
$\epsilon$ denotes the wave amplitude.
In particular we prove the existence of two disjoint instability
intervals, separated by order $\epsilon^3$ and each with width of $\epsilon^{9/2}$.
As the Floquet exponent varies in these intervals,
the unstable eigenvalues trace two distinct, strictly
convex spectral isolas in the complex plane, with
maximal growth rate of order $\epsilon^{7/2}$. On the negative branch, all six eigenvalues near the origin remain purely imaginary throughout the interval $\e^2\leq\mu\leq 2\e^2$ for all sufficiently small $\e>0$.

The proof of the instability result relies on an exact symplectic decomposition
of the six-dimensional spectral problem into three
two-dimensional Hamiltonian blocks, two of which
generate the unstable isolas, while the third has
purely imaginary eigenvalues.
In contrast with the classical Benjamin--Feir
instability of Stokes waves, characterized by a
figure-eight spectrum through the origin, the
resonance produces detached spectral isolas.
\end{abstract}

\small
\tableofcontents
\normalsize

\section{Introduction}\label{sec1}

Wilton ripples are a family of spatially periodic
traveling-wave solutions of the gravity--capillary
water-wave equations. Their origin goes back to
Wilton's work of 1915 \cite{Wilton}, which extended
the classical perturbative construction of Stokes
waves to include surface tension.
Wilton observed that, for certain resonant values
of the capillarity parameter $\kappa$, the 
Stokes expansion becomes singular due to a resonance
between the fundamental mode and a higher harmonic.
This resonance prevents the higher harmonic from
being determined as a successive correction in the
standard Stokes expansion and requires a modified
perturbative construction, in which the resonant
harmonic appears already at low-perturbative  orders. 
In particular, for the second-order resonance, the
fundamental mode and its second harmonic are both
present at leading order.

The resonance has a particularly transparent
interpretation in infinite depth.
Normalizing gravity to $g=1$ and the fundamental
wave number to $1$, the linear gravity--capillary
dispersion relation is
\[
\omega(k):=\sqrt{|k|+\kappa |k|^3},
\qquad
c_{\mathrm{ph}}(k):=\frac{\omega(k)}{|k|}
=\sqrt{\frac{1}{|k|}+\kappa |k|}.
\]
The fundamental wave number and its $n$th harmonic
have the same phase speed precisely when
\[
c_{\mathrm{ph}}(1)=c_{\mathrm{ph}}(n)
\quad\Longleftrightarrow\quad
\kappa=\frac1n,\qquad n\geq 2.
\]
Thus, whenever the capillarity $\kappa$ belongs to
the resonant set
$\left\{\frac1n \colon n\geq 2\right\}$,
the linearized traveling-wave equation at the zero
state, with speed equal to the common $c_{\mathrm{ph}}(1)$ speed,
has a four-dimensional kernel, instead of the
two-dimensional kernel arising in the non-resonant
case (reduced to one-dimensional exploiting symmetries). 
This degeneracy makes the construction of
Wilton ripples a  higher-dimensional
bifurcation problem.

The existence and bifurcation of Wilton ripples
have been investigated through several approaches.
Following the formal studies of Pierson and Fife
\cite{PiersonFife1961} and Chen and Saffman
\cite{ChenSaffman1979}, Reeder and Shinbrot
\cite{reeder1981wiltonI,reeder1981wiltonII}
established rigorous existence results through
an iterative method in Banach spaces, accounting
for the nonlinear solvability conditions associated
with the resonant kernel (see also the recent \cite{akers2021wilton, langer2024wilton}).
Toland and Jones \cite{TolandJones1985}
subsequently investigated the existence and
multiplicity of resonant gravity--capillary waves,
combining Lyapunov--Schmidt reduction with 
secondary-bifurcation theory.
A different approach, started by Craig and Nicholls \cite{CraigNicholls2000},   exploits the variational
structure of the water-wave equations by combining variational and
topological methods.
More recently, Barbieri, Berti, Maspero and
Mazzucchelli
\cite{BBMM}
developed a variational Lyapunov--Schmidt
reduction and applied critical-point theory
to prove the existence of small-amplitude
traveling waves for all values of the physical
parameters, including resonant configurations
with a four-dimensional kernel.
More general variational 
methods have  been developed by Barbieri, Berti and Mazzucchelli
\cite{BarbieriBertiMazzucchelli2026_3D}
for the bifurcation  of three-dimensional gravity--capillary traveling waves.
Finally, very recently  Shavit, Wang and Buckmaster
\cite{ShavitWangBuckmaster2026} developed a
neural-network-based approach to compute
large-amplitude Wilton ripples, achieving
near-machine-precision accuracy and complementing
the numerical construction with rigorous
computer-assisted  certification.

\paragraph{Linear instability of Wilton ripples.}
The aim of this paper is to investigate the linear
stability and instability of small-amplitude Wilton
ripples under long-wave perturbations.
Precisely we consider the two-dimensional gravity--capillary
water-wave equations in infinite depth, at the resonant
value 
\begin{equation}\label{def:fR}
\kappa=\tfrac12
\end{equation}
corresponding to Wilton ripples
of order $n=2$.
Reeder and Shinbrot
\cite{reeder1981wiltonI,reeder1981wiltonII}
established the existence of two families of
small-amplitude, $2\pi$-periodic traveling waves
\[
\big(
\eta_{\e,\tts}(x-c_\e t),
\psi_{\e,\tts}(x-c_\e t)
\big),
\qquad \tts\in\{\pm1\},
\]
whose profiles depend analytically on the amplitude
parameter $\e$, see Theorem \ref{Thm: Wilton expansion}, and whose 
wave speed satisfies
\[
c_\e=\ck+\mathcal O(\e),
\qquad \ck:=\sqrt{\tfrac32}.
\]
As a consequence of the resonance, the second harmonic
appears already at leading order in their small-amplitude
expansions, see \eqref{exp:Wil2}.

We focus on the positive Wilton branch $\tts=1$, and, 
in a reference frame moving with the wave,
the linearized water-wave equations at 
$\big( \eta_{\e,1}, \psi_{\e,1}\big)$, take the form
$h_t=\mathcal L_{\e, 1} h$,
where $\mathcal L_{\e, 1}$ is a linear operator with
$2\pi$-periodic coefficients; see \eqref{mathcal L e}.
This operator is obtained from the Zakharov formulation
through invertible changes of variables involving
Alinhac's good unknown \eqref{alinhac} and the
Levi--Civita transformation \eqref{LC}.

To prove existence of long-wave unstable solutions, we  seek exponentially growing solutions of the form
\[
h(t,x)=\operatorname{Re}
\bigl(e^{\lambda t}e^{\mathrm i\mu x}v(x)\bigr),
\qquad
v(x+2\pi)=v(x),
\qquad \operatorname{Re}\lambda>0,
\]
where $\mu\in\mathbb R$ is the Floquet exponent.
By Bloch-Floquet theory, such $\lambda$ is an eigenvalue of the
operator 
$\mathcal L_{\mu,\e, 1}
:=e^{-\mathrm i\mu x}\mathcal L_{\e,1}
e^{\mathrm i\mu x}$
acting on $2\pi$-periodic functions.

Jones \cite{Jones1992} used coupled amplitude equations to study the long-wave stability of capillary-gravity waves in second-harmonic resonance. Unstable eigenvalues have been
investigated numerically by Deconinck, Trichtchenko \cite{DT} and Trichtchenko, Deconinck, Wilkening \cite{trichtchenko2016instability}.
More recently, for small-amplitude Wilton ripples in finite depth, Hur and Yang \cite{HY2} established stability and instability criteria near nonzero spectral resonances and derived conditional criteria for modulational instability near the origin.

For second-order Wilton ripples in infinite depth, the six
eigenvalues of $\mathcal L_{\mu,\e,\tts}$ near the origin
behave differently on the two branches. On the positive branch,
we identify two disjoint instability intervals for the Floquet
exponent $\mu$ in a parabolic neighborhood of $\mu=\frac98\e^2$.
As $\mu$ varies over each interval, the corresponding pair of
eigenvalues traces a spectral isola detached from the origin.
On the negative branch, all six eigenvalues remain purely
imaginary throughout $\e^2\leq\mu\leq2\e^2$. 

\begin{itemize}
    \item[$\bullet$] {\bf Notation.}
    We use the notation
$\cO_{\mathrm{u}}(a_1,\ldots,a_p)$
to denote a 
remainder function $R(a_1, \ldots, a_p)$ satisfying
$|R(a_1,\ldots,a_p)|\leq C\sum_{j=1}^p |a_j|$
for some constant $C>0$ independent of the
parameters $a_1,\ldots,a_p$.
\end{itemize}

\begin{figure}[htbp]
\centering
\includegraphics[width=\linewidth]{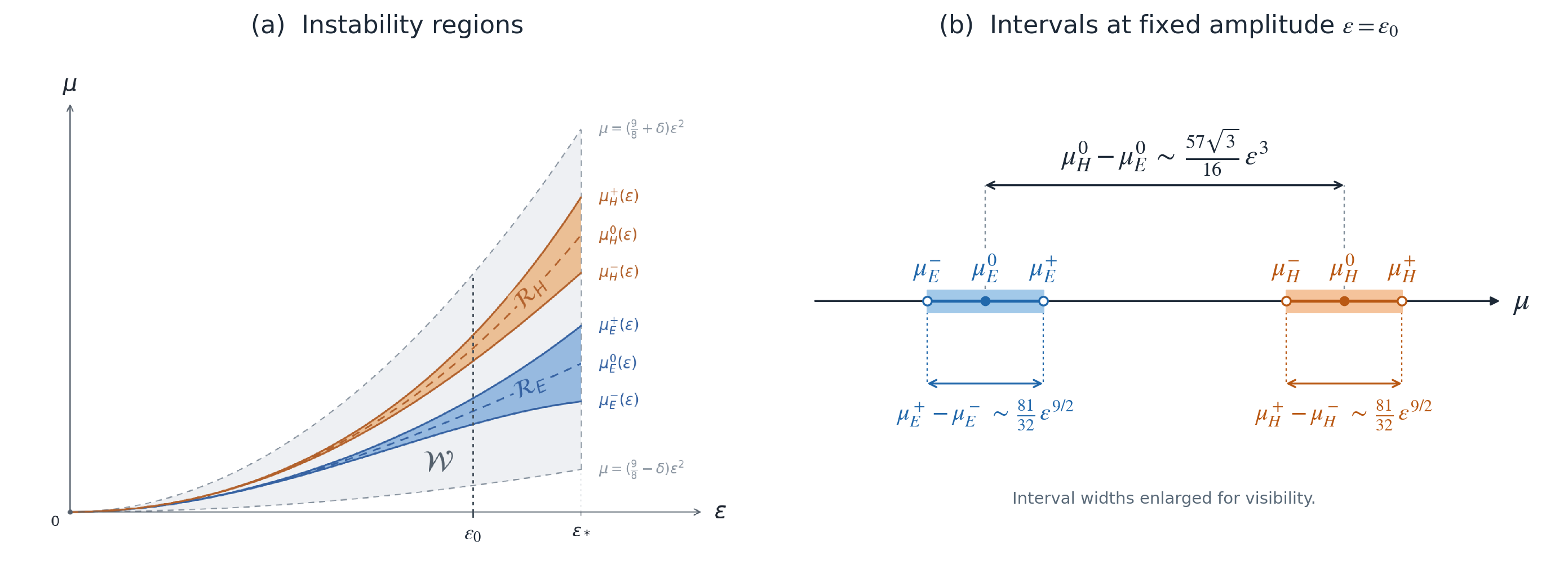}
\caption{(a) The region $\cW$ of Theorem~\ref{thm:main}.
In the blue region $\cR_E$, the pair $\lambda_E^\pm$ leaves
the imaginary axis. In the orange region $\cR_H$, the pair
$\lambda_H^\pm$ leaves the imaginary axis. In the gray region,
all six eigenvalues near the origin are purely imaginary.
(b) A slice at fixed amplitude $\e=\e_0$, indicated by the
vertical dotted line in (a). The upper arrow shows the distance
between the two centers, of order $\e^3$. The lower arrows show
the widths of the two instability intervals, each of order
$\e^{9/2}$, as given by \eqref{eq:intro-original-boundaries}.
The interval widths are enlarged for visibility.
}
\label{fig:thm1}
\end{figure}

\begin{figure}[h]
\includegraphics[width=12cm]{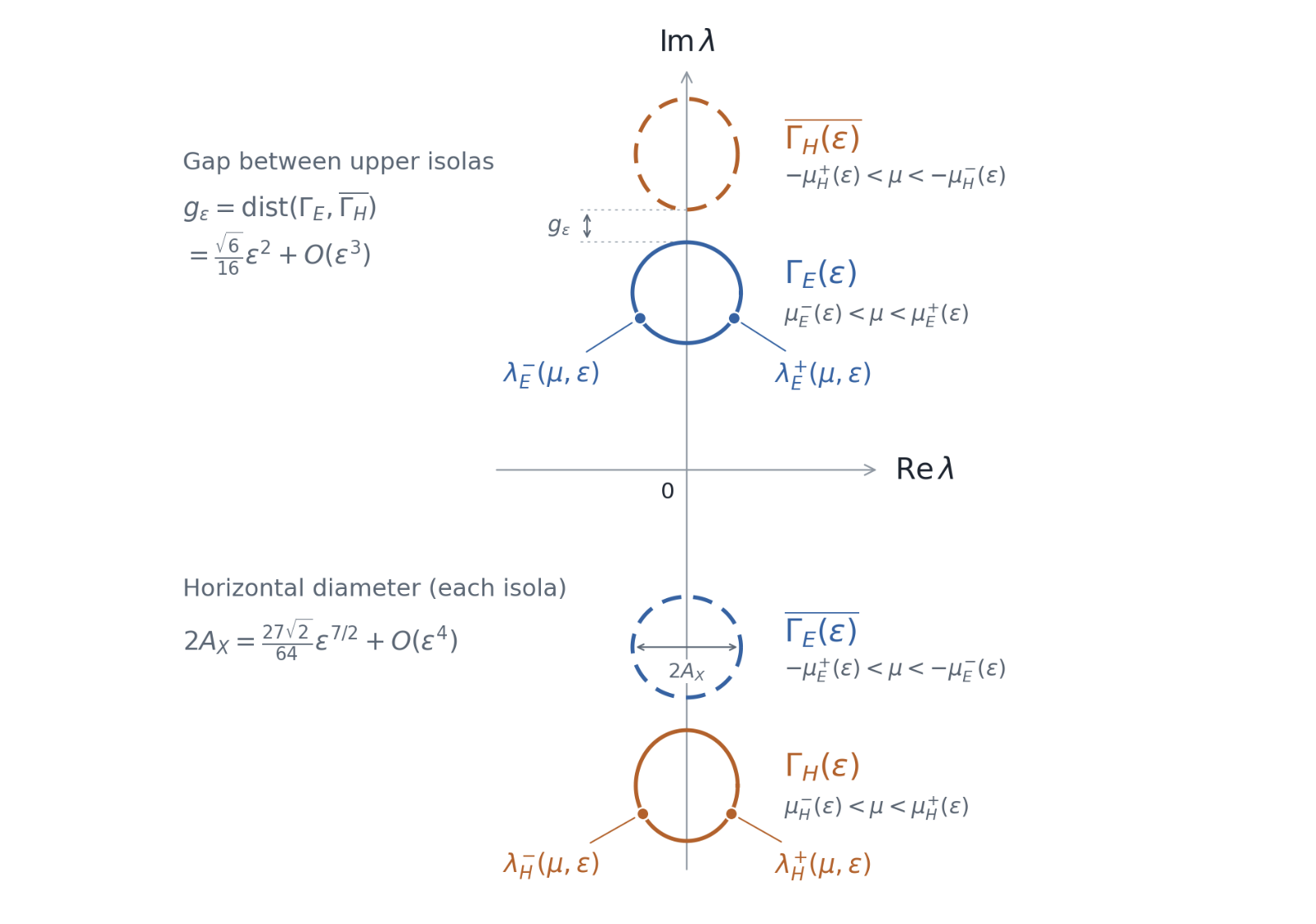}
\centering
\caption{
The spectral isolas $\Gamma_E(\e)$ and $\Gamma_H(\e)$
traced by the eigenvalue pairs
$\lambda_E^\pm(\mu,\e)$ and
$\lambda_H^\pm(\mu,\e)$
as the positive Floquet exponent varies over the
corresponding instability intervals.
The reflected isolas are generated by negative
Floquet exponents.
The sizes of the isolas relative to their
mutual distances are enlarged for readability.
}
\label{fig:isolas}
\end{figure}

We first state a result concerning the instability of the positive branch of Wilton ripples near the origin.
\begin{theorem}
\label{thm:main}
There exist $\epsilon_*>0$, $0<\delta<\frac18$, and continuous functions
$\mu_X^{0},\mu_X^{-},\mu_X^{+}:[0,\epsilon_*)\to\mathbb R$,
$X\in\{E,H\}$, vanishing at $\epsilon=0$, fulfilling 
\[
 \mu_{E}^-<\mu_E^{0}<\mu_E^{+}
 <\mu_H^{-}<\mu_H^{0}<\mu_H^{+} \ , \quad \mbox{ for }   0<\epsilon<\epsilon_* ,
\]
and expanding, for $X \in \{E, H\}$, as 
\begin{equation}\label{eq:intro-original-boundaries}
\begin{aligned}
 \mu_{X}^0(\epsilon)=\frac98\epsilon^2-
 \left(\frac{153}{64}+\sigma_X\frac{57\sqrt3}{32}\right)
 \epsilon^3+\cO_{\mathrm{u}}(\epsilon^4),  \qquad 
 \mu_X^{\pm}(\epsilon)=\mu_X^{0}(\epsilon)
 \pm\frac{81}{64}\epsilon^{9/2}+\cO_{\mathrm{u}}(\epsilon^5),
\end{aligned}
\end{equation}
where $\sigma_E = 1$, $ \sigma_H = -1$, 
such that the following holds true. For any $(\mu, \e)$ in the region 
\begin{equation}\label{def:W}
 \mathcal W:= \cW(\e_*, \delta):=\left\{(\mu,\epsilon):0<\epsilon<\epsilon_*,\quad
 \left|\mu-\frac98\epsilon^2\right|<\delta\epsilon^2\right\}
\end{equation}
the operator $\cL_{\mu, \e, 1}$ has six eigenvalues  near the origin, counted with multiplicity, separated from the rest of the spectrum, that we denote as   $\lambda_E^\pm(\mu, \e)$, $\lambda_H^\pm(\mu, \e)$, $\lambda_G^\pm(\mu, \e)$.
The eigenvalues $ \lambda_G^\pm(\mu, \e)$ are always purely imaginary in $\mathcal W$, whereas $\lambda_X^\pm(\mu, \e)$, $X \in \{ E, H\}$, have non-zero real part if and only if $(\mu, \e)$ belongs to the unstable region
\begin{equation}\label{def:RX}
   \mathcal R_X:=\left\{(\mu,\epsilon)\in\mathcal W:
 \mu_{X}^-(\epsilon)<\mu<\mu_{X}^+(\epsilon)\right\},
 \quad X\in\{E,H\} \ , 
\end{equation}
see Figure \ref{fig:thm1}.
At each boundary $\mu=\mu_{X}^\pm(\epsilon)$, the corresponding
pair of eigenvalues $\lambda_X^{\pm}(\mu, \e)$ collide  on the imaginary axis.
\end{theorem}
The next theorem provides asymptotic expansions of the
unstable eigenvalues and proves that, for each fixed
$\e>0$, they trace closed 
spectral isolas in the complex plane,
as the Floquet parameter $\mu$ varies in the intervals $[\mu_X^-(\e), \mu_X^+(\e)]$.
\begin{theorem}\label{prop:eigenvalue-isolas}
With the notation of Theorem~\ref{thm:main}, the following holds
for each $X\in\{E,H\}$.
\smallskip

\noindent{\bf Eigenvalue asymptotics.} 
There exist real-valued functions $\omega_X,D_X$ defined on  $\mathcal W$, such that, denoting 
$\nu_X^\pm(\e):=\mu_X^{\pm}(\e)-\mu_{X}^0(\e)$, one has
\begin{equation}\label{eig.thm}
 \lambda_X^\pm(\mu_{X}^0(\epsilon) + \nu,\epsilon)=
 \begin{cases}
 \im \omega_X(\mu_{X}^0(\epsilon) + \nu,\epsilon)\pm\sqrt{D_X(\mu_{X}^0(\epsilon) + \nu,\epsilon)},
 &\nu \in (\nu_X^-(\e), \nu_X^+(\e)),\\[1mm]
 \im \omega_X(\mu_{X}^0(\epsilon) + \nu,\epsilon)\pm \im \sqrt{|D_X(\mu_{X}^0(\epsilon) + \nu,\epsilon)|},
 &\nu \not\in (\nu_X^-(\e), \nu_X^+(\e)) ,
 \end{cases}
\end{equation}
where 
\begin{equation}\label{eq:intro-original-normal-form}
\begin{aligned}
 D_X(\mu_{X}^0(\epsilon)+\nu,\epsilon)
 &=\frac{729}{8192}\epsilon^7-\frac{\nu^2}{18\epsilon^2}
 +\cO_{\mathrm{u}}\left(\epsilon^{15/2}, \epsilon^5|\nu|, \frac{\nu^2}{\epsilon}, \frac{|\nu|^3}{\epsilon^4}\right).
\end{aligned}
\end{equation}

\noindent{\bf The  isolas.}
For each fixed $0<\epsilon<\epsilon_*$ and as $\mu$ varies
in $[\mu_X^-(\e), \mu_X^+(\e)]$,
the pair of unstable eigenvalues $\lambda_X^\pm(\mu, \e)$ depicts  a regular real analytic simple closed curve $\Gamma_X(\e)$ bounding a strictly convex domain in the complex plane, intersecting orthogonally the imaginary axis and encircling a convex region. The curve $\Gamma_E(\e)$ is in the upper imaginary semiplane, whereas $\Gamma_H(\e)$ is in the lower one,   see Figure \ref{fig:isolas}.
 In particular,
\[
 \max_{\lambda\in\Gamma_X(\epsilon)}|\operatorname{Re}\lambda|
 =\frac{27\sqrt2}{128}\epsilon^{7/2}+\cO_{\mathrm{u}}(\epsilon^4).
\]
\end{theorem}

Let us make some comments.

\smallskip
\noindent
\emph{1. The instability regions.}
Theorems~\ref{thm:main} and~\ref{prop:eigenvalue-isolas}
provide a complete description
of the six eigenvalues near the origin for
$(\mu,\e)\in\cW$, and identifies two distinct
regions of modulational instability.
The  curve $\mu=\frac98\e^2$
gives the common leading-order location of two
spectral collisions, which separate at order
$\e^3$.
The centers of the two instability
intervals satisfy
$\mu_H^0(\e)-\mu_E^0(\e)
=\frac{57\sqrt3}{16}\e^3+\cO_{\mathrm{u}}(\e^4),$
whereas each interval has width
$\mu_X^+(\e)-\mu_X^-(\e)
=\frac{81}{32}\e^{9/2}+\cO_{\mathrm{u}}(\e^5)$, 
$X\in\{E,H\}$; see Figure~\ref{fig:thm1}.
Thus, the instability is confined to two narrow,
well-separated intervals of Floquet exponents.
Inside each interval, a pair of eigenvalues
acquires nonzero real parts of opposite signs,
whereas outside these intervals all six eigenvalues
 remain purely imaginary.

\smallskip
\noindent
\emph{2. The spectral isolas and their symmetries.}
For each fixed $0<\e<\e_*$, as the positive Floquet
exponent $\mu$ varies in
$[\mu_X^-(\e),\mu_X^+(\e)]$, the two eigenvalues
$\lambda_X^\pm(\mu,\e)$ trace a regular,
real-analytic, simple closed curve
$\Gamma_X(\e)$, $X\in\{E,H\}$.
These \emph{spectral isolas}
bound strictly convex domains and intersect
the imaginary axis orthogonally at the
endpoints of the corresponding instability intervals.
The isola $\Gamma_E(\e)$ lies in the upper
half-plane, whereas $\Gamma_H(\e)$ lies
in the lower half-plane.
Each isola is symmetric with respect to the
imaginary axis, as a consequence of the
Hamiltonian symmetry
$\lambda\mapsto-\overline{\lambda}$
at fixed $\mu$.
On the other hand, the relation
$\sigma(\cL_{-\mu,\e,1})
=\overline{\sigma(\cL_{\mu,\e,1})}$
produces the conjugate isolas
$\overline{\Gamma_E(\e)}$ and
$\overline{\Gamma_H(\e)}$
for negative Floquet exponents.
Consequently, the near-origin unstable spectrum
contains four isolas, whose union is symmetric
with respect to both coordinate axes; see Figure \ref{fig:isolas}.

The two isolas in the upper half-plane,
$\Gamma_E(\e)$ and $\overline{\Gamma_H(\e)}$,
are well separated, as
$\operatorname{dist}
\bigl(\Gamma_E(\e),\overline{\Gamma_H(\e)}\bigr)
=\frac{\sqrt6}{16}\e^2+\cO_{\mathrm{u}}(\e^3),$
whereas their horizontal and vertical diameters
are of order 
$\frac{27\sqrt2}{64}\e^{7/2}$.
Thus, the distance between neighboring isolas
is asymptotically much larger than their size, 
and the  maximal unstable eigenvalues  have real part 
of order $\e^{7/2}$.

\smallskip
\noindent
{\bf Strategy of the proof.}
At $(\mu,\e)=(0,0)$, the Bloch operator
$\mathcal L_{0,0}$ has a zero eigenvalue of
algebraic multiplicity six.
For $\mu=0$ and sufficiently small $\e>0$, four of these eigenvalues
remain at zero, as a consequence of the translation invariance, the invariance under the addition of a constant to
the velocity potential, the existence of a smooth family of traveling waves. The other two bifurcate along the imaginary axis with leading
terms $\pm\im\sqrt{9/8}\,\e$.
The emergence of this pair is a consequence
of the harmonic resonance. Indeed, since the fundamental mode and
its second harmonic both belong to the kernel of the unperturbed
operator, their coupling through the order-$\e$ terms of the
linearized operator already produces a nontrivial splitting of the
zero eigenvalue at order $\e$, see Remark \ref{remark im omega0}. For $\e=0$ and $\mu\to0^+$, the two long-wave eigenvalues have
leading terms $\pm\im\sqrt{\mu}$, while the other four are
$\cO(\mu)$; see \eqref{eigenvalues of Lmu0}.
Comparing these expansions suggests two collisions near
$\mu=\frac98\e^2$, one on the positive imaginary axis and
the other on the negative imaginary axis.

The leading-order calculation predicts the same two collisions on both Wilton branches, yet only the positive branch develops instability near them.
For the positive branch $\tts=1$, each collision
generates an instability interval in which
a pair of eigenvalues leaves the imaginary axis.
These intervals are disjoint and give rise to the two spectral isolas described Theorem~\ref{prop:eigenvalue-isolas}; see Figure \ref{fig:isolas}.
In contrast, for the negative branch $\tts=-1$,
all six eigenvalues remain purely imaginary
throughout the region
$\e^2\leq\mu\leq2\e^2$,
as established in Proposition \ref{final proposition}.
This difference motivates the separate symplectic reductions developed for the two branches.

Following \cite{BMV1,BMV3,BMV_ed}, we use Kato's transformation
to represent the Bloch operator on its six-dimensional spectral
subspace near the origin by a Hamiltonian matrix. We set
$\e=\ve^2$, $\mu=(\tau\ve)^4$ and perform the singular symplectic
change of basis in Lemma~\ref{lem:first-symplectic-change}.
On the positive Wilton branch, the expansions
\eqref{eq:E1-first-change} and \eqref{eq:H1-first-change} show
that the eigenvalues of $\mathsf J_2 E^{(1)}$ have leading terms
$\im\sqrt{\mu}$ and $\im\sqrt{9/8}\,\e$, while those of
$\mathsf J_2 H^{(1)}$ have the corresponding negative imaginary
leading terms. The eigenvalues of $\mathsf J_2 G^{(1)}$ are
$\cO(\e^2)$. Thus, at $\tau=(9/8)^{1/4}$, the two eigenvalues
within each of the first two blocks agree to leading order.
The spectra of distinct blocks remain separated by a distance
of order $\e$, which allows us to solve the Sylvester equations
and complete the exact symplectic block diagonalization.

The main difficulty is to determine whether the eigenvalues of
the two critical blocks leave the imaginary axis. For $X\in\{E,H\}$,
\eqref{0911:2314} gives
\[
  \lambda_X^\pm=\im b_X\pm\sqrt{d_X^2-a_X^2},
\]
where $a_X$, $b_X$ and $d_X$ are real-valued. On the center curve
defined by $a_X=0$, we have $|\operatorname{Re}\lambda_X^\pm|=|d_X|$.
The coefficients of $\e^{5/2}$ and $\e^3$ in $d_X$ vanish. Along the center curve, the first nonzero term of $d_X$ is of order $\e^{7/2}$. Its coefficient includes the correction
from the off-diagonal couplings computed in
Proposition~\ref{prop:first-block-decoupling}.

By \eqref{eq:intro-original-boundaries}, the distance
$\mu_H^0(\e)-\mu_E^0(\e)$ between the centers of the two
instability intervals is of order $\e^3$, while each interval
has width of order $\e^{9/2}$; see Figure~\ref{fig:thm1}. The $\cO(\e^4)$ error bound
in the cubic approximation of each center is therefore too coarse
to resolve the interval. In Lemma~\ref{lem:mu-plus-center-endpoints},
we locate the endpoints relative to the exact center by solving
$a_X=-d_X$ and $a_X=d_X$. As $\mu$ varies between these endpoints,
the two eigenvalues trace a closed curve. After translation and
division by $\e^{7/2}$, its parametrization converges in $C^2$
to that of a circle, yielding the strict convexity of each isola
(Proposition~\ref{prop:isola}).

On the negative Wilton branch, the leading spectra of the $E$ and
$H$ blocks coincide at $\tau=(9/8)^{1/4}$, so we keep the four
modes in a single block. The proof of
Proposition~\ref{final proposition} uses the Hamiltonian identity
\[
2\operatorname{Re}\lambda\,w^*\mathsf B w=0
\]
for an eigenpair $(\lambda,w)$ of $\mathsf L=\mathsf J\mathsf B$.
The leading quadratic form on the four-dimensional block is
positive definite. An estimate on the coupling to the remaining
two modes gives $w^*\mathsf B w>0$ for the eigenvectors associated
with the four eigenvalues of order $\e$. For the remaining pair,
a purely imaginary shift of the matrix makes the corresponding
quadratic form negative on its eigenvectors. The identity then
gives $\operatorname{Re}\lambda=0$ for all six eigenvalues.

The different spectral behaviors of the two Wilton branches can be
understood in terms of the Krein signatures of the colliding modes.
Indeed, the leading quadratic forms of the $E$ and $H$ blocks have
eigenvalues $\e\tau^2$ and $-\s\e\sqrt{9/8}$. Thus, on the positive
Wilton branch, the quadratic form restricted to each leading-order
collision eigenspace is indefinite, with one positive and one
negative direction, allowing the eigenvalues to leave the imaginary
axis. On the negative branch, the corresponding quadratic form
is positive definite, which means that the Krein signatures are identical, preventing such an instability. 

\paragraph{Comparison with the Benjamin--Feir instability
of Stokes waves.}

The stability of periodic traveling water waves
has been a central question since Stokes's
pioneering work of 1847 \cite{stokes};
see also the survey \cite{ebb}.
In 1960s, experimental and theoretical work by Benjamin and Feir  \cite{BF}, Lighthill \cite{Li} and Zakharov \cite{Zak68}
showed that small-amplitude
gravity Stokes waves in deep water are unstable
under long-wave perturbations.
The first rigorous proofs were obtained by
Bridges and Mielke \cite{BrM} in finite depth
and Nguyen and Strauss
\cite{NS} in infinite depth.
A complete description of the spectrum near
the origin was later obtained by Berti,
Maspero and Ventura for gravity Stokes waves
in deep water \cite{BMV1}, finite depth
\cite{BMV3}, and at the critical
Whitham--Benjamin depth \cite{BMV_ed}, 
establishing the  \emph{figure-eight}
configuration traced by the unstable eigenvalues
as the Floquet exponent varies. In deep water, the unstable interval of Floquet exponents has width of order $\e$. The modulational stability of non-resonant gravity-capillary Stokes waves has also been studied in finite depth \cite{HY2,SW,HM} and in deep water \cite{HW}.

The instability of Wilton ripples 
in Theorems~\ref{thm:main} and~\ref{prop:eigenvalue-isolas} exhibits a
qualitatively different spectral picture.
The second-order harmonic
resonance of Wilton ripples leads to six
eigenvalues near the origin --instead of the four of Stokes waves-- whose instability
occurs in two narrow intervals of positive
Floquet exponents centered at
$\mu=\frac98\e^2$ to leading order.
The corresponding unstable eigenvalues trace
two distinct spectral isolas, detached from
the origin, rather than a figure eight.
On the positive Wilton branch, the underlying six-dimensional spectral
problem admits an exact symplectic
decomposition into three two-dimensional
blocks, two of which generate the isolas,
while the third has purely imaginary
eigenvalues in the region under consideration.
We emphasize that this comparison concerns
the spectrum near the origin, since spectral
isolas may also arise near nonzero spectral
collisions for non-resonant Stokes waves.

Finally, we mention related results on
transverse perturbations
\cite{CNS,CNS2,JRSY,HTW, BMR}, the stability of Stokes waves with constant vorticity
\cite{HsiaoThesis,HHY1,CN1,HHY2,FS,CN2}, other aspects of the stability of gravity
Stokes waves
\cite{HY,ChenSu,BMV4,BCMV,
Mc1,Mc2,CD,BDS,CDT},
and instability phenomena in shallow-water
and approximate dispersive models
\cite{HP,HJ,BHJ,BJ,GH,HK,BBHZ,BHR,MaRa, BMP}.

\section{Wilton ripples and their linearized operator}

First, let us introduce the gravity-capillary water waves equations and the  Wilton ripples solution.      

\smallskip
\noindent{\bf The water waves equations.}
We consider the Euler equations of hydrodynamics for a 2-dimensional perfect and incompressible fluid under the action of gravity and capillary forces at the free surface.
The fluid fills the region
    \begin{align*}
       \mathcal{D}_\eta:=\{(x,y)\in \mathbb{R}\times \mathbb{R}: -\infty < y <\eta(t,x)\}\ ,
    \end{align*}
with infinite depth. The irrotational velocity
field is the gradient of a harmonic scalar potential  $\Phi = \Phi(t,x,y)$ determined by  its trace $\psi(t,x):= \Phi(t, x,\eta(t,x))$ at the free surface $y=\eta(t,x)$.
Actually $\Phi$ is
the unique solution of the elliptic equation $\Delta \Phi  = 0$ in $\mathcal{D}_\eta$  with Dirichlet datum
$\Phi(t, x,\eta(t, x))=\psi(t, x)$ and $\Phi_y(t,x,y) = 0$ at $y \to - \infty$.
 Imposing that the fluid particles
at the free surface remain on it along the evolution (kinematic boundary condition) and that the pressure of
the fluid plus the capillary forces at the free surface is equal to the constant atmospheric pressure (dynamic boundary condition), the time evolution of the fluid is determined by the non-local quasi-linear equations \cite{Zak68, CS}
\begin{equation} \label{Craig-Sulem formula}
\left\{\begin{aligned}
    \eta_t&=G(\eta)\psi\\
    \psi_t&=-g\eta-\frac{\psi_x^2}{2}+\frac{1}{2(1+\eta^2_x)}\left(G(\eta)\psi+\eta_x\psi_x\right)^2+\kappa\left(\frac{\eta_x}{(1+\eta_x^2)^{1/2}}\right)_x,
\end{aligned}\right.
\end{equation}
where $g>0$ is the gravity constant and $G(\eta)$ denotes the Dirichlet-Neumann operator $[G(\eta)\psi](x):=\Phi_y(x,\eta(x))-\Phi_x(x,\eta(x))\eta_x(x)$.
 In the sequel,
with no loss of generality, we set the gravity constant $g=1$.

System \eqref{Craig-Sulem formula} is Hamiltonian and can be written as
\begin{align} \label{pa_t eta psi}
     \partial_t
    \begin{bmatrix}
    \eta\\
    \psi
    \end{bmatrix} &=
    \cJ
    \begin{bmatrix}
    \nabla_{\eta} \mathcal{H}\\
    \nabla_\psi \mathcal{H}
    \end{bmatrix}
    , \quad \cJ :=
    \begin{bmatrix}
    0& \mathrm{Id}\\
    -\mathrm{Id} &0
    \end{bmatrix},
\end{align}
where $\nabla$ denote the $L^2$-gradient, and the Hamiltonian
\begin{align} \label{mathcal H}
\mathcal{H}(\eta,\psi):=\int_{\mathbb{T}} \left(\frac{1}{2}\psi G(\eta)\psi+\frac{1}{2}\eta^2+\kappa(\sqrt{1+\eta^2_x}-1)\right) \mathrm{d}x , \quad \mathbb{T} := \mathbb{R}/(2\pi\mathbb{Z}),
\end{align}
is the sum of the kinetic and potential energy of the fluid.
In addition,  the water waves system \eqref{Craig-Sulem formula} is reversible with respect to the involution
\begin{equation} \label{rho involution}
    \begin{aligned}
        \rho \begin{bmatrix}
            \eta(x)\\
            \psi(x)
        \end{bmatrix}:=
        \begin{bmatrix}
            \eta(-x)\\
            -\psi(-x)
        \end{bmatrix},~~\mathrm{i.e.}~\mathcal{H}\circ \rho=\mathcal{H},
    \end{aligned}
\end{equation}
and it is space invariant.

 \paragraph{Wilton ripples.} As mentioned in the introduction, we assume that $\kappa = \frac{1}{2}$ (see \eqref{def:fR}). In a moving reference frame with constant speed $c$ (and with normalized gravity $g=1)$, the water waves system \eqref{Craig-Sulem formula} becomes
 \begin{equation} \label{in the reference frame}
\left\{\begin{aligned}
    \eta_t&=c\,\eta_x+G(\eta)\psi,\\
    \psi_t&=c\,\psi_x-\eta-\frac{\psi_x^2}{2}+\frac{1}{2(1+\eta^2_x)}\left(G(\eta)\psi+\eta_x\psi_x\right)^2+\kappa \frac{\eta_{xx}}{(1+\eta_x^2)^{3/2}} .
\end{aligned}\right.
\end{equation}
By \cite{reeder1981wiltonII} (see in particular  Lemma 5.1, Theorem 5.2 and Corollary 5.3), there exist 
small amplitude {\em Wilton ripples of order $2$}, namely   stationary solutions, $2\pi$-periodic of \eqref{in the reference frame} bifurcating from a double resonance.
To state the precise result, 
we denote by $B(r):=\{x\in\mathbb{R}: |x|<r\}$ the open ball of radius $r$ centered at zero, and by 
$ H^{\sigma,s}_{\mathtt{ev}}(\T) $  (respectively $ H^{\sigma,s}_{\mathtt{odd}}(\T) $) the space of even (respectively odd), real valued $ 2 \pi $-periodic functions
$ u(x) = \sum_{k \in \mathbb{Z}} u_k e^{\im k x} $
such that $ \| u \|_{\sigma,s}^2 := \sum_{k \in \mathbb{Z}} |u_k|^2 \langle k \rangle^{2s} 
e^{2 \sigma |k|} < + \infty$. 

\begin{theorem} [Wilton ripples of order $2$] \label{Thm: Wilton expansion} 
Let $\kappa = \frac{1}{2}$. There exist $\e_{\star}>0$ and a family  of real analytic solutions $(\eta_{\e,\tts}(x), \psi_{\e,\tts}(x), c_{\e,\tts})$, parameterized by the amplitude $|\e| \leq \e_{\star}$ and by $\tts \in \{ \pm 1 \}$, of 
\begin{equation}\label{travelingWWstokes}
c \, \eta_x+G(\eta)\psi = 0 \, , \quad 
c \, \psi_x -  \eta - \dfrac{\psi_x^2}{2} + 
\dfrac{1}{2(1+\eta_x^2)} \big( G(\eta) \psi + \eta_x \psi_x \big)^2  +
\kappa \frac{\eta_{xx}}{(1+\eta_x^2)^{3/2}}
= 0 \, , 
\end{equation}
  such that
 $ \eta_{\e,\tts} (x), \psi_{\e,\tts} (x) $ are $2\pi$-periodic;  $\eta_{\e,\tts} (x) $ is even
and $\psi_{\e,\tts} (x) $ is odd, of the form 
  \begin{align}
  \eta_{\e,\tts} (x) &= \e \left[  \cos(x) + \frac{\tts}{2} \cos(2x) \right] +\e^2 \bigg[ - \cos(2x) -\frac32 \tts \cos(3 x) - \frac14 \cos(4x) \bigg] +\cO(\e^3)  , \label{exp:Wil2} \\ 
  \notag
  \psi_{\e,\tts} (x)  &=  \e \left[  \ck \sin(x) +  \frac{3 \tts}{4 \ck} \sin(2 x) \right]  +\e^2 
  \bigg[ 
- \frac{\ck \tts}{4} \sin(x) -\frac{7\ck}{8} \sin(2 x) -\frac{3 \ck \tts}{2}\sin(3x) - \frac{\ck}{4} \sin(4x) \bigg] +\cO(\e^3), \\
  c_{\e,\tts} &= \ck + \frac{\ck \tts}{4}\e  -\frac{29\ck}{32}\e^2 +\cO(\e^3)\quad \text{where} \quad \ck := (1+\kappa)^{\frac{1}{2}}= \sqrt{\frac{3}{2}} \, , \label{exp:c}
  \end{align}
Moreover, for any $ \sigma \geq  0 $ and $ {s > \frac72} $, there exists $ \e_*>0 $ such that the map $\e \mapsto (\eta_{\e,\tts}, \psi_{\e,\tts}, c_{\e,\tts})$ is analytic from $B(\e_*) \to H^{\sigma,s}_{\mathtt{ev}} (\T)\times H^{\sigma,s}_{\mathtt{odd}}(\T)\times \R$.
\end{theorem}

The expansions in 
\eqref{exp:Wil2}--\eqref{exp:c} are proved in Appendix \ref{sec:App2}.

\paragraph{Linearization.}
We linearize system \eqref{in the reference frame} at the Wilton ripples $(\eta_{\e,\tts}(x),\psi_{\e,\tts}(x))$ given in Theorem \ref{Thm: Wilton expansion} and evaluate $c$ at $c_{\e,\tts}$. 
Using the shape derivative formula \cite{LD book} 
$
\mathrm{d}_\eta G(\eta)[\hat \eta][\psi] = - G(\eta)(B\hat \eta) - \pa_x( V \hat \eta), 
$
where
the functions $(V(x),B(x))$ are the horizontal and vertical components of the velocity field $(\Phi_x,\Phi_y)$ at the free surface and are given by
\begin{align} \label{espV}
   V:= V(x)&:=-B(\eta_{\e,\tts})_x+(\psi_{\e,\tts})_x,\\ \label{espB}
  B:=  B(x)&:=\frac{G(\eta_{\e,\tts})\psi_{\e,\tts}+(\psi_{\e,\tts})_x(\eta_{\e,\tts})_x}{1+(\eta_{\e,\tts})^2_x}=\frac{(\psi_{\e,\tts})_x-c_{\e,\tts}}{1+(\eta_{\e,\tts})^2_x} (\eta_{\e,\tts})_x,
\end{align}
one obtains the real, autonomous linearized  system
\small
\begin{equation} \label{first linear eq}
    \begin{aligned}
        \begin{bmatrix}
            \hat{\eta}_t\\
            \hat{\psi}_t
        \end{bmatrix}=
        \left[\begin{array}{c|c} 
	 -G(\eta_{\e,\tts})B-\partial_x\circ(V-c_{\e,\tts}) & G(\eta_{\e,\tts})  \\ 
	\hline 
	-1+B(V-c_{\e,\tts})\partial_x-B\partial_x\circ(V-c_{\e,\tts})-BG(\eta_{\e,\tts})\circ B+\kappa \partial_x \circ l \circ\partial_x & -(V-c_{\e,\tts})\partial_x+B G(\eta_{\e,\tts})
\end{array}\right] 
        \begin{bmatrix}
            \hat{\eta}\\
            \hat{\psi}
        \end{bmatrix},
    \end{aligned}
\end{equation}
\normalsize
where
\begin{align} \label{esp l}
   l:= l(x):=\frac{1}{(1+(\eta_{\e,\tts})_x^2)^{3/2}}.
\end{align}
The map  $\e \to (V, B, l)$ is analytic as a map $B(\e_0) \to H^{\sigma, s-1}_{\mathtt{ev}}(\T) \times H^{\sigma, s-1}_{\mathtt{odd}}(\T) \times H^{\sigma, s-1}_{\mathtt{ev}}(\T) $. 
The real system \eqref{first linear eq} is Hamiltonian, i.e. of the form
$\cJ \cA$ with $\cA$ symmetric, $\cA = \cA^\top$, where the transpose is taken with respect to the real scalar product of $L^2(\T, \R) \times L^2(\T, \R)$.
Moreover the linear operator in \eqref{first linear eq} is reversible, namely it anti-commutes with the involution $\rho$ in \eqref{rho involution}. \\
Next, we conjugate \eqref{first linear eq} by using the time-independent ``good unknown of Alinhac'' linear transformation 
\begin{align}\label{alinhac}
    \begin{bmatrix}
        \hat{\eta}\\
        \hat{\psi}
    \end{bmatrix}:= Z
    \begin{bmatrix}
        u\\
        v
    \end{bmatrix},~~Z:=Z_{\e,\tts}=\begin{bmatrix}
        1&0\\
        B&1
    \end{bmatrix},
    ~~Z^{-1}=\begin{bmatrix}
        1&0\\
        -B&1
    \end{bmatrix},
\end{align}
yielding the linear system 
\begin{align} \label{second linear eq}
    \begin{bmatrix}
        u_t\\
        v_t
    \end{bmatrix}=  \widetilde{\cL}_{0,\e,\tts} 
    \begin{bmatrix}
        u\\
        v
    \end{bmatrix} \ , \qquad 
    \widetilde{\cL}_{0,\e,\tts} := \begin{bmatrix}
        -\partial_x\circ (V-c_{\e,\tts}) & G(\eta_{\e,\tts})\\
        -1-(V-c_{\e,\tts})B_x+\kappa \partial_x\circ l\circ \partial_x & -(V-c_{\e,\tts})\partial_x
    \end{bmatrix} \ , 
\end{align}
which  is Hamiltonian and reversible since the transformation $Z$ is symplectic, $\mathrm{i.e.}$ $Z^T\mathcal{J} Z=\mathcal{J}$, and satisfies $Z\circ \rho=\rho\circ Z$.

Next, 
we perform a conformal change of variables to flatten 
the water surface. 
Following \cite{NS, BMV1} 
 there exists a diffeomorphism of $\mathbb{T}$,
 $ x\mapsto x+\mathfrak{p}(x)$, with a small $2\pi$-periodic  function $\mathfrak{p}(x) := \mathfrak{p}_{\e,\tts}(x)$, 
 such that, by defining the associated composition operator $ (\mathfrak{P}u)(x) := u(x+\mathfrak{p}(x))$, the Dirichlet-Neumann operator writes as \cite[Prop. 3.3]{NS}
\begin{equation*}
 G(\eta_{\e,\tts}) = \pa_x \circ \mathfrak{P}^{-1} \circ {\mathfrak H} \circ \mathfrak{P} \, , 
\end{equation*}
where $ {\mathfrak H} $ is the Hilbert transform, i.e. the  Fourier multiplier operator
$$
 \mathfrak{H}(e^{\im j x}):= - \im\, \textup{sign}(j) e^{\im j x} \, , 
 \quad  \forall j \in \Z \setminus \{0\} \, , 
 \quad \mathfrak{H}(1) := 0 \, . 
$$
The function $\mathfrak p(x)$ is determined as a fixed point of 
\begin{equation} \label{def:ttp}
\mathfrak{p}  =  \mathfrak{H}[\eta_{\e,\tts} ( x + \mathfrak{p}(x))] \, . 
  \end{equation}
The map $\e \to \mathfrak{p}$ is analytic as a map $B(\e_0) \to H^s_{\mathtt{odd}}(\T)$. Indeed, by Lemma \ref{thm:wilton-spatial-analyticity} the map $(\e,p)\mapsto p-\mathfrak{H}[\eta_{\e,\tts}\circ(\mathrm{Id}+p)]$ is analytic near $(0,0)$ in $\mathbb{R}\times H^s_{\mathtt{odd}}(\mathbb{T})$, and its derivative in $p$ there is the identity since $\eta_{0,\tts}=0$. The analytic implicit function theorem applies as in \cite[p.~660, footnote~2]{BMV1}. In addition, in  Appendix \ref{sec:App2} we prove the following expansion.

\begin{lemma} \label{lem:pfexp}

There exists $\e_0>0$ such that for any $\e \in (0,\e_0)$ the function $\mathfrak{p}$ defined in \eqref{def:ttp} admits the following expansions
\begin{align}
\mathfrak{p}(x)   &= \e \mathfrak{p}_1(x) + \e^2 \mathfrak{p}_2(x) + \cO(\e^3) , \label{pexp}
\end{align}
where
\begin{align}
\mathfrak{p}_1(x) &=  \sin (x) + \frac{\tts}{2}  \sin(2 x) , \qquad 
\mathfrak{p}_2(x) = - \frac{3\tts}{4} \sin(x) - \frac12 \sin( 2x) - \frac{3\tts}{4}\sin( 3x)  . \label{fexp}
\end{align}
\end{lemma}
Under the symplectic and reversibility-preserving change of variables 
\begin{align}\label{LC}
    h=\mathcal{P}\begin{bmatrix}
        u\\
        v
\end{bmatrix},\qquad \mathcal{P}:=\mathcal{P}_{\e,\tts}=\begin{bmatrix}
    (1+\mathfrak{p}_x)\mathfrak{P} & 0\\
    0 & \mathfrak{P}
\end{bmatrix} \ , 
\end{align}
one  transforms the system \eqref{second linear eq} into the linear system $h_t=\mathcal{L}_{0,\e,\tts} h$ where $\mathcal{L}_{0,\e,\tts}$ is the Hamiltonian and reversible real operator
\begin{equation} \label{mathcal L e}
\begin{aligned}
    \mathcal{L}_{0,\e,\tts} := \mathcal{L}_{0,\e,\tts}(\kappa) :=\mathcal{P} \widetilde{\cL}_{0,\e,\tts} \mathcal{P}^{-1} = &\begin{bmatrix}
\partial_x\circ(\ck+p_{\e,\tts}(x)) &  |D| \\
        -(1+a_{\e,\tts}(x))+\kappa \Sigma_{\e,\tts} & ~~(\ck+p_{\e,\tts}(x))\partial_x
    \end{bmatrix}\\
    =&\begin{bmatrix}
        0& \mathrm{Id}\\
        -\mathrm{Id}& 0
    \end{bmatrix}\begin{bmatrix}
 (1+a_{\e,\tts}(x))-\kappa \Sigma_{\e,\tts}& ~~-(\ck+p_{\e,\tts}(x))\partial_x\\
 \partial_x\circ(\ck+p_{\e,\tts}(x))& ~~|D| 
    \end{bmatrix},
\end{aligned}    
\end{equation}
where the functions $p_{\e,\tts}(x)$ and $a_{\e,\tts}(x)$ are given by
\begin{equation}\label{def:pa}
\ck+p_{\e,\tts}(x) :=  \displaystyle{\frac{ c_{\e,\tts}-V(x+\mathfrak{p}(x))}{ 1+\mathfrak{p}_x(x)}} \, , \quad 1+a_{\e,\tts}(x):=   \displaystyle{\frac{1+ (V(x + \mathfrak{p}(x)) - c_{\e,\tts})
 B_x(x + \mathfrak{p}(x))  }{1+\mathfrak{p}_x(x)}} \,, 
\end{equation}
and the operator $\Sigma_{\e,\tts}$  by
\begin{equation} \label{def:Sigma g}
\begin{aligned}
  &   \Sigma_{\e,\tts}: =\frac{1}{1+\mathfrak{p}_x}\partial_x \circ g_{\e,\tts}(x)\circ \partial_x \circ \frac{1}{1+\mathfrak{p}_x},\qquad 
    g_{\e,\tts}(x): =\frac{1}{\big(1+(\eta_{\e,\tts})_x^2(x+\mathfrak{p}(x))\big)^{3/2}}\frac{1}{1+\mathfrak{p}_x}.
\end{aligned}    
\end{equation}
By the analyticity result of the map $\e \mapsto (V, B, l)$ given above,  the map
$\e \to (p_{\e,\tts}, a_{\e,\tts}, g_{\e,\tts})$ is analytic as a map $B(\e_0) \to H^s_{\mathtt{ev}}(\T)\times H^s_{\mathtt{ev}}(\T) \times H^s_{\mathtt{ev}}(\T)$.
In Appendix \ref{sec:App2} we provide their Taylor expansions, that we collect in the following lemma:
\begin{lemma}\label{lem:pa.exp}
The analytic functions $p_{\e,\tts} (x) $ and $a_{\e,\tts} (x) $  in \eqref{def:pa} 
are even in $ x $, and
\begin{equation}\label{SN1}
p_{\e,\tts} (x)  
= \e p_1 (x) + \e^2 p_2 (x)  + \cO(\e^3) \, , \qquad   
a_{\e,\tts} (x)  
= \e a_1(x) +\e^2 a_2 (x) + \cO(\e^3) \, , 
\end{equation}
where
\begin{align}
     p_1(x) & = \frac{ \ck}{4} \, \tts  - \sqrt{6} \cos(x) - \sqrt{6} \, \tts \cos(2 x) ,  \label{pino1fd} \\
     p_2(x) & = \frac{35\ck}{32}  +
      3 \ck\tts \cos(x) + \frac{7\ck}{2} \cos(2 x) +\frac{17\ck \tts}{2}\cos(3x) +\sqrt{6} \cos(4 x) ,  \label{pino2fd} 
\end{align}
and 
 \begin{align} 
 a_1(x) & =  - \frac{5}{2} \cos(x) -4 \, \tts \cos(2 x) , \label{aino1fd}  \\
     a_2(x) & =  \frac{31}{4} + \frac{71\tts}{8}\cos(x) + \frac{33}{4} \cos(2 x) + \frac{233\tts }{8} \cos(3x) + 11\cos(4 x) ,  \label{aino2fd} 
\end{align}
The function $g_{\e,\tts}(x)$ in \eqref{def:Sigma g} is even in $x$ and expands as
\begin{equation} \label{SN1 g}
    g_{\e,\tts}(x) = 1+ \e g_1(x)  + \e^2 g_2(x) + \cO(\e^3),
\end{equation}
where
\begin{align}\label{exp:g1}
    g_1(x) &= -\cos(x) -\tts\cos(2 x), \\
    \label{exp:g2}
    g_2(x) &= -\frac{1}{2}  +\frac{\tts}{4}\cos(x)+\frac{9}{4}\cos(2 x)+ \frac{19\tts}{4}\cos(3 x) + \frac{5}{4}\cos(4 x) .
\end{align}
Finally, the self-adjoint operator $\Sigma_{\e,\tts}$ in \eqref{def:Sigma g} expands as
\begin{equation} \label{Sigma epsilon}
    \Sigma_{\e,\tts} = \pa_{xx} + \e \Sigma_1 + \e^2 \Sigma_2 + \cO(\e^3), 
\end{equation}
where
\begin{align}\label{Sigma 1}
    \Sigma_j:=& d_j(x)\,\pa_{xx}+e_j(x)\pa_{x} + h_j(x), \quad j=1,2 
\end{align}
and
\begin{align}
\label{d1x}
    d_1(x)&=-3 \cos(x) -3 \tts  \cos(2 x) , \quad 
     e_1(x) = 3 \sin(x) + 6 \tts  \sin(2 x) \ , \quad  
     h_1(x) = \cos(x) + 4 \tts \cos(2 x) , \\
     d_2(x)&= \frac{9}{2} + \frac{27\tts}{4}\cos(x) + \frac{27}{4}\cos(2x) + \frac{57\tts}{4}\cos(3x) + \frac{15}{4}\cos(4x), \label{d2x}  \\
     e_2(x)& =   -\frac{27\tts}{4}\sin(x)  -\frac{27}{2}\sin(2x)  -\frac{171\tts}{4}\sin(3x) -15\sin(4x) ,  \label{e2x} \\
     h_2(x)& = -\frac{5}{2}  -\frac{19\tts}{4}\cos(x) -\frac{15}{2}\cos(2x) -\frac{145\tts}{4}\cos(3x)  -14\cos(4x) . \label{h2 h02} 
\end{align}
\end{lemma}

\paragraph{Bloch-Floquet expansions.} 
Since the operator $\mathcal{L}_{\e,\tts}$ in \eqref{mathcal L e} has $2\pi$-periodic coefficients, Bloch-Floquet theory   implies that $\lambda\in\mathbb{C}$ belongs to the $L^2(\mathbb{R})-$spectrum of $\mathcal{L}_{\e,\tts}$ if and only if there exists a nontrivial Floquet–Bloch mode $\tilde{h}(x)=e^{\im \mu x} v(x)$, where $v$ is $2\pi-$periodic and $\mu\in[-\frac{1}{2},\frac{1}{2})$, such that $\lambda \tilde{h}=\mathcal{L}_{\e,\tts} \tilde{h}$, or equivalently, $\lambda v=e^{-\im \mu x}\mathcal{L}_{\e,\tts} e^{\im \mu x}v$. Therefore, we obtain
\begin{align*}
    \sigma_{L^2(\mathbb{R})}(\mathcal{L}_{\e,\tts})=\bigcup_{\mu\in[-\frac{1}{2},\frac{1}{2})} \sigma_{L^2(\mathbb{T})}(\mathcal{L}_{\mu,\e,\tts})~~\mbox{where}~~~\mathcal{L}_{\mu,\e,\tts}:=e^{-\im\mu x}\mathcal{L}_{\e,\tts} e^{\im\mu x}.
\end{align*}
In particular, if $\lambda$ is an eigenvalue of $\mathcal{L}_{\mu,\e,\tts}$ on $L^2(\mathbb{T},\mathbb{C}^2)$ with eigenvector $v(x)$, then $h(t,x)=e^{\lambda t}e^{\im \mu x}v(x)$ solves $h_t=\mathcal{L}_{\e,\tts} h$. We remark that: \\
$(i)$ if $A=Op(a)$ is a pseudo-differential operator with symbol $a(x,\xi)$, which is $2\pi$-periodic in $x$, then $A_\mu:=e^{-\im\mu x}A e^{\im\mu x}=Op(a(x,\xi+\mu))$. \\
$(ii)$ If $A$ is a real operator then $\overline{A_\mu}=A_{-\mu}$. As a consequence the spectrum $\sigma(A_{-\mu})=\overline{\sigma(A_\mu)}$ and we can study $\sigma(A_\mu)$ just for $\mu\geq 0$.\\
$(iii)$  $\sigma(A_\mu)$ is a $1$-periodic set with respect to $\mu$, so one can restrict to $\mu\in[0,\frac{1}{2}]$.

By the previous remarks the Floquet operator associated with the real operator $\mathcal{L}_{\e,\tts}$ in \eqref{mathcal L e} is the complex Hamiltonian and reversible operator
\begin{equation} \label{mathcal L mu e}
\begin{aligned}
    \mathcal{L}_{\mu,\e,\tts}:=&\begin{bmatrix}
(\partial_x+\im\mu)\circ(\ck+p_{\e,\tts}(x)) &  |D+\mu| \\
        -(1+a_{\e,\tts}(x))+\kappa \Sigma_{\mu,\e,\tts} & ~~(\ck+p_{\e,\tts}(x))(\partial_x+\im\mu)
    \end{bmatrix}\\
    =&\underbrace{\begin{bmatrix}
        0& \mathrm{Id}\\
        -\mathrm{Id}& 0
    \end{bmatrix}}_{=:\mathcal{J}}\underbrace{\begin{bmatrix}
 (1+a_{\e,\tts}(x))-\kappa \Sigma_{\mu,\e,\tts}& ~~-(\ck+p_{\e,\tts}(x))(\partial_x+\im\mu)\\
 (\partial_x+\im\mu)\circ(\ck+p_{\e,\tts}(x))& ~~|D+\mu| 
    \end{bmatrix}}_{=:\mathfrak{B}_{\mu,\e,\tts}},
\end{aligned}    
\end{equation}
where $\Sigma_{\mu,\e,\tts}:=e^{-\im\mu x}\Sigma_{\e,\tts} e^{\im\mu x}$ is given by the selfadjoint operator
\begin{align}\label{Sigma}
    \Sigma_{\mu,\e,\tts}=\frac{1}{1+\mathfrak{p}_x}(\partial_x + \im \mu) \circ g_{\e,\tts}(x)\circ (\partial_x+\im\mu) \circ \frac{1}{1+\mathfrak{p}_x} \ . 
\end{align}

We regard $\mathcal{L}_{\mu,\e,\tts}$ as an operator with domain $Y:=H^2(\mathbb{T},\mathbb{C})\times H^1(\mathbb{T},\mathbb{C})$ and codomain $X:=L^2(\mathbb{T},\mathbb{C})\times L^2(\mathbb{T},\mathbb{C})$, equipped with the complex scalar product
\begin{align} \label{complex product}
    (f,g):=\frac{1}{2\pi}\int_0^{2\pi} \left(f_1\overline{g_1}+f_2\overline{g_2} \right)\,\mathrm{d}x,~~\forall~f=\vet{f_1}{f_2},~~g=\vet{g_1}{g_2}\in L^2(\mathbb{T},\mathbb{C}^2).
\end{align}
We also denote $\|f\|^2=(f,f)$.

The complex operator $\mathcal{L}_{\mu,\e,\tts}$ in \eqref{mathcal L mu e} is complex Hamiltonian and reversible. Recall that if $\mathcal{L} : Y \to X$ is a complex linear operator, we say that it is 
\begin{itemize}
    \item \textbf{Complex Hamiltonian:} if there exists a self-adjoint operator, namely $\mathcal{B}=\mathcal{B}^*$, where $\mathcal{B}^*$ (with domain $Y$) is the adjoint with respect to the complex scalar product \eqref{complex product} such that $\mathcal{L} = \mathcal{J} \mathcal{B}$.
    \item \textbf{Reversible:} if 
    \begin{align} \label{reversible}
    \mathcal{L} \circ \overline{\rho} = -\overline{\rho} \circ \mathcal{L}, \quad \text{where} \quad 
    \overline{\rho} \begin{bmatrix} \eta(x) \\ \psi(x) \end{bmatrix} := 
    \begin{bmatrix} \overline{\eta}(-x) \\ -\overline{\psi}(-x) \end{bmatrix} \ . 
    \end{align}
\end{itemize}
The property \eqref{reversible} for $\mathcal{L}_{\mu,\e,\tts}$ follows because $\mathcal{L}_\e$ is a real operator which is reversible with respect to the involution $\rho$ in \eqref{rho involution}. Equivalently, since $\mathcal{J}\circ \overline{\rho}=-\overline{\rho}\circ \mathcal{J}$, the self-adjoint operator $\mathfrak{B}_{\mu,\e,\tts}$ is reversibility-preserving, $\mathrm{i.e.}$
\begin{align}\label{B rho=rho B}
    \mathfrak{B}_{\mu,\e,\tts}\circ\overline{\rho}=\overline{\rho}\circ \mathfrak{B}_{\mu,\e,\tts}.
\end{align}
In addition $(\mu,\e)\rightarrow \mathcal{L}_{\mu,\e,\tts}\in\mathcal{L}( Y,X)$ is analytic, since the functions $\e \mapsto a_{\e,\tts}, ~p_{\e,\tts}$ and $g_{\e,\tts}$ defined in \eqref{SN1}, \eqref{SN1 g} are analytic and $\mathcal{L}_{\mu,\e,\tts}$ is analytic with respect to $\mu$.
Recall also the identity \cite[Section 5.1]{NS}
\begin{align} \label{D+mu}
    |D+\mu|=|D|+\mu(\mathrm{sgn}(D)+\Pi_0),~~{\forall \mu\in[0,\frac{1}{2})},
\end{align}
where $\mathrm{sgn}(D)$ is the Fourier multiplier operator, acting on $2\pi$-periodic functions, with symbol
\begin{equation} \label{sgn D}
\mathrm{sgn}(k):=1 \ \ \forall k >0 \ , \ \
\mathrm{sgn}(0):=0  , \ \ \ 
\mathrm{sgn}(k):= -1 \ \ \forall k < 0 \  , 
\end{equation}
and $\Pi_0$ is the projector operator on the zero mode, $\Pi_0 f(x):=\frac{1}{2\pi}\int_\mathbb{T} f(x) dx$.

\begin{itemize}
\item
{\bf Notation:} We denote by $\mathrm{even}(x)$ a real $2\pi$-periodic function which is even in $x$, and by $\mathrm{odd}(x)$ a real $2\pi$-periodic function which is odd in $x$.
\end{itemize}

Our goal is to prove the existence of eigenvalues of $\mathcal{L}_{\mu,\e,\tts}$ in \eqref{mathcal L mu e} with non zero real part. We remark that the Hamiltonian structure of $\mathcal{L}_{\mu,\e,\tts}$ implies that eigenvalues with non zero real part may arise only from multiple eigenvalues of $\mathcal{L}_{\mu,0,\tts}$ (``Krein criterion''), because if $\lambda$ is an eigenvalue of $\mathcal{L}_{\mu,\e,\tts}$ then also $-\overline{\lambda}$ is, and the total algebraic multiplicity of the eigenvalues is conserved under small perturbation. We now describe the spectrum of $\mathcal{L}_{\mu,0,\tts}$: since the operator $\mathcal{L}_{\mu,0,\tts}$ does not depend on $\tts$, we denote it by $\mathcal{L}_{\mu,0}$ for simplicity.

\paragraph{The spectrum of $\mathcal{L}_{\mu,0}$.} The spectrum of the Fourier multiplier matrix operator
\begin{equation*} 
\begin{aligned}
\mathcal{L}_{\mu,0}:=\mathcal{L}_{\mu,0,\tts}=&\begin{bmatrix}
\ck(\partial_x+\im\mu)&  |D+\mu| \\
        -1+\kappa (\pa_x + \im \mu)^2 & ~~\ck(\partial_x+\im\mu)
    \end{bmatrix}
\end{aligned}
\end{equation*}
consists of the purely imaginary eigenvalues $\{\lambda^{\pm}_k(\mu),~k\in\mathbb{Z}\}$, where
\begin{align} \label{eigenvalues of Lmu0}
    \lambda^{\pm}_k(\mu):=\im\left(\ck(\pm k+\mu)\mp \sqrt{\left(1+\kappa(k\pm\mu)^2\right)\, |k\pm \mu|}\right) \ .
\end{align}
For $\kappa = \frac{1}{2}$, at $\mu=0$ the real operator $\mathcal{L}_{0,0}$ possesses the eigenvalue $0$ with algebraic multiplicity $6$,
\begin{align*}
    \lambda^+_0(0)=\lambda^-_0(0)=\lambda^+_1(0)=\lambda^-_1(0)=\lambda^+_2(0)=\lambda^-_2(0)=0,
\end{align*}
and geometric multiplicity $5$. A real basis of the Kernel of $\mathcal{L}_{0,0}$ is
\small
\begin{align}
    f^-_0:=\vet{0}{1}, ~~ f^+_1 &:=\vet{ \ck^{-1/2} \cos(x)}{\ck^{1/2} \sin(x)}, ~~f^-_1:=\vet{- \ck^{-1/2} \sin(x)}{ \ck^{1/2} \cos(x)}, \label{eigenfunc of mathcall L00 -part1}\\
    f^+_2 &:=\vet{\ck^{-1/2} \cos(2 x)}{\ck^{1/2} \sin(2 x)}, ~~ f^-_2 :=\vet{- \ck^{-1/2} \sin(2 x)}{ \ck^{1/2} \cos(2 x)}, \label{eigenfunc of mathcall L00 - part2}
\end{align}
\normalsize
together with the generalized eigenvector
\begin{align} \label{eigenfunc f+0}
    f^+_0:=\vet{1}{0}, ~~\mathcal{L}_{0,0} f^+_0=-f^-_0.
\end{align}
Notice that the above basis of $\mathrm{Ker}(\cL_{0,0}^{2})$ is symplectic, namely it satisfies
\begin{align} \label{eq:sympKerL00}
    \left( \cJ f^{-}_{k} , f^{+}_{k} \right) &= 1 \ , \quad
    \left( \cJ f^{\sigma}_k , f^{\sigma}_{k} \right)  = 0 \ ,
    \quad \mbox{ if } k \neq k' \mbox{ then } \left( \cJ f^{\sigma}_k , f^{\sigma'}_{k'} \right)  = 0
    \quad\forall \sigma,\sigma^{\prime} = \pm .
\end{align}

Furthermore, $0$ is an isolated eigenvalue for $\mathcal{L}_{0,0}$, namely the spectrum $\sigma(\mathcal{L}_{0,0})$ decomposes into two separated parts,
\begin{align} \label{sigma decomposition}
    \sigma(\mathcal{L}_{0,0})=\sigma'(\mathcal{L}_{0,0})\cup\sigma''(\mathcal{L}_{0,0}), ~~\mbox{where}~~\sigma'(\mathcal{L}_{0,0}):=\{0\},
\end{align}
and $\sigma''(\mathcal{L}_{0,0}):=\{\lambda^{\sigma}_k(0),~k\neq 0,1,2, ~\sigma=\pm\}$.
Notice that from \eqref{eigenvalues of Lmu0} we have that as $\mu \to 0^{+}$
\begin{align*}
\lambda^{\sigma}_0(\mu) = - \sigma \mathrm{i} \mu^{1/2} + \mathcal{O}(\mu) , \quad 
\lambda^{\sigma}_1(\mu) = \mathrm{i} \frac{ 1 }{2 \sqrt{6}} \mu - \mathrm{i} \frac{11}{24 \sqrt{6}} \sigma \mu^2 + \mathcal{O}(\mu^3) , \quad 
\lambda^{\sigma}_2(\mu) = - \mathrm{i} \frac{1}{2\sqrt{6}} \mu - \mathrm{i} \frac{23}{48 \sqrt{6}} \sigma \mu^2 + \mathcal{O}(\mu^3) .
\end{align*}

By Kato's perturbation theory, for any $\mu,\e \neq 0$ sufficiently small, the perturbed spectrum $\sigma(\cL_{\mu,\e})$ admits a disjoint decomposition
\begin{align*} 
    \sigma(\mathcal{L}_{\mu,\e,\tts})=\sigma'(\mathcal{L}_{\mu,\e,\tts})\cup \sigma''(\mathcal{L}_{\mu,\e,\tts}),
\end{align*}
where $\sigma'(\mathcal{L}_{\mu,\e,\tts})$ consists of $6$ eigenvalues close to $0$. We denote by $\mathcal{V}_{\mu,\e}$ the spectral subspace associated with $\sigma'(\mathcal{L}_{\mu,\e,\tts})$, which has dimension $6$ and it is invariant by $\mathcal{L}_{\mu,\e,\tts}$. Our goal is to prove that, for $\e$ small, for values of the Floquet exponent $\mu$ close to $9\e^2/8$, the $6 \times 6$ matrix which represents the operator $\mathcal{L}_{\mu,\e,\tts}:\mathcal{V}_{\mu,\e,\tts}\rightarrow \mathcal{V}_{\mu,\e,\tts}$ possesses a pair of eigenvalues close to zero with opposite non zero real parts.

Now let us introduce a notation that we shall use throughout the paper.

\begin{itemize}
\item
\textbf{Notation:} we denote by $\mathcal{O}(\mu^{m_1} \epsilon^{n_1}, \dots, \mu^{m_p} \epsilon^{n_p})$, $m_j, n_j \in \mathbb{N}_{0}$ (for us $\mathbb{N}_{0} := \{0,1,2, \dots \}$), analytic functions of $(\mu, \epsilon)$ with values in a Banach space $X$ which satisfy, for some $C > 0$, the bound
\[
\|\mathcal{O}(\mu^{m_1}\epsilon^{n_1},\ldots,\mu^{m_p}\epsilon^{n_p})\|_X \leq C \sum_{j=1}^{p} |\mu|^{m_j} |\epsilon|^{n_j}
\]
for small values of $(\mu,\epsilon)$. Similarly, $r_k(\cdots)$ denotes a real-valued analytic scalar remainder with the same bound.
\end{itemize}

Our complete spectral result is the following one:
\begin{theorem}\label{thm:main.mat}
There exists $\e_0 >0 $ such that for any $(\mu, \e) \in \cW$, cf. \eqref{def:W}, the operator $\cL_{\mu, \e, 1}\colon \cV_{\mu, \e} \to \cV_{\mu, \e}$ can be represented by a $6\times 6$ Hamiltonian and reversible matrix of the form 
\begin{equation}\label{U}
\begin{pmatrix}
 \mathsf{U}_E&0&0\\
 0& \mathsf{U}_H &0\\
 0&0&\mathsf{S}
 \end{pmatrix},
\end{equation}
where  $\mathsf{U}_E$, $ \mathsf{U}_H$, $\mathsf{S}$ are $2 \times 2$ matrices, with identical diagonal entries each.
Such matrices are real analytic in the variables $(\tau, \ve)$ parametrizing $\cW$ as 
\begin{equation*} 
\e = \varepsilon^2 \ , \quad \mu = (\tau \varepsilon)^4  \ , \quad (\tau, \ve) \in \left[\left(\frac98-\delta\right)^{\frac14},\, \left(\frac98+\delta\right)^{\frac14}\right]\times  (0 , \e_0^{1/2} ) 
\end{equation*}
and they have the asymptotic expansion, for $X \in \{E, H\} $
\small
\begin{equation}\label{UX}
\begin{aligned}
\mathsf{U}_X
&=
\begin{pmatrix}
\im b_X(\tau, \ve)
&
a_X(\tau, \ve) + d_X(\tau, \ve)  \\
- a_X(\tau, \ve) + d_X(\tau, \ve)
&
\im b_X(\tau, \ve)
\end{pmatrix}, \qquad 
\mathsf{S}
=
\begin{pmatrix}
\im r(\ve^4)
&
\dfrac{4\ck}{81}\tau^8\ve^6+r(\ve^8)\\[1ex]
-\dfrac{\ck}{4}\ve^2 + r(\ve^4) & \im r(\ve^4)
\end{pmatrix}
\end{aligned}
\end{equation}
\normalsize
where \(a_X,d_X,b_X\) are real-valued real-analytic function expanding as
\begin{equation}\label{functionsB2-theorem}
\begin{aligned}
& a_X(\tau, \ve)=\varepsilon^2\big[ 
\alpha(\tau) + \gamma_X(\tau)\ve^2 + r(\ve^4)\big],
\qquad
d_X(\tau, \ve)=\varepsilon^7\bigl(-\delta_X(\tau)+r(\varepsilon)\bigr),  \quad X \in \{E, H\}\\
& b_E(\tau, \ve)=\beta(\tau) \varepsilon^2+\big( \ck \tau^4 -\gamma_E(\tau)\big) \varepsilon^4 +r(\varepsilon^6),\qquad 
b_H(\tau,\ve)=-\beta(\tau) \varepsilon^2+ \big( \ck \tau^4+\gamma_H(\tau)\big)\varepsilon^4+r(\varepsilon^6) 
\end{aligned}
\end{equation}
and
\begin{equation}\label{final.coeff}
    \begin{aligned}
&        \alpha(\tau) = \frac12 \left(\tau^2 - \sqrt{\frac98}\right) \ , \quad \beta(\tau) = \frac12\left(\tau^2 +\sqrt{\frac98}\right) \ , \\ 
\gamma_E(\tau)
&:=\frac{51\sqrt2}{128}+\frac{19\ck}{36}\tau^4,
\qquad
\gamma_H(\tau)
:=\frac{51\sqrt2}{128}-\frac{19\ck}{36}\tau^4,
\\[1ex]
\delta_E(\tau)
&:=\frac{2^{1/4}\tau}{576}
\Big[
\ck\big(-45\sqrt2+24\tau^2-24\sqrt2\tau^4-32\tau^6\big)
+27\big(3-2\sqrt2\tau^2\big)
\Big],
\\[1ex]
\delta_H(\tau)
&:=\frac{2^{1/4}\tau}{576}
\Big[
\ck\big(-45\sqrt2+24\tau^2-24\sqrt2\tau^4-32\tau^6\big)
-27\big(3-2\sqrt2\tau^2\big)
\Big].
\end{aligned}
\end{equation}
\end{theorem}

\section{Perturbative Approach to the Separated Eigenvalues} \label{sec:Kato}

In this section, we analyze the splitting of the eigenvalues of $\mathcal{L}_{\mu, \e,\tts}$ close to 0 for small values of $\mu$ and $\epsilon$, using Kato’s similarity transformation theory \cite[I-§4-6, II-§4]{Kato1966} and \cite{BMV1, BMV3, BMV_ed}. To this end, it is convenient to rewrite the operator $\mathcal{L}_{\mu,\e,\tts}$ in \eqref{mathcal L mu e} as
\begin{equation} \label{mathcal L= ichmu+mathscr L}
\mathcal{L}_{\mu,\e,\tts} = \im \ck \mu + \mathscr{L}_{\mu,\e,\tts}, \quad \mu > 0,
\end{equation}
where, using also \eqref{D+mu}, $\mathscr{L}_{\mu,\e,\tts}$ is the Hamiltonian operator
\begin{equation} \label{mathscr L mu e}
\mathscr{L}_{\mu,\e,\tts} = \mathcal{J} \mathcal{B}_{\mu,\e,\tts},
\end{equation}
with $\mathcal{B}_{\mu,\e,\tts}$ the self-adjoint operator
\begin{equation} \label{mathcal B mu e}
\mathcal{B}_{\mu,\e,\tts} \coloneqq 
\begin{bmatrix} 
1+ a_{\e,\tts}(x)-\kappa\Sigma_{\mu,\e,\tts} & -(\ck + p_{\e,\tts}(x)) \partial_x - \im \mu p_{\e,\tts}(x) \\ \partial_x \circ(\ck + p_{\e,\tts}(x)) + \im \mu p_{\e,\tts}(x) & |D + \mu| 
\end{bmatrix} \ ,
\quad \Sigma_{\mu,\e,\tts} \mbox{ in } \eqref{Sigma} \ . 
\end{equation}
In addition  $\mathscr{L}_{\mu,\e,\tts}$   is also complex-reversible, namely it satisfies, by \eqref{reversible},
\begin{equation} \label{mathscr rho=-rho mathscr}
\mathscr{L}_{\mu,\e,\tts} \circ \bar{\rho} = - \bar{\rho} \circ \mathscr{L}_{\mu,\e,\tts},
\end{equation}
whereas $\mathcal{B}_{\mu,\e,\tts}$ is reversibility-preserving, i.e. fulfills \eqref{B rho=rho B}. Note also that $\mathcal{B}_{0,\e,\tts}$ is a real operator.

The scalar operator $\im \ck \mu \equiv \im \ck \mu \,\text{Id}$ just translates the spectrum of $\mathcal{L}_{\mu,\e,\tts}$ along the imaginary axis of the quantity $\im \ck \mu$, that is, in view of \eqref{mathcal L= ichmu+mathscr L},
\begin{equation*}
\sigma(\mathcal{L}_{\mu,\e,\tts}) = \im \ck \mu + \sigma(\mathscr{L}_{\mu,\e,\tts}).
\end{equation*}
Thus in the sequel we focus on studying the spectrum of $\mathscr{L}_{\mu,\e,\tts}$.

Note also that $\mathscr{L}_{0,\e,\tts} = \mathcal{L}_{0,\e,\tts}$ for any $\epsilon \geq 0$. In particular $\mathscr{L}_{0,0}$ has zero as an isolated eigenvalue with algebraic multiplicity $6$ and generalized kernel spanned by the vectors $\{f_2^{+}, f_2^{-}, f_1^{+}, f_1^{-}, f_0^{+}, f_0^{-}\}$ in \eqref{eigenfunc of mathcall L00 -part1}, \eqref{eigenfunc of mathcall L00 - part2} and \eqref{eigenfunc f+0}; furthermore, its spectrum is separated as in \eqref{sigma decomposition}.

We remark that, in view of \eqref{D+mu}, the operator $\mathscr{L}_{\mu,\e,\tts}$ is analytic with respect to $\mu$. Moreover, the operator $\mathscr{L}_{\mu,\e,\tts}: Y \subset X \to X$ has domain $Y := H^2(\mathbb{T},\mathbb{C})\times H^1(\mathbb{T},\mathbb{C})$ and codomain $X := L^2(\mathbb{T},\mathbb{C})\times L^2(\mathbb{T},\mathbb{C})$.

\begin{lemma} \label{kato thm}
Let $\tts \in \{ \pm 1 \}$. Let $\Gamma$ be a simple, closed, counterclockwise-oriented curve around $0$ in the complex plane separating $\sigma' (\mathscr{L}_{0,0}) = \{0\}$ and the other part of the spectrum $\sigma'' (\mathscr{L}_{0,0})$ in \eqref{sigma decomposition}. There exist $\mu_0,\,\epsilon_0,  > 0$ such that for any $(\mu, \epsilon) \in B(\mu_0) \times B(\epsilon_0)$ the following statements hold:

\begin{enumerate}
    \item \textit{The curve $\Gamma$ belongs to the resolvent set of the operator $\mathscr{L}_{\mu,\e,\tts} : Y \subset X \to X$ defined in \eqref{mathscr L mu e}.}
    \item \textit{The operators}
    \begin{equation} \label{Projection P mu e}
        P_{\mu,\e,\tts} := - \frac{1}{2\pi \im} \oint_\Gamma (\mathscr{L}_{\mu,\e,\tts} - \lambda)^{-1} \ \mathrm{d}\lambda : X \to Y
    \end{equation}
    \textit{are well-defined projectors commuting with $\mathscr{L}_{\mu,\e,\tts}$, i.e., $P_{\mu,\e,\tts}^2 = P_{\mu,\e,\tts}$ and $P_{\mu,\e,\tts} \mathscr{L}_{\mu,\e,\tts} = \mathscr{L}_{\mu,\e,\tts} P_{\mu,\e,\tts}$. The map $(\mu, \epsilon) \mapsto P_{\mu,\e,\tts}$ is analytic from $B(\mu_0) \times B(\epsilon_0)$ to $\mathcal{L}(X,Y)$.}
    \item \textit{The domain $Y$ of the operator $\mathscr{L}_{\mu,\e,\tts}$ decomposes as the direct sum}
    \begin{equation*} 
        Y = \mathcal{V}_{\mu,\e,\tts} \oplus \bigl(Y\cap{}\ker(P_{\mu,\e,\tts})\bigr), \quad \mathcal{V}_{\mu,\e,\tts} \coloneqq \mathrm{Rg}(P_{\mu,\e,\tts}) = \ker(\mathrm{Id} - P_{\mu,\e,\tts}),
    \end{equation*}
    of closed invariant subspaces, namely $\mathscr{L}_{\mu,\e,\tts} : \mathcal{V}_{\mu,\e,\tts} \to \mathcal{V}_{\mu,\e,\tts}$, $\mathscr{L}_{\mu,\e,\tts} : Y\cap{}\ker(P_{\mu,\e,\tts}) \to \ker(P_{\mu,\e,\tts})$. Moreover
    \begin{equation} \label{spectrum separated by Gamma}
    \begin{aligned}
        \sigma(\mathscr{L}_{\mu,\e,\tts}) \cap \{z \in \mathbb{C} \text{ inside } \Gamma\} &= \sigma(\mathscr{L}_{\mu,\e,\tts} |_{\mathcal{V}_{\mu,\e,\tts}}) = \sigma'(\mathscr{L}_{\mu,\e,\tts}), \\
        \sigma(\mathscr{L}_{\mu,\e,\tts}) \cap \{z \in \mathbb{C} \text{ outside } \Gamma\} &= \sigma(\mathscr{L}_{\mu,\e,\tts} |_{\ker(P_{\mu,\e,\tts})}) = \sigma''(\mathscr{L}_{\mu,\e,\tts}).
    \end{aligned}
    \end{equation}
    \item The projectors $P_{\mu,\e,\tts}$ are similar to each other; the transformation operators
    \begin{equation} \label{U transformation operators}
        U_{\mu,\e,\tts} \coloneqq (\mathrm{Id} - (P_{\mu,\e,\tts} - P_{0,0})^2)^{-1/2} \big[P_{\mu,\e,\tts} P_{0,0} + (\mathrm{Id} - P_{\mu,\e,\tts})(\mathrm{Id} - P_{0,0}) \big]
    \end{equation}
    \textit{are bounded and invertible in $Y$ and in $X$, with inverse}
    \begin{equation*} 
        U_{\mu,\e,\tts}^{-1} = \big[P_{0,0} P_{\mu,\e,\tts} + (\mathrm{Id} - P_{0,0})(\mathrm{Id} - P_{\mu,\e,\tts})\big](\mathrm{Id} - (P_{\mu,\e,\tts} - P_{0,0})^2)^{-1/2},
    \end{equation*}
    and $U_{\mu,\e,\tts} P_{0,0} U_{\mu,\e,\tts}^{-1} = P_{\mu,\e,\tts}$ as well as $U_{\mu,\e,\tts}^{-1} P_{\mu,\e,\tts} U_{\mu,\e,\tts} = P_{0,0}$. The map $(\mu, \epsilon) \mapsto U_{\mu,\e,\tts}$ is analytic from $B(\mu_0) \times B(\epsilon_0)$ to $\mathcal{L}(Y)$.
    \item The subspaces $\mathcal{V}_{\mu,\e,\tts} = \mathrm{Rg}(P_{\mu,\e,\tts})$ are isomorphic to each other: $\mathcal{V}_{\mu,\e,\tts} = U_{\mu,\e,\tts} \mathcal{V}_{0,0}$. In particular $\dim \mathcal{V}_{\mu,\e,\tts} = \dim \mathcal{V}_{0,0} = 6$, for any $(\mu, \epsilon) \in B(\mu_0) \times B(\epsilon_0)$.
\end{enumerate}
\end{lemma}

The proof of Lemma \ref{kato thm} is similar to the one of \cite[Lemma 3.1]{BMV1} and we skip it.  Recalling \eqref{mathscr L mu e}-\eqref{mathscr rho=-rho mathscr}, the Hamiltonian and reversible nature of the operator $\mathscr{L}_{\mu,\e}$ imply additional algebraic properties for spectral projectors $P_{\mu,\e}$ and the transformation operators $U_{\mu,\e}$ as follows.  

\begin{lemma} \label{properties of U and P}
For any $(\mu, \epsilon) \in B(\mu_0) \times B(\epsilon_0)$, the following holds true:

\begin{itemize}
    \item[(i)] The projectors $P_{\mu,\e,\tts}$ defined in \eqref{Projection P mu e} are skew-Hamiltonian, namely $\mathcal{J} P_{\mu,\e,\tts} = P_{\mu,\e,\tts}^* \mathcal{J}$, and reversibility preserving, i.e. $\bar{\rho} P_{\mu,\e,\tts} = P_{\mu,\e,\tts} \bar{\rho}$.
    
    \item[(ii)] The transformation operators $U_{\mu,\e,\tts}$ in \eqref{U transformation operators} are symplectic, namely $U_{\mu,\e,\tts}^* \mathcal{J} U_{\mu,\epsilon,\tts} = \mathcal{J}$, and reversibility preserving.
    
    \item[(iii)] $P_{0,\e,\tts}$ and $U_{0,\e,\tts}$ are real operators, i.e. $\bar{P}_{0,\e,\tts} = P_{0,\e,\tts}$ and $\bar{U}_{0,\e,\tts} = U_{0,\e,\tts}$.
\end{itemize}
    
\end{lemma} 

See \cite[Lemma 3.2]{BMV1} for details. By the previous lemma, the antilinear involution $\bar{\rho}$ commutes with the spectral projectors $P_{\mu,\e,\tts}$ and then $\bar{\rho}$ leaves invariant the subspace $\mathcal{V}_{\mu,\e,\tts} = \mathrm{Rg}(P_{\mu,\e,\tts})$.

\paragraph{Symplectic and reversible basis of $\mathcal{V}_{\mu,\e,\tts}$.} It is convenient to represent the Hamiltonian and reversible operator $\mathscr{L}_{\mu,\e,\tts}: \mathcal{V}_{\mu,\e,\tts} \to \mathcal{V}_{\mu,\e,\tts}$ in a basis which is symplectic and reversible, according to the following definition:

\begin{definition}[Symplectic and reversible basis]  \label{Symplectic and reversible basis}
    A basis $\mathsf{F} \coloneqq \{ \mathsf{f}_2^+, \mathsf{f}_2^-, \mathsf{f}_1^+, \mathsf{f}_1^-, \mathsf{f}_0^+, \mathsf{f}_0^- \}$ of $\mathcal{V}_{\mu,\e,\tts}$ is \textit{symplectic} if, for any $k, k' = 0,1,2$,
\begin{equation} \label{basis is symplectic}
    \begin{aligned}
    (\mathcal{J} \mathsf{f}_k^\mp, \mathsf{f}_k^\pm) = \pm 1, ~~(\mathcal{J} \mathsf{f}_k^\sigma, \mathsf{f}_k^\sigma) = 0, \quad \forall \,\sigma = \pm; \\
    \text{if } k \neq k', \text{ then } (\mathcal{J} \mathsf{f}_k^\sigma, \mathsf{f}_{k'}^{\sigma'}) = 0, \quad \forall\, \sigma, \sigma' = \pm.
\end{aligned}
\end{equation}   
It is \textit{reversible} if
\begin{equation} \label{basis is reversible}
    \begin{aligned}
    \bar{\rho} \mathsf{f}_k^\sigma &= \sigma \mathsf{f}_k^\sigma, \quad \forall\, \sigma = \pm, k = 0,1,2.
\end{aligned}
\end{equation}
\end{definition}

\begin{remark}[Parity structure of the reversible basis \eqref{basis is reversible}]
The elements of a reversible basis \( \mathsf{F}=\{\mathsf{f}^+_2,\mathsf{f}^-_2,\mathsf{f}^+_1,\mathsf{f}^-_1,\mathsf{f}^+_0,\mathsf{f}^-_0\} \) enjoys specific parity properties. Specifically,
\begin{align} \label{Parity structure}
\mathsf{f}_k^+(x) = 
\begin{bmatrix}
\mathit{even}(x) + \im\,\mathit{odd}(x) \\
\mathit{odd}(x) + \im\,\mathit{even}(x)
\end{bmatrix}, \quad
\mathsf{f}_k^-(x) = 
\begin{bmatrix}
\mathit{odd}(x) + \im\,\mathit{even}(x) \\
\mathit{even}(x) + \im\,\mathit{odd}(x)
\end{bmatrix}.
\end{align}
This structure follows from the reversibility of the problem, specifically from the involution \( \overline{\rho} \) defined in equation \eqref{reversible}, which implies that the real and imaginary parts of each component satisfy definite parity conditions.
\end{remark}

\begin{remark}[Symplectic expansion using the basis in \eqref{basis is symplectic}]
We can express any vector \( \mathsf{f} \in \mathcal{V}_{\mu,\epsilon} \) as a linear combination of the symplectic basis:
\begin{align} \label{remark f expansion}
    \mathsf{f} = \alpha_2^+ \mathsf{f}_2^+ + \alpha_2^- \mathsf{f}_2^- + \alpha_1^+ \mathsf{f}_1^+ + \alpha_1^- \mathsf{f}_1^- + \alpha_0^+ \mathsf{f}_0^+ + \alpha_0^- \mathsf{f}_0^-,
\end{align}
for suitable coefficients \( \alpha_k^\sigma \in \mathbb{C} \). These coefficients are computed by applying the symplectic form \( \mathcal{J} \), taking \( L^2 \)-scalar products with the basis elements, and using the symplecticity \eqref{basis is symplectic}. Therefore, we may rewrite \eqref{remark f expansion} as
\begin{align} \label{expasion of f by using symplectic basis}
    \mathsf{f} &= -(\mathcal{J}\mathsf{f},\mathsf{f}^-_2)\mathsf{f}^+_2 + (\mathcal{J}\mathsf{f},\mathsf{f}^+_2)\mathsf{f}^-_2 - (\mathcal{J}\mathsf{f},\mathsf{f}^-_1)\mathsf{f}^+_1 + (\mathcal{J}\mathsf{f},\mathsf{f}^+_1)\mathsf{f}^-_1 - (\mathcal{J}\mathsf{f},\mathsf{f}^-_0)\mathsf{f}^+_0 + (\mathcal{J}\mathsf{f},\mathsf{f}^+_0)\mathsf{f}^-_0.
\end{align}
\end{remark}

We now compute the matrix representation of $\mathscr{L}_{\mu,\e}$ in a symplectic and reversible basis of $\mathcal{V}_{\mu,\e}$.

\begin{lemma} \label{matrix representation mathsf L}
The $6 \times 6$ matrix that represents the Hamiltonian and reversible operator $\mathscr{L}_{\mu,\e,\tts} = \mathcal{J} \mathcal{B}_{\mu,\e,\tts} : \mathcal{V}_{\mu,\e,\tts} \to \mathcal{V}_{\mu,\e,\tts}$ with respect to a symplectic and reversible basis $\mathsf{F} = \{ \mathsf{f}_2^{+}, \mathsf{f}_2^{-}, \mathsf{f}_1^{+}, \mathsf{f}_1^{-}, \mathsf{f}_0^{+}, \mathsf{f}_0^{-} \}$ of $\mathcal{V}_{\mu,\e,\tts}$ is
\begin{align*} 
  \mathsf{L}_{\mu,\e,\tts} &\coloneqq  \mathsf{J}_6 \mathsf{B}_{\mu,\e,\tts}, \quad \mathsf{J}_6 \coloneqq 
\begin{pmatrix}
\mathsf{J}_2 & 0 & 0\\
0 & \mathsf{J}_2 & 0\\
0 & 0 & \mathsf{J}_2 
\end{pmatrix}
, \quad
    \mathsf{J}_2 := 
    \begin{pmatrix}
        0 & 1 \\
        -1 & 0
    \end{pmatrix}, \quad \text{where} \quad \mathsf{B}_{\mu,\e,\tts} = \mathsf{B}_{\mu,\e,\tts}^{*}.
\end{align*}
The self-adjoint matrix
\begin{align} \label{Bmuepsilon matrix representation}
    \mathsf{B}_{\mu,\e,\tts} = 
    \begin{pmatrix}
        (\mathcal{B}_{\mu,\e,\tts} \mathsf{f}_2^{+}, \mathsf{f}_2^{+}) & (\mathcal{B}_{\mu,\e,\tts} \mathsf{f}_2^{-}, \mathsf{f}_2^{+}) & (\mathcal{B}_{\mu,\e,\tts} \mathsf{f}_1^{+}, \mathsf{f}_2^{+}) & (\mathcal{B}_{\mu,\e,\tts} \mathsf{f}_1^{-}, \mathsf{f}_2^{+}) & (\mathcal{B}_{\mu,\e,\tts} \mathsf{f}_0^{+}, \mathsf{f}_2^{+}) & (\mathcal{B}_{\mu,\e,\tts} \mathsf{f}_0^{-}, \mathsf{f}_2^{+}) \\
        (\mathcal{B}_{\mu,\e,\tts} \mathsf{f}_2^{+}, \mathsf{f}_2^{-}) & (\mathcal{B}_{\mu,\e,\tts} \mathsf{f}_2^{-}, \mathsf{f}_2^{-}) & (\mathcal{B}_{\mu,\e,\tts} \mathsf{f}_1^{+}, \mathsf{f}_2^{-}) & (\mathcal{B}_{\mu,\e,\tts} \mathsf{f}_1^{-}, \mathsf{f}_2^{-}) & (\mathcal{B}_{\mu,\e,\tts} \mathsf{f}_0^{+}, \mathsf{f}_2^{-}) & (\mathcal{B}_{\mu,\e,\tts} \mathsf{f}_0^{-}, \mathsf{f}_2^{-}) \\
        (\mathcal{B}_{\mu,\e,\tts} \mathsf{f}_2^{+}, \mathsf{f}_1^{+}) & (\mathcal{B}_{\mu,\e,\tts} \mathsf{f}_2^{-}, \mathsf{f}_1^{+}) & (\mathcal{B}_{\mu,\e,\tts} \mathsf{f}_1^{+}, \mathsf{f}_1^{+}) & (\mathcal{B}_{\mu,\e,\tts} \mathsf{f}_1^{-}, \mathsf{f}_1^{+}) & (\mathcal{B}_{\mu,\e,\tts} \mathsf{f}_0^{+}, \mathsf{f}_1^{+}) & (\mathcal{B}_{\mu,\e,\tts} \mathsf{f}_0^{-}, \mathsf{f}_1^{+}) \\
        (\mathcal{B}_{\mu,\e,\tts} \mathsf{f}_2^{+}, \mathsf{f}_1^{-}) & (\mathcal{B}_{\mu,\e,\tts} \mathsf{f}_2^{-}, \mathsf{f}_1^{-}) & (\mathcal{B}_{\mu,\e,\tts} \mathsf{f}_1^{+}, \mathsf{f}_1^{-}) & (\mathcal{B}_{\mu,\e,\tts} \mathsf{f}_1^{-}, \mathsf{f}_1^{-}) & (\mathcal{B}_{\mu,\e,\tts} \mathsf{f}_0^{+}, \mathsf{f}_1^{-}) & (\mathcal{B}_{\mu,\e,\tts} \mathsf{f}_0^{-}, \mathsf{f}_1^{-}) \\
        (\mathcal{B}_{\mu,\e,\tts} \mathsf{f}_2^{+}, \mathsf{f}_0^{+}) & (\mathcal{B}_{\mu,\e,\tts} \mathsf{f}_2^{-}, \mathsf{f}_0^{+}) & (\mathcal{B}_{\mu,\e,\tts} \mathsf{f}_1^{+}, \mathsf{f}_0^{+}) & (\mathcal{B}_{\mu,\e,\tts} \mathsf{f}_1^{-}, \mathsf{f}_0^{+}) & (\mathcal{B}_{\mu,\e,\tts} \mathsf{f}_0^{+}, \mathsf{f}_0^{+}) & (\mathcal{B}_{\mu,\e,\tts} \mathsf{f}_0^{-}, \mathsf{f}_0^{+}) \\
        (\mathcal{B}_{\mu,\e,\tts} \mathsf{f}_2^{+}, \mathsf{f}_0^{-}) & (\mathcal{B}_{\mu,\e,\tts} \mathsf{f}_2^{-}, \mathsf{f}_0^{-}) & (\mathcal{B}_{\mu,\e,\tts} \mathsf{f}_1^{+}, \mathsf{f}_0^{-}) & (\mathcal{B}_{\mu,\e,\tts} \mathsf{f}_1^{-}, \mathsf{f}_0^{-}) & (\mathcal{B}_{\mu,\e,\tts} \mathsf{f}_0^{+}, \mathsf{f}_0^{-}) & (\mathcal{B}_{\mu,\e,\tts} \mathsf{f}_0^{-}, \mathsf{f}_0^{-})
    \end{pmatrix}
    .
\end{align}
The entries of the matrix $\mathsf{B}_{\mu,\epsilon}$ are alternatively real or purely imaginary: for any $\sigma = \pm, k = 0, 1, 2$,
\small
\begin{align} \label{B are alternatively real or imaginary}
    (\mathcal{B}_{\mu,\epsilon} \mathsf{f}_k^{\sigma}, \mathsf{f}_{k'}^{\sigma}) \text{ is real}, \quad
    (\mathcal{B}_{\mu,\epsilon} \mathsf{f}_k^{\sigma}, \mathsf{f}_{k'}^{-\sigma}) \text{ is purely imaginary}.
\end{align}  
\normalsize
\end{lemma} 
\begin{proof}
    Recalling \eqref{mathscr L mu e} and \eqref{expasion of f by using symplectic basis}, for $\sigma=\pm$, $k=0,1,2$, we obtain
    \begin{align*}
        &\mathscr{L}_{\mu,\epsilon}\mathsf{f}^{\sigma}_k \\
        &=-(\mathcal{J}\mathscr{L}_{\mu,\epsilon}\mathsf{f}^{\sigma}_k,\mathsf{f}^-_2)\mathsf{f}^+_2+(\mathcal{J}\mathscr{L}_{\mu,\epsilon}\mathsf{f}^{\sigma}_k,\mathsf{f}^+_2)\mathsf{f}^-_2-(\mathcal{J}\mathscr{L}_{\mu,\epsilon}\mathsf{f}^{\sigma}_k,\mathsf{f}^-_1)\mathsf{f}^+_1+(\mathcal{J}\mathscr{L}_{\mu,\epsilon}\mathsf{f}^{\sigma}_k,\mathsf{f}^+_1)\mathsf{f}^-_1-(\mathcal{J}\mathscr{L}_{\mu,\epsilon}\mathsf{f}^{\sigma}_k,\mathsf{f}^-_0)\mathsf{f}^+_0+(\mathcal{J}\mathscr{L}_{\mu,\epsilon}\mathsf{f}^{\sigma}_k,\mathsf{f}^+_0)\mathsf{f}^-_0\\
        &=(\mathcal{B}_{\mu,\epsilon}\mathsf{f}^{\sigma}_k,\mathsf{f}^-_2)\mathsf{f}^+_2-(\mathcal{B}_{\mu,\epsilon}\mathsf{f}^{\sigma}_k,\mathsf{f}^+_2)\mathsf{f}^-_2+(\mathcal{B}_{\mu,\epsilon}\mathsf{f}^{\sigma}_k,\mathsf{f}^-_1)\mathsf{f}^+_1-(\mathcal{B}_{\mu,\epsilon}\mathsf{f}^{\sigma}_k,\mathsf{f}^+_1)\mathsf{f}^-_1+(\mathcal{B}_{\mu,\epsilon}\mathsf{f}^{\sigma}_k,\mathsf{f}^-_0)\mathsf{f}^+_0-(\mathcal{B}_{\mu,\epsilon}\mathsf{f}^{\sigma}_k,\mathsf{f}^+_0)\mathsf{f}^-_0.
    \end{align*}
This verifies that the matrix representation of $\mathscr{L}_{\mu,\e,\tts}$ with respect to $\mathsf{F}$ is $\mathsf{J}_6\mathsf{B}_{\mu,\e,\tts}$. Also, the matrix $\mathsf{B}_{\mu,\e,\tts}$ is self-adjoint because $\mathcal{B}_{\mu,\e,\tts}$ is self-adjoint. Next, recall from \eqref{complex product} and \eqref{reversible} that the inner product satisfies
\begin{align*} 
    (\mathsf{f},\mathsf{g})=\overline{(\overline{\rho}\mathsf{f},\overline{\rho}\mathsf{g})}.
\end{align*}
Then, since $\mathcal{B}_{\mu,\e,\tts}$ is both self-adjoint and reversibility-preserving \eqref{B rho=rho B} and \eqref{basis is reversible}, we compute:
\begin{align*}
    (\mathcal{B}_{\mu,\e,\tts} \mathsf{f}_k^{\sigma}, \mathsf{f}_{k'}^{\sigma'})=\overline{(\overline{\rho}\mathcal{B}_{\mu,\e,\tts} \mathsf{f}^{\sigma}_k,\overline{\rho} \mathsf{f}^{\sigma'}_{k'})}=\overline{(\mathcal{B}_{\mu,\e,\tts}\overline{\rho} \mathsf{f}^{\sigma}_k,\overline{\rho} \mathsf{f}^{\sigma'}_{k'})}=\sigma\sigma'\overline{(\mathcal{B}_{\mu,\e,\tts}\mathsf{f}^{\sigma}_k,\mathsf{f}^{\sigma'}_{k'})},
\end{align*}
which proves \eqref{B are alternatively real or imaginary}.
\end{proof}

We conclude this section recalling some notation. A $2n \times 2n$, $n = 1, 2,3$, matrix of the form $\mathsf{L} = \mathsf{J}_{2n} \mathsf{B}$ is \textit{Hamiltonian} if $\mathsf{B}$ is a self-adjoint matrix, i.e. $\mathsf{B} = \mathsf{B}^*$. It is \textit{reversible} if $\mathsf{B}$ is reversibility-preserving, i.e.
\[
\rho_{2n} \circ \mathsf{B} = \mathsf{B} \circ \rho_{2n},
\]
where
\[
\rho_6 := \begin{pmatrix} \rho_2 & 0 & 0 \\ 0 & \rho_2 & 0 \\ 0 & 0 & \rho_2 \end{pmatrix}, \quad
\rho_4 := \begin{pmatrix} \rho_2 & 0 \\ 0 & \rho_2 \end{pmatrix}, \quad
\rho_2 := \begin{pmatrix} \mathfrak{c} & 0 \\ 0 & -\mathfrak{c} \end{pmatrix}, 
\]
and $\mathfrak{c} : z \mapsto \bar{z}$ is the conjugation of the complex plane.
Equivalently, $\rho_{2n} \circ \mathsf{L} = - \mathsf{L} \circ \rho_{2n}.$

The transformations preserving the Hamiltonian structure are called \textit{symplectic}, and satisfy
\begin{align*} 
Y^* \mathsf{J}_6 Y = \mathsf{J}_6.
\end{align*}
If $Y$ is symplectic then $Y^*$ and $Y^{-1}$ are symplectic as well. A Hamiltonian matrix $\mathsf{L} = \mathsf{J}_6 \mathsf{B}$, with $\mathsf{B} = \mathsf{B}^*$, is conjugated through a symplectic matrix $Y$ in a new Hamiltonian matrix
\begin{align*} 
    \mathsf{L}_1=Y^{-1}\mathsf{L}Y=Y^{-1}\mathsf{J}_6(Y^{-1})^*Y^*\mathsf{B}Y=\mathsf{J}_6 \mathsf{B}_1\,\, \ \ \mbox{ where }\,\,\mathsf{B}_1:=Y^*\mathsf{B}Y=\mathsf{B}^*_1.
\end{align*}
A $6 \times 6$ matrix $\mathsf{B} = (\mathsf{B}_{ij})_{i,j=1,\ldots,6}$ is reversibility-preserving if and only if its entries are alternatively real and purely imaginary, namely $\mathsf{B}_{ij}$ is real when $i + j$ is even and purely imaginary otherwise, as in \eqref{B are alternatively real or imaginary}. A $6 \times 6$ complex matrix $\mathsf{L} = (\mathsf{L}_{ij})_{i,j=1,\ldots,6}$ is reversible if and only if $\mathsf{L}_{ij}$ is purely imaginary when $i + j$ is even and real otherwise.

We finally mention that the flow of a Hamiltonian reversibility-preserving matrix is symplectic and reversibility-preserving (see \cite[Lemma 3.8]{BMV1}).

\section{Matrix Representation of $\mathscr{L}_{\mu,\e,\tts}$ on $\mathcal{V}_{\mu,\e,\tts}$} \label{sec:MatrRepr}

Using the transformation operator $U_{\mu,\e,\tts}$ in \eqref{U transformation operators}, we construct the basis of $\mathcal{V}_{\mu,\e,\tts}$
\begin{equation} \label{F basis set and f}
    \begin{aligned}
        \mathcal{F}&:=\{f^+_2(\mu,\e,\tts), 
        \,  f^-_2(\mu,\e,\tts), \, f^+_1(\mu,\e,\tts),
        \, f^-_1(\mu,\e,\tts),\, f^+_0(\mu,\e,\tts),
        \, f^-_0(\mu,\e,\tts)\}, \\ 
        f^\sigma_k(\mu,\e,\tts) &\coloneqq U_{\mu,\e,\tts} f^{\sigma}_k,~~\sigma=\pm,~k=0,1,2,
    \end{aligned}
\end{equation}
where 
\begin{equation} \label{eigenfunc of mathcall L00 2}
    f^+_j=\vet{\ck^{-1/2}\cos(j x)}{\ck^{1/2}\sin(j x)},\, \, 
    f^-_j=\vet{-\ck^{-1/2}\sin(j x)}{\ck^{1/2}\cos(j x)}, \quad j=1,2, \quad f^+_0=\vet{1}{0},~~f^-_0=\vet{0}{1},
\end{equation}
form a basis of $\mathcal{V}_{0,0}=\mathrm{Rg}(P_{0,0})$, see \eqref{eigenfunc of mathcall L00 -part1}, \eqref{eigenfunc of mathcall L00 - part2} and \eqref{eigenfunc f+0}. Note that the real valued vectors $\{f^\pm_2,f^\pm_1,f^\pm_0\}$ form a symplectic and reversible basis for $\mathcal{V}_{0,0}$, according to Definition \ref{Symplectic and reversible basis}. Then, by Lemma \ref{kato thm} and Lemma \ref{properties of U and P} we deduce that:
\begin{lemma} \label{F is symplectic and reversible}
    The basis $\mathcal{F}$ of $\mathcal{V}_{\mu,\e,\tts}$ defined in \eqref{F basis set and f} is symplectic and reversible, i.e. satisfies \eqref{basis is symplectic} and \eqref{basis is reversible}. Each map $(\mu,\e) \mapsto f^\sigma_k(\mu,\e,\tts)$ is analytic as a map $B(\mu_0)\times B(\e_0)\rightarrow H^2(\mathbb{T})\times H^1(\mathbb{T})$.
\end{lemma}
\begin{proof}
By Lemma \ref{properties of U and P}-(ii), the operators $U_{\mu,\e,\tts}$ are symplectic and reversibility preserving, i.e. 
$U_{\mu,\e,\tts}^* \mathcal{J} U_{\mu,\e,\tts} = \mathcal{J}$ and $U_{\mu,\e,\tts}\circ\overline{\rho}=\overline{\rho}\circ U_{\mu,\e,\tts}$. Recall \eqref{F basis set and f}, i.e. $f_k^\sigma(\mu,\e,\tts) \coloneqq U_{\mu,\e,\tts} f_k^\sigma, \sigma=\pm, k=0,1,2.$
Then for any $k,k',\sigma,\sigma'$,
$(\mathcal{J} f_k^\sigma(\mu,\e,\tts), f_{k'}^{\sigma'}(\mu,\e,\tts))
= (\mathcal{J} f_k^\sigma, f_{k'}^{\sigma'}),$
so the symplectic relations \eqref{basis is symplectic} are preserved. Moreover,
$\overline{\rho} f_k^\sigma(\mu,\epsilon,\tts)
= U_{\mu,\e,\tts}\circ \overline{\rho} f_k^\sigma
= \sigma f_k^\sigma(\mu,\e,\tts),$
which shows that the reversibility conditions \eqref{basis is reversible} hold as well. 
Finally, the analyticity of $f_k^\sigma(\mu,\e,\tts)$ follows from the analyticity of 
$U_{\mu,\epsilon, \tts}$ (see Lemma \ref{kato thm}). 
\end{proof}

We then expand the vectors $f^\sigma_k(\mu,\e,\tts)$ in $(\mu,\e)$. We denote by $\mathrm{even}_0(x)$ a real, even, $2\pi$-periodic function with zero space average. In the sequel $\cO(\mu^m\e^n)\vet{\mathrm{even}(x)}{\mathrm{odd}(x)}$ denotes an analytic map in $(\mu,\e)$ with values in $Y=H^2(\mathbb{T},\mathbb{C})\times H^1(\mathbb{T},\mathbb{C})$, whose first component is $\mathrm{even}(x)$ and the second one $\mathrm{odd}(x)$; we have a similar meaning for $\cO(\mu^m\e^n)\vet{\mathrm{odd}(x)}{\mathrm{even}(x)}$, etc $\ldots$.

\begin{lemma}[Expansion of the basis $\mathcal{F}$]
\label{expansion of the basis F}
For small values of $(\mu,\e)$, the basis $\mathcal{F}$
in \eqref{F basis set and f} has the expansion
\small
\begin{equation}\label{f+2}
\begin{aligned}
f_2^+(\mu,\e,\tts)
&= f_2^+-\frac{\im\mu}{24}f^-_{-2}
+\e\left[
\frac14 f^+_{-1}-2f^+_3+\frac12 f^+_{-3}
-\frac{\tts}{2}f^+_4+\frac{\tts}{2}f^+_{-4}
\right]\\
&\quad+\mathcal{O}(\mu^2)
\begin{bmatrix}
\mathit{even}_0(x)+\im\mathit{odd}(x)\\
\mathit{odd}(x)+\im\mathit{even}_0(x)
\end{bmatrix}
+\mathcal{O}(\e^2)
\begin{bmatrix}
\mathit{even}_0(x)\\
\mathit{odd}(x)
\end{bmatrix}
+\im\mu\e
\begin{bmatrix}
\mathit{odd}(x)\\
\mathit{even}(x)
\end{bmatrix}
+\mathcal{O}(\mu^2\e,\mu\e^2).
\end{aligned}
\end{equation}

\begin{equation}\label{f-2}
\begin{aligned}
f_2^-(\mu,\e,\tts)
&= f_2^--\frac{\im\mu}{24}f^+_{-2}
+\e\left[
-\frac14 f^-_{-1}-2f^-_3-\frac12 f^-_{-3}
-\frac{\tts}{2}f^-_4-\frac{\tts}{2}f^-_{-4}
\right]\\
&\quad+\mathcal{O}(\mu^2)
\begin{bmatrix}
\mathit{odd}(x)+\im\mathit{even}_0(x)\\
\mathit{even}_0(x)+\im\mathit{odd}(x)
\end{bmatrix}
+\mathcal{O}(\e^2)
\begin{bmatrix}
\mathit{odd}(x)\\
\mathit{even}(x)
\end{bmatrix}+\im\mu\e
\begin{bmatrix}
{\mathit{even}_0}(x)\\
\mathit{odd}(x)
\end{bmatrix}
+\mathcal{O}(\mu^2\e,\mu\e^2).
\end{aligned}
\end{equation}

\begin{equation}\label{f+1}
\begin{aligned}
f_1^+(\mu,\e,\tts)
&= f_1^++\frac{\im\mu}{12}f^-_{-1}
+\e\left[
\frac{3\tts}{8}f^+_{-1}+\frac38 f^+_{-2}
-\frac{5\tts}{4}f^+_3+\frac{\tts}{2}f^+_{-3}
\right]\\
&\quad+\mathcal{O}(\mu^2)
\begin{bmatrix}
\mathit{even}_0(x)+\im\mathit{odd}(x)\\
\mathit{odd}(x)+\im\mathit{even}_0(x)
\end{bmatrix}
+\mathcal{O}(\e^2)
\begin{bmatrix}
\mathit{even}_0(x)\\
\mathit{odd}(x)
\end{bmatrix}+\im\mu\e
\begin{bmatrix}
\mathit{odd}(x)\\
\mathit{even}(x)
\end{bmatrix}
+\mathcal{O}(\mu^2\e,\mu\e^2).
\end{aligned}
\end{equation}

\begin{equation}\label{f-1}
\begin{aligned}
f_1^-(\mu,\e,\tts)
&= f_1^-+\frac{\im\mu}{12}f^+_{-1}
+\e\left[
\frac{3\tts}{8}f^-_{-1}-\frac38 f^-_{-2}
-\frac{5\tts}{4}f^-_3-\frac{\tts}{2}f^-_{-3}
\right]\\
&\quad+\mathcal{O}(\mu^2)
\begin{bmatrix}
\mathit{odd}(x)+\im\mathit{even}_0(x)\\
\mathit{even}_0(x)+\im\mathit{odd}(x)
\end{bmatrix}
+\mathcal{O}(\e^2)
\begin{bmatrix}
\mathit{odd}(x)\\
\mathit{even}(x)
\end{bmatrix}+\im\mu\e
\begin{bmatrix}
{\mathit{even}_0}(x)\\
\mathit{odd}(x)
\end{bmatrix}
+\mathcal{O}(\mu^2\e,\mu\e^2).
\end{aligned}
\end{equation}

\begin{equation}\label{f+0}
\begin{aligned}
f_0^+(\mu,\e,\tts)
&=\begin{bmatrix}1\\0\end{bmatrix}
+\e\ck^{1/2}\left[f^+_{-1}+\tts f^+_{-2}\right]
+\mathcal{O}(\e^2)
\begin{bmatrix}
\mathit{even}_0(x)\\
\mathit{odd}(x)
\end{bmatrix}\\
&\quad+\mu\e\left(
\frac14\ck^{-3/2}f^+_{-1}
+\frac{\tts}{16}\ck^{-3/2}f^+_{-2}
+\im
\begin{bmatrix}
\mathit{odd}(x)\\
\mathit{even}_0(x)
\end{bmatrix}
\right)
+\mathcal{O}(\mu^2\e,\mu\e^2).
\end{aligned}
\end{equation}

\begin{equation}\label{f-0}
\begin{aligned}
f_0^-(\mu,\e,\tts)
&=\begin{bmatrix}0\\1\end{bmatrix}
+\mu\e\left(
\frac12\ck^{-1/2}f^-_{-1}
+\frac{\tts}{4}\ck^{-1/2}f^-_{-2}
+\im
\begin{bmatrix}
\mathit{even}_0(x)\\
\mathit{odd}(x)
\end{bmatrix}
\right)+\mathcal{O}(\mu^2\e,\mu\e^2).
\end{aligned}
\end{equation}
\normalsize
For $\mu=0$, the basis
$\{f_j^\pm(0,\e,\tts):j=0,1,2\}$ is real and
\begin{equation}\label{f-0=01}
\begin{aligned}
f_j^+(0,\e,\tts)
&=\begin{bmatrix}
\mathit{even}_0(x)\\
\mathit{odd}(x)
\end{bmatrix},
\qquad
f_j^-(0,\e,\tts)
=\begin{bmatrix}
\mathit{odd}(x)\\
\mathit{even}(x)
\end{bmatrix},
\quad j=1,2,\\
f_0^+(0,\e,\tts)
&=\vet{1}{0}
+\begin{bmatrix}
\mathit{even}_0(x)\\
\mathit{odd}(x)
\end{bmatrix},
\qquad
f_0^-(0,\e,\tts)
=\begin{bmatrix}0\\1\end{bmatrix}.
\end{aligned}
\end{equation}
\end{lemma}

The proof of Lemma \ref{expansion of the basis F} is deferred to Appendix \ref{secA1}.
We now state the main result of this section.

\begin{proposition} \label{B decomposition}
The matrix that represents the Hamiltonian and reversible operator $\mathscr{L}_{\mu,\e,\tts} : \mathcal{V}_{\mu,\e,\tts} \to \mathcal{V}_{\mu,\e,\tts}$ in the symplectic and reversible basis $ \mathcal{F}$ of $\mathcal{V}_{\mu,\e,\tts}$ defined in  \eqref{F basis set and f}, is a Hamiltonian matrix $\mathsf{L}_{\mu,\e,\tts} = \mathsf{J}_6 \mathsf{B}_{\mu,\e,\tts}$, where $\mathsf{B}_{\mu,\e,\tts}$ is a self-adjoint and reversibility preserving (i.e., satisfying \eqref{B are alternatively real or imaginary}) $6 \times 6$ matrix of the form
\begin{align} \label{B mu epsilon}
    \mathsf{B}_{\mu,\e,\tts} =
\begin{pmatrix}
E & I & F \\
I^* & H & L \\
F^* & L^* & G
\end{pmatrix}, \quad E = E^*, \quad H = H^*, \quad G = G^*,
\end{align}
where $E, H, G $ are the $2 \times 2$ selfadjoint and reversible matrices
\begin{equation} \label{E}
\begin{aligned}
E &:=
\begin{pmatrix}
\mathtt{e}_{1} \e + \mathtt{e}_{2} \e^2  + r(\e^3,  \mu^2)  & \mathrm{i}\mathtt{e}_{4}  \mu     + \im r(\mu \e,  \mu^2)  \\
-\mathrm{i} \mathtt{e}_{4}  \mu    -\im r(\mu \e,  \mu^2) &  \mathtt{e}_{1} \e + \mathtt{e}_{6} \e^2  + r(\e^3, \mu^2) 
\end{pmatrix} \ , 
\qquad 
G :=
\begin{pmatrix}
1+ r(\e^3, \mu \e^2, \mu^2)    & \im r(\mu\e, \mu^2)   \\
 -\im r(\mu\e, \mu^2)   & \mu+r(\mu^2)
\end{pmatrix} \ , 
\end{aligned}
\end{equation}
\begin{equation} \label{H}
\begin{aligned}
H &:=
\begin{pmatrix}
\mathtt{h}_{1} \e^2  +r(\e^3, \mu \e^2, \mu^2)   & \im   \mathtt{h}_{3} \mu  + \im r(\mu \e, \mu^2) \\
- \im \mathtt{h}_{3} \mu   - \im r(\mu \e, \mu^2)  & \mathtt{e}_{1} \e + \mathtt{h}_{5} \e^2 +r(\e^3,  \mu \e^2, \mu^2)
\end{pmatrix} \ , 
\end{aligned}
\end{equation}
whereas $F, I, L$ are the $2\times 2$ reversible matrices
\begin{equation} \label{F}
\begin{aligned}
F &:= 
\begin{pmatrix}
\mathtt{f}_{1}\e^3  + r(\e^4, \mu \e^2, \mu^2)  & \mathrm{i} \mu \e \left(\mathtt{f}_{2}  + \mathtt{f}_{3}\e\right)  +\im r(\mu \e^3, \mu^2\e) \\
 \mathrm{i} \mathtt{f}_{4}  \mu \e +\im r(\mu \e^2, \mu^2) & \mathtt{f}_5\mu\e^2+r(\mu \e^3, \mu^2\e)
\end{pmatrix} ,
\end{aligned}
\end{equation}

\begin{equation} \label{I}
\begin{aligned}
I &:=
\begin{pmatrix}
\mathtt{i}_{1} \e + \mathtt{i}_{2} \e^2 + r(\e^3, \mu^2)  &  \im r(\mu \e, \mu^2) \\
 \im r(\mu \e, \mu^2)   & \mathtt{i}_{1} \e + \mathtt{i}_{4} \e^2 +  r(\e^3, \mu^2)
\end{pmatrix} \ , 
\end{aligned}
\end{equation}

\begin{equation} \label{L}
\begin{aligned}
L &:= 
\begin{pmatrix}
\mathtt{l}_1 \e^3  + r(\e^4, \mu\e^2, \mu^2) &\im\mu \e  \left(\mathtt{l}_{2} + \mathtt{l}_{3}\e\right)  + \im r(\mu \e^3, \mu^2\e)\\
\im  \mathtt{l}_{4} \mu \e  + \im r(\mu \e^2, \mu^2 )&  r(\mu \e^3, \mu^2\e)
\end{pmatrix} \ , 
\end{aligned}
\end{equation}
where 
\begin{equation}\label{coeff.B}
\begin{aligned}
& \mathtt{e}_{1} := - \frac{1}{2} \ck \tts ,\qquad  \mathtt{e}_{2} := - \frac{13}{16} \ck ,
\qquad 
\mathtt{e}_{4}:= \frac{7}{6}\ck \ , 
\qquad \mathtt{e}_{6} := - \frac{5}{16} \ck\,,
\\
& \mathtt{f}_{1}:=  \frac{3\tts}{16} \ck^{3/2} \ , 
\qquad 
\mathtt{f}_{2}:=\tts \ck^{1/2} \ , \qquad 
\mathtt{f}_{3}:=-\frac32 \ck^{1/2}  \ , 
\qquad 
\mathtt{f}_{4}:= - \frac{\tts}{6} \ck^{3/2} \ ,  \qquad 
\mathtt{f}_{5}:=\frac{\ck^{1/2}}8\\
& \mathtt{h}_{1} := -\frac{3}{16} \ck , \qquad
\mathtt{h}_{3}:= \frac{5}{6} \ck\ , \qquad 
\mathtt{h}_{5} := \frac{7}{16} \ck\,, \\
& \mathtt{i}_{1} := \frac{1}{2} \ck ,
\qquad 
\mathtt{i}_{2} := - \frac{35}{16} \tts \ck  ,
\qquad \mathtt{i}_{4} := -\frac{\tts}{16} \ck \ , \\
& \mathtt{l}_{1}:=  -\frac{\ck^{3/2}}{8} \ , 
\qquad \mathtt{l}_{2}:=  \ck^{1/2}\ , 
\qquad
\mathtt{l}_{3}:= -\frac34 \tts \ck^{1/2}\ , 
\qquad
\mathtt{l}_{4}:= - \frac12 \ck^{-\frac12}\, \  . 
\end{aligned}
\end{equation}
The entries $E_{22}$, $H_{22}$
and $I_{22}$ in \eqref{E}, \eqref{H} and \eqref{I} satisfy 
\begin{align} \label{lemma C3}
 E_{22}(0,\e,\tts) H_{22}(0,\e,\tts)=(I_{22}(0,\e,\tts))^2.
\end{align}
\end{proposition}

The rest of the section is devoted to the proof of such result. 
We decompose $\cB_{\mu,\e,\tts}$ in \eqref{mathcal B mu e} as
\begin{equation*}
    \cB_{\mu,\e,\tts} = \cB_{0,\e,\tts}  +\mathcal{B}^{\flat} + \mathcal{B}^\sharp+\mathcal{B}^{cap},
\end{equation*}
where $\mathcal{B}_{0,\e,\tts}$, $\mathcal{B}^{\flat}$, $\mathcal{B}^\sharp$, $\mathcal{B}^{cap}$ are the self-adjoint and reversibility preserving operators:
\begin{align} \label{B epsilon}
\mathcal{B}_{0,\e,\tts} :=
\begin{bmatrix}
1 + a_{\e,\tts}(x)-\kappa \Sigma_{0,\e,\tts} & -(\ck + p_{\e,\tts}(x))\partial_x \\
\partial_x \circ (\ck + p_{\e,\tts}(x)) & |D|
\end{bmatrix},
\qquad 
\mathcal{B}^{\flat} :=
\begin{bmatrix}
0 & 0 \\
0 & \mu\left(\mathrm{sgn}(D)+\Pi_0\right)
\end{bmatrix},
\end{align}
\begin{align} \label{B sharp}
\mathcal{B}^\sharp :=
\begin{bmatrix}
0 & -\im\, \mu p_{\e,\tts} \\
\im\, \mu p_{\e,\tts} & 0
\end{bmatrix},
\qquad
\mathcal{B}^{cap} :=
\begin{bmatrix}
-\kappa \left(\mu\, \partial_{\mu}\Sigma_{0,\e,\tts}+\frac{1}{2}\mu^2 \, \partial_{\mu}^2 \Sigma_{0,\e,\tts}\right) & 0 \\
0 & 0
\end{bmatrix},
\end{align}
where $a_{\e,\tts}(x)$ and $p_{\e,\tts}(x)$ are given by \eqref{SN1}, and where
\begin{align}\label{Bcap.op}
\partial_{\mu}\Sigma_{0,\e,\tts}:=\frac{\im}{1+\mathfrak{p}_x}\left(g_{\e,\tts}(x)\circ\pa_x+\pa_x\circ g_{\e,\tts}(x)\right)\frac{1}{1+\mathfrak{p}_x}, \quad \partial_{\mu}^2 \Sigma_{0,\e,\tts}:=-\frac{2g_{\e,\tts}(x)}{(1+\mathfrak{p}_x)^2}.
\end{align}
with $g_{\e, \tts}(x)$ in \eqref{SN1 g}.
We note that $\mathcal{B}^{\flat}$ and $\mathcal{B}^{\sharp}$ are linear in $\mu$.

\begin{lemma} [Expansion of $\mathsf{B}_{0,\e,\tts}$] \label{Expansion of B epsilon}
The self-adjoint and reversibility preserving matrix
$\mathsf{B}_{\e,\tts} := \mathsf{B}_{\e,\tts}(\mu)$ associated, as in \eqref{Bmuepsilon matrix representation},
with the self-adjoint and reversibility preserving operator
$\mathcal{B}_{0,\e,\tts}$ defined in \eqref{B epsilon}, with respect to the basis $\mathcal{F}$ of $\mathcal{V}_{\mu,\e,\tts}$ in \eqref{F basis set and f}, expands as
\small
\begin{align}
\label{BeExp}
& \mathsf B_{\e,\tts}
=\\
\notag
&\left(
\begin{NiceArray}{cc|cc|cc}[code-for-first-col=\scriptstyle]
\mathtt{e}_{1}\e+\mathtt{e}_{2}\e^2+ r(\e^3)
 & \im r(\mu\e^2)
 &\mathtt{i}_{1}\e+\mathtt{i}_{2}\e^2+ r(\e^3)
 & \im r(\mu\e^2)
 &\mathtt{f}_1\e^3 + r(\e^4, \mu \e^2)
 & \im r(\mu \e^3) \\
- \im r(\mu\e^2)
&\mathtt{e}_{1}\e+\mathtt{e}_{6}\e^2 + r(\e^3)
 & \im r(\mu\e^2)
 &\mathtt{i}_{1}\e+\mathtt{i}_{4}\e^2 + r(\e^3)
 & \im r(\mu \e^2)
 & r(\mu \e^3)\\
\hline
\mathtt{i}_{1}\e+\mathtt{i}_{2}\e^2 + r(\e^3)
 & -\im r(\mu \e^2)
 &\mathtt{h}_{1}\e^2+ r(\e^3)
 & \im r(\mu \e^2)
 &\mathtt{l}_1\e^3 + r(\e^4, \mu \e^2) 
 &\im r(\mu \e^3)\\
- \im r(\mu \e^2)
 &\mathtt{i}_{1}\e+\mathtt{i}_{4}\e^2 + r(\e^3)
 & - \im r(\mu \e^2)
 &\mathtt{e}_{1}\e+\mathtt{h}_{5}\e^2+r(\e^3)
 & \im r(\mu \e^2) &  r(\mu \e^3)\\
\hline
\mathtt{f}_{1}\e^3 + r(\e^4, \mu \e^2)
 & - \im r(\mu \e^2)
 &\mathtt{l}_{1}\e^3 + r(\e^4, \mu \e^2)
 & - \im r(\mu \e^2) & 1 + r(\e^3,\mu \e^2) 
 & \im r(\mu \e^2)\\
 - \im r(\mu \e^3)
& r(\mu \e^3) 
& - \im r(\mu \e^3)
& r(\mu \e^3)
& - \im r(\mu \e^2)
& 0
\end{NiceArray}
\right)\\
\notag
&
+ \cO(\mu^2), 
\end{align}
\normalsize
with the coefficients defined in \eqref{coeff.B}.
\end{lemma}

\begin{remark} \label{remark im omega0}
By  Lemma~\ref{Expansion of B epsilon}, the characteristic polynomial of the  matrix $\mathsf{L}_{0,\e,\tts}= \mathsf{J}_6 \mathsf{B}_{\e,\tts}(0)$ is given by 
   $ \det(\mathsf{L}_{0,\e,\tts} - \lambda I_6) = \lambda^2 \left( \lambda^4 + T_\e \lambda^2 + D_\e \right)$
where
\begin{align}
    T_\e &\coloneqq E_{11}(0,\e,\tts) E_{22}(0,\e,\tts) + 2 I_{11}(0,\e,\tts) I_{22}(0,\e,\tts) + H_{11}(0,\e,\tts) H_{22}(0,\e,\tts) = \frac{9}{8} \e^2 + r(\e^3) , \label{eq:Texp} \\
    D_\e &\coloneqq \left( E_{11}(0,\e,\tts) H_{11}(0,\e,\tts) - I_{11}(0,\e,\tts)^2 \right) \left( E_{22}(0,\e,\tts) H_{22}(0,\e,\tts) - I_{22}(0,\e,\tts)^2\right) .\label{eq:Dexp}
\end{align}
By \eqref{lemma C3}, one has $D_\e\equiv 0$. Therefore, for sufficiently small $\e\neq0$, $\mathsf{L}_{0,\e,\tts}$ has two nonzero purely imaginary eigenvalues $\pm\im\sqrt{T_\e}$.
\end{remark}

\begin{proof}

We expand the matrix $\mathsf{B}_{\e,\tts}(\mu)$ as
\begin{equation*} 
\mathsf{B}_{\e,\tts}(\mu) = \mathsf{B}_{\e,\tts}(0) + \mu(\partial_\mu \mathsf{B}_{\e,\tts})(0)  + \mathcal{O}(\mu^2).
\end{equation*}

\noindent\textbf{The matrix $\mathsf{B}_{\e,\tts}(0)$.}
We start recalling that the self-adjoint matrix $\mathsf{B}_{\e,\tts}(0)$ has the following form 
\small
\begin{align} \label{B epsilon 0}
    \mathsf{B}_{\e,\tts}(0) =
\begin{pmatrix}
 \left( \mathcal{B}_{\e,\tts} f_2^{+},\ f_{2}^{+} \right) &
 \left( \mathcal{B}_{\e,\tts} f_2^{-},\ f_{2}^{+} \right) &
 \left( \mathcal{B}_{\e,\tts} f_1^{+},\ f_{2}^{+} \right) &
 \left( \mathcal{B}_{\e,\tts} f_1^{-},\ f_{2}^{+} \right) &
 \left( \mathcal{B}_{\e,\tts} f_0^{+},\ f_{2}^{+} \right) &
 \left( \mathcal{B}_{\e,\tts} f_0^{-},\ f_{2}^{+} \right) \\
 \left( \mathcal{B}_{\e,\tts} f_2^{+},\ f_{2}^{-} \right) &
 \left( \mathcal{B}_{\e,\tts} f_2^{-},\ f_{2}^{-} \right) &
 \left( \mathcal{B}_{\e,\tts} f_1^{+},\ f_{2}^{-} \right) &
 \left( \mathcal{B}_{\e,\tts} f_1^{-},\ f_{2}^{-} \right) &
 \left( \mathcal{B}_{\e,\tts} f_0^{+},\ f_{2}^{-} \right) &
 \left( \mathcal{B}_{\e,\tts} f_0^{-},\ f_{2}^{-} \right) \\
 \left( \mathcal{B}_{\e,\tts} f_2^{+},\ f_{1}^{+} \right) &
 \left( \mathcal{B}_{\e,\tts} f_2^{-},\ f_{1}^{+} \right) &
 \left( \mathcal{B}_{\e,\tts} f_1^{+},\ f_{1}^{+} \right) &
 \left( \mathcal{B}_{\e,\tts} f_1^{-},\ f_{1}^{+} \right) &
 \left( \mathcal{B}_{\e,\tts} f_0^{+},\ f_{1}^{+} \right) &
 \left( \mathcal{B}_{\e,\tts} f_0^{-},\ f_{1}^{+} \right) \\
 \left( \mathcal{B}_{\e,\tts} f_2^{+},\ f_{1}^{-} \right) &
 \left( \mathcal{B}_{\e,\tts} f_2^{-},\ f_{1}^{-} \right) &
 \left( \mathcal{B}_{\e,\tts} f_1^{+},\ f_{1}^{-} \right) &
 \left( \mathcal{B}_{\e,\tts} f_1^{-},\ f_{1}^{-} \right) &
 \left( \mathcal{B}_{\e,\tts} f_0^{+},\ f_{1}^{-} \right) &
 \left( \mathcal{B}_{\e,\tts} f_0^{-},\ f_{1}^{-} \right) \\
 \left( \mathcal{B}_{\e,\tts} f_2^{+},\ f_{0}^{+} \right) &
 \left( \mathcal{B}_{\e,\tts} f_2^{-},\ f_{0}^{+} \right) &
 \left( \mathcal{B}_{\e,\tts} f_1^{+},\ f_{0}^{+} \right) &
 \left( \mathcal{B}_{\e,\tts} f_1^{-},\ f_{0}^{+} \right) &
 \left( \mathcal{B}_{\e,\tts} f_0^{+},\ f_{0}^{+} \right) &
 \left( \mathcal{B}_{\e,\tts} f_0^{-},\ f_{0}^{+} \right)\\
 \left( \mathcal{B}_{\e,\tts} f_2^{+},\ f_{0}^{-} \right) &
 \left( \mathcal{B}_{\e,\tts} f_2^{-},\ f_{0}^{-} \right) &
 \left( \mathcal{B}_{\e,\tts} f_1^{+},\ f_{0}^{-} \right) &
 \left( \mathcal{B}_{\e,\tts} f_1^{-},\ f_{0}^{-} \right) &
 \left( \mathcal{B}_{\e,\tts} f_0^{+},\ f_{0}^{-} \right) &
 \left( \mathcal{B}_{\e,\tts} f_0^{-},\ f_{0}^{-} \right)
\end{pmatrix}
,
\end{align}
\normalsize
where in the above matrix $f^{\sigma}_k \equiv  f^{\sigma}_k(0,\e,\tts)$, for $\sigma=\pm$, $k=0,1,2$. The matrix $\mathsf{B}_{\e,\tts}(0)$ is real, because the operator $\mathcal{B}_{\e,\tts}$ is real and the basis $\{ f_k^{\pm}(0, \e,\tts)\}_{k=0,1,2}$ is real.
On the other hand, by \eqref{B are alternatively real or imaginary}, its matrix elements $(\mathsf{B}_{\e,\tts}(0))_{i,j}$ vanish for $i + j$ odd. In addition, $f_0^{-}(0,\e,\tts) =
\begin{bmatrix}
0 \\ 1
\end{bmatrix}
$ by \eqref{f-0}, and, by \eqref{mathcal B mu e}, we get $\mathcal{B}_{\e,\tts} f_0^{-}(0,\e,\tts) = 0$, for any $\e$.
We deduce that the self-adjoint matrix $\mathsf{B}_{\e,\tts}(0)$ in \eqref{B epsilon 0} has the form
\small
\begin{equation} \label{B_e 0}
\mathsf{B}_{\e,\tts}(0) =
\left(
\begin{array}{cc|cc|cc}
E_{11}(0,\e,\tts)
  & 0
  & I_{11}(0,\e,\tts)
  & 0
  &  F_{11}(0,\e,\tts)
  & 0\\
    0
  & E_{22}(0,\e,\tts)
  & 0
  & I_{22}(0,\e,\tts)
  & 0
  & 0 \\
\hline
I_{11}(0,\e,\tts)
  & 0
  & H_{11}(0,\e,\tts)
  & 0
  & L_{11}(0,\e,\tts)
  & 0 \\
0
  & I_{22}(0,\e,\tts)
  & 0
  & H_{22}(0,\e,\tts)
  & 0
  & 0 \\
  \hline
  F_{11}(0,\e,\tts)
  & 0
  & L_{11}(0,\e,\tts)
  & 0
  &  G_{11}(0,\e,\tts)
  & 0 \\
  0
  & 0
  & 0
  & 0
  & 0
  & 0
\end{array}
\right)
,
\end{equation}
\normalsize
with real entries. 
We now compute the Taylor expansions of these entries.
First, from \eqref{mathcal B mu e}, we get 
\begin{align} \label{exp B0e}
\mathcal{B}_{0,\e,\tts} & = \mathcal{B}_{0} + \e \mathcal{B}_1 + \e^2 \mathcal{B}_2 + \cO(\e^3), \\
\notag
& \mathcal{B}_{0} \coloneqq
\begin{bmatrix}
    1-\kappa\partial_{xx} & -\ck  \partial_x  \\
    \ck \partial_x  & |D|
\end{bmatrix} , \quad 
\mathcal{B}_{j} \coloneqq
\begin{bmatrix}
    a_j(x) - \frac12 \Sigma_j & - p_j(x) \partial_x \\
    \partial_x \circ p_j(x) & 0
\end{bmatrix} , 
\ j= 1, 2, 
\end{align}
where $a_1$, $a_2$, $p_1$, $p_2$, $\Sigma_1$ and $\Sigma_2$ are given by Lemma \ref{lem:pa.exp}. 
We also need the first derivatives of the basis vectors at $(\mu,\e)=(0,0)$, which, by 
Lemma \ref{expansion of the basis F},
are given by
\small
\begin{align}
 \partial_\e f_2^+(0,0)
 &=\begin{bmatrix}
 \ck^{-1/2}\left(\frac14\cos (x)-\frac32\cos(3x)\right)\\
 \ck^{1/2}\left(-\frac14\sin (x)-\frac52\sin(3x)-\tts\sin(4x)\right)
 \end{bmatrix}, \quad 
 \partial_\e f_2^-(0,0)
 =\begin{bmatrix}
 \ck^{-1/2}\left(-\frac14\sin (x)+\frac32\sin(3x)\right)\\
 \ck^{1/2}\left(-\frac14\cos (x)-\frac52\cos(3x)-\tts\cos(4x)\right)
 \end{bmatrix}, \nonumber\\
 \partial_\e  f_1^+(0,0)
 &=\begin{bmatrix}
 \ck^{-1/2}\left(\frac38 \tts \cos (x)+\frac38\cos(2x)-\frac34 \tts \cos(3x)\right)\\
 \ck^{1/2}\left(-\frac38 \tts \sin (x)-\frac38\sin(2x)-\frac74 \tts \sin(3x)\right)
 \end{bmatrix}, \quad 
 \partial_\e  f_1^-(0,0)
 =\begin{bmatrix}
 \ck^{-1/2}\left(\frac38 \tts \sin (x)-\frac38\sin(2x)+\frac34 \tts \sin(3x)\right)\\
 \ck^{1/2}\left(\frac38 \tts \cos (x)-\frac38\cos(2x)-\frac74 \tts \cos(3x)\right)
 \end{bmatrix}, \nonumber\\
 \label{basis first derivatives quadratic calculation}
 \partial_\e f_0^+(0,0)
 &=\begin{bmatrix}
 \cos (x)+\tts\cos(2x)\\
 -\ck\left(\sin (x)+\tts\sin(2x)\right)
 \end{bmatrix},
 \qquad \partial_\e  f_0^-(0,0)=\begin{bmatrix}0\\0\end{bmatrix}.
\end{align}
\normalsize
Then, using also 
 $\cB_{0}f_j^\pm=0$, $j=1,2$, 
$ \cB_{0}f_0^-=0$, 
$\cB_{0}f_0^+=f_0^+$ and that 
 \eqref{f-0=01} imply that the first component of every $\partial_\e^2 f_j^\pm(0,0)$ has zero space average,
 by  \eqref{exp B0e} and \eqref{basis first derivatives quadratic calculation} we obtain
\begin{align*}
 E_{11}(0,\e,\tts)
 =&-\frac{\ck \tts }{2}\e
 +\e^2\big[ \underbrace{(\mathcal B_2f_2^+,f_2^+)}_{ = \frac{57 \ck}{16} }
 +2\underbrace{(\mathcal B_1\partial_\e  f_2^+(0,0),f_2^+)}_{=-\frac{35\ck}{8}}
 +\underbrace{(\mathcal B_{0}\partial_\e f_2^+(0,0),
 \partial_\e  f_2^+(0,0))}_{=\frac{35\ck}{8}}\big]
 +r(\e^3)\\
 & = -\frac{\ck \tts }{2}\e -\frac{13 \ck}{16}\e^2+r(\e^3).
\end{align*}
where we used \eqref{pino2fd}, \eqref{aino2fd}, \eqref{d2x}--\eqref{h2 h02} for the first quadratic term and 
\eqref{pino1fd}, \eqref{aino1fd} and
\eqref{basis first derivatives quadratic calculation} for the second and third. 
The other six entries are obtained in the same way:
\begin{align*}
 E_{22}(0,\e,\tts)
 &=-\frac12\ck\tts\e-\frac{5}{16}\ck\e^2+r(\e^3),\quad  G_{11}(0,\e,\tts)
 =1+r(\e^3) \ , \\
 I_{11}(0,\e,\tts)
 &=\frac12\ck\e-\frac{35}{16}\tts\ck\e^2+r(\e^3),\quad 
 I_{22}(0,\e,\tts)
 =\frac12\ck\e-\frac{1}{16}\tts\ck\e^2+r(\e^3),
 \\
 H_{11}(0,\e,\tts)
 &=-\frac{3}{16}\ck\e^2+r(\e^3),\quad 
 H_{22}(0,\e,\tts)
 =-\frac12\ck\tts\e+\frac{7}{16}\ck\e^2+r(\e^3) , 
 \nonumber
\end{align*}
The entries displayed as zero in \eqref{B_e 0} vanish by
the parity relations and
$\mathcal B_{0,\e,\tts}f_0^-(0,\e,\tts)=0$.

Regarding the functions  $F_{11}(0,\e,\tts)$ and $L_{11}(0,\e,\tts)$, we show in Lemma~\ref{lem:appendix-pure-cubic-cross} that  they have  cubic Taylor expansion
\[
 F_{11}(0,\e,\tts)
 = \frac{3\ck^{3/2} \tts}{16} \e^3+r(\e^4),
 \qquad
 L_{11}(0,\e,\tts)
 =- \frac{\ck^{3/2}}{8}\e^3+r(\e^4).
\]
\textbf{Terms linear in $\mu$ of  $\mathsf{B}_{\e,\tts}(\mu)$.}
We have
\begin{align} \label{X+X}
    \partial_\mu \mathsf{B}_{\e,\tts}(0) = X + X^* \quad \text{where} \quad
X :=  \left( \left( \mathcal{B}_{0,\e,\tts} f_k^\sigma(0,\e,\tts), \, (\partial_\mu f_{k'}^{\sigma'})(0,\e,\tts) \right) \right)_{\sigma,\sigma'=\pm;\,k,k'=0,1,2}.
\end{align}
By Lemma \ref{expansion of the basis F},  the operator $\sL_{0,\e,\tts}: \cV_{0,\e,\tts} \to \cV_{0,\e,\tts}$ given by \eqref{mathscr L mu e} satisfies
$ \sL_{0,\e,\tts}  f_0^{-}(0,\e,\tts) = 0$
so that $\cB_{0,\e,\tts}  f_0^{-}(0,\e,\tts) = 0$.
Thus the last column of $X$ is zero. To compute the other entries, we Taylor expand
the terms $(\partial_\mu  f_{k'}^{\sigma'})(0,\e,\tts)$ getting, by 
 Lemma \ref{expansion of the basis F},  
$$
\partial_{\mu}  f_j^+(0,\e,\tts) =
 \im  \,
\rho^+_{j} f^-_{-j}
+  \im \epsilon
\begin{bmatrix}
\mathit{odd}(x) \\
\mathit{even}(x)
\end{bmatrix}
+ \mathcal{O}( \epsilon^2)\,,
\quad 
\partial_{\mu} f_j^-(0,\e,\tts) =
 \im  \, \rho^-_{j} f^+_{-j}
+   \im \epsilon 
\begin{bmatrix}
{\mathit{even}_0}(x) \\
\mathit{odd}(x)
\end{bmatrix}
 + \mathcal{O}( \epsilon^2)\,, \ \ j = 1, 2
$$
where  $\rho_2^+ = \rho_2^- = - \frac{1}{24}$, $\rho_1^+ = \rho_1^- =\frac{1}{12}$, whereas
\begin{equation*}
\begin{aligned}
\partial_{\mu} f_0^+(0,\e,\tts) &=
 \epsilon \left(
\frac{1}{4} \ck^{-3/2} f^+_{-1} + \tts \frac{1}{16} \ck^{-3/2} f^+_{-2}+ \im
\begin{bmatrix}
\mathit{odd}(x) \\
\mathit{even}_0(x)
\end{bmatrix}
\right)  + \mathcal{O}( \epsilon^2)\,,
\end{aligned}
\end{equation*}

\begin{equation*} 
\begin{aligned}
\partial_{\mu}  f_0^-(0,\e,\tts) &=
\epsilon\left(
\frac{1}{2} \ck^{-1/2} f^-_{-1} + \frac{\tts}{4} \ck^{-1/2} f^-_{-2}
+ \im
\begin{bmatrix}
\mathit{even}_0(x) \\
\mathit{odd}(x)
\end{bmatrix}\right)
+ \mathcal{O}(\epsilon^2)\,.
\end{aligned}
\end{equation*}
Using the above expansions, \eqref{exp B0e},
\eqref{X+X}, and Lemma~\ref{expansion of the basis F}, a direct computation yields
\begin{align} \label{mat X}
X &=
\begin{pmatrix}
    \cO(\e^3) & \cO(\e^2) & \cO(\e^3) & \cO(\e^2) & \cO(\e^2) & 0 \\
    \cO(\e^2) & \cO(\e^3) & \cO(\e^2) & \cO(\e^3) & \cO(\e^2) & 0 \\
    \cO(\e^3) & \cO(\e^2) & \cO(\e^3) & \cO(\e^2) & \cO(\e^2) & 0 \\
    \cO(\e^2) & \cO(\e^3) & \cO(\e^2) & \cO(\e^3) & \cO(\e^2) & 0 \\
    \cO(\e^2) & \cO(\e^2) & \cO(\e^2) & \cO(\e^2) & \cO(\e^2) & 0 \\
     X_{61} & \cO(\e^3) & X_{63} & \cO(\e^3) & \cO(\e^2) & 0 \\
\end{pmatrix}.
\end{align}
We are left to show that  both $X_{61}$ and $X_{63}$ are $\cO(\e^3)$.
This follows since,  for $j=1,2$,
\begin{align*}
(\mathcal{B}_{0,\e,\tts}f^+_j(0,\e,\tts),\pa_\mu f^-_0(0,\e,\tts))
 =&\e^2 \left(\mathcal{B}_1 f^+_j+\mathcal{B}_0 \pa_\e f^+_j(0,0),\pa_\mu\pa_\e f^-_0(0,0)\right)
 +\cO(\e^3)\overset{\eqref{lemma C1}}{=}\cO(\e^3), 
\end{align*}
where we also used that, by \eqref{f-0}, $\pa_\mu f_0^-(0,0) = 0$.
Combining \eqref{B_e 0} and \eqref{mat X}  we obtain \eqref{BeExp}.
\end{proof}

\begin{lemma}
[Expansion of $\mathsf{B}^{\flat}$] \label{Expansion of B flat}
The self-adjoint and reversibility-preserving matrix $\mathsf{B}^{\flat}$ associated, as in \eqref{Bmuepsilon matrix representation}, to the self-adjoint and reversibility-preserving operator $\mathcal{B}^{\flat}$, defined in \eqref{B epsilon}, with respect to the basis $\mathcal{F}$ of $\mathcal{V}_{\mu, \epsilon}$ in \eqref{F basis set and f}, admits the expansion
\small
\begin{align*}
&\mathsf{B}^{\flat} =  \\
&
\left(
\begin{array}{cc|cc|cc}
0 & \mathrm{i} \frac{\ck}{2} \mu + \im r(\mu \e^2) & 0 & 
\im r(\mu \e)  & 0 & 0 \\
-\mathrm{i} \frac{\ck}{2} \mu - \im r(\mu \e^2) & r(\mu \e^4)  & \im r( \mu \e) & r(\mu \e^4) &  \mathrm{i} \tts \frac{\ck^{3/2}}{2} \mu \e + \im r(\mu \e^2) & 
\frac{\ck^{1/2}}8\mu\e^2 +r(\mu\e^3)
\\
\hline
0 & -\im r(\mu \e) & 0  & \mathrm{i} \frac{\ck}{2} \mu + \im r(\mu \e^2) & 0 & 0 \\
- \im r(\mu \e) & r(\mu \e^4) & -\mathrm{i} \frac{\ck}{2} \mu - \im r(\mu \e^2) & r(\mu \e^4)  &  \im \frac{\ck^{3/2}}{2} \mu \e +\im r(\mu \e^2)  & r(\mu\e^3) \\
\hline
0 & - \im \tts \frac{\ck^{3/2}}{2} \mu \e -\im r(\mu \e^2)  & 0 & - \mathrm{i}\frac{\ck^{3/2}}{2}  \mu \e +\im r(\mu \e^2)  & 0 & 0 \\
0 & \frac{\ck^{1/2}}8\mu\e^2 +r(\mu\e^3) & 0 & r(\mu\e^3) & 0 & \mu 
\end{array}
\right) \nonumber \\
&\quad+  \mathcal{O}( \mu^2). \nonumber
\end{align*}
\normalsize
\end{lemma}

\begin{proof}
We compute the matrix
\begin{align*}
 \left( \left( \mathcal{B}^{\flat}f_k^\sigma(\mu,\e,\tts), \, f_{k'}^{\sigma'}(\mu,\e,\tts) \right) \right)_{\sigma,\sigma'=\pm;\,k,k'=0,1,2}
\end{align*}
with respect to the basis $\mathcal{F}$. Since $\mathcal{B}^{\flat}$ is self-adjoint and reversibility-preserving, it suffices to determine the entries in the upper triangular part and then use the corresponding symmetry relations.

We begin with the last column and row. By \eqref{f-0} and \eqref{B epsilon}, we obtain
$
\mathcal{B}^{\flat}  f_0^{-}(\mu,\e,\tts)
=
\begin{bmatrix}
0\\
\mu
\end{bmatrix}
+\mathcal{O}(\mu^2\epsilon).
$
Using \eqref{f+2}--\eqref{f-0}  we infer that
\begin{equation*}
\begin{aligned}
 \left( \mathcal{B}^{\flat} f_0^{-}(\mu,\e,\tts), f_0^{-}(\mu,\e,\tts) \right)
&= \mu + \mathcal{O}(\mu^3\e, \mu^2\epsilon^{2}),\quad 
\bigl( \mathcal{B}^{\flat}  f_0^{-}(\mu,\e,\tts), f_k^{+}(\mu,\e,\tts) \bigr)= \mathcal{O}(\mu^2\epsilon),
\qquad
 k =0,1,2.
\end{aligned}
\end{equation*}
For the remaining entries,  
formula \eqref{B epsilon} for the operator $\cB^{\flat}$, together with the above formula for $\cB^{\flat} f_0^{-}(\mu,\e,\tts)$, gives for $j=1,2$
\begin{align*}
  \left( \mathcal{B}^{\flat}f_0^{-}(\mu,\e,\tts), \, f_{j}^{-}(\mu,\e,\tts) \right)  &= \mu \left( f_0^{-}, \, f_{j}^{-}(0,\e,\tts) \right) + \cO(\mu^2)\\
  & = \mu \underbrace{\left( f_0^{-}, \, f_{j}^{-} \right) }_{=0}
  + \mu  \e\underbrace{\left( f_0^{-}, \, \pa_{\e}f_{j}^{-}(0,0) \right)}_{=0 \small \mbox{ by Lemma } \ref{expansion of the basis F}}
  + \mu \e^2\left( f_0^{-}, \, \frac{1}{2}\pa_{\e}^2f_{j}^{-}(0,0) \right)
  +\cO(\mu^2, \mu \e^3)
\end{align*}
To compute the last term, use that $P_{\mu,\e}^2=P_{\mu,\e}$ implies that $P_{0,0} P_{0,2} P_{0,0} = - P_{0,0}P_{0,1}^{2}P_{0,0}$ and then \eqref{A2} gives $  P_{0,0}U_{0,2}P_{0,0}
=-\frac12 P_{0,0}P_{0,1}^{2}P_{0,0}$. Also, using Lemma~\ref{properties of U and P} (i), we have the following identity
\begin{align} \label{20261117}
P_{0,0}^* f_0^- = -P_{0,0}^*\mathcal J f^+_0=-\mathcal J P_{0,0} f^+_0=-\mathcal J f^+_0= f_0^-.
\end{align}
Then, we obtain for $j=1,2$ 
\begin{align*}
   \frac{1}{2} \left(f_0^-,\pa_{\e}^2f_{j}^{-}(0,0)\right)&=\bigl(f_0^-,U_{0,2}f_{j}^{-}(0,0)\bigr)=\bigl(f_0^-,U_{0,2}P_{0,0}f_{j}^{-}(0,0)\bigr)\\ &\overset{\eqref{20261117}}{=} \bigl(P^*_{0,0}f_0^-,U_{0,2}P_{0,0}f_{j}^{-}(0,0)\bigr)=\bigl(f_0^-,P_{0,0}U_{0,2}P_{0,0}f_{j}^{-}(0,0)\bigr)\\
    &=\left(f_0^-,-\frac12 P_{0,0}P_{0,1}^{2}P_{0,0}f_{j}^{-}(0,0)\right)=-\frac{1}{2}\left(f^-_0,P_{0,1}^2 f^-_j\right)\\
    &=\frac{1}{2}\left(\mathcal{J} f^+_0, P_{0,1}^2 f^-_j\right)=\frac{1}{2}\left(P_{0,1}^*\mathcal{J} f^+_0, P_{0,1} f^-_j\right)=\frac{1}{2}\left(\mathcal{J} P_{0,1} f^+_0, P_{0,1} f^-_j\right)=
\begin{cases}
\displaystyle \frac{\ck^{1/2}}8,
   & j=2,\\[1ex]
0, & j=1.
\end{cases}
\end{align*}
This implies
\begin{align*}
    \left( \mathcal{B}^{\flat}f_0^{-}(\mu,\e,\tts), \, f_{j}^{-}(\mu,\e,\tts) \right)=\begin{cases}
\displaystyle (\mathsf B^{\flat})_{26}=\frac{\ck^{1/2}}8 \mu\e^2+\cO(\mu\e^3,\mu^2),\\[1ex]
(\mathsf B^{\flat})_{46}=\cO(\mu\e^3,\mu^2).
\end{cases}
\end{align*}
This gives the expansion of the last column and, by self-adjointness, of the last row as well.

It remains to determine the upper-left $5\times5$ block
of $\mathsf B^\flat$. Since $\mathcal B^\flat$ is linear in $\mu$ and the basis vectors
are analytic in $\mu$, we have
\begin{align*}
\bigl(\mathcal B^\flat f_k^\sigma(\mu,\e,\tts),
f_{k'}^{\sigma'}(\mu,\e,\tts)\bigr)=
\bigl(\mathcal B^\flat f_k^\sigma(0,\e,\tts),
f_{k'}^{\sigma'}(0,\e,\tts)\bigr)
+\cO(\mu^2),
\qquad k,k'=0,1,2,\quad \sigma,\sigma'=\pm.
\end{align*}
For entries $(1,4)$ and $(2,3)$, using
$f_k^\sigma(0,\e,\tts)=f_k^\sigma+\cO(\e)$
and Fourier orthogonality, we obtain
\begin{align*}
(\mathsf B^\flat)_{14}=\underbrace{\bigl(\mathcal B^\flat f_1^-,f_2^+\bigr)}_{=0}
+\cO(\mu\e,\mu^2),\qquad (\mathsf B^\flat)_{23}
=\underbrace{\bigl(\mathcal B^\flat f_1^+,f_2^-\bigr)}_{=0}
+\cO(\mu\e,\mu^2).
\end{align*}
Their purely imaginary character follows from reversibility.

Next, we consider all pairings with $\sigma=\sigma'=+$.
By Lemma~\ref{expansion of the basis F}, the second components $\bigl(f_k^+(0,\e,\tts)\bigr)_2$, $k=0,1,2$, are real odd functions. Hence $\Pi_0 \bigl(f_k^+(0,\e,\tts)\bigr)_2=0$ and
$\operatorname{sgn}(D)\bigl(f_k^+(0,\e,\tts)\bigr)_2$ is even. Therefore,
\begin{align*}
\bigl(\mathcal B^\flat f_k^+(0,\e,\tts),
f_{k'}^+(0,\e,\tts)\bigr)
&=\frac{\mu}{2\pi}\int_0^{2\pi}
\bigl(\operatorname{sgn}(D)\bigl(f_k^+(0,\e,\tts)\bigr)_2\bigr)
\overline{\bigl(f_{k'}^+(0,\e,\tts)\bigr)_2}\,\mathrm{d}x=0,
\qquad k,k'=0,1,2.
\end{align*}
Consequently, the entries $(1,1)$, $(1,3)$, $(1,5)$,
$(3,3)$, $(3,5)$ and $(5,5)$ are all $\cO(\mu^2)$.

By Lemma~\ref{expansion of the basis F}, $\Pi_0\bigl(f_k^-(0,\e,\tts)\bigr)_2=\cO(\e^2)$ and $\operatorname{sgn}(D)\bigl(f_k^-(0,\e,\tts)\bigr)_2$ is odd for $k=1,2$. A straightforward calculation gives, for $k,k'=1,2$,
\begin{align*}
\bigl(\mathcal B^\flat f_k^-(0,\e,\tts),
f_{k'}^-(0,\e,\tts)\bigr)
&=\frac{\mu}{2\pi}\int_0^{2\pi}
\bigl((\operatorname{sgn}(D)+\Pi_0)\bigl(f_k^-(0,\e,\tts)\bigr)_2\bigr)
\overline{\bigl(f_{k'}^-(0,\e,\tts)\bigr)_2}\,\mathrm{d}x=\cO(\mu\e^4).
\end{align*}
Therefore, the entries $(2,2)$, $(2,4)$, $(4,4)$ are all $\cO(\mu\e^4,\mu^2)$.

Finally, using Lemma~\ref{expansion of the basis F},
a direct calculation yields $(\mathsf B^{\flat})_{12}$, $(\mathsf B^{\flat})_{34}$, $(\mathsf B^{\flat})_{25}$, $(\mathsf B^{\flat})_{45}$. The other entries are determined by self-adjointness of $\mathcal{B}^{\flat}$.
\end{proof}

\begin{lemma}
[Expansion of $\mathsf{B}^{\sharp}$] \label{Expansion of B sharp}
The self-adjoint and reversibility-preserving matrix $\mathsf{B}^{\sharp}$ associated, as in \eqref{Bmuepsilon matrix representation}, to the self-adjoint and reversibility-preserving operators $\mathcal{B}^{\sharp}$, defined in \eqref{B sharp}, with respect to the basis $ \mathcal{F}$ of $\mathcal{V}_{\mu,\epsilon}$ in  \eqref{F basis set and f}, admits the expansion
\setlength{\arraycolsep}{3.5pt}
\small
\begin{align*} 
\mathsf{B}^{\sharp} &=
\left(
\begin{array}{cc|cc|cc}
0 &\im r(\mu \e)  & 0 & \im r(\mu \e)  & 0 &   \im c^\sharp \\
-\im r(\mu \e)  & 0 &  \im r(\mu \e)  & 0 & -\im e^\sharp & 0 \\
\hline
0 & -\im r(\mu\e) & 0 & \im r(\mu\e) & 0 & \mathrm{i} g^\sharp \\
-\mathrm{i}r(\mu \e) & 0 & -\im r(\mu \e) & 0 & -\mathrm{i}h^\sharp & 0 \\
\hline
0 & \mathrm{i} e^\sharp & 0 & \mathrm{i} h^\sharp  & 0 & \im r(\mu \e) \\
-\mathrm{i} c^\sharp  & 0 & -\mathrm{i} g^\sharp & 0 & -\im r(\mu \e) & 0
\end{array}
\right) +  \mathcal{O}(\mu^2 \e) , 
\end{align*}
\normalsize
where 
\begin{equation}\label{Bsharp.coeff}
\begin{aligned}
c^\sharp :=  \ck^{1/2} \mu \e\big(\tts -\frac32 \e + r(\e^2) \big), \quad e^\sharp := \tts \ck^{3/2} \mu (\e + r(\e^2)), \\
g^\sharp := \ck^{1/2} \mu \e\big(1-\frac34 \tts \e + r(\e^2) \big), \quad h^\sharp:= \ck^{3/2} \mu (\e + r(\e^2)).
\end{aligned}
\end{equation}
\end{lemma}

\begin{proof}

  Since $\mathcal{B}^{\sharp} = -\mathrm{i} \mu p_{\e,\tts} \mathcal{J}$ and $p_{\e,\tts} = \mathcal{O}(\epsilon)$ by \eqref{SN1}, we have the expansion
 \begin{equation*}
 \left( \mathcal{B}^{\sharp} f_k^\sigma (\mu,\e,\tts),\, f_{k'}^{\sigma'}(\mu,\e,\tts) \right)_{\sigma,\sigma'=\pm;\,k,k'=0,1,2}
 = \left( \mathcal{B}^{\sharp} f_k^\sigma (0,\e,\tts),\, f_{k'}^{\sigma'}(0,\e,\tts) \right)_{\sigma,\sigma'=\pm;\,k,k'=0,1,2}
 + \mathcal{O}(\mu^2 \epsilon).
 \end{equation*}
The matrix entries $\left( \mathcal{B}^{\sharp} f_k^\sigma (0,\e,\tts),\, f_{k'}^{\sigma}(0,\e,\tts) \right),\,\sigma = \pm,\,k, k' = 0, 1,2$ are zero, because they are simultaneously real by \eqref{B are alternatively real or imaginary}, and purely imaginary, being the operator $\mathcal{B}^{\sharp}$ purely imaginary and the basis $\{ f_k^{\pm}(0,\e,\tts) \}_{k=0,1,2}$ real. Hence $\mathsf{B}^{\sharp}$ has the form
\small
 \begin{equation*}
 \mathsf{B}^{\sharp} =
 \left(
 \begin{array}{cc|cc|cc}
 0 & \mathrm{i} a^{\sharp} & 0 & \mathrm{i} b^{\sharp} &  0 & \mathrm{i} c^{\sharp}\\
 -\mathrm{i} a^{\sharp} & 0 & -\mathrm{i} d^{\sharp} & 0 & -\mathrm{i} e^{\sharp}  & 0 \\
 \hline
 0 & \mathrm{i} d^{\sharp} & 0 & \mathrm{i} f^{\sharp} & 0 & \mathrm{i} g^{\sharp} \\
 -\mathrm{i} b^{\sharp} & 0 & -\mathrm{i} f^{\sharp} & 0 & -\mathrm{i} h^{\sharp} & 0 \\
 \hline
 0 & \mathrm{i} e^{\sharp} & 0 & \mathrm{i} h^{\sharp} & 0 & \mathrm{i} j^{\sharp} \\
 - \mathrm{i} c^{\sharp} & 0 & -\mathrm{i} g^{\sharp} & 0 & -\mathrm{i} j^{\sharp} & 0
 \end{array}
 \right)
 + \mathcal{O}(\mu^2 \epsilon),
 \end{equation*}
 where $a^{\sharp}, b^{\sharp}, c^{\sharp}, d^{\sharp}, e^{\sharp}, f^{\sharp}, g^{\sharp}, h^{\sharp}, j^{\sharp}$ are real numbers given by
\begin{align*}
a^\sharp &:=  -\im (\mathcal B^{\sharp} f_2^-(0, \e,\tts),f_2^+(0,\e,\tts)) , \quad 
b^\sharp :=  -\im (\mathcal B^{\sharp} f_1^-(0, \e,\tts),f_2^+(0,\e,\tts)) , \quad
c^\sharp :=  -\im (\mathcal B^{\sharp} f_0^-(0, \e,\tts),f_2^+(0,\e,\tts)) , \\
d^\sharp &:=  -\im (\mathcal B^{\sharp} f_2^-(0, \e,\tts),f_1^+(0,\e,\tts)) , \quad
e^\sharp :=  \im (\mathcal B^{\sharp} f_0^+(0, \e,\tts),f_2^-(0,\e,\tts)) , \quad
f^\sharp := -\im
\bigl(\mathcal B^\sharp f_1^-(0,\e,\tts),
f_1^+(0,\e,\tts)\bigr),\\
g^\sharp &:=  -\im (\mathcal B^{\sharp} f_0^-(0, \e,\tts),f_1^+(0,\e,\tts)), \quad
h^\sharp :=  \im (\mathcal B^{\sharp} f_0^+(0, \e,\tts),f_1^-(0,\e,\tts)) , \quad
j^\sharp :=  -\im (\mathcal B^{\sharp} f_0^-(0, \e,\tts),f_0^+(0,\e,\tts)) .
\end{align*}
\normalsize
Observe that, since $\mathcal{B}^{\sharp} = \mathcal{O}(\mu \epsilon)$ in $\mathcal{L}(Y,X)$, the quantities $a^{\sharp}, b^{\sharp}, c^{\sharp}, d^{\sharp}, e^{\sharp}, f^{\sharp}, g^{\sharp}, h^{\sharp}, j^{\sharp}$ are of order $\cO(\mu \epsilon)$.

Using $f_0^-(0,\e,\tts)=f_0^-$, we obtain
\begin{align*}
\bigl(\mathcal B^\sharp f_0^-(0,\e,\tts),
f_j^+(0,\e,\tts)\bigr)
=-\frac{\im\mu}{2\pi}\int_0^{2\pi}
p_{\e,\tts}(x)
\bigl(f_j^+(0,\e,\tts)\bigr)_1(x)\,\mathrm{d}x,
\qquad j=1,2.
\end{align*}
Using \eqref{pino1fd}--\eqref{pino2fd} and the first-order
basis expansions \eqref{f+2} and \eqref{f+1}, we obtain
the expansions of $c^\sharp$ and $g^\sharp$ in
\eqref{Bsharp.coeff}. Moreover, using $p_{\e,\tts}=\e p_1+\cO(\e^2)$ and
$f_k^\sigma(0,\e,\tts)=f_k^\sigma+\cO(\e)$, a direct
computation gives
\begin{align*}
\bigl(\mathcal B^\sharp f_0^+(0,\e,\tts),
f_2^-(0,\e,\tts)\bigr)=-\im\tts\ck^{3/2}\mu\e+\cO(\mu\e^2),\qquad \bigl(\mathcal B^\sharp f_0^+(0,\e,\tts),
f_1^-(0,\e,\tts)\bigr)=-\im\ck^{3/2}\mu\e+\cO(\mu\e^2).
\end{align*}
These identities yield their expansions in \eqref{Bsharp.coeff}.
\end{proof}

We consider the term  $\mathcal{B}^{cap}$ that takes into account the capillary term.

\begin{lemma}
[Expansion of $\mathsf{B}^{cap}$] \label{Expansion of B cap}
The self-adjoint and reversibility-preserving matrix $\mathsf{B}^{cap}$ associated, as in \eqref{Bmuepsilon matrix representation}, to the self-adjoint and reversibility-preserving operator $\mathcal{B}^{cap}$, defined in \eqref{B sharp}, with respect to the basis $\mathcal{F}$ of $\mathcal{V}_{\mu,\epsilon}$ in  \eqref{F basis set and f}, admits the expansion
\small
\begin{align*}
\notag
& \mathsf{B}^{cap} =
\frac{1}{2 \ck}
 \left(
\begin{array}{cc|cc|cc}
0 & \mathrm{i} a^{cap} & 0 &  \im r(\mu \e) & 0 & 0 \\
- \im  a^{cap} & 0 & \im r(\mu\e) & 0 & - \mathrm{i} e^{cap}_1  & 0 \\
\hline
0 & - \im r( \mu \e) & 0 & \mathrm{i}f^{cap}_1 & 0 & 0 \\
- \im r(\mu \e) & 0 & -\mathrm{i} f^{cap}_1 & 0  & -\mathrm{i} h^{cap}_1& 0 \\
\hline
0 & \mathrm{i}  e^{cap}_1 & 0 & \mathrm{i} h^{cap}_1 & 0 & 0 \\
0 & 0 & 0 & 0 & 0 & 0
\end{array}
\right) +\mathcal{O}(\mu^2) , 
\end{align*}
\normalsize
where
\begin{equation}\label{B.cap.coeff}
\begin{aligned}
    a^{cap}:= \mu (2  + r( \e)), \quad 
   e^{cap}_1:= - \tts \ck^{1/2} \mu \e ( 1+ r(\e)), \\ 
   f^{cap}_1 := \mu (1 +r(\e)), 
   \quad 
   h^{cap}_1:= -\frac{\ck^{1/2}}{2}  \mu \e (1 + r(\e)).
   \end{aligned}
\end{equation}
\end{lemma}
\begin{proof}
    Recalling \eqref{SN1 g} and \eqref{exp:g1}, we write the operator $\mathcal{B}^{cap}$ in \eqref{B sharp}--\eqref{Bcap.op} as
    \begin{align*}
        \mathcal{B}^{cap} = -\frac12\mu\mathcal{B}^{cap}_1 + \cO(\mu^2) \ , \quad 
          \mathcal{B}^{cap}_1:=\frac{\im}{1+\mathfrak{p}_x}\left(g_{\e,\tts}(x)\circ\pa_x+\pa_x\circ g_{\e,\tts}(x)\right)\frac{1}{1+\mathfrak{p}_x} \Lambda_{\tilde{\Pi}}
    \end{align*}
    where $\cO(\mu^2) \in\mathcal{L}(Y, X)$,  $g_{\e, \tts}$ is given by \eqref{SN1 g} and $\Lambda_{\tilde{\Pi}}:=\begin{bmatrix}
\mathrm{Id} & 0 \\
0 & 0
\end{bmatrix}$.
To determine the matrix representation of $-\frac12\mu\mathcal{B}^{cap}_{1}$, we use the expansion
\begin{equation*} 
 \left( -\frac12\mu\mathcal{B}^{cap}_{1} f_k^\sigma (\mu,\e,\tts),\, f_{k'}^{\sigma'}(\mu,\e,\tts) \right)
= \left( -\frac12\mu\mathcal{B}^{cap}_{1} f_k^\sigma (0,\e,\tts),\, f_{k'}^{\sigma'}(0,\e,\tts) \right)
+ \mathcal{O}(\mu^2).
\end{equation*}
The matrix entries $\left( -\frac12\mu\mathcal{B}^{cap}_{1} f_k^\sigma (0, \e,\tts),\, f_{k'}^{\sigma}(0,\e,\tts) \right),\, \sigma = \pm,\,k, k' = 0, 1,2$ are zero, because they are simultaneously real by \eqref{B are alternatively real or imaginary}, and purely imaginary, being the operator $-\frac12\mu\mathcal{B}^{cap}_{1}$ purely imaginary and the basis $\{ f_k^{\pm}(0,\e,\tts) \}_{k=0,1,2}$ real. Hence the matrix $-\frac12\mu\mathsf{B}^{cap}_{1}$
 representing the action of
the operator $-\frac12\mu\mathcal{B}^{cap}_{1}$ on the basis
$\mathcal{F}$ of $\mathcal{V}_{\mu,\e}$
has the form 
\small
\begin{equation*}
-\frac12\mu\mathsf{B}^{cap}_{1} =
\frac{1}{2\ck}
\left(
\begin{array}{cc|cc|cc}
 0 & \mathrm{i} a^{cap} & 0 & \mathrm{i} b^{cap}_{1} &  0 & \mathrm{i} c^{cap}_{1} \\
 -\mathrm{i} a^{cap} & 0 & -\mathrm{i} d^{cap}_{1} & 0 & -\mathrm{i} e^{cap}_{1}  & 0 \\
 \hline
 0 & \mathrm{i} d^{cap}_{1} & 0 & \mathrm{i} f^{cap}_{1} & 0 & \mathrm{i} g^{cap}_{1} \\
 -\mathrm{i} b^{cap}_{1} & 0 & -\mathrm{i} f^{cap}_{1} & 0 & -\mathrm{i} h^{cap}_{1} & 0 \\
 \hline
 0 & \mathrm{i} e^{cap}_{1} & 0 & \mathrm{i} h^{cap}_{1} & 0 & \mathrm{i} j^{cap}_{1} \\
 - \mathrm{i} c^{cap}_{1} & 0 & -\mathrm{i} g^{cap}_{1} & 0 & -\mathrm{i} j^{cap}_{1} & 0
\end{array}
\right)
+ \mathcal{O}(\mu^2),
\end{equation*}
\normalsize
where $a^{cap},b^{cap}_1,\ldots,j^{cap}_1$
are real coefficients.

Using that
$\mathcal{B}^{cap}_{1}
= \im (2 \tilde{g}_{\e,\tts}(x) \pa_x + \tilde{g}_{\e,\tts}') \Lambda_{\tilde \Pi}$ with 
$\tilde{g}_{\e,\tts}:=g_{\e,\tts}/(1+\frak{p}_x)^2= 1 + \e(-3 \cos(x) - 3 \tts \cos(2x)) + \cO(\e^2)$, we get
the expressions in \eqref{B.cap.coeff}, and that 
$b^{cap}_1, \, d^{cap}_{1} = r(\mu \e)$. 

Finally, since by Lemma \ref{expansion of the basis F} we have $f_0^-(0,\epsilon,\tts)=\vet{0}{1}$ for small $\e$, using that $\mathcal B^{cap}_1\vet{0}{1}=0$ we get 
that $c^{cap}_1 = g^{cap}_1  = j^{cap}_1 = 0$.
\end{proof}

\begin{lemma}
\label{lem:E22H22 is I22I22}
The entries $E_{22}$, $H_{22}$
and $I_{22}$ in \eqref{E}, \eqref{H} and \eqref{I} satisfy \eqref{lemma C3}.
\end{lemma}
\begin{proof}
For $\e\ne0$, define the normalized translation vector by
\begin{equation}\label{eq:U3}
 \tilde{U}_3:=\ck^{-1/2}\e^{-1}\mathscr T_{\e,\tts}\begin{bmatrix}
 \eta_{\e,\tts,x}\\
 \psi_{\e,\tts,x}
 \end{bmatrix}=f^{-}_1 + \tts f^{-}_2 + \cO(\e),
 \qquad \mathscr T_{\e,\tts}=\mathcal P_{\e,\tts} Z_{\e,\tts}^{-1},
\end{equation}
where $(\eta_{\e,\tts},\psi_{\e,\tts})$ is an analytic solution in Theorem~\ref{Thm: Wilton expansion}, for $\tts \in \{ \pm 1 \}$. For each $a\in\R$, the translated wave
$(\eta_{\e,\tts}(\cdot+a),\psi_{\e,\tts}(\cdot+a))$ solves
\eqref{travelingWWstokes} at the same speed $c_{\e,\tts}$.
Differentiating in $a$ at zero puts
$(\eta_{\e,\tts,x},\psi_{\e,\tts,x})^\top$ in the kernel of the
original linearized operator. Conjugating by $\mathscr T_{\e,\tts}$ gives
$\mathscr L_{0,\e,\tts}\tilde U_3=0$, hence
$\tilde U_3\in\mathcal V_{0,\e,\tts}$ by \eqref{Projection P mu e}.
The leading term in \eqref{eq:U3} follows from \eqref{exp:Wil2} and
$\mathscr T_{0,\tts}=\mathrm{Id}$, and defines $\tilde U_3$ at $\e=0$ by continuity.
The vector $\tilde{U}_3$ is real $\vet{\mbox{odd}(x)}{\mbox{even}(x)}$. By \begingroup \eqref{f-0=01}\endgroup{}, one may expand it as
\begin{equation*}
 \tilde{U}_3=a_2(\e,\tts) f_2^-(0,\e,\tts)+a_1(\e,\tts) f_1^-(0,\e,\tts)+a_0(\e,\tts) f_0^-(0,\e,\tts).
\end{equation*}
Due to the leading order in \eqref{eq:U3} 
\begin{align} \label{alpha2 alpha1}
    a_2(\e,\tts)=\tts+\cO(\e),\qquad\text{and} \qquad a_1(\e,\tts)=1+\cO(\e).
\end{align} 
On the other hand, $-\mathcal{J}\mathscr{L}_{0,\e,\tts} \tilde{U}_3=\mathcal{B}_{0,\e,\tts}\tilde{U}_3=0$, so, recalling \eqref{B_e 0}, we obtain
\small
\begin{equation} \label{B_e 0 U3}
\left(
\begin{array}{cc|cc|cc}
E_{11}(0,\e,\tts)
  & 0
  & I_{11}(0,\e,\tts)
  & 0
  &  F_{11}(0,\e,\tts)
  & 0\\
    0
  & E_{22}(0,\e,\tts)
  & 0
  & I_{22}(0,\e,\tts)
  & 0
  & 0 \\
\hline
I_{11}(0,\e,\tts)
  & 0
  & H_{11}(0,\e,\tts)
  & 0
  & L_{11}(0,\e,\tts)
  & 0 \\
0
  & I_{22}(0,\e,\tts)
  & 0
  & H_{22}(0,\e,\tts)
  & 0
  & 0 \\
  \hline
  F_{11}(0,\e,\tts)
  & 0
  & L_{11}(0,\e,\tts)
  & 0
  &  G_{11}(0,\e,\tts)
  & 0 \\
  0
  & 0
  & 0
  & 0
  & 0
  & 0
\end{array}
\right) \begin{pmatrix}
0\\
a_2(\e,\tts)\\
0\\
a_1(\e,\tts)\\
0\\
a_0(\e,\tts)
\end{pmatrix}=\begin{pmatrix}
0\\
0\\
0\\
0\\
0\\
0
\end{pmatrix}.
\end{equation}
\normalsize
The second and the fourth rows in \eqref{B_e 0 U3} imply that 
\begin{equation} \label{E22 I22 H22}
\begin{pmatrix}
E_{22}(0,\e,\tts) & I_{22}(0,\e,\tts) \\
I_{22}(0,\e,\tts) & H_{22}(0,\e,\tts)
\end{pmatrix} \begin{pmatrix}
a_2(\e,\tts)  \\
a_1(\e,\tts) 
\end{pmatrix} =\begin{pmatrix}
0  \\
0
\end{pmatrix}.
\end{equation} 
Combining \eqref{alpha2 alpha1} and \eqref{E22 I22 H22}, we get the identity \eqref{lemma C3}.
\end{proof}

\begin{proof}[Proof of Proposition  \ref{B decomposition}]
 It follows from lemmata  \ref{Expansion of B epsilon}, \ref{Expansion of B flat},  \ref{Expansion of B sharp}, \ref{Expansion of B cap}, \ref{lem:E22H22 is I22I22}.
 Finally we improve the remainders of the entries (16), (26), (36), (46) (hence also of their corresponding symmetric ones), that in view of these lemmata are of the form $r(\mu \e^3, \mu^2)$.
From \eqref{mathcal B mu e} and Lemma \ref{fmu0 eq f0}, for any $\mu$ sufficiently small
 $$
\mathcal{B}_{\mu,0,\tts}=
\begin{bmatrix}
 1-\kappa (\pa_x + \im \mu)^2 & 
-\ck\partial_x&  \\
~~\ck\partial_x       &  |D+\mu|
    \end{bmatrix} \ , \ \ f_0^-(\mu, 0) = f_0^- \ , \quad \int_\T f_{j}^{\pm}(\mu, 0) \di x = \vet{0}{0}  \ , \  j=1,2 \  . 
 $$
Hence at $\e =0$, the matrix elements 
$(\mathsf{B}_{\mu,0,\tts})_{16}$, 
$(\mathsf{B}_{\mu,0,\tts})_{26}$, $(\mathsf{B}_{\mu,0,\tts})_{36}$, $(\mathsf{B}_{\mu,0,\tts})_{46}$, are respectively 
 \begin{align*}
 \left( \mathcal{B}_{\mu, 0, \tts} f_0^{-},\ f_{2}^{+}(\mu,0) \right) = 
 \left( \mathcal{B}_{\mu, 0, \tts} f_0^{-},\ f_{2}^{-}(\mu,0) \right) = 
 \left( \mathcal{B}_{\mu, 0, \tts} f_0^{-},\ f_{1}^{+}(\mu,0) \right) = 
 \left( \mathcal{B}_{\mu, 0, \tts} f_0^{-},\ f_{1}^{-}(\mu,0) \right) = 0  \ . 
 \end{align*}
Hence  their analytic remainders, so far of order
$r(\mu \e^3, \mu^2)$,  must vanish at $\e =0$, hence they must be of the form $r(\mu \e^3, \mu^2\e)$.
\end{proof}

\section{Block-Decoupling}\label{sec:block-decoupling}
The goal of this section is to conjugate the Hamiltonian
matrix $\mathsf{L}_{\mu,\e,\tts}
=\mathsf{J}_6\mathsf{B}_{\mu,\e,\tts}$
of Proposition \ref{B decomposition} to a
block-diagonal matrix with three $2\times2$ blocks,
for $(\mu,\e)$ in the parabolic region $\cW$
defined in \eqref{def:W}.
To this end, we introduce the rescaled parameters
\begin{equation}\label{eq:sylvester-band-new}
    \e=\ve^2,\qquad
    \mu=(\tau\ve)^4,\qquad
    (\tau,\ve)\in[1,2^{1/4}]\times(0,\ve_0),
    \quad \ve_0:=\e_0^{1/2}.
\end{equation}
This change of variables is adapted to the scaling
$\mu\sim\e^2$ and maps the distinguished curve
$\mu=\frac98\e^2$ to the fixed value
$\tau=(9/8)^{1/4}$.

\medskip 
$\bullet $ {\bf Notation:} From now on we denote by  $r(\ve^a)$  real-valued, real-analytic functions in $(\tau, \ve)$, fulfilling the estimate
$|r(\ve^a)|\leq C |\ve^a|$ uniformly for $\tau \in [1,2^{1/4}]$.

\medskip
From Proposition \ref{B decomposition} we obtain immediately
\begin{lemma}\label{lem:L0}
    In the variables $(\tau, \ve)$ of \eqref{eq:sylvester-band-new}, one has 
    \begin{equation}\label{L^0}
         \mathsf L^{(0)}_{\tau,\ve,\tts}:=  \mathsf L_{(\tau\ve)^4, \ve^2,\tts}  
         = 
      \mathsf{J}_6    \mathsf B^{(0)}_{\tau,\ve,\tts} =
         \mathsf{J}_6 \begin{pmatrix}
E^{(0)} & I^{(0)} & F^{(0)} \\
(I^{(0)})^* & H^{(0)} & L^{(0)} \\
(F^{(0)})^* & (L^{(0)})^* & G^{(0)}
\end{pmatrix} \ , 
    \end{equation}
    where $E^{(0)}, H^{(0)}, G^{(0)}$ are  selfadjoint and reversibility preserving, 
$F^{(0)}, I^{(0)}, L^{(0)}$ reversibility preserving, and
\small
\begin{equation*}
E^{(0)}=
\begin{pmatrix}
\mathtt e_1\ve^2+\mathtt e_2\ve^4+r(\ve^6)
&\im\mathtt e_4\tau^4\ve^4+\im r(\ve^6)\\[3pt]
-\im\mathtt e_4\tau^4\ve^4-\im r(\ve^6)
&\mathtt e_1\ve^2+\mathtt e_6\ve^4+r(\ve^6)
\end{pmatrix} \ , 
\quad 
F^{(0)}=
\begin{pmatrix}
\mathtt f_1\ve^6+r(\ve^8)
&\im\tau^4\ve^6\bigl(\mathtt f_2+\mathtt f_3\ve^2\bigr)
+\im {r(\ve^{10})}\\[3pt]
\im\mathtt f_4\tau^4\ve^6+\im r(\ve^8)
&
\mathtt{f}_5 \tau^4\ve^8+ {r(\ve^{10})}
\end{pmatrix} \ ,
\end{equation*}

\begin{equation*}
G^{(0)}=
\begin{pmatrix}
1+r(\ve^4)&\im r(\ve^6)\\[3pt]
-\im r(\ve^6)&\tau^4\ve^4+r(\ve^{8})
\end{pmatrix}, 
\qquad 
H^{(0)}=
\begin{pmatrix}
\mathtt h_1\ve^4+r(\ve^6)
&\im\mathtt h_3\tau^4\ve^4+\im r(\ve^6)\\[3pt]
-\im\mathtt h_3\tau^4\ve^4-\im r(\ve^6)
&\mathtt e_1\ve^2+\mathtt h_5\ve^4+r(\ve^6)
\end{pmatrix}, 
\end{equation*}

\begin{equation*}
I^{(0)}=
\begin{pmatrix}
\mathtt i_1\ve^2+\mathtt i_2\ve^4+r(\ve^6)
&\im r(\ve^6)\\[3pt]
\im r(\ve^6)
&\mathtt i_1\ve^2+\mathtt i_4\ve^4+r(\ve^6)
\end{pmatrix}, 
\quad 
L^{(0)}=
\begin{pmatrix}
\mathtt l_1\ve^6+r(\ve^8)
&\im\tau^4 \ve^6\bigl(\mathtt l_2+\mathtt l_3\ve^2\bigr)
+\im {r(\ve^{10})}\\[3pt]
\im\mathtt l_4\tau^4\ve^6+\im r(\ve^8)
&{r(\ve^{10})}
\end{pmatrix},
\end{equation*}
\normalsize
with the coefficients defined in \eqref{coeff.B}.
\end{lemma}

We now perform a series of transformations to block diagonalize 
$\mathsf L^{(0)}_{\tau,\ve,\tts}$. 
The first one is a reordering of the basis elements.

\begin{lemma}
\label{lem:first-symplectic-change}
There exists $\ve_0>0$ such that for any $(\tau, \ve) \in [1,2^{1/4}]\times(0,\ve_0)$, the singular, 
symplectic and reversibility preserving matrix 
\begin{equation}\label{eq:Y-first-change}
Y_{\tts}(\tau, \ve) :=\frac12
\begin{pmatrix}
-{2^{-\frac14}\ck^{-\frac12}}
&-\im{2^{-\frac14}\ck^{-\frac12}}
&-{2^{-\frac14}\ck^{-\frac12}}
&\im{2^{-\frac14}\ck^{-\frac12}}
&{\ck^{-1}}&0
\\[2pt]
\im 2^{\frac14}{\ck^{-\frac32}}
&-{2^{\frac14}}{\ck^{-\frac32}}
&-{2^{\frac14}\im}{\ck^{-\frac32}}
&-{2^{\frac14}}{\ck^{-\frac32}}
&0&{2}{\ck^{-1}}
\\[2pt]
\tts{ 2^{-\frac14}\ck^{-\frac12}}
&{\tts \im}{2^{-\frac14}\ck^{-\frac12}}
&{\tts}{2^{-\frac14}\ck^{-\frac12}}
&-{\tts\im }{2^{-\frac14}\ck^{-\frac12}}
&2\tts\ck^{-1}&0
\\[2pt]
-{\im 2^{-\frac34}\tts}{\ck^{-\frac32}}
&{2^{-\frac34}\tts}{\ck^{-\frac32}}
&{\im 2^{-\frac34}\tts}{\ck^{-\frac32}}
&{2^{-\frac34}\tts}{\ck^{-\frac32}}
&0&{2\tts}{\ck^{-1}}
\\[2pt]
\tau \ve&-\im \tau \ve 
&\tau \ve&\im \tau \ve&0&0
\\[2pt]
\im (\tau \ve)^{-1}&(\tau \ve)^{-1}
&-\im (\tau \ve)^{-1}&(\tau \ve)^{-1}&0&0
\end{pmatrix} \ , \quad \ve > 0 \ . 
\end{equation}
 conjugates the 
the $6\times 6$ matrix 
$\mathsf L_{\tau,\ve,\tts}^{(0)}$ in \eqref{L^0} to the 
$6 \times 6$ Hamiltonian and reversible matrix
\begin{equation}\label{eq:L1-first-change}
\begin{aligned}
\mathsf L^{(1)}_{\tau,\ve,\tts}
&:=Y_{\tts}\bigl(\tau,\ve \bigr)^{-1}
\mathsf L_{\tau, \ve,\tts}^{(0)}
Y_{\tts}\bigl(\tau,\ve\bigr)=\mathsf J_6\mathsf B^{(1)}_{\tau,\ve,\tts} 
= \mathsf J_6
\begin{pmatrix}
E^{(1)}&I^{(1)}&F^{(1)}\\
[I^{(1)}]^*&H^{(1)}&L^{(1)}\\
[F^{(1)}]^*&[L^{(1)}]^*&G^{(1)}
\end{pmatrix}
\end{aligned}
\end{equation}
which admits a real-analytic extension to $\ve=0$.
Moreover, the $2\times 2$ selfadjoint and reversibility preserving  matrices  $E^{(1)}, H^{(1)}, G^{(1)}$ 
expand as
\begin{equation}\label{eq:E1-first-change}
\begin{aligned}
E^{(1)} &= \ve^2
\begin{pmatrix}
{\widehat\alpha_{\tts}}+\gamma_E\ve^2+\mathtt q_E^{(1)}\ve^4-\delta_E^{(1)}\ve^5+ r(\ve^6)
&
\im\left(-{\widehat\beta_{\tts}}+\gamma_E\ve^2+ r(\ve^4)\right)\\[1ex]
-\im\left(-{\widehat\beta_{\tts}}+\gamma_E\ve^2+ r(\ve^4)\right)
&
{\widehat\alpha_{\tts}}+\gamma_E\ve^2+\mathtt q_E^{(1)}\ve^4+\delta_E^{(1)}\ve^5+r(\ve^6)
\end{pmatrix} \ , 
\end{aligned}
\end{equation}
\begin{equation}\label{eq:H1-first-change}
\begin{aligned}
H^{(1)} &=
\ve^2
\begin{pmatrix}
{\widehat\alpha_{\tts}}+\gamma_H\ve^2+\mathtt q_H^{(1)}\ve^4+\delta_{H}^{(1)}\ve^5 + r(\ve^6)
&
\im\left({\widehat\beta_{\tts}}-\gamma_H\ve^2+ r(\ve^4)\right) \\[1ex]
- \im\left({\widehat\beta_{\tts}}-\gamma_H\ve^2+ r(\ve^4)\right)
&
{\widehat\alpha_{\tts}}+\gamma_H\ve^2+\mathtt q_H^{(1)}\ve^4-\delta_{H}^{(1)}\ve^5 + r(\ve^6)
\end{pmatrix}\  , 
\end{aligned}
\end{equation}
\begin{equation}\label{eq:G1-first-change}
G^{(1)}
=
\ve^2
\begin{pmatrix}
\dfrac{\ck}{4}\tts+r(\ve^2)
&
\im r(\ve^2)
\\[6pt]
-\im r(\ve^2)
&
-\dfrac{3}{16}\ck\tts\ve^4+r(\ve^6)
\end{pmatrix} \ , 
\end{equation}
whereas the $2\times 2$ reversibility preserving matrices
$I^{(1)}, F^{(1)}, L^{(1)}$ as
\begin{equation}\label{eq:Lblock1-first-change}
\begin{aligned}
F^{(1)} &=
\ve^4
\begin{pmatrix}
\varphi_1-\psi_5\ve + r(\ve^2) & \im\,\varphi_2 + \im r(\ve^2)\\[1ex]
-\im\left(\varphi_1+\psi_5\ve\right) + \im r(\ve^2) & \varphi_2 + r(\ve^2)
\end{pmatrix}, \quad 
L^{(1)} =
\ve^4
\begin{pmatrix}
\phi_1+\psi_5\ve + r(\ve^2)& -\im\,\phi_2 + \im r(\ve^2)\\[1ex]
\im\,\phi_1-\im\,\psi_5\ve + \im r(\ve^2) & \phi_2+ r(\ve^2)
\end{pmatrix}
,\\
I^{(1)} & =
\ve^4
\begin{pmatrix}
\iota_I^{(1)} + r(\ve^2) & -\im\,\iota_I^{(1)} + \im r(\ve^2)\\[1ex]
-\im\,\iota_I^{(1)} + \im r(\ve^2)& -\iota_I^{(1)} + r(\ve^2)
\end{pmatrix},
\end{aligned}
\end{equation}
where $\gamma_E, \gamma_H$ in \eqref{final.coeff}, and
\small
\begin{equation}
\begin{aligned}
{\widehat\alpha_{\tts}(\tau)}&=\frac12\left(\tau^2-\tts\sqrt{\frac98}\right),
\qquad
{\widehat\beta_{\tts}(\tau)
:=\frac12\left(\tau^2+\tts\sqrt{\frac98}\right),}\\
\delta_E^{(1)}(\tau)
&=\frac{2^{1/4}\tau}{192}
 \Big[(36\sqrt2-8\ck)\tau^2+15\sqrt2\ck\,\mathtt s\Big], \qquad 
\delta_H^{(1)}(\tau)
=\frac{2^{1/4}\tau}{192}
 \Big[(36\sqrt2+8\ck)\tau^2-15\sqrt2\ck\,\mathtt s\Big],
 \\[1ex] 
\varphi_1(\tau)&=2^{1/4}\ck^{1/2}
 \left(\frac{7\sqrt2\ck}{32}-\frac{\tau^4}{18}\right), 
 \qquad
 \varphi_2(\tau)=2^{1/4}\ck^{1/2}
 \left(\frac3{16}-\frac{\sqrt2\ck}{18}\tau^4\right), 
 \qquad 
 \psi_5(\tau)=\frac{\ck^{3/2}}2\,\mathtt s\,\tau^3  , \\  
\phi_1(\tau)&=2^{1/4}\ck^{1/2}
 \left(\frac{7\sqrt2\ck}{32}+\frac{\tau^4}{18}\right), 
 \qquad 
\phi_2(\tau)=2^{1/4}\ck^{1/2}
 \left(\frac3{16}+\frac{\sqrt2\ck}{18}\tau^4\right) \ , \qquad 
\iota_I^{(1)}:=\frac{57\sqrt2}{128}.
\label{eq:coeff-first-change}
\end{aligned}
\end{equation}
\normalsize
and $\mathtt{q}_{E}^{(1)}$,  $\mathtt{q}_H^{(1)}$, $ \iota_I^{(1)}$ are real analytic functions of $\tau$ for each fixed $\tts$. 
In particular, $(\widehat\alpha_1,\widehat\beta_1)=(\alpha,\beta)$ and  $(\widehat\alpha_{-1},\widehat\beta_{-1})=(\beta,\alpha)$,
where $\alpha,\beta$ are defined in \eqref{final.coeff}.
\end{lemma}

\begin{proof}
A direct computation shows that 
$Y_\tts^*\mathsf J_6Y_\tts=\mathsf J_6$, and that $Y_\tts$ preserves reversibility. Hence
$\mathsf B^{(1)}_{\tau, \ve, \tts}=Y_\tts^*\mathsf B^{(0)}_{\tau, \ve, \tts}Y_\tts$
is self-adjoint and reversibility preserving, and the conjugation has
the stated Hamiltonian form.

We prove that the sixth-order diagonal coefficients of $E^{(1)}$ and $H^{(1)}$ within each block coincide without computing them. The matrix $\mathsf B^{(0)}_{\tau, \ve, \tts}$ in \eqref{L^0} is even in $\ve$, since it is obtained by 
evaluating  $ \mathsf{B}_{\mu,\e,\tts}$  in \eqref{B mu epsilon} at  $(\mu,\e)=(\tau^4\ve^4,\ve^2)$.
Note that 
\[
 Y_\tts(\tau,-\ve)=Y_\tts(\tau,\ve)\mathsf Q,\qquad
 \mathsf Q:=\operatorname{diag}(\im\mathsf J_2,-\im\mathsf J_2,I_2), 
\]
consequently, $\mathsf B^{(1)}_{\tau,-\ve, \tts}
 =\mathsf Q^*\mathsf B^{(1)}_{\tau,\ve, \tts}\mathsf Q $ and thus
\begin{equation}\label{eq:first-change-exchange-symmetry}
X^{(1)}(\tau, -\ve) =  (\im\mathsf J_2)^*X^{(1)}(\tau, \ve)  (\im\mathsf J_2) \ , \quad X \in \{E, H\}.
\end{equation}
Since for a self-adjoint reversible $2\times2$ matrix,
$ (\im\mathsf J_2)^*
 \begin{pmatrix}a&\im b\\-\im b&d\end{pmatrix}
 (\im\mathsf J_2)
 =\begin{pmatrix}d&\im b\\-\im b&a\end{pmatrix}$, 
it follows  that
\begin{equation}\label{eq:first-change-even-diagonals}
 X_{11}^{(1)}(\tau,-\ve)=X_{22}^{(1)}(\tau,\ve),\qquad
 X_{12}^{(1)}(\tau,-\ve)=X_{12}^{(1)}(\tau,\ve) \ , \quad X \in \{E, H\} \ . 
\end{equation}
Thus every even Taylor coefficient of the  diagonal entries of
each block coincide, so in particular at order 6.

The remaining  coefficients follow by direct computation. In particular, in terms of the coefficients in \eqref{coeff.B},
\small
\[
 \begin{aligned}
 \gamma_{E}={}&
 \frac{\mathtt e_2+\mathtt h_1-2\tts\mathtt i_2}{4\sqrt2\ck}
 +\frac{4\mathtt e_6+\mathtt h_5-4\tts\mathtt i_4}{8\sqrt2\ck^3}
+\frac{2\mathtt e_4+\mathtt h_3}{4\ck^2}\tau^4 , 
 \quad 
 \gamma_{H}={}
 \frac{\mathtt e_2+\mathtt h_1-2\tts\mathtt i_2}{4\sqrt2\ck}
 +\frac{4\mathtt e_6+\mathtt h_5-4\tts\mathtt i_4}{8\sqrt2\ck^3}
-\frac{2\mathtt e_4+\mathtt h_3}{4\ck^2}\tau^4 ,
 \end{aligned}
\]
\normalsize
and 
\small
\[
 \begin{aligned}
 \delta_E^{(1)} = -  W_+ \ , \quad 
 \delta_H^{(1)} =   W_-  \ , \quad 
 W_\sigma=\left[
 -\frac{(\mathtt f_1-\tts\mathtt l_1)\tau
       +\sigma(\tts\mathtt l_3-\mathtt f_3)\tau^3}{2^{5/4}\ck^{1/2}}
 +\frac{\sigma(2\mathtt f_4-\tts\mathtt l_4)\tau^5}{2^{7/4}\ck^{3/2}}
 +\frac{\mathtt f_5\tau^3}{2^{3/4}\ck^{3/2}}
 \right].
 \end{aligned}
\]
\normalsize
The quintic terms in the diagonal of $E^{(1)}$, $H^{(1)}$ and in the out-diagonal of $I^{(1)}$ vanish due to an exact cancellation, indeed they have respectively the form 
$$
[\ve^5]E^{(1)}=\zeta_5 \begin{pmatrix}1& 0\\ 0&-1\end{pmatrix},
\qquad
[\ve^5]H^{(1)}=\zeta_5 
\begin{pmatrix}-1&0\\ 0 & 1 \end{pmatrix}
\ , \qquad 
[\ve^5]I^{(1)}  = \zeta_5
\begin{pmatrix}0&-\im\\ \im&0\end{pmatrix}, \, \quad 
$$
with
$$
\zeta_5
=\frac{2^{-1/4}\ck^{3/2}}{3}\,
\tau^3(\mathtt f_2-\mathtt s\mathtt l_2) \stackrel{\eqref{coeff.B}}{=} 0.
$$
 Finally 
\small
\[
\begin{aligned}
\varphi_1
&=2^{1/4}\ck^{1/2}
\left[
\frac{-\mathtt e_2-\tts\mathtt i_2+2\mathtt h_1}{6\sqrt2}
+\frac{\ck(\mathtt h_3-\mathtt e_4)}9\,\tau^4
\right], \quad 
\varphi_2
=2^{1/4}\ck^{1/2}
\left[
\frac{\ck(-2\mathtt e_6-\tts\mathtt i_4+\mathtt h_5)}9
+\frac{\mathtt h_3-\mathtt e_4}{3\sqrt2}\,\tau^4
\right], \\
\phi_1
&=2^{1/4}\ck^{1/2}
\left[
\frac{-\mathtt e_2-\tts\mathtt i_2+2\mathtt h_1}{6\sqrt2}
-\frac{\ck(\mathtt h_3-\mathtt e_4)}9\,\tau^4
\right],\quad 
\phi_2
=2^{1/4}\ck^{1/2}
\left[
\frac{\ck(-2\mathtt e_6-\tts\mathtt i_4+\mathtt h_5)}9
-\frac{\mathtt h_3-\mathtt e_4}{3\sqrt2}\,\tau^4
\right] \ , \\
 \psi_5& =\frac{\mathtt f_2+2\tts\mathtt l_2}{4\ck}\tau^3 , \qquad 
 \iota_I^{(1)}
=\frac{\ck}{18\sqrt2}
\left[
3(\mathtt e_2+\mathtt h_1-2\tts \mathtt i_2)
-(4\mathtt e_6+\mathtt h_5-4\tts \mathtt i_4)
\right] \ . 
\end{aligned}
\]
\normalsize
Substituting the values of \eqref{coeff.B} gives the explicit formulas of  \eqref{final.coeff}, \eqref{eq:coeff-first-change}.

Finally consider 
$ G_{22}^{(1)}=\ck^{-2}
 \big(E_{22}+2\tts I_{22}+H_{22}\big).$ 
The identity $E_{22}H_{22}=I_{22}^2$ at $\mu = 0$ in \eqref{lemma C3} yields
\begin{equation}\label{1811:1652}
 E_{22}(0, \e)+2\tts I_{22}(0, \e)+H_{22}(0, \e)
 =\frac{(E_{22}(0, \e)+\tts I_{22}(0, \e))^2}{E_{22}(0, \e)}
\stackrel{\eqref{E}, \eqref{I}}{=}\frac{(\mathtt e_6+\tts \mathtt i_4)^2}{\mathtt e_1}\e^3
 +r(\e^4).
\end{equation}
Then, using \eqref{1811:1652} and the expansions of Proposition \ref{B decomposition}, 
\begin{align*}
     G_{22}^{(1)}(\mu, \e) & =   G_{22}^{(1)}(0, \e) + \mu \partial_\mu G_{22}^{(1)}(0, \e) + r(\mu^2)  = 
 \frac{(\mathtt e_6+\tts \mathtt i_4)^2}{\ck^2 \mathtt e_1}\e^3
 +r(\e^4, \mu \e^2, \mu^2)
\end{align*}
where we used that $\partial_\mu G_{22}^{(1)}(0, \e) = r(\e^2)$ in view of the expansions in Proposition \ref{B decomposition}. 
Thus, in the variables $\mu = \tau^4 \ve^4$, $\epsilon = \ve^2$, one gets
$$
 G_{22}^{(1)}(\tau^4 \ve^4, \ve^2)
 =\frac{(\mathtt e_6+\tts\mathtt i_4)^2}{\ck^2\mathtt e_1}\ve^6
 +r(\ve^8)
 =-\frac{3\tts \ck}{16}\ve^6+r(\ve^8).
$$
The remaining entries of $G^{(1)}$ follow from a direct computation.
\end{proof}

\begin{remark}
The goal of the first transformation is to reorganize and rescale the columns and rows of the matrix $\mathsf L^{(0)}_{\tau, \ve, \tts}$ so that the diagonal blocks of $\mathsf L^{(1)}_{\tau, \ve, \tts}$ have separated spectrum for $\tts=1$. Indeed, at principal order 
\[
 \begin{aligned}
 \sigma(\mathsf J_2E^{(1)})
 &=\{\im\tau^2\ve^2+\cO(\ve^4),\;
    \im{\tau^2_0}\ve^2+\cO(\ve^4)\},\quad  \sigma(\mathsf J_2G^{(1)})=\cO(\ve^4),\\
 \sigma(\mathsf J_2H^{(1)})
&  =\{-\im\tau^2\ve^2+\cO(\ve^4),\;
       -\im{\tau^2_0}\ve^2+\cO(\ve^4)\}, \quad {\tau_0}:= \left(\frac98 \right)^{\frac14} \ . 
 \end{aligned}
\]
Thus the three spectra are pairwise separated, including at $\tau=\tau_0$, where each critical
block has a double leading eigenvalue. This is the separation needed
for the subsequent block decoupling.
For $\tts=-1$, the leading spectra of the $\mathsf J_2E^{(1)}$ and $\mathsf J_2H^{(1)}$ blocks coincide
at $\tau_0$; these blocks are therefore kept together in the
$4\times4$ part of the reduction.
 \end{remark}

\begin{remark}[The four eigenvalues of order $\ve^2$ for $\tts=-1$]
\label{rem:-s}
Fix $\rho>0$. We claim that on the regime \eqref{eq:sylvester-band-new}, every eigenvalue $\lambda$ of $\mathsf L^{(1)}_{\tau,\ve,-1}$ with
$|\lambda|\geq \rho\ve^2$ is purely imaginary for sufficiently
small $\ve>0$. Indeed, we write the matrix in \eqref{eq:L1-first-change} as
$$
\mathsf B^{(1)}
=\begin{pmatrix} A&C\\ C^*&G^{(1)} \end{pmatrix},
\qquad
C=\begin{pmatrix} F^{(1)}\\ L^{(1)} \end{pmatrix}.
$$
By \eqref{eq:E1-first-change}--\eqref{eq:coeff-first-change},
$\ve^{-2}A=A_0+\cO(\ve^2)$ with $A_0\geq I_4$, since the
eigenvalues of $A_0$ are $\tau^2$ and $\sqrt{9/8}$, each of multiplicity two.
Moreover, $C=\cO(\ve^4)$ and
$\ve^{-2}\mathsf J_2G^{(1)}=N+\cO(\ve^2)$, where
$N=\mathsf J_2\operatorname{diag}(-\ck/4,0)$ satisfies $N^2=0$.

Let $z=\ve^{-2}\lambda$ and let $w=(x,y)$ be a corresponding
eigenvector. The lower block equation is
\[
(zI_2-\ve^{-2}\mathsf J_2G^{(1)})y
=\ve^{-2}\mathsf J_2C^*x.
\]
Since $(zI_2-N)^{-1}=z^{-1}I_2+z^{-2}N$ is uniformly bounded for
$|z|\geq\rho$, the inverse on the left is uniformly bounded for
small $\ve$. Hence $\|y\|\lesssim\ve^2\|x\|$ and
\[
\ve^{-2}w^*\mathsf B^{(1)}w
=x^*A_0x+\cO(\ve^2)\|x\|^2>0.
\]
The Hamiltonian identity then implies
$2\operatorname{Re}\lambda\,w^*\mathsf B^{(1)}w=0$.
Thus $\operatorname{Re}\lambda=0$.
\end{remark}

\subsection{Non-perturbative Step of Block-Decoupling}
\label{subsec:first-block-decoupling}

The goal of this section is to perform a first
block-decoupling step that reduces the order of the
off-diagonal blocks of the matrix
$\mathsf B^{(1)}_{\tau,\ve,\tts}$ in
\eqref{eq:L1-first-change}.
The resulting second-order corrections to the
diagonal blocks $E^{(1)}, H^{(1)}, G^{(1)}$
play a crucial role in determining the
eigenvalue asymptotics.

\medskip
\noindent
{\bf From now on, we focus on the positive Wilton
branch $\tts=1$.}
The modifications required for the negative branch
$\tts=-1$ are discussed in Appendix \ref{app:s-1}.
The two branches must be treated separately because,
for $\tts=1$, the block-decoupling amounts to solving
three independent Sylvester equations, whereas for
$\tts=-1$ only two of these equations are solvable.
Consequently, in the negative branch, the transformed
matrix $\mathsf B^{(2)}_{\tau,\ve,-1}$ retains a
$4\times4$ block coupled at leading order, alongside
a separate $2\times2$ block.

\begin{proposition}[{\bf Non perturbative block-decoupling, $\tts =1$}] \label{prop:first-block-decoupling}
There exists $\ve_1>0$ such that for any $(\tau, \ve) \in [1,2^{1/4}]\times(0,\ve_1)$,
the following holds true. 
There exist three $2\times2$
reversibility-preserving matrices
\begin{equation}\label{eq:three-X-form}
 X_k := X_{k,1} := \begin{pmatrix}x_{11}^{(k)}&\im x_{12}^{(k)}\\
                   \im x_{21}^{(k)}&x_{22}^{(k)}\end{pmatrix},
 \qquad x_{ij}^{(k)}\in\R,\quad k=1,2,3,
\end{equation}
real analytic in  $(\tau,\ve)$ and expanding as
\begin{equation}
\begin{aligned}
\label{X1}
X_1 = 
\begin{pmatrix}
\fd \ve^2  + r(\ve^4)& -\im \fd \ve^2 + \im  r(\ve^4)\\
-\im \fd\ve^2 + \im r(\ve^4) & - \fd \ve^2+ r(\ve^4)
\end{pmatrix} \ , \quad 
X_2 &=
\begin{pmatrix}
\fa_2 \ve^2+ \fc \ve^3  + r(\ve^4)
&
\im \fb_2 \, \ve^2 
+ \im r(\ve^4)
\\
- \im \fa_2 \ve^2 + \im \fc \ve^3+\im r(\ve^4)
&
\fb_2 \ve^2 + r(\ve^4) 
\end{pmatrix} \ , 
\\[1ex]
X_3&=
 \begin{pmatrix}\fa_3 \ve^2 -\fc\ve^3 + r(\ve^4) &-\im \fb_3 \ve^2 + \im r(\ve^4)\\
 \im \fa_3 \ve^2 + \im \fc \ve^3 + \im r(\ve^4) & \fb_3\ve^2 + r(\ve^4)
 \end{pmatrix},
\end{aligned}
\end{equation}
where 
\begin{equation}\label{X.coeff}
\begin{aligned}
\fa_2&=2^{1/4}\ck^{1/2}
\left(\frac{\ck}{4}-\frac{\sqrt2}{54}\tau^4\right),
&\fb_2&=2^{1/4}\ck^{1/2}
\left(\frac{\sqrt2}{8}-\frac{2\ck}{27}\tau^4\right),
& \fd &=\frac{19}{64},
\\[1ex]
\fa_3&=2^{1/4}\ck^{1/2}
\left(\frac{\ck}{4}+\frac{\sqrt2}{54}\tau^4\right),
&\fb_3&=2^{1/4}\ck^{1/2}
\left(\frac{\sqrt2}{8}+\frac{2\ck}{27}\tau^4\right), 
&\fc &=\frac{\ck^{3/2}}{2}\tau , 
\end{aligned}
\end{equation}
such that by conjugating the 
Hamiltonian and reversible
matrix $\mathsf L^{(1)}_{\tau,\ve,1}$, defined in \eqref{eq:L1-first-change}, with the symplectic and reversibility-preserving $6\times 6$
matrix
\begin{equation}\label{eq:S1-three-blocks}
\exp(S_{1}^{(1)}) \quad \mbox{ where } \quad 
 S_1^{(1)} =\mathsf J_6
 \begin{pmatrix}
 0&M_1&M_2\\M_1^*&0&M_3\\M_2^*&M_3^*&0
 \end{pmatrix},\qquad M_k:=\mathsf J_2X_k,\quad k=1,2,3,
\end{equation}
we get the Hamiltonian and reversible
matrix 
\begin{align}\label{L2}
 \mathsf L^{(2)}_{\tau,\ve,1}
 &:=\exp(S_{1}^{(1)})\mathsf L^{(1)}_{\tau,\ve,1}\exp(-S_{1}^{(1)})
 =\mathsf J_6\mathsf B^{(2)}_{\tau,\ve,1} , \quad  \mathsf B^{(2)}_{\tau,\ve,1} = \begin{pmatrix}
 E^{(2)}&I^{(2)}&F^{(2)}\\
 [I^{(2)}]^*&H^{(2)}&L^{(2)}\\
 [F^{(2)}]^*&[L^{(2)}]^*&G^{(2)}
 \end{pmatrix}
 ,
 \end{align}
where the diagonal blocks are self-adjoint and expand as
\begin{align}
E^{(2)} &=
\ve^2
\begin{pmatrix}
\alpha+\gamma_E\ve^2+ \mathtt q_E^{(2)} \ve^4 + \delta_{E}\ve^5 + r(\ve^6)
& \im\left(-\beta+\gamma_E\ve^2 + r(\ve^4) \right) \\[1ex]
-\im\left(-\beta+\gamma_E\ve^2 +  r(\ve^4) \right)
& \alpha+\gamma_E\ve^2+ \mathtt q_E^{(2)} \ve^4 -\delta_{E}\ve^5 +r(\ve^6)
\end{pmatrix}\ ,
\label{eq:E2-final}\\[2ex]
H^{(2)} &=\ve^2
\begin{pmatrix}
\alpha +\gamma_H\ve^2+ \mathtt q_H^{(2)} \ve^4 +\delta_{H}\ve^5 + r(\ve^6)
& \im\big(\beta-\gamma_H\ve^2 +  r(\ve^4) \big)\\[1ex]
-\im\big(\beta-\gamma_H\ve^2   + r(\ve^4) \big)
& \alpha+\gamma_H\ve^2+ \mathtt q_H^{(2)} \ve^4 -\delta_{H}\ve^5 + r(\ve^6)
\end{pmatrix}\ ,
\label{eq:H2-final}\\[2ex]
G^{(2)} &=
\ve^2
\begin{pmatrix}
\frac{\ck}{4}+ r(\ve^2) & \im r( \ve^2) \\[1ex]
-\im r(\ve^2) & \frac{4\ck}{81}\tau^8\ve^4 + r(\ve^6)
\end{pmatrix}\ ,
\label{eq:G2-final}
\end{align}
where $\alpha, \beta, \gamma_E, \gamma_H, \delta_E, \delta_H$ as in \eqref{final.coeff}
and $\mathtt{q}_{E}^{(2)}$,  $\mathtt{q}_H^{(2)}$, are real analytic functions of $\tau$ for each fixed $\tts$.

The off-diagonal blocks are
\begin{align*}
    I^{(2)},F^{(2)},L^{(2)}=\begin{pmatrix}
        r(\ve^6) & \im r(\ve^6) \\
         \im r(\ve^6) & r(\ve^6)
    \end{pmatrix}.
\end{align*}
\end{proposition}

The rest of the subsection is devoted to the proof.
For simplicity, we write $S_{1}:=S_{1}^{(1)}$. \\
We split the matrix
\eqref{eq:L1-first-change} into its $2\times2$ block-diagonal and
off-diagonal parts:
\begin{align}
\label{eq:D1-R1-three-blocks}
 \mathsf L^{(1)}_{\tau,\ve,1}&=D^{(1)}+R^{(1)},\\
 \notag
 D^{(1)}&:=\operatorname{diag}(D_1,D_2,D_0)
 =\operatorname{diag}(\mathsf J_2E^{(1)},
                     \mathsf J_2H^{(1)},\mathsf J_2G^{(1)}),
                     \quad 
 R^{(1)}:=\mathsf J_6
 \begin{pmatrix}0&I^{(1)}&F^{(1)}\\
 [I^{(1)}]^*&0&L^{(1)}\\
 [F^{(1)}]^*&[L^{(1)}]^*&0\end{pmatrix}.
\end{align}
We expand the conjugation in Lie series, with $\ad_S(A):=[S,A] = SA - AS$:
\begin{align}\label{eq:lie-expansion-three-blocks}
 \mathsf L^{(2)}_{\tau,\ve,1} = \exp(S_1)\mathsf L^{(1)}_{\tau,\ve,\tts} & \exp(-S_1)
 ={}D^{(1)}+[S_1,D^{(1)}]+\frac12[S_1,[S_1,D^{(1)}]]
       +R^{(1)}+[S_1,R^{(1)}]\\
       \notag
&  +\frac12\int_0^1(1-t)^2
 e^{t S_1}\ad_{S_1}^3(D^{(1)})e^{-t S_1}\,\de t
 +\int_0^1(1-t)
 e^{t S_1}\ad_{S_1}^2(R^{(1)})e^{-t S_1}\,\de t.
\end{align}
We choose $S_1$ by solving the homological equation
\begin{equation}\label{eq:homological-three-blocks}
 R^{(1)}+[S_1,D^{(1)}]=0
\end{equation}
which, in view of \eqref{eq:D1-R1-three-blocks}, \eqref{eq:S1-three-blocks}, read as
\small
\begin{equation*}
\mathsf J_6
\begin{pmatrix}
0 & I^{(1)}+M_1\mathsf J_2H^{(1)}-E^{(1)}\,\mathsf J_2M_1 & F^{(1)}+M_2\mathsf J_2G^{(1)}-E^{(1)}\,\mathsf J_2M_2\\[1.5ex]
(I^{(1)})^*+M_1^*\mathsf J_2E^{(1)}-H^{(1)}\,\mathsf J_2M_1^* & 0 & L^{(1)}+M_3\mathsf J_2G^{(1)}-H^{(1)}\,\mathsf J_2M_3\\[1.5ex]
(F^{(1)})^*+M_2^*\mathsf J_2E^{(1)}-G^{(1)}\,\mathsf J_2M_2^* & (L^{(1)})^*+M_3^*\mathsf J_2H^{(1)}-G^{(1)}\,\mathsf J_2M_3^* & 0
\end{pmatrix} = 0.
\end{equation*}
\normalsize
Note that the vanishing of the three equations in the upper off-diagonal blocks implies the vanishing of the 
 lower off-diagonal blocks and viceversa by  self-adjointness.
Since $M_k=\mathsf J_2X_k$ (cf. \eqref{eq:S1-three-blocks}), we obtain the three decoupled
Sylvester equations
\begin{equation}\label{eq:three-sylvester-equations}
\begin{aligned}
 D_1X_1-X_1D_2=-\mathsf J_2I^{(1)},\qquad D_1X_2-X_2D_0=-\mathsf J_2F^{(1)},\qquad
 D_2X_3-X_3D_0=-\mathsf J_2L^{(1)}, 
\end{aligned}
\end{equation}
which we now solve  explicitly.
Write
\begin{equation}\label{eq:real-entries-three-blocks}
\begin{gathered}
 E^{(1)}=\begin{pmatrix}E_{11}^{(1)}&\im E_{12}^{(1)}\\
 -\im E_{12}^{(1)}&E_{22}^{(1)}\end{pmatrix},\quad
 H^{(1)}=\begin{pmatrix}H_{11}^{(1)}&\im H_{12}^{(1)}\\
 -\im H_{12}^{(1)}&H_{22}^{(1)}\end{pmatrix},\quad
 G^{(1)}=\begin{pmatrix}G_{11}^{(1)}&\im G_{12}^{(1)}\\
 -\im G_{12}^{(1)}&G_{22}^{(1)}\end{pmatrix},\\
 I^{(1)}=\begin{pmatrix}I_{11}^{(1)}&\im I_{12}^{(1)}\\
 \im I_{21}^{(1)}&I_{22}^{(1)}\end{pmatrix},\quad
 F^{(1)}=\begin{pmatrix}F_{11}^{(1)}&\im F_{12}^{(1)}\\
 \im F_{21}^{(1)}&F_{22}^{(1)}\end{pmatrix},\quad
 L^{(1)}=\begin{pmatrix}L_{11}^{(1)}&\im L_{12}^{(1)}\\
 \im L_{21}^{(1)}&L_{22}^{(1)}\end{pmatrix},
\end{gathered}
\end{equation}
where all the displayed scalar entries are real.
Then \eqref{eq:three-sylvester-equations} is equivalent to the three decoupled systems
\begin{equation}\label{eq:three-real-systems}
\begin{gathered}
 A_1\begin{pmatrix}x_{11}^{(1)}\\x_{12}^{(1)}\\x_{21}^{(1)}\\x_{22}^{(1)}\end{pmatrix}
 =\begin{pmatrix}-I_{21}^{(1)}\\I_{22}^{(1)}\\-I_{11}^{(1)}\\I_{12}^{(1)}\end{pmatrix},\qquad
 A_2\begin{pmatrix}x_{11}^{(2)}\\x_{12}^{(2)}\\x_{21}^{(2)}\\x_{22}^{(2)}\end{pmatrix}
 =\begin{pmatrix}-F_{21}^{(1)}\\F_{22}^{(1)}\\-F_{11}^{(1)}\\F_{12}^{(1)}\end{pmatrix},
 \qquad A_3\begin{pmatrix}x_{11}^{(3)}\\x_{12}^{(3)}\\x_{21}^{(3)}\\x_{22}^{(3)}\end{pmatrix}
 =\begin{pmatrix}-L_{21}^{(1)}\\L_{22}^{(1)}\\-L_{11}^{(1)}\\L_{12}^{(1)}\end{pmatrix},
\end{gathered}
\end{equation}
where the $4\times 4$ matrices $A_k$, $k=1,2,3$ have the form 
\small
\begin{equation}\label{eq:Ak-three-blocks}
A_k = \begin{pmatrix}
  a_k & b_k & c_k & 0\\
  d_k & a_k & 0   & -c_k\\
  e_k & 0   & a_k & -b_k\\
  0   & -e_k& -d_k& a_k
\end{pmatrix},
\qquad
\begin{array}{@{}l@{\;\;}l@{\;\;}l@{\;\;}l@{\;\;}l@{}}
  a_1 = H_{12}^{(1)}-E_{12}^{(1)}, & b_1 = H_{11}^{(1)}, & c_1 = E_{22}^{(1)}, & d_1 = H_{22}^{(1)}, & e_1 = E_{11}^{(1)},\\[2pt]
  a_2 = G_{12}^{(1)}-E_{12}^{(1)}, & b_2 = G_{11}^{(1)}, & c_2 = E_{22}^{(1)}, & d_2 = G_{22}^{(1)}, & e_2 = E_{11}^{(1)},\\[2pt]
  a_3 = G_{12}^{(1)}-H_{12}^{(1)}, & b_3 = G_{11}^{(1)}, & c_3 = H_{22}^{(1)}, & d_3 = G_{22}^{(1)}, & e_3 = H_{11}^{(1)}.
\end{array}
\end{equation}
\normalsize
Using also  \cite[Lemma 5.5]{BMV3}, one has 
\begin{equation}\label{eq:Ak-determinant}
 \det A_k=a_k^4-2a_k^2(b_kd_k+c_ke_k)
                 +(b_kd_k-c_ke_k)^2
\end{equation}
and, in case $\det A_k \neq 0$, 
the inverse is
\small
\begin{equation}\label{eq:Ak-inverse}
\begin{aligned}
 A_k^{-1}&=\frac1{\det A_k}
 \begin{pmatrix}
 a_k(a_k^2-b_kd_k-c_ke_k)&b_k(-a_k^2+b_kd_k-c_ke_k)&-c_k(a_k^2+b_kd_k-c_ke_k)&-2a_kb_kc_k\\
 d_k(-a_k^2+b_kd_k-c_ke_k)&a_k(a_k^2-b_kd_k-c_ke_k)&2a_kc_kd_k&c_k(a_k^2+b_kd_k-c_ke_k)\\
 -e_k(a_k^2+b_kd_k-c_ke_k)&2a_kb_ke_k&a_k(a_k^2-b_kd_k-c_ke_k)&b_k(a_k^2-b_kd_k+c_ke_k)\\
 -2a_kd_ke_k&e_k(a_k^2+b_kd_k-c_ke_k)&d_k(a_k^2-b_kd_k+c_ke_k)&a_k(a_k^2-b_kd_k-c_ke_k)
 \end{pmatrix}.
\end{aligned}
\end{equation}
\normalsize
Using the expansions in \eqref{eq:E1-first-change}--\eqref{eq:G1-first-change},
the entries of $A_k$ have the Taylor expansions
\begin{align*}\notag
a_1 &= 2 \beta \ve^2 + r(\ve^4), \qquad
b_1,\ c_1,\ d_1,\ e_1 = \alpha \ve^2 + r(\ve^4),\\
\notag
a_2 &= \beta \ve^2 + r(\ve^4), \qquad
b_2 =\tfrac{\ck}{4}\,\ve^2 + r(\ve^4), \qquad
c_2,\ e_2 = \alpha\ve^2 + r(\ve^4), \qquad
d_2 = r(\ve^4),\\
\notag
a_3 &= -\beta \ve^2 + r(\ve^4), \qquad
b_3 = \tfrac{\ck}{4}\,\ve^{2} + r(\ve^4), \quad
c_3,\ e_3 = \alpha\ve^2 + r(\ve^4), \quad
d_3 = r(\ve^4).
\end{align*}
Hence
\begin{equation}\label{eq:three-determinant-expansions}
\begin{aligned}
 \det A_1
 &=2^{5/2} 3 \,\tau^2   \beta^2 \ve^8
       \big(1+r(\ve^2)\big), \quad 
 \det A_2 =\frac98 \tau^4 \ve^8 \big(1+r(\ve^2)\big), \quad 
 \det A_3=\frac98 \tau^4 \ve^8 \big(1+r(\ve^2)\big).
\end{aligned}
\end{equation}
Then for $\ve >0$, the matrices $A_k$ are invertible and 
 \eqref{eq:three-real-systems} has a unique real solution. Using formula \eqref{eq:Lblock1-first-change} for the right-hand side  of \eqref{eq:three-real-systems}, we obtain  $X_1$, $X_2$, $X_3$ as in 
\eqref{eq:three-X-form} with
\[
\begin{aligned}
\fa_2&=\frac{\varphi_1}{\beta - \alpha}
       -\frac{\ck\varphi_2}{4(\beta - \alpha)^2},
& \fb_2 &=\frac{\varphi_2}{\beta - \alpha}, 
& \fd&=\frac{\iota_I^{(1)}}{2(\beta - \alpha)},
\\
\fa_3 &=\frac{\phi_1}{\beta - \alpha}
       -\frac{\ck\phi_2}{4(\beta - \alpha)^2},
&\fb_3 &=\frac{\phi_2}{\beta - \alpha}, &\fc&=\frac{\psi_5}{\alpha + \beta},
\end{aligned}
\]
which, using the values in \eqref{final.coeff}, \eqref{eq:coeff-first-change}, give the coefficients in \eqref{X.coeff}.

Using \eqref{eq:homological-three-blocks}, the Lie expansion
\eqref{eq:lie-expansion-three-blocks} simplifies to the exact identity
\begin{equation}\label{eq:lie-simplified-three-blocks}
 \mathsf L^{(2)}_{\tau,\ve,1}
 =D^{(1)}+\frac12[S_1,R^{(1)}]
 +\frac12\int_0^1(1-t^2)
 e^{t S_1}\ad_{S_1}^2(R^{(1)})e^{-t S_1}\,\de t.
\end{equation}
The new  term 
\begin{equation}\label{eq:commutator-six-blocks}
 \frac12[S_1,R^{(1)}]=\mathsf J_6
 \begin{pmatrix}
 \widetilde E&\widetilde I&\widetilde F\\
 \widetilde I^*&\widetilde H&\widetilde L\\
 \widetilde F^*&\widetilde L^*&\widetilde G
 \end{pmatrix}
\end{equation}
is given by 
\begin{equation}\label{eq:tilde-block-identities}
{\small
\begin{aligned}
  \widetilde E &= \operatorname{Sym}\big(
    \mathsf J_2X_1\mathsf J_2[I^{(1)}]^*
    + \mathsf J_2X_2\mathsf J_2[F^{(1)}]^*\big),
  &\qquad
  \widetilde H &= \operatorname{Sym}\big(
    X_1^*I^{(1)} + \mathsf J_2X_3\mathsf J_2[L^{(1)}]^*\big), \\[4pt]
  \widetilde G &= \operatorname{Sym}\big(X_2^*F^{(1)} + X_3^*L^{(1)}\big),
  &\qquad
  \widetilde I &= \tfrac12\big(
    \mathsf J_2X_2\mathsf J_2[L^{(1)}]^*
    + F^{(1)}\mathsf J_2X_3^*\mathsf J_2\big), \\[4pt]
  \widetilde F &= \tfrac12\big(
    \mathsf J_2X_1\mathsf J_2L^{(1)} + I^{(1)}X_3\big),
  &\qquad
  \widetilde L &= \tfrac12\big(X_1^*F^{(1)} + [I^{(1)}]^*X_2\big),
\end{aligned}
}
\end{equation}
 where we wrote $\operatorname{Sym}(A)=\frac12(A+A^*)$.
 The matrices $\widetilde E, \widetilde H, \widetilde G$ are selfadjoints.
The next lemma computes the Taylor expansions of such blocks.
\begin{lemma}\label{lem:cubic-commutator}
Uniformly on \eqref{eq:sylvester-band-new}, the six blocks in
\eqref{eq:commutator-six-blocks} have the expansions
\begin{align}\notag
\widetilde E &= \ve^2
\begin{pmatrix}
e_1\,\ve^4+e_2\,\ve^5  + r(\ve^6)& \im r(\ve^4)\\[1ex]
-\im  r(\ve^4) & e_1\,\ve^4-e_2 \,\ve^5 + r(\ve^6)
\end{pmatrix} , \; \quad
\widetilde G =
\begin{pmatrix}
r(\ve^6) & \im r(\ve^6)\\[1ex]
-\im r(\ve^6) & g_2\ve^6 + r(\ve^8)
\end{pmatrix},  \\
\label{tildeEHG}
\widetilde H & =\ve^2
\begin{pmatrix}
h_1\,\ve^4+h_2\,\ve^5  + r(\ve^6)& \im r(\ve^4)\\[1ex]
-\im r(\ve^4)& h_1\,\ve^4-h_2\,\ve^5+ r(\ve^6)
\end{pmatrix} , 
\end{align}
where 
\begin{equation}\label{coeff.tilde.mat}
\begin{aligned}
e_1&:=\sqrt2\left(\frac{459}{4096}-\frac{\tau^8}{108}\right)
      +\frac{\ck}{9}\tau^4,
\qquad 
e_2:=\frac{2^{1/4}\tau}{576}
\Big[
27\big(3+2\sqrt2\tau^2\big)
-8\ck\tau^4\big(3\sqrt2+4\tau^2\big)
\Big],\\
h_1&:=\sqrt2\left(\frac{459}{4096}-\frac{\tau^8}{108}\right)
      -\frac{\ck}{9}\tau^4,\qquad 
h_2:=-\frac{2^{1/4}\tau}{576}
\Big[
27\big(3+2\sqrt2\tau^2\big)
+8\ck\tau^4\big(3\sqrt2+4\tau^2\big)
\Big] \ , \\
g_2& := \frac{3 \ck}{16} +\frac{4\ck}{81}\tau^8,
\end{aligned}
\end{equation}
whereas 
$\widetilde I , \ \widetilde F, \ \widetilde L = \cO(\ve^6)$.
\end{lemma}
\begin{proof}
Substitute \eqref{X1} 
and \eqref{eq:Lblock1-first-change}
into \eqref{eq:tilde-block-identities}, getting the expansion in \eqref{tildeEHG} with  coefficients 
\[
\begin{aligned}
e_1&=2\fd\iota_{I}^{(1)}-\fa_2 \varphi_2-\fb_2 \varphi_1, 
\qquad 
e_2 =\fc\varphi_2+\psi_5 \fb_2, 
\qquad 
g_2 = 2 (\fb_2\varphi_2 + \fb_3 \phi_2) \ , 
\\
h_1& =2\fd\iota_{I}^{(1)}-\fa_3\phi_2-\fb_3 \phi_1,\qquad 
h_2=-\fc\phi_2-\psi_5 \fb_3.
\end{aligned}
\]
Then \eqref{coeff.tilde.mat} follows from 
\eqref{X.coeff} and \eqref{final.coeff}, \eqref{eq:coeff-first-change}.
\end{proof}

We finally estimate the integral term in
\eqref{eq:lie-simplified-three-blocks}.

\begin{lemma}\label{lem:lie-remainder-three-blocks}
The integral in \eqref{eq:lie-simplified-three-blocks} has the form
\begin{equation}\label{eq:hat-blocks-three-blocks}
 \frac12\int_0^1(1-t^2)
 e^{t S_1}\ad_{S_1}^2(R^{(1)})e^{-t S_1}\,\de t
 =\mathsf J_6
 \begin{pmatrix}
 \widehat E&\widehat I&\widehat F\\
 \widehat I^*&\widehat H&\widehat L\\
 \widehat F^*&\widehat L^*&\widehat G
 \end{pmatrix},
\end{equation}
where $\widehat E,\widehat H,\widehat G$ are self-adjoint,
all six blocks are reversibility preserving, and
\begin{equation}\label{eq:hat-blocks-estimates}
 \widehat E,\widehat H,\widehat G,
 \widehat I,\widehat F,\widehat L=\cO(\ve^8).
\end{equation}
\end{lemma}

\begin{proof}
The bounds $S_1=\cO(\ve^2)$ and $R^{(1)}=\cO(\ve^4)$ imply
\[
 \|\ad_{S_1}^2(R^{(1)})\|
 \leq4\|S_1\|^2\|R^{(1)}\|\leq C\ve^8,\qquad
 \|e^{\pm t S_1}\|\leq e^{C\ve^2},\quad 0\leq t \leq1.
\]
This proves \eqref{eq:hat-blocks-estimates}.
The integrand is Hamiltonian and reversible, which gives the
block structure in \eqref{eq:hat-blocks-three-blocks}.
\end{proof}

\begin{proof}[Proof of Proposition~\ref{prop:first-block-decoupling}]
The unique solutions of \eqref{eq:three-real-systems} give the matrices
$X_1,X_2,X_3$, with the expansions already proved.
Their real analyticity for positive $(\tau,\ve)$ follows from
\eqref{eq:Ak-inverse} and the nonvanishing determinants
\eqref{eq:three-determinant-expansions}.
Equations \eqref{eq:lie-simplified-three-blocks},
\eqref{eq:commutator-six-blocks}, and
\eqref{eq:hat-blocks-three-blocks} give the exact identities
\begin{equation*}
\begin{aligned}
 E^{(2)}&=E^{(1)}+\widetilde E+\widehat E,&
 I^{(2)}&=\widetilde I+\widehat I,\\
 H^{(2)}&=H^{(1)}+\widetilde H+\widehat H,&
 F^{(2)}&=\widetilde F+\widehat F,\\
 G^{(2)}&=G^{(1)}+\widetilde G+\widehat G,&
 L^{(2)}&=\widetilde L+\widehat L.
\end{aligned}
\end{equation*}
Lemmata~\ref{lem:cubic-commutator} and
\ref{lem:lie-remainder-three-blocks} now yield the stated off-diagonal estimates and all the asserted expansions.
In particular, the previously computed $\ve^4$ coefficients in
$E^{(1)},H^{(1)},G^{(1)}$ are retained.
\end{proof}

\subsection{Complete Block-Decoupling}
\label{subsec:complete-block-decoupling}
In this section we complete block-diagonalize the $6\times6$ Hamiltonian and reversible
matrix $\mathsf L^{(2)}_{\tau,\ve,1}$ in \eqref{L2}.

In this case we first split $\mathsf L^{(2)}_{\tau,\ve,1}$ into its $2\times2$ block-diagonal and off-diagonal Hamiltonian and reversible matrices
\begin{equation*}
\begin{aligned}
 \mathsf L^{(2)}_{\tau,\ve,1}&=D_1^{(2)}+R_1^{(2)},\\
 D_1^{(2)}&:=\begin{pmatrix}
 \mathsf J_2E^{(2)}&0&0\\
 0&\mathsf J_2H^{(2)}&0\\
 0&0&\mathsf J_2G^{(2)}
 \end{pmatrix},\qquad
 R_1^{(2)}:=\begin{pmatrix}
 0&\mathsf J_2I^{(2)}&\mathsf J_2F^{(2)}\\
 \mathsf J_2[I^{(2)}]^*&0&\mathsf J_2L^{(2)}\\
 \mathsf J_2[F^{(2)}]^*&\mathsf J_2[L^{(2)}]^*&0
 \end{pmatrix}.
\end{aligned}
\end{equation*}
By Proposition~\ref{prop:first-block-decoupling},
$D_1^{(2)}=\cO(\ve^2)$ and $R_1^{(2)}=\cO(\ve^6)$.

\begin{lemma}\label{lem:complete-block-decoupling}
There exist a $6\times6$ reversibility-preserving Hamiltonian matrix
$S_1^{(2)}:=S_1^{(2)}(\tau,\ve)$ of the form \eqref{eq:S1-three-blocks} and a $6\times6$ block-diagonal reversible Hamiltonian matrix $P_1:=P_1(\tau,\ve)$ such that
\begin{equation*}
   S_1^{(2)} = \cO(\ve^4) , \qquad P_1 = \cO(\ve^{10}) ,
\end{equation*}
and
\begin{equation}\label{eq:complete-conjugation}
 \mathsf L^{(3)}_{\tau,\ve,1}
 :=\exp(S_1^{(2)})\mathsf L^{(2)}_{\tau,\ve,1}\exp(-S_1^{(2)})=D_1^{(2)}+P_1.
\end{equation}
Both matrices are real analytic in $(\tau,\ve)$ for positive parameters in a neighborhood of the band \eqref{eq:sylvester-band-new}. In particular,
\begin{equation*}
\begin{aligned}
 \mathsf L^{(3)}_{\tau,\ve,1}
   =\mathsf J_6\mathsf B^{(3)}_{\tau,\ve,1}
=\begin{pmatrix}
 \mathsf J_2 E^{(3)}&0&0\\
 0&\mathsf J_2 H^{(3)}&0\\
 0&0&\mathsf J_2 G^{(3)}
 \end{pmatrix},
\end{aligned}
\end{equation*}
where $E^{(3)}$, $H^{(3)} $, $G^{(3)}$ are self-adjoint and reversibility
preserving, and expand as 
\eqref{eq:E2-final}--\eqref{eq:G2-final}.
\end{lemma}

\begin{proof}
We set for brevity $S=S_1^{(2)}$, $D^{(2)}=D_1^{(2)}$, $R^{(2)}=R_1^{(2)}$.
Equation \eqref{eq:complete-conjugation} is equivalent to the system
\begin{equation}\label{eq:complete-projected-system}
\begin{cases}
 \Pi_D\big(e^S(D^{(2)}+R^{(2)})e^{-S}-D^{(2)}\big)&=P_1,\\
 \Pi_\emptyset\big(e^S(D^{(2)}+R^{(2)})e^{-S}\big)&=0,
\end{cases}
\end{equation}
where $\Pi_D$ is the projector onto the block-diagonal matrices and
$\Pi_\emptyset$ onto the block-off-diagonal ones, with block sizes
$2,2,2$. The second equation in \eqref{eq:complete-projected-system}
is equivalent, by a Lie expansion, to
\begin{equation}\label{eq:complete-nonlinear-equation}
\begin{aligned}
 &R^{(2)}+[S,D^{(2)}]+R(S)=0,\\
 &R(S):=\Pi_\emptyset[S,R^{(2)}]
 +\Pi_\emptyset\int_0^1(1-t)e^{t S}
       \ad_S^2(D^{(2)}+R^{(2)})e^{-t S}\,\de t.
\end{aligned}
\end{equation}
Here $[S,D^{(2)}]$ has zero diagonal blocks.
With three blocks, $[S,R^{(2)}]$ may have nonzero off-diagonal blocks,
which are included in $R(S)$.
The nonlinear homological equation
\begin{equation}\label{eq:complete-homological}
 [S,D^{(2)}]=-R^{(2)}-R(S)
\end{equation}
is equivalent to solve a $12\times12$ real system
\begin{equation}\label{eq:complete-real-system}
 A\vec x=f(\tau,\ve,\vec x),\qquad
 f(\tau,\ve,\vec x)=\ve^2 \vec\nu(\tau,\ve)
                         +\ve^2 \vec g(\tau,\ve,\vec x)
\end{equation}
associated, as in \eqref{eq:three-real-systems},  to \eqref{eq:complete-homological}.
Precisely $\vec x\in\R^{12}$ consists of the four real entries of each
of the three matrices $X_1,X_2,X_3$ representing $S$ as in
\eqref{eq:S1-three-blocks}, in the order used in
\eqref{eq:three-real-systems}; the matrix
$A=\operatorname{diag}(A_1,A_2,A_3)$
is obtained from \eqref{eq:Ak-three-blocks} by replacing
$E^{(1)},H^{(1)},G^{(1)}$ with $E^{(2)},H^{(2)},G^{(2)}$.
Since these blocks differs  by $\cO(\ve^6)$, by \eqref{eq:Ak-inverse} and \eqref{eq:three-determinant-expansions}  we obtain the bound
\begin{equation}\label{eq:complete-inverse}
 A^{-1}=\ve^{-2} B(\tau,\ve),\qquad \|B(\tau,\ve)\|\leq C.
\end{equation}
The vector $\ve^2 \vec\nu$ is associated with $-R^{(2)}$, namely
\begin{equation*}
\begin{aligned}
 \ve^2 \vec\nu:=\big(&-I_{21}^{(2)},I_{22}^{(2)},-I_{11}^{(2)},I_{12}^{(2)}, -F_{21}^{(2)},F_{22}^{(2)},-F_{11}^{(2)},F_{12}^{(2)}, -L_{21}^{(2)},L_{22}^{(2)},-L_{11}^{(2)},L_{12}^{(2)}\big)^{\mathsf T},
 \qquad |\vec\nu|\leq C\ve^4.
\end{aligned}
\end{equation*}
As in \eqref{eq:real-entries-three-blocks}, these scalar entries are
real, with the factors $\im$ removed from the off-diagonal entries.
The vector $\ve^2 \vec g$ is associated with $-R(S)$, which is
Hamiltonian, reversible, and block off-diagonal.
Since $S$ depends linearly on $\vec x$, the bounds for $D^{(2)}$ and
$R^{(2)}$ give, for $|\vec x|$ in a fixed small ball,
\begin{equation}\label{eq:complete-nonlinear-bounds}
\begin{aligned}
 |\vec g(\tau,\ve,\vec x)|
 &\leq C\big(\ve^4|\vec x|+|\vec x|^2\big), \qquad
 \|\partial_{\vec x}\vec g(\tau,\ve,\vec x)\|
 &\leq C\big(\ve^4+|\vec x|\big).
\end{aligned}
\end{equation}
The term linear in $\vec x$ comes from
$\Pi_\emptyset[S,R^{(2)}]$; the integral in
\eqref{eq:complete-nonlinear-equation} vanishes to second order in $\vec x$.

By \eqref{eq:complete-inverse}, system
\eqref{eq:complete-real-system} is equivalent to
\begin{equation}\label{eq:complete-fixed-point}
 \vec x=B(\tau,\ve)\vec\nu(\tau,\ve)
       +B(\tau,\ve)\vec g(\tau,\ve,\vec x).
\end{equation}
For a sufficiently large fixed $C_0$, the right-hand side maps
$|\vec x|\leq C_0\ve^4$ into itself and has derivative of size
$\cO(\ve^4)$ there. It is therefore a contraction for small $\ve$,
uniformly on \eqref{eq:sylvester-band-new}, and gives a unique solution
$\vec x=\cO(\ve^4)$. The analytic implicit function theorem applied
to \eqref{eq:complete-fixed-point} gives real analyticity in the
positive parameters $(\tau,\ve)$. Uniqueness makes the local analytic
solutions agree on overlaps. This proves the assertion for $S_1^{(2)}$.

The first equation in \eqref{eq:complete-projected-system} now gives
\begin{equation}\label{eq:complete-diagonal-correction}
 P_1=\Pi_D[S,R^{(2)}]
   +\Pi_D\int_0^1(1-t)e^{t S}
       \ad_S^2(D^{(2)}+R^{(2)})e^{-t S}\,\de t.
\end{equation}
Consequently
\[
 \|P_1\|\leq C\big(\|S\|\|R^{(2)}\|
                 +\|S\|^2(\|D^{(2)}\|+\|R^{(2)}\|)\big)
          \leq C\ve^{10}.
\]
Finally, $\exp(S)$ is symplectic and reversibility preserving, so
$P_1$ is Hamiltonian and reversible. Its block-diagonal structure
follows from \eqref{eq:complete-diagonal-correction}.
\end{proof}

\begin{proof}[Proof of Theorem \ref{thm:main.mat}]
Recalling \eqref{mathcal L= ichmu+mathscr L}, the matrix
representing the full operator is
$ \im\ck \tau^4 \ve^4 I_6+\mathsf L^{(3)}_{\tau,\ve,1}$. Then \eqref{U} 
   follows by Lemma~\ref{lem:complete-block-decoupling} putting
 $\mathsf{U}_E:= \im\ck \tau^4 \ve^4 I_2 +\mathsf J_2 E^{(3)}$,   $\mathsf{U}_H:=\im\ck \tau^4 \ve^4 I_2+ \mathsf J_2 H^{(3)}$ and
  $\mathsf{S}:= \im\ck \tau^4 \ve^4 I_2+\mathsf J_2 G^{(3)}$, and defining for $X \in \{E, H\}$
  \begin{align} \label{aX dX bX}
a_X:= \frac{X_{11}^{(3)}+X_{22}^{(3)}}{2} \ , 
\quad 
d_X:= \frac{X_{22}^{(3)}-X_{11}^{(3)}}{2} \ , \quad   
b_X := \ck \tau^4 \ve^4 - X_{12}^{(3)} \ .
\end{align}
  The expansions in \eqref{final.coeff} follows from those in 
 \eqref{eq:E2-final}, \eqref{eq:H2-final}. 
\end{proof}

\section{Spectral Analysis}
\label{sec:resonant-instability}
\begin{figure}[h]
\includegraphics[width=10cm]{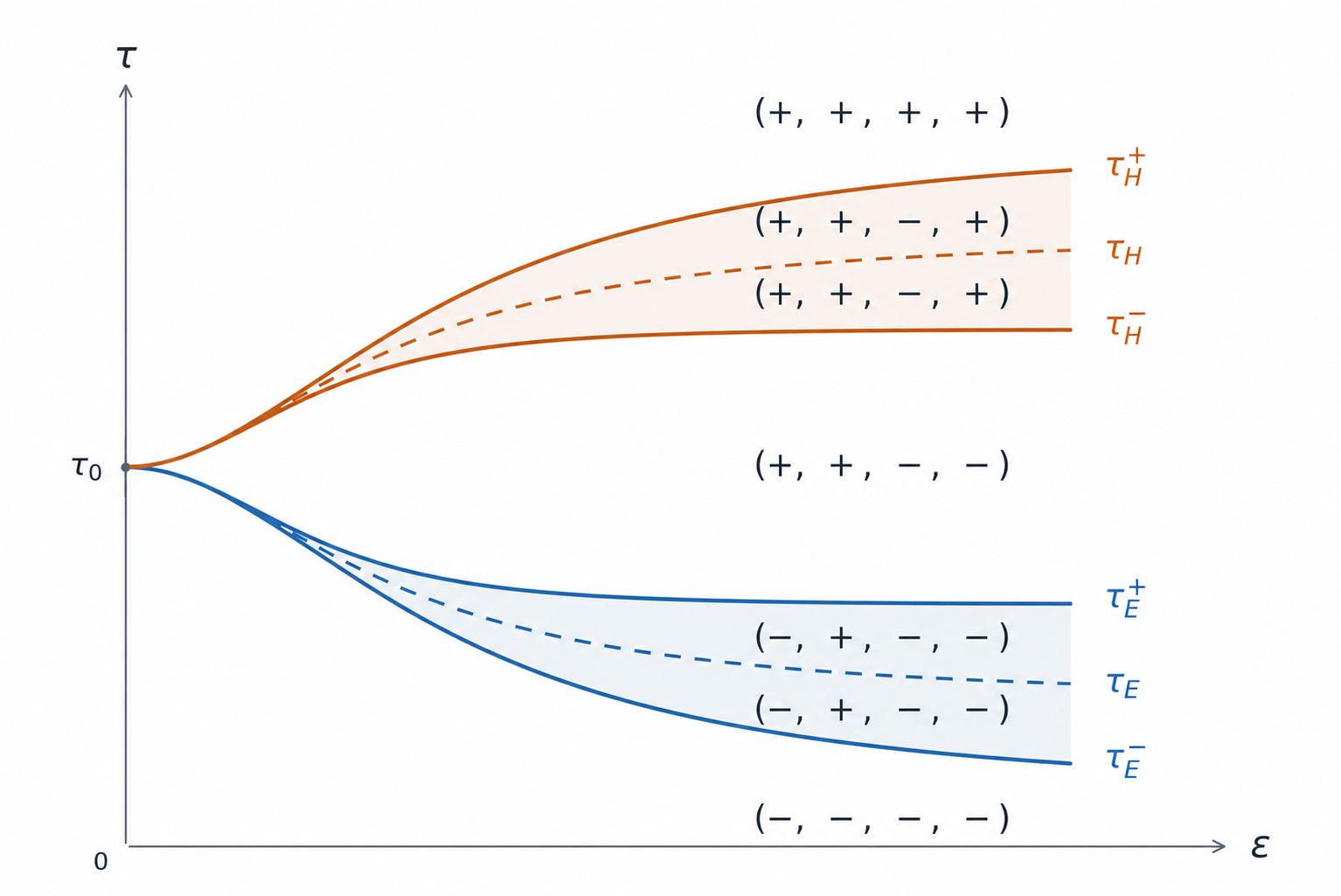}
\centering
\caption{The real analytic curves $\tau_E, \tau_{E}^\pm,  \tau_H, \tau_{H}^\pm$. The signs between the six curves are the signs of the functions $(a_E - d_E, a_E + d_E, a_H-d_H, a_H + d_H)$ in this order. In the light blue region between $\tau_{E}^-$ and $\tau_{E}^+$, the functions $a_E - d_E, a_E + d_E$ have opposite signs so $\mathsf U_E$ has two eigenvalues with {nonzero} real part. In the light orange region, the two eigenvalues of $\mathsf U_H$ have nonzero real part.}
\label{fig:regions}
\end{figure}

In this section we show that the matrix in \eqref{U} has unstable eigenvalues. 
In view of \eqref{UX}, the eigenvalues of $\mathsf U_X$, $X \in \{E, H\}$ are given by 
\begin{equation} \label{0911:2314}
    \begin{aligned}
    \lambda_X^\pm(\tau, \ve) & = \im b_X(\tau, \ve) \pm\sqrt{(d_X(\tau, \ve) - a_X(\tau, \ve))(d_X(\tau, \ve) + a_X(\tau, \ve))} \\
    & =  \im b_X(\tau, \ve) \pm \sqrt{ d_X^2(\tau, \ve)  - a_X^2(\tau, \ve)} .
\end{aligned}
\end{equation}
To describe the instability region, we shall prove the existence of analytic curves  $\tau_X(\ve)$ along which 
$$
a_X(\tau_X(\ve), \ve)\equiv 0 \ , \quad d_X(\tau_X(\ve), \ve) \not \equiv 0 \ . 
$$
\begin{lemma}
\label{lem:mu-plus-center-endpoints}
There exists $\ve_1>0$ and six  real analytic functions $\tau_E, \tau_{E}^+, \tau_{E}^-, \tau_H, \tau_{H}^+, \tau_{H}^-\colon [0, \ve_1)\to (1, 2^\frac14)$ with 
$$
\tau_{E}^-(\ve)<\tau_E(\ve)<\tau_{E}^+(\ve) <\tau_{H}^-(\ve)<\tau_H(\ve)<\tau_{H}^+(\ve) \ , \quad \mbox{ for }\ve >0
$$ 
and 
\begin{align} \label{tau0}
\tau_{E}^-(0)=\tau_E(0)=\tau_{E}^+(0) = \tau_{H}^-(0) =\tau_{H}(0)= \tau_{H}^+(0) = \tau_0= \left(\tfrac98 \right)^{\frac14}
\end{align}
satisfying
\begin{equation}\label{eq:mu-plus-center-definition}
\begin{aligned}
& a_X(\tau_X(\ve),\ve)  \equiv 0 ,  \qquad 
(a_X + d_X)(\tau_{X}^-(\ve),\ve) \equiv 0,
\qquad
(a_X - d_X)(\tau_{X}^+(\ve),\ve)\equiv 0
\end{aligned}
\end{equation}
and Taylor expanding, for  $ X \in \{E, H\}$, as 
\begin{equation}\label{eq:mu-plus-center-expansion}
\begin{aligned}
& \tau_X(\ve)=\tau_0
- \frac{\gamma_X(\tau_0)}{\tau_0} \ve^2 
+r(\ve^4) ,\qquad 
& \tau_{X}^\pm(\ve)=\tau_X(\ve) \mp \frac{\delta_X(\tau_0)}{\tau_0}\ve^{5}+r(\ve^6),\qquad 
\end{aligned}
\end{equation}
For $(\tau, \ve) \in [1, 2^\frac14] \times (0, \ve_1)$, and $X \in \{E, H\}$ one has 
\begin{equation}\label{0911:2207}
    (a_X - d_X)(\tau, \ve) > 0 \mbox{ iff } \tau > {\tau_{X}^+(\ve)}  \ , 
    \qquad 
      (a_X + d_X)(\tau, \ve) > 0 \mbox{ iff } \tau > {\tau_{X}^-(\ve)}.
\end{equation} 
See Figure \ref{fig:regions}.
\end{lemma}

\begin{proof} 
Fix $X\in\{E,H\}$. In view of the \eqref{functionsB2-theorem} and \eqref{eq:E2-final}, the functions
\begin{align*}
    \mathsf A_X(\tau,\ve):=\ve^{-2} a_X(\tau,\ve),\qquad \mathsf D_X(\tau,\ve):=\ve^{-2} d_X(\tau,\ve),\qquad \mathsf Z^{\pm}_X(\tau,\ve)=\mathsf A_X(\tau,\ve)\mp \mathsf D_X(\tau,\ve),
\end{align*}
extends analytically at $\ve = 0$. 
For each \(\mathtt{g}\in\{\mathsf A_X,\mathsf Z^{+}_X, \mathsf Z^{-}_X\}\), again using \eqref{functionsB2-theorem}, \eqref{final.coeff} and recalling \eqref{tau0}, 
\begin{equation}\label{0811:2318}
\mathtt{g}(\tau_0,0)=\alpha(\tau_0)=0,
\qquad
\partial_\tau \mathtt{g}(\tau_0,0)=\alpha'(\tau_0)=\tau_0\ne0  \ , \qquad 
\partial^2_\tau \mathtt{g}(\tau_0,0)=\alpha''(\tau_0)=1\ne0.
\end{equation}
The analytic implicit function theorem therefore gives analytic zero curves
\(\tau_X,\tau_{X}^+,\tau_{X}^-\), respectively fulfilling \eqref{eq:mu-plus-center-definition}. Since
$\mathsf A_X(\tau,\varepsilon)
=\alpha(\tau)+\gamma_X(\tau)\varepsilon^2
+r(\varepsilon^4)$,
the Taylor expansion of 
\(\mathsf A_X(\tau_X(\varepsilon),\varepsilon)=0\) gives
\begin{equation}\label{tauE.exp}
\tau_X(\varepsilon)
=\tau_0-\frac{\gamma_X(\tau_0)}{\tau_0}\varepsilon^2
+\cO(\varepsilon^4)
 \stackrel{\eqref{final.coeff}}{=}  \tau_0-\frac{1}{4\tau_0^3}
\left(\frac{153}{64}+\sigma_X\frac{57\sqrt3}{32}\right)\varepsilon^2
+\cO(\varepsilon^4),\quad\sigma_E=1,\quad\sigma_H=-1
\end{equation}
which is the first of \eqref{eq:mu-plus-center-expansion}.
Moreover, since $\mathsf A_X(\tau_X(\ve),\ve)=0$,
using \eqref{functionsB2-theorem} and \eqref{tauE.exp},
we obtain
\begin{equation}\label{0811:2320}
\begin{aligned}
 \mathsf Z_X^\pm(\tau_X(\ve),\ve)=\mp \mathsf D_X(\tau_X(\ve),\ve)=\pm\delta_X(\tau_X(\ve))\ve^5+\cO(\ve^6)=\pm\delta_X(\tau_0)\ve^5+\cO(\ve^6).
\end{aligned}
\end{equation}
Then, using  $\mathsf Z_X^{\pm}(\tau_{X}^{\pm}(\ve), \ve) \equiv  0$, by Taylor's formula and \eqref{0811:2318} we get
\begin{equation}\label{1111:1049}
    \begin{aligned}
        & \mathsf Z_X^{\pm}(\tau, \ve)  = \mathsf Z_X^{\pm}(\tau, \ve)-\mathsf Z_X^{\pm}(\tau_{X}^{\pm}(\ve), \ve)=(\tau - \tau_X^{\pm}(\ve)) e_{X}^{\pm}(\tau, \ve) \ , \\ &e_{X}^{\pm}(\tau, \ve):= \int_0^1 \pa_\tau \mathsf Z_X^{\pm}\big((1-s)\tau + s \tau_{X}^{\pm}(\ve), \ve\big) \, \di s 
      \mbox{ with }  e_{X}^{\pm}(\tau, \ve) = \frac{\tau + \tau_0}{2} + \cO(\ve^2) > 0 
    \end{aligned}
\end{equation}
uniformly for $\tau\in[1,2^{1/4}]$ and sufficiently small $\ve>0$. Substituting $\tau=\tau_X(\ve)$ into \eqref{1111:1049} and using \eqref{0811:2320}, together with
$e_X^\pm(\tau_X(\ve),\ve)=\tau_0+\cO(\ve^2)$, yields
\begin{align} \label{20261120 delta tauX}
    \tau_X^{\pm}(\ve)-\tau_X=-\frac{\mathsf Z^{\pm}_X(\tau_X,\ve)}{e^{\pm}_X(\tau_X,\ve)}=\mp \frac{\delta_X(\tau_0)}{\tau_0}\ve^5+\cO(\ve^6).
\end{align}
This proves the second expansion in \eqref{eq:mu-plus-center-expansion}. 

By \eqref{final.coeff}, we have $\delta_E(\tau_0)= - \frac{27 \sqrt2}{128} = \delta_H(\tau_0)$. Thus, \eqref{20261120 delta tauX} implies
$\tau_X^-(\ve)<\tau_X(\ve)<\tau_X^+(\ve)$
for sufficiently small $\ve>0$.
Moreover, \eqref{tauE.exp} and
\eqref{20261120 delta tauX} give
\begin{align*}
 \tau_H^-(\ve)-\tau_E^+(\ve)
 =\frac{\gamma_E(\tau_0)-\gamma_H(\tau_0)}
        {\tau_0}\ve^2+\cO(\ve^4)
 \overset{\eqref{final.coeff}}{=}\frac{19\ck}{18}\tau_0^3\ve^2+\cO(\ve^4)>0.
\end{align*}
This proves the full ordering of the six curves. 

Finally, \eqref{1111:1049} implies $(a_X\mp d_X)(\tau,\ve)
 =\ve^2\bigl(\tau-\tau_X^\pm(\ve)\bigr)
 e_X^\pm(\tau,\ve)$. Since $\ve^2e_X^\pm(\tau,\ve)>0$, we prove \eqref{0911:2207}.
\end{proof}

By Lemma \ref{lem:mu-plus-center-endpoints} and
\eqref{0911:2314}, for
$(\tau,\ve)\in[1,2^{1/4}]\times(0,\ve_1)$ we have
(cf. Figure \ref{fig:regions}):
\begin{itemize}
    \item If $\tau\le\tau_E^-(\ve)$,
$\tau_E^+(\ve)\le\tau\le\tau_H^-(\ve)$,
or $\tau\ge\tau_H^+(\ve)$, then $d_X^2-a_X^2\le0$ for $X\in\{E,H\}$.
    Hence both $\mathsf U_E$ and $\mathsf U_H$
    have purely imaginary eigenvalues.

    \item If $\tau_E^-(\ve)<\tau<\tau_E^+(\ve)$,
    then $d_E^2-a_E^2>0$ and $d_H^2-a_H^2<0$.
    Hence $\mathsf U_E$ has two eigenvalues with nonzero real parts,
    whereas $\mathsf U_H$ has purely imaginary eigenvalues.

    \item If $\tau_H^-(\ve)<\tau<\tau_H^+(\ve)$,
    then $d_H^2-a_H^2>0$ and $d_E^2-a_E^2<0$.
    Hence $\mathsf U_H$ has two eigenvalues with nonzero real parts,
    whereas $\mathsf U_E$ has purely imaginary eigenvalues.
\end{itemize}
The following proposition gives the asymptotic expansions
of the eigenvalues.

\begin{proposition}[{\bf Asymptotic expansions of the eigenvalues and isolas}]
\label{prop:isola}
There exist $\varepsilon_1>0$ and an open interval
$I \subset (1,2^{1/4})$  containing $\tau_0=\big(\frac98\big)^{\frac14}$  such that the following hold for any $(\tau,\varepsilon)\in I\times(0,\varepsilon_1)$.

\begin{enumerate}
\item[$\mathrm{1.}$] {\bf Eigenvalue expansions:} The matrix  $\mathsf U_X$, $X \in \{E, H\}$ in \eqref{U} has eigenvalues
\begin{equation}\label{1111:1253}
 \lambda^\pm_X(\tau,\varepsilon)=
 \begin{cases}
 \im b_X(\tau,\varepsilon)
 \pm\sqrt{\mathscr D_X(\tau,\varepsilon)},
 &\tau\in (\tau_{X}^-(\varepsilon),\tau_{X}^+(\varepsilon)),\\[1mm]
 \im b_X(\tau,\varepsilon)
 \pm \im \sqrt{|\mathscr D_X(\tau,\varepsilon)|},
 &\tau\in I\setminus (\tau_{X}^-(\varepsilon),\tau_{X}^+(\varepsilon)), 
 \end{cases}
\end{equation}
with $\tau_{X}^\pm(\varepsilon)$ as in Lemma \ref{lem:mu-plus-center-endpoints}. 
The discriminant $\mathscr D_X(\tau,\varepsilon):= (d_X^2- a_X^2)(\tau, \ve)$ is positive for  $\tau \in (\tau_{X}^-(\varepsilon),\tau_{X}^+(\varepsilon))$, vanishes at its two
endpoints, and is negative for  $\tau \in I\setminus\left[\tau_{X}^-(\varepsilon),\tau_{X}^+(\varepsilon)\right]$.
The functions $b_X$ (cf. \eqref{aX dX bX}) and $\mathscr D_X$ are real analytic, and have the expansions
\begin{equation}\label{1011:2056}
\begin{aligned}
 \mathscr D_X(\tau_X(\varepsilon)+\nu,\varepsilon)
 &=\varepsilon^4\Bigl[\big(\frac{27\sqrt2}{128}\big)^2\varepsilon^{10}-\tau_0^2\nu^2
 +\cO\bigl(\varepsilon^{11}, \varepsilon^{10}|\nu|, 
 \varepsilon^2\nu^2, |\nu|^3\bigr)\Bigr],\\
 b_X(\tau_X(\varepsilon)+\nu,\varepsilon)
 &=b_X(\tau_X(\varepsilon), \ve)+ \sigma_X\tau_0\ve^2 \nu 
+\cO\bigl(\varepsilon^4|\nu|, \varepsilon^2\nu^2\bigr) , \qquad 
\end{aligned}
\end{equation}
where ${\sigma_E}:=1$, ${\sigma_H} :=-1$.

The matrix  $\mathsf  S$ has purely imaginary eigenvalues for all $(\tau, \ve) \in I\times (0, \ve_1)$. 

\item[$\mathrm{2.}$] {\bf Geometry and size of the isolas:}
Let  $X\in\{E,H\}$ and $\ve \in (0, \ve_1)$. The   set
\begin{equation}\label{1111:1256}
 \Gamma_X(\varepsilon):=
 \bigcup_{\tau\in \left[\tau_{X}^-(\varepsilon),\tau_{X}^+(\varepsilon)\right]}
 \{\lambda_X^+(\tau,\varepsilon),\lambda_X^-(\tau,\varepsilon)\}
\end{equation}
 is a regular real analytic simple closed
curve bounding a strictly convex domain and invariant under
$\lambda\mapsto-\overline\lambda$. Its two intersections with the
imaginary axis have horizontal tangents. Moreover,
\[
 \Gamma_E(\varepsilon)\subset\{\operatorname{Im}\lambda>0\},\qquad
 \Gamma_H(\varepsilon)\subset\{\operatorname{Im}\lambda<0\},
\]
so the two curves are disjoint. They admit real analytic
parameterizations of the form
\begin{equation}\label{1111:1325}
 \Lambda_X(\theta,\varepsilon)
 =\im b_X(\tau_X(\ve), \ve)+\frac{27\sqrt2}{128}\varepsilon^7e^{\im \sigma_X\theta}
 +\cO(\varepsilon^8),
 \qquad 0\leq\theta\leq2\pi.
\end{equation}
In particular,
\begin{equation}\label{1111:1434}
 \max_{\lambda\in\Gamma_X(\varepsilon)}|\operatorname{Re}\lambda|
 =\frac{27\sqrt2}{128}\varepsilon^7+\cO(\ve^8).
\end{equation}
\end{enumerate}
\end{proposition}

\begin{proof}
1. The formulas for the eigenvalues follows from \eqref{0911:2314}.
We prove the expansion for  $\mathscr{D}_X(\tau_X(\ve) + \nu , \ve)$ in \eqref{1011:2056}. 
Using \eqref{functionsB2-theorem} --\eqref{final.coeff},  
we get 
\begin{align*}
& a_X(\tau_X(\ve)+\nu,\ve)
 =\ve^2 \left(\tau_0\nu+\cO(\nu^2, \ve^2|\nu|) \right), \quad 
 d_X(\tau_X(\ve)+\nu,\ve)
 =\ve^2\left(-\delta_X(\tau_0) \ve^5 + \cO(\ve^6, \ve^5 |\nu|) \right) \ . 
 \end{align*}
Inserting these  into $\mathscr D_X= d_X^2-a_X^2$ gives
the stated expansion, since $\delta_E(\tau_0)= - \frac{27 \sqrt2}{128} = \delta_H(\tau_0)$.

The expansion for $ b_X(\tau_X(\ve)+\nu,\ve)$ follows from  \eqref{functionsB2-theorem}, \eqref{final.coeff} and \eqref{eq:mu-plus-center-expansion}.

Finally,  $\mathsf S$ in \eqref{UX} has clearly  purely imaginary eigenvalues in $I\times (0, \ve_1)$ since the out-diagonal entries have opposite signs.

\smallskip
2. We parametrize the set $\Gamma_X(\ve)$ for a fixed  $\ve \in (0, \ve_1)$. 
First, using \eqref{1111:1049}, for $\tau \in (\tau_{X}^-(\ve), \tau_{X}^+(\ve))$
we get
\begin{equation}\label{1111:1254}
\begin{aligned}
    \mathscr D_X(\tau, \ve)
& =\varepsilon^4 K^2(\tau, \ve)  
(\tau_{X}^+(\ve)-\tau)(\tau-\tau_{X}^-(\ve)),\mbox{ with }K(\tau, \ve):= \sqrt{e_{X}^+(\tau,\varepsilon) e_{X}^-(\tau,\varepsilon) } \ , 
\, e_{X}^\pm(\tau, \ve)  \mbox{ in }  \eqref{1111:1049} \ 
\end{aligned}
\end{equation}
where we used that, in view of the second of \eqref{1111:1049}, 
for $\tau$ sufficiently close to $\tau_0$, the product $e_{X}^+(\tau,\varepsilon) e_{X}^-(\tau,\varepsilon) > 0 $.
Next introduce 
\begin{equation}\label{1111:1309}
\begin{aligned}
&  \tau(\theta, \ve ):=m(\ve) +h(\ve) \sin\theta, \quad \theta \in \T \quad \mbox{ where }\\
& m(\ve):=\frac{\tau_{X}^+(\ve) + \tau_{X}^-(\ve)}{2}
 \stackrel{\eqref{eq:mu-plus-center-expansion}}{=}\tau_X(\ve) +r(\ve^6),\quad
 h(\ve):=\frac{\tau_{X}^+(\ve) - \tau_{X}^-(\ve)}{2}
\stackrel{\eqref{eq:mu-plus-center-expansion}}{=}- \frac{\delta_X(\tau_0)}{\tau_0}\ve^5+r(\ve^6) \ ,
\end{aligned}
\end{equation}
so that as $\theta$ runs in $\T$, the variable $\tau(\theta, \ve)$ parametrizes the interval $[\tau_{X}^-(\ve), \tau_{X}^+(\ve)]$ forth and back.
In view of  \eqref{1111:1253}, \eqref{1111:1254}, and writing $\tau_{X}^+(\ve)= m(\ve) + h(\ve)$, $\tau_{X}^-(\ve)= m(\ve) - h(\ve)$, the set 
$\Gamma_X(\varepsilon)$ in \eqref{1111:1256} is then parametrized by 
\begin{equation}\label{1111:1328}
\begin{aligned}
 \Lambda_X(\theta,\ve)
 &:=\ve^2 h(\ve) \, K(\tau(\theta, \ve),\ve )\cos(\theta)
 +\im b_X(\tau(\theta, \ve),\ve),\qquad 0\leq\theta\leq2\pi, 
\end{aligned}
\end{equation}
with the  sign of $\cos(\theta)$ selecting one of the two eigenvalue branches.
The  expansions  in \eqref{1011:2056},  \eqref{1111:1049}, \eqref{1111:1309} then give
$$
 \Lambda_X(\theta,\ve)
 =\im b_X(\tau_X(\ve),\ve)
 +\frac{27\sqrt2}{128 }\varepsilon^7e^{\im \sigma_X \theta}
 +\cO_{C^2}(\ve^8),
$$
where the reminders are uniform in $\theta$ in $C^2$ thanks to the analyticity
of all factors.

{In order to prove that $\Gamma_X(\ve)$ is a simple curve, we observe that since $(\partial_{\tau}b_X)(\tau_X(\ve),\ve) = \sigma_X \tau_0\ve^2 + r(\ve^4)$ for $\tau\in[\tau_X^-(\ve),\tau_X^+(\ve)]$, the function $b_X$ is strictly monotone on this interval
for sufficiently small $\ve$. Hence the right and left branches have no self-intersections, and meet only at the two endpoints. The above expansion for $\Lambda_X$ implies that $\Gamma_X(\ve)$ bounds a strictly convex domain.}
The exact parameterization in \eqref{1111:1328} also gives
$\Lambda_X(\pi-\theta)=-\overline{\Lambda_X(\theta)}$.
At $\cos(\theta)=0$, its imaginary derivative vanishes and its real
derivative equals
$-\varepsilon^2h(\ve) K(\tau(\theta, \ve),\ve)\sin(\theta)\ne0$,
so both imaginary-axis intersections have horizontal tangents.

Finally  
$b_X(\tau_X(\ve), \ve) =\sigma_X\tau_0^2\ve^2+\cO(\ve^4)$, which places
the $E$ and $H$ isolas in the upper and lower half-planes, respectively.
The bound in \eqref{1111:1434}  follows immediately from the expansion \eqref{1111:1325}.
\end{proof}

\begin{proof}[Proof of Theorem~\ref{thm:main} and Theorem~\ref{prop:eigenvalue-isolas}]
 We start from  Proposition~\ref{prop:isola} 
 put $\ve:=\sqrt\epsilon$ and $\tau=(\mu/\epsilon^2)^{1/4}$.
Choose $\delta>0$ so small that 
$\left[\left(\frac98-\delta\right)^{1/4},
\left(\frac98+\delta\right)^{1/4}\right]\subset I$ 
and then shrink $\epsilon_*$ as needed.
Define, for $X\in\{E,H\}$,
\[
 \mu_X^0(\epsilon):=\epsilon^2\tau_X(\sqrt\epsilon)^4,
 \qquad
 \mu_X^\pm(\epsilon):=\epsilon^2\tau_X^\pm(\sqrt\epsilon)^4.
\]
These functions extend continuously by zero at $\epsilon=0$ and are real
analytic for $\epsilon>0$. In view of the expansions in 
\eqref{eq:mu-plus-center-expansion}
\[
\begin{aligned}
 \mu_X^0(\e)
 &=\frac98\epsilon^2-4\tau_0^2\gamma_X(\tau_0)\epsilon^3
 +\cO_{\mathrm{u}}(\epsilon^4),\quad 
 \mu_X^\pm(\e)-\mu_X^0(\e)
 =\pm4\tau_0^2\frac{27\sqrt2}{128}\epsilon^{9/2} +\cO_{\mathrm{u}}(\epsilon^5).
\end{aligned}
\]
Then \eqref{eq:intro-original-boundaries} follows using  
$ 4\tau_0^2\gamma_X(\tau_0)
 =\frac{153}{64}+\sigma_X\frac{57\sqrt3}{32},$
and 
$ 4\tau_0^2\frac{27\sqrt2}{128}=\frac{81}{64}$.
 All six curves lie in
$\mathcal W$ for small $\epsilon$, and the strictly increasing map
$\tau\mapsto\epsilon^2\tau^4$ preserves their ordering and maps the
instability intervals onto those in \eqref{def:RX}.
Set
\[
 \omega_X(\mu,\epsilon)
:=b_X\left((\mu/\epsilon^2)^{1/4},\sqrt\epsilon\right),\qquad
 D_X(\mu,\epsilon)
 :=\mathscr D_X\left((\mu/\epsilon^2)^{1/4},\sqrt\epsilon\right).
\]
The eigenvalue formulas, signs, collisions and purely imaginary stable
pair follow directly from the  Proposition~\ref{prop:isola}.

For the discriminant expansion in \eqref{eq:intro-original-normal-form}, fix $X$ and $\epsilon>0$ and put
$\mu_0:=\mu_X^0(\epsilon)$.
Write 
$a(\mu, \e) := a_X\big((\mu/\epsilon^2)^{1/4},\sqrt\epsilon \big)$, $d(\mu, \e) := d_X\big((\mu/\epsilon^2)^{1/4},\sqrt\epsilon \big)$ and 
$D_X(\mu, \e)=d^2(\mu, \e)- a^2(\mu, \e)$.
Then $a(\mu_0, \e)=0$, 
$d(\mu_0, \e)=\frac{27\sqrt2}{128}\epsilon^{7/2}+\cO_{\mathrm{u}}(\epsilon^4)$ and using the  chain rule $\partial_\mu=\frac{1}{4\epsilon^2\tau^3}\partial_\tau$ one gets
\[
 \partial_\mu a(\mu_0, \e)=\frac{1}{4\tau_0^2\epsilon}+\cO_{\mathrm{u}}(1), 
 \quad 
 \partial_\mu d=\cO_{\mathrm{u}}(\epsilon^{3/2}) \ , \quad
 \partial_\mu^2d=\cO_{\mathrm{u}}(\epsilon^{-1/2})
\]
 Hence
\[
\begin{aligned}
 D_X(\mu_0,\epsilon)
 &=\big(\frac{27\sqrt2}{128}\big)^2\epsilon^7+\cO_{\mathrm{u}}(\epsilon^{15/2}),\qquad 
 \partial_\mu D_X(\mu_0,\epsilon)
 =2d(\mu_0,\epsilon)\,\partial_\mu d(\mu_0,\epsilon)=\cO_{\mathrm{u}}(\epsilon^5),\\
 \frac12\partial_\mu^2D_X(\mu_0,\epsilon)
 &=(\partial_\mu d(\mu_0,\epsilon))^2+d(\mu_0,\epsilon)\,\partial_\mu^2d(\mu_0,\epsilon)
   -(\partial_\mu a(\mu_0,\epsilon))^2
 =-\frac{1}{16\tau_0^4\epsilon^2}
   +\cO_{\mathrm{u}}(\epsilon^{-1}).
\end{aligned}
\]
Moreover, $D_X=\epsilon^2F_X(\tau,\sqrt\epsilon)$ with $F_X$ analytic
and its first three $\tau$ derivatives uniformly bounded. The chain rule
therefore gives $\partial_\mu^3D_X=\cO_{\mathrm{u}}(\epsilon^{-4})$ on $\mathcal W$.
Taylor's formula in $\mu$, using $16\tau_0^4=18$  now yields
\eqref{eq:intro-original-normal-form}
uniformly whenever $(\mu_0+\nu,\epsilon)\in\mathcal W$.

Finally, the monotone change $\mu=\epsilon^2\tau^4$ maps the closed
endpoint intervals bijectively onto one another. Hence the isola is the
same subset of the complex plane as in Proposition~\ref{prop:isola} at
$\ve=\sqrt\epsilon$. Its geometric properties and half-plane location
therefore carry over unchanged. The maximal-real-part estimate follows
by substituting $\ve=\sqrt\epsilon$ into
$\frac{27\sqrt2}{128}\ve^7+\cO(\ve^8)$.
\end{proof}

\appendix
\section{Expansion and analyticity of Wilton ripples in infinite depth}\label{sec:App2}

In this section we prove the analyticity of the Wilton ripples and provide their Taylor expansion. We begin with the following result due to Reeder and  Shinbrot \cite{reeder1981wiltonII}:
\begin{theorem}[\cite{reeder1981wiltonII}]\label{thm:rs-wilton-amplitude}
For $g=1$, $\kappa=\frac12$ and infinite depth, there are two branches of $2\pi$-periodic Wilton ripples, parametrized by $\tts \in \{ \pm 1 \}$. For each $\tts \in \{\pm 1\}$,  $\e \mapsto (\eta_{\e,\tts}(x), \psi_{\e,\tts}(x), c_{\e,\tts})$ is real analytic as a map $B(\e_0) \to H^2_{\mathtt{ev}}(\T)\times H^2_{\mathtt{odd}}(\T)\times\R$, and one has the Taylor expansion 
$\eta_{\e,\tts}(x)=\e\bigl(\cos (x)+\tfrac{\tts}{2}\cos(2x)\bigr)+\cO(\e^2)$ in $H^2_{\mathtt{ev}}(\T)$
\end{theorem}

\begin{proof}
Take $\widehat\tau=0$ and $A=1$ in Reeder--Shinbrot \cite{reeder1981wiltonII}. Applying Theorem 5.2 of \cite{reeder1981wiltonII} to the two initial profiles of Corollary 4.2 in \cite{reeder1981wiltonII} yields the branches $\tts=\pm 1$. Corollary 5.3 in \cite{reeder1981wiltonII} gives their analytic dependence on $\e$ and the Taylor expansion.
\end{proof}

We now  prove the expansions \eqref{exp:Wil2}-\eqref{exp:c}, \eqref{pexp}--\eqref{fexp}, 
\eqref{SN1}-\eqref{h2 h02}.\\

\noindent{\bf Proof  of \eqref{exp:Wil2}-\eqref{exp:c}.}
Writing 
 \begin{equation}\label{etapsic}
\begin{aligned}
 & \eta_{\e,\tts}(x) = \e \eta_1(x) + \e^2 \eta_2(x) + \cO(\e^3) \, , \\
 &  \psi_{\e,\tts}(x) = \e \psi_1(x) + \e^2 \psi_2(x) + \cO(\e^3) \, ,  
  \end{aligned}
\qquad \quad c_{\e,\tts} = \ck + \e c_1 + \e^2 c_2+ \cO(\e^3) \, ,  
\end{equation}
where  $\eta_i$ is $\mathit{even}(x)$ and $\psi_i$ is $\mathit{odd}(x)$ for $i=1,2$ (we omit the dependence on $\tts$ in order to simplify the notation), 
we solve order by order in $ \e $ the equations \eqref{travelingWWstokes},  
that we rewrite  as
\begin{equation}
\label{Sts}
\begin{cases}
-c \, \psi_x +  \eta   + \dfrac{\psi_x^2}{2} - 
\dfrac{\eta_x^2}{2(1+\eta_x^2)} ( c  -  \psi_x )^2-\kappa\left(\frac{\eta_x}{(1+\eta_x^2)^{1/2}}\right)_x  = 0\,, \\
c \, \eta_x+G(\eta)\psi = 0 	\, ,
\end{cases}
\end{equation}
having substituted $G(\eta)\psi $ with $-c \, \eta_x $ in the first equation.
 We expand the Dirichlet-Neumann operator 
$ G(\eta)=  G_0+ G_1(\eta) + G_2(\eta) + G_3(\eta) + G_4(\eta) + \cO(\eta^5)  $
where, according to \cite[formula (2.13)] {CS}
\begin{equation}\label{expDiriNeu}
    \begin{aligned}
        G_0 \psi=|D|\psi,\qquad G_m(\eta)\psi=\frac{1}{m !}|D|^{m-1}D\left(\eta^m D\psi\right)-\sum_{q=1}^m \frac{1}{q !}|D|^{q}\left(\eta^q G_{m-q}(\eta)\psi\right),\quad m\geq 1.
    \end{aligned}
\end{equation}
{\bf First order in $ \e $.}  Substituting the expansions in \eqref{etapsic} into \eqref{Sts}, we get the linear system 
\begin{equation}\label{cB0}
 \left\{\begin{matrix} -\ck (\psi_1)_x + \eta_1 -\kappa (\eta_1)_{xx}= 0\,, \\
 \ck (\eta_1)_x + G_0\psi_1 =0 \, ,  \end{matrix}\right.
 \quad  \text{i.e.} \, \vet{\eta_1}{\psi_1} \in \text{Ker }\cB_0  \text{ with } \cB_0 := \begin{bmatrix} 1-\kappa\pa_{xx} & -\ck\pa_x \\ \ck\pa_x & G_0  \end{bmatrix},  
\end{equation}
where $\eta_1$ is $\mathit{even}(x)$ and $\psi_1$ is $\mathit{odd}(x)$.
\begin{lemma}\label{lem:B0}
If $\kappa = \frac{1}{2}$,  the kernel of 
 the linear operator $\cB_0$ in \eqref{cB0} is two-dimensional and given by
\begin{equation}\label{chk}
\text{Ker }\cB_0= \text{span}\,\{ v_1, v_2 \} , \quad
v_1 \coloneqq \vet{\cos(x)}{\ck \sin(x)} , \;
v_2 \coloneqq \vet{\cos(2 x)}{ \ck \sin(2 x)} .
\end{equation}
\end{lemma}

\begin{proof} The action of $\cB_0$ on each subspace spanned by $\footnotesize{\,\Big\{\vet{\cos(kx)}{0}, \vet{0}{\sin(kx)}\Big\}} $, $k\in \N$, is represented by the $2\times 2$ matrix $\footnotesize{ \begin{bmatrix} 1+\kappa k^2 & -\ck k \\ -\ck k & k  \end{bmatrix}}$. Its determinant is given by
 $$ (1+\kappa k^2)k  - \ck^2 k^2\stackrel{\eqref{exp:c}}{=} k^2 \Big(\frac{(1+\kappa k^2)}{k} -(1+\kappa) \Big)$$
 and, provided $\kappa = \frac{1}{2}$, 
 vanishes  if and only if $k \in \{ 1,2 \}$.
  For $k=1$ and for $k=2$ we obtain the kernel of $\cB_0$ given in \eqref{chk}. For $k=0$ it has no kernel since $\psi_1(x)$ is odd.
\end{proof}
We set
\begin{align}\label{2411:1210}
    \begin{bmatrix}
        \eta_1(x) \\ \psi_1(x)
    \end{bmatrix}
    &= \mathtt{a}_1 v_1 + \mathtt{a}_2 v_2 , \quad \mathtt{a}_1, \mathtt{a}_2 \in \mathbb{R},     \ \ v_1, v_2 \mbox{ in } \eqref{chk}, 
\end{align}
in agreement with \eqref{exp:Wil2}. Without loss of generality, we can choose $\mathtt{a}_1 = 1$.
\\[1mm]
{\bf Second order in $ \e $.} 
By \eqref{Sts}, and since 
$ \ck^2 (\eta_1)_x^2 = (G_0\psi_1)^2  $, we get  the linear system
\begin{equation}\label{syslin2}
 \cB_0 \vet{\eta_2}{\psi_2} = \vet{c_1(\psi_1)_x-\frac12 (\psi_1)_x^2 + \frac12 (G_0\psi_1)^2 }{-c_1(\eta_1)_x - G_1(\eta_1)\psi_1} \, , 
\end{equation}
where  $\cB_0$ is the self-adjoint operator in \eqref{cB0}. 
The system \eqref{syslin2} admits a solution if and only if its right-hand term is orthogonal to the Kernel  of $\cB_0$ in \eqref{chk}, namely
\begin{equation}\label{orth1}
 \Big(\vet{c_1(\psi_1)_x-\frac12 (\psi_1)_x^2 + \frac12 (G_0\psi_1)^2 }{-c_1(\eta_1)_x - G_1(\eta_1)\psi_1} \; , \; v_j \Big)=0 , \quad j=1,2 . 
\end{equation}
Equations \eqref{orth1} lead to a system of 2 equations in the unknowns $\mathtt{a}_2$ and $c_1$, which we want to determine. Using \eqref{orth1}, together with the first order expansion \eqref{exp:Wil2} and \eqref{expDiriNeu}, we obtain
\begin{align*} 
    \mathtt{a}_2 = \tts \left( \frac{1}{4}  \right)^{\frac{1}{2}} , &\quad c_1 = \tts \left( \frac{ 3  }{ 32  } \right)^{\frac{1}{2}} ,
\end{align*}
in agreement with \eqref{exp:Wil2}-\eqref{exp:c}. Setting $\eta_2=\eta_2^{[0]} + \eta_2^{[2]} \cos(2x) + \eta_2^{[3]} \cos(3 x) + \eta_2^{[4]} \cos(4x)$ and $\psi_2=\psi_2^{[1]} \sin(x) + \psi_2^{[2]} \sin(2 x) + \psi_2^{[3]} \sin(3x) + \psi_2^{[4]} \sin(4x)$, by \eqref{syslin2} we obtain $\eta_2^{[0]}=0$, together with the explicit formulae for $\eta_2^{[3]}$, $\eta_2^{[4]}$, $\psi_2^{[1]}$, $\psi_2^{[3]}$ and $\psi_2^{[4]}$ in \eqref{exp:Wil2}-\eqref{exp:c}. Moreover, we obtain that 
\begin{align} 
    \psi_2^{[2]} &= \frac{1}{8} \ck (8 \eta_2^{[2]}+1) .\label{eq:psi22formula}
\end{align}

\noindent
{\bf Third order in $ \e $.} 
It remains to determine 
$ c_2 $ and $\eta_2^{[2]}$ in \eqref{exp:Wil2}-\eqref{exp:c}.
 We get the linear system 
\begin{equation}\label{syslin3}
 \cB_0 \vet{\eta_3}{\psi_3} = \vet{c_1 (\psi_2)_x + c_2(\psi_1)_x 
 - (\psi_1)_x (\psi_2)_x - (\eta_1)_x^2 (\psi_1)_x \ck + \ck c_1 (\eta_1)_x^2 +  (\eta_1)_x (\eta_2)_x \ck^2-\frac{3}{2}\kappa (\eta_1)^2_x (\eta_1)_{xx}
 }{-c_1 (\eta_2)_x -c_2(\eta_1)_x - G_1(\eta_1)\psi_2- G_1(\eta_2)\psi_1 - G_2(\eta_1)\psi_1} \, . 
\end{equation}
System \eqref{syslin3} has 
a solution if and only if the right hand side is orthogonal to the Kernel of
$ \cB_0 $ given in \eqref{chk}. We use this condition in order to determine $c_2$ and $\eta_2^{[2]}$: namely, if we define by $\Pi_j$, $j=1,2$, the $L^2$-orthogonal projector on span$\, \{\cos(j x),\sin(j x)\} $, we have
\begin{align*}
    \Pi_1[ G_1(\eta_1)\psi_2] &= \psi_2^{[1]} \, \mathtt{a}_2 \, \sin(x), \qquad 
    \Pi_2[ G_1(\eta_1)\psi_2] = 0, \\
    \Pi_1[ G_1(\eta_2)\psi_1] &= \frac{(1+\kappa)\eta_2^{[2]}+(1+4\kappa) \eta_2^{[3]} \mathtt{a}_2 }{\sqrt{1+\kappa }} \sin(x) , \qquad 
    \Pi_2[ G_1(\eta_2)\psi_1] = 2\frac{(1+\kappa)\eta_2^{[3]}+(1+4\kappa) \eta_2^{[4]} \mathtt{a}_2 }{\sqrt{1+\kappa }} \sin(2x), \\
    \Pi_1[ G_2(\eta_1)\psi_1] &= \frac{1+\kappa +2 \mathtt{a}_2^2 (2+5 \kappa )}{4 \sqrt{1+\kappa} } \sin(x) , \qquad 
    \Pi_2[ G_2(\eta_1)\psi_1] = \frac{ \mathtt{a}_2 \left(1+\kappa +\mathtt{a}_2^2 (1+4 \kappa )\right)}{\sqrt{1+\kappa }} \sin(2x).
\end{align*}
By applying the projectors $\Pi_j$ to the right-hand side of \eqref{syslin3} and using the formulae for $\eta_2^{[0]}$, $\eta_2^{[3]}$, $\eta_2^{[4]}$, $\psi_2^{[1]}$, $\psi_2^{[3]}$ and $\psi_2^{[4]}$, we obtain

\begin{align*}
  &  \Pi_1 \vet{c_1 (\psi_2)_x + c_2(\psi_1)_x 
 - (\psi_1)_x (\psi_2)_x - (\eta_1)_x^2 (\psi_1)_x \ck + \ck c_1 (\eta_1)_x^2  +  (\eta_1)_x (\eta_2)_x \ck^2-\frac{3}{2}\kappa (\eta_1)^2_x (\eta_1)_{xx}
 }{-c_1 (\eta_2)_x -c_2(\eta_1)_x - G_1(\eta_1)\psi_2- G_1(\eta_2)\psi_1 - G_2(\eta_1)\psi_1} \\
 &= \vet{ \left[ -\frac{9}{32}+\ck c_2 \right] \cos(x) }{ \left[ c_2 - \ck (-1+\eta_2^{[2]}) \right] \sin(x) }, 
\end{align*}

\begin{align*}
  &  \Pi_2 \vet{c_1 (\psi_2)_x + c_2(\psi_1)_x 
 - (\psi_1)_x (\psi_2)_x - (\eta_1)_x^2 (\psi_1)_x \ck + \ck c_1 (\eta_1)_x^2 +  (\eta_1)_x (\eta_2)_x \ck^2-\frac{3}{2}\kappa (\eta_1)^2_x (\eta_1)_{xx}
 }{-c_1 (\eta_2)_x -c_2(\eta_1)_x - G_1(\eta_1)\psi_2- G_1(\eta_2)\psi_1 - G_2(\eta_1)\psi_1} \\
 &= \tts \vet{ \frac{1}{32} \left[ 3+16 \sqrt{6} c_2+24 \eta_2^{[2]} \right] \cos(2x)  }{\frac{1}{8} \left[ 8 c_2 +\sqrt{6} (11+2 \eta_2^{[2]})\right] \sin(2x)} .
\end{align*}

Therefore, the orthogonality conditions give the formulae for $c_2$ and $\eta_2^{[2]}$ respectively. Moreover, by using formula \eqref{eq:psi22formula}, we can also deduce the formula for $\psi_2^{[2]}$.

\medskip 
We now show that the map $\e \to (\eta_{\e,\tts},\psi_{\e,\tts})$ is analytic as a map $B(\e_\sigma) \to H^{\sigma,s}_{\mathtt{ev}}(\T)\times H^{\sigma,s}_{\mathtt{odd}}(\T)$ (defined above Theorem \ref{Thm: Wilton expansion}).
\begin{lemma}[Spatial analyticity]\label{thm:wilton-spatial-analyticity}
For every $\sigma>0$ and $s>7/2$, there exists $\e_{\sigma}>0$ such that the map
$ \e\mapsto  (\eta_{\e,\tts},\psi_{\e,\tts})\in H^{\sigma,s}_{\mathtt{ev}}(\T)\times H^{\sigma,s}_{\mathtt{odd}}(\T)$
is real analytic for $|\e|<\e_{\sigma}$.
\end{lemma}

\begin{proof}
We fix a branch $\tts\in\{\pm1\}$ and omit
the dependence on $\tts$ from the notation.
By Theorem \ref{thm:rs-wilton-amplitude},
the maps
\[
B(\e_0)\to H^2_{\mathtt{ev}}(\T)\times H^2_{\mathtt{odd}}(\T)\times \R \ , \quad 
 \e\mapsto (U_\e, c_\e) \ ,\ \  U_\e:=(\eta_\e,\psi_\e),
\]
are real analytic. Hence they admit convergent Taylor expansions
\[
U_\e=\sum_{n\ge1}\e^n U_n,
\qquad
c_\e=\ck+\sum_{n\ge1}\e^n c_n,
\qquad |\e|<\e_0 ,
\]
where $U_n=(\eta_n,\psi_n)$.
Fix $0<R<\e_0$. By the convergence of the
Taylor series, there exists $C>0$
such that
\begin{equation}\label{C R negative n}
\|U_n\|_{H^2(\T)\times H^2(\T)}
\le C R^{-n},
\qquad n\ge1.
\end{equation}
We now  prove that the Fourier support of
$U_n$ satisfies
\begin{equation}\label{supp}
\operatorname{supp}\widehat U_n
\subseteq
\{k\in\Z:|k|\le 2n\},
\qquad n\ge1.
\end{equation}
The proof is by induction on $n$.
For $n=1$, the claim follows from
\eqref{2411:1210}: the Fourier support of $U_1$
is contained in $\{\pm1,\pm2\}$.
Let $n\ge2$ and assume that \eqref{supp}
holds for all $1\le j<n$.
Writing the traveling-wave equations
\eqref{Sts} in the form
\[
\mathcal B_0 U_\e
=\mathcal N(U_\e,c_\e),
\]
where $\mathcal B_0$ is defined in
\eqref{cB0}, we have
\[
\mathcal N(U,c)=
\begin{pmatrix}
(c-\ck)\psi_x-\dfrac12\psi_x^2
+\dfrac{\eta_x^2}{2(1+\eta_x^2)}
(c-\psi_x)^2
+\kappa\partial_x
\left(
\dfrac{\eta_x}{\sqrt{1+\eta_x^2}}
-\eta_x
\right)
\\[2mm]
-(c-\ck)\eta_x
-\bigl(G(\eta)-G_0\bigr)\psi
\end{pmatrix}.
\]
The nonlinear map $\mathcal N$ has no
terms of degree less than two in
$(U,c-\ck)$.
Collecting the terms of order $\e^n$
therefore yields
\begin{equation}\label{eq:spatial-analyticity-recursion}
\mathcal B_0 U_n=F_n,
\end{equation}
where $F_n$ depends only on
$U_1,\ldots,U_{n-1}$ and
$c_1,\ldots,c_{n-1}$.
We claim that
\begin{equation}\label{eq:spatial-analyticity-Fn}
\operatorname{supp}\widehat F_n
\subseteq\{k\in\Z:|k|\le2n\}.
\end{equation}
To prove this, recall the expansion of
the Dirichlet--Neumann operator
in \eqref{expDiriNeu}. For every $m\ge1$,
$G_m(\eta)\psi$ is homogeneous of degree
$m$ in $\eta$ and linear in $\psi$,
and is obtained by compositions of
Fourier multipliers and multiplication
operators.
Consequently, each contribution to the
coefficient of order $\e^n$ in
$G_m(\eta_\e)\psi_\e$ contains
$m$ factors
$\eta_{j_1},\ldots,\eta_{j_m}$
and one factor $\psi_\ell$, where
\[
j_1+\cdots+j_m+\ell=n,
\qquad
j_1,\ldots,j_m,\ell\ge1.
\]
All these indices are strictly smaller
than $n$. By the induction hypothesis,
the Fourier support of such a product
is contained in
\[
\left\{
k\in\Z:
|k|\le
2(j_1+\cdots+j_m+\ell)
\right\}
=
\{k\in\Z:|k|\le2n\}.
\]
Here we use that multiplication adds
Fourier frequencies, whereas
differentiation and Fourier multipliers
preserve Fourier support.
The same argument applies to all the
remaining nonlinear terms in
$\mathcal N$. Indeed, their coefficients
at order $\e^n$ are finite sums of
products of lower-order profile
coefficients, whose amplitude orders
add up to at most $n$.
The coefficients of $c_\e-\ck$ are
independent of $x$ and therefore do
not affect Fourier support.
This proves \eqref{eq:spatial-analyticity-Fn}.

Taking Fourier coefficients in
\eqref{eq:spatial-analyticity-recursion},
we obtain
\[
(\mathcal B_0 U_n)^\wedge (k)=0,
\qquad |k|>2n.
\]
By Lemma \ref{lem:B0}, the only
resonant Fourier frequencies of the Fourier multiplier 
$\mathcal B_0$ are $k=\pm1,\pm2$, so it is 
invertible for $|k|>2n$.
Thus $\widehat U_n(k)=0$ for
$|k|>2n$, completing the induction
and proving \eqref{supp}.

We now use the Fourier support estimate
to establish convergence in analytic
Sobolev spaces.
For every $\sigma>0$ and $s>7/2$,
we have
\[
\begin{aligned}
\|U_n\|_{H^{\sigma,s}(\T)
\times H^{\sigma,s}(\T)}^2
&
=\sum_{|k|\le2n} \langle k\rangle^{2s} e^{2\sigma|k|} |\widehat U_n(k)|^2 \quad\le \langle2n\rangle^{2(s-2)} e^{4\sigma n} \sum_{|k|\le2n} \langle k\rangle^4 |\widehat U_n(k)|^2
\\
& =
\langle2n\rangle^{2(s-2)}
e^{4\sigma n}
\|U_n\|_{H^2(\T)\times H^2(\T)}^2
\stackrel{\eqref{C R negative n}}{\leq } C^2\langle2n\rangle^{2(s-2)}
\left(\frac{e^{2\sigma}}{R}\right)^{2n}  \ . 
\end{aligned}
\]
Hence
\[
\begin{aligned}
&\sum_{n\ge1}|\e|^n
\|U_n\|_{H^{\sigma,s}(\T)
\times H^{\sigma,s}(\T)} \le
C_R\sum_{n\ge1}
\langle2n\rangle^{s-2}
\left(
\frac{|\e|e^{2\sigma}}{R}
\right)^n
\end{aligned}
\]
and the series  converges
whenever $|\e|<\e_\sigma:=Re^{-2\sigma}$.
Therefore, the Taylor series of $U_\e$
converges absolutely in
$H^{\sigma,s}_{\mathtt{ev}}(\T)
\times H^{\sigma,s}_{\mathtt{odd}}(\T)$
for $|\e|<\e_\sigma$, proving the claim.
\end{proof}

\smallskip

\noindent{\bf Proof of Lemma \ref{lem:pfexp}.} 
We expand  the function $\mathfrak{p}(x)  = \e\mathfrak{p}_1(x) + \e^2 \mathfrak{p}_2(x) + \cO(\e^3)  $ defined by the fixed point equation \eqref{def:ttp}. 
Then 
$ \mathfrak{p}(x)
 = \mathfrak{H} \big[\e\eta_1 +\e^2\big(\eta_2 + 
 (\eta_1)_x \mathfrak{p}_1 \big)+\cO(\e^3)\big] $, 
and, using that $\mathfrak H \cos(kx) = \sin(kx)$, for any  $ k \in \N$, we get 
\begin{align}
 \mathfrak{p}_1(x) 
 &  = \sin(x) + \mathtt{a}_2 \sin(2 x) \,  \ , \quad 
 \mathfrak{p}_2(x) = \mathfrak{H}((\eta_1)_x \mathfrak{p}_1 +\eta_2 ) , \label{pfra2}
\end{align}
where $\mathfrak{p}_2$ satisfies \eqref{fexp}. Hence, we proved Lemma \ref{lem:pfexp}.
\\[1mm]
{\bf Proof of Lemma \ref{lem:pa.exp}.} 
In view of \eqref{exp:Wil2}, the expansions of the functions $B$, $V$  in \eqref{espV} and \eqref{espB} are
\begin{align}
 B&= : \e B_1(x) + \e^2 B_2(x) + \cO(\e^3),
 \qquad B_1(x) =  \ck \sin(x) + \ck \, \tts  \sin(2x) ,\label{espB1} \\
 \nonumber 
B_2(x) &= \frac{\ck  \tts}{4}  \sin(x) - \frac{9\ck}{4} \sin(2x)-\frac{11\ck  \tts}{2}\sin(3x) -\frac{3\ck}{2} \sin(4x) , 
\end{align}
  and
\begin{align}
V &= : \e V_1(x) + \e^2 V_2(x)  + \cO(\e^3) , 
\qquad V_1(x) = \ck \cos(x) + \ck \, \tts \cos(2x) , \label{espV1} \\
V_2(x) &= \ck +  \frac{3\ck  \tts }{4} \cos(x) -\frac{9\ck}{4} \cos(2x) -\frac{11\ck  \tts }{2}  \cos(3x) -\frac{3\ck}{2}  \cos(4x) . \nonumber 
\end{align}
In  view of  \eqref{def:pa}, denoting derivatives w.r.t $x$ with  a prime and suppressing dependence on $x$ when trivial, we have
\begin{align}
\ck+p_{\e,\tts}(x)  &= (\ck+\e c_1+\e^2 c_2 - V(x) - V'(x)\mathfrak{p}(x)+\cO(\e^3)) (1-\mathfrak{p}'(x)+(\mathfrak{p}'(x))^2+\cO(\e^3))\notag \\
&= \ck + \e \underbrace{(c_1 - V_1 -\ck \mathfrak{p}_1')}_{=: p_1}
 +\e^2 \underbrace{\big( c_2 + V_1\mathfrak{p}_1' - V_2 -c_1 \mathfrak{p}_1' - V_1'\mathfrak{p}_1 -\ck \mathfrak{p}_2' +\ck (\mathfrak{p}_1')^2 \big)}_{=:p_2} + \, \cO(\e^3) \, .\label{pino12imp}
\end{align} 
Similarly by \eqref{def:pa}
\begin{align}
 1+a_{\e,\tts}(&x) : = \frac{1}{1+\mathfrak{p}_x(x)} - (\ck+p_{\e,\tts}(x))B_x(x+\mathfrak{p}(x)) \notag \\
   = &  1+\e\underbrace{\big(-\mathfrak{p}_1'- \ck B_1'\big)}_{=:a_1} +\e^2\underbrace{\big((\mathfrak{p}_1')^2-\mathfrak{p}_2'- \ck B_2'-\ck B_1''\mathfrak{p}_1(x)-c_1 B_1' + B_1' V_1 + \ck  B_1'\mathfrak{p}_1' \big)}_{=:a_2}+\cO(\e^3)\, . \label{aino12imp}
 \end{align}
By \eqref{espV1}, \eqref{fexp}, \eqref{exp:Wil2}, \eqref{pfra2}, \eqref{espB1} we 
deduce that the functions $p_1 $, $p_2 $, $a_1 $, $a_2 $ in \eqref{pino12imp} and \eqref{aino12imp} have an expansion as in \eqref{pino1fd}-\eqref{aino2fd}.

In view of \eqref{exp:Wil2}-\eqref{exp:c}, the expansion of the function $l(x)$ in \eqref{esp l} is 
\begin{align}
 l(x)&= 1 + \e^2 l_2(x)  + \cO(\e^3)= 1+ \e^2  \left[ l_2^{[0]} + \sum_{j=1}^4 l_2^{[j]} \cos(j x) \right]+\cO(\e^3) , \label{espl1} \\
 l_2^{[0]} &= -\frac{3}{4} - 3 \mathtt{a}_2^2 =-\frac{3}{2}, \quad l_2^{[1]} = -3 \mathtt{a}_2=-\frac{3}{2}\tts , \quad l_2^{[2]} = \frac{3}{4} , \quad l_2^{[3]} = 3 \mathtt{a}_2 =\frac{3}{2}\tts, \quad l_2^{[4]} = 3 \mathtt{a}_2^2=\frac{3}{4} . \nonumber
\end{align}
In  view of  \eqref{def:Sigma g}, denoting derivatives w.r.t $x$ with  a prime and suppressing dependence on $x$ when trivial, we have
\begin{align*}
g_{\e,\tts}(x)  &= (1-\e \, \mathfrak{p}_1'+\e^2(-\mathfrak{p}_2'+(\mathfrak{p}_1')^2)+\cO(\e^3))(1+\e^2 l_2+\cO(\e^3))\notag \\
&= 1 + \e \underbrace{(-\mathfrak{p}_1')}_{=: g_1}+\e^2 \underbrace{\big(-\mathfrak{p}_2'+(\mathfrak{p}_1')^2+l_2)}_{=:g_2} + \, \cO(\e^3) \, .
\end{align*} 
By \eqref{espl1}, the expansions \eqref{exp:g1}--\eqref{exp:g2} follow.
Finally, the operator $\Sigma_{\e,\tts}$ in \eqref{def:Sigma g} expands as
$\Sigma_{\e,\tts} = \pa_{xx} + \e \Sigma_1 + \e^2 \Sigma_2 + \cO(\e^3)$
where
\begin{align*}
    \Sigma_1:=& (g_{1}-2\mathfrak{p}_1')\,\pa_{xx}+(g_1'-2\,\mathfrak{p}_1'')\pa_{x}-\mathfrak{p}_1''',\\ \notag
    \Sigma_2:=&(g_2-2\mathfrak{p}_2'+3(\mathfrak{p}_1')^2-2g_1\mathfrak{p}_1') \pa_{xx}+(g_2'-2\mathfrak{p}_2''-2g_1\mathfrak{p}_1''-2g_1'\mathfrak{p}_1'+6\mathfrak{p}_1''\mathfrak{p}_1')\pa_{x}\\ 
    &+2(\mathfrak{p}_1'')^2-\mathfrak{p}_2'''-\mathfrak{p}_1''g_1'+3\mathfrak{p}_1'''\mathfrak{p}_1'-g_1\mathfrak{p}_1'''.
\end{align*}
Inserting the expansions of $g_1,g_2$ in \eqref{exp:g1}, \eqref{exp:g2} and those of $\mathfrak{p}_1$, $\mathfrak{p}_2$ in \eqref{fexp}, \eqref{pfra2} proves the claimed expansions 
in \eqref{Sigma 1} with coefficients in \eqref{d1x}--\eqref{h2 h02}.

\qed

\section{Expansion of the Kato Basis}\label{secA1}

In this appendix we prove Lemma \ref{expansion of the basis F}.

\paragraph{Decomposition of $H^2(\mathbb{T},\mathbb{C})\times H^1(\mathbb{T},\mathbb{C})$ into invariant subspaces under $\sL_{0,0}$.}

Let $\kappa = \frac{1}{2}$. In addition to the basis of $\mathrm{ker}(\sL_{0,0})$, we also introduce the following vectors for $j \in \N$
\begin{align}
f_j^+ &:=\vet{ \ck^{-1/2} \cos(jx)}{ \ck^{1/2} \sin(jx)}, \quad f_j^-:=\vet{- \ck^{-1/2} \sin(jx)}{ \ck^{1/2} \cos(jx)} ,  \quad j \geq 3, \label{func fj+ fj-}\\
f_{-j}^+ &:=\vet{ \ck^{-1/2} \cos(jx)}{- \ck^{1/2} \sin(jx)}, \quad f_{-j}^- :=\vet{ \ck^{-1/2} \sin(jx)}{ \ck^{1/2} \cos(jx)} , \quad j \geq 1. \notag
\end{align}
Notice that for $\kappa = \frac{1}{2}$
\begin{align*}
    \sL_{0,0} f^{-}_{-1} = 2 \ck f^{+}_{-1} , \quad \sL_{0,0} f^{+}_{-1} = - 2 \ck f^{-}_{-1} , \quad
    \sL_{0,0} f^{-}_{-2} =  4 \ck f^{+}_{-2} , \quad \sL_{0,0} f^{+}_{-2} = - 4 \ck  f^{-}_{-2} .
\end{align*}
Moreover, recalling \eqref{eigenfunc of mathcall L00 -part1}-\eqref{eigenfunc f+0} and \eqref{func fj+ fj-}, we notice that
\begin{align*} 
    \cJ f^{-}_j = \vet{\ck^{1/2} \cos(jx)}{\ck^{-1/2} \sin(jx) } , &\quad \cJ f^{+}_j = \vet{ \ck^{1/2} \sin(jx) }{ - \ck^{-1/2} \cos(jx) }  , \quad j \in \N .
\end{align*}
One can directly check that the space $Y=H^2(\mathbb{T},\mathbb{C})\times H^1(\mathbb{T},\mathbb{C})$ admits the following decomposition (see Lemma B.2 in \cite{HM}).

\begin{lemma} \label{lem:Dec}
    The space $Y$ decomposes as $Y = \mathcal{V}_{0,0} \oplus \mathcal{U} \oplus \mathcal{W}_{Y}$, with
    $\mathcal{W}_{Y} = \overline{\bigoplus_{k=3}^{\infty} \mathcal{W}_k}^{Y}$
where the subspaces $\mathcal{V}_{0,0}, \mathcal{U}$, and $\mathcal{W}_k$, defined below, are invariant under $\sL_{0,0}$ and the following properties hold:

\begin{itemize}
    \item[(i)] $\mathcal{V}_{0,0} = \text{span}\{f^{-}_{0}, f^{+}_{0}, f^{-}_{1}, f^{+}_{1}, f^{-}_2, f^{+}_{2} \}$ is the generalized kernel of $\sL_{0,0}$. For any $\lambda \neq 0$ the operator $\sL_{0,0}-\lambda: \mathcal{V}_{0,0} \to \mathcal{V}_{0,0}$  is invertible, and
    \begin{align*}
        (\sL_{0,0}-\lambda)^{-1} f^{\sigma}_{j} &= -\frac{1}{\lambda} f^{\sigma}_{j} , \quad j =1,2, \quad \sigma = \pm, \\
        (\sL_{0,0}-\lambda)^{-1} f^{-}_{0} &= -\frac{1}{\lambda} f^{-}_{0} ,  \quad
        (\sL_{0,0}-\lambda)^{-1} f^{+}_{0} = -\frac{1}{\lambda} f^{+}_{0} + \frac{1}{\lambda^2} f^{-}_{0} .
    \end{align*}

    \item[(ii)] $\mathcal{U} := \text{span}\{ f_{-1}^{+}, f_{-1}^{-}, f_{-2}^{+}, f_{-2}^{-} \}$. For any $\lambda \notin \{ \pm \mathrm{i} 2 \ck , \pm \mathrm{i} 4 \ck \}$ the operator $\sL_{0,0}-\lambda: \mathcal{U} \to \mathcal{U}$  is invertible, and
    \begin{align*}
        (\sL_{0,0}-\lambda)^{-1} f^{+}_{-1} &= \frac{1}{\lambda^2 + 4 \ck^2} ( 2\ck f^{-}_{-1} -\lambda f^{+}_{-1} ) , \qquad (\sL_{0,0}-\lambda)^{-1} f^{-}_{-1} = \frac{1}{\lambda^2 + 4 \ck^2} ( -2\ck f^{+}_{-1} -\lambda f^{-}_{-1} )  , \\
        (\sL_{0,0}-\lambda)^{-1} f^{+}_{-2} &= \frac{1}{\lambda^2 + 16 \ck^2} ( 4\ck f^{-}_{-2} -\lambda f^{+}_{-2} ) , \qquad (\sL_{0,0}-\lambda)^{-1} f^{-}_{-2} = \frac{1}{\lambda^2 + 16 \ck^2} ( -4\ck f^{+}_{-2} -\lambda f^{-}_{-2} ) .
    \end{align*}

    \item[(iii)] Each subspace $\mathcal{W}_k := \text{span}\{ f_k^+, f_k^-, f_{-k}^+, f_{-k}^- \}$ is invariant under $\mathscr{L}_{0,0}$. Let
$ \mathcal{W}_{L^2} = \overline{\bigoplus_{k=3}^{\infty} \mathcal{W}_k}^{ L^2}$.
    There exists a sufficiently small $\delta>0$ such that for $|\lambda| < \delta$ the operator $\sL_{0,0}-\lambda:\cW_{Y} \to \cW_{L^2}$ is invertible, and for any $f \in \cW_{L^2}$
    \begin{align*}
        (\sL_{0,0}-\lambda)^{-1} f &= \left[ ( \ck^2 \partial_x^2 + |D| (1-\kappa \partial^2_{x}) )\right]^{-1}
        \begin{bmatrix}
            \ck \partial_x & - |D| \\
            1-\kappa \partial^2_{x} & \ck \partial_x
        \end{bmatrix}
        f + \lambda \varphi_{f}(\lambda,x) ,
    \end{align*}
    for some analytic function $\lambda \mapsto \varphi_{f}(\lambda,\cdot) \in Y$.

\end{itemize}
\end{lemma}

We provide the expansion of the basis $f_k^{\pm}(\mu,\e,\tts) = U_{\mu,\e,\tts} f_k^{\pm}$, $k = 0,1,2$, in \eqref{F basis set and f}, where $f_k^{\pm}$ defined in \eqref{eigenfunc of mathcall L00 2} belong to the subspace $\mathcal{V}_{0,0} := \operatorname{Rg}(P_{0,0})$. 

Following \cite[Sec. 3 and Appendix B]{BMV_ed}, for an operator $A=A(\mu,\e)$ we denote its Taylor coefficients as
\begin{align} \label{Taylor}
    A_{i,j} &\coloneqq \frac{1}{i! j!} (\partial_{\mu}^i \partial_{\e}^j A)(0,0) , \quad A_k \coloneqq A_k(\mu,\e) = \sum_{\substack{i+j=k \\ i,j \geq 0}} A_{i,j} \mu^i \e^j .
\end{align}
The projectors $P_{\mu,\e,\tts}$ in \eqref{Projection P mu e} admit the expansion (for the sake of notation, we omit the dependence on $\tts$ throughout this section)
\begin{equation}\label{expP}
P_{\mu,\e,\tts} = P_0 + P_1 + P_2 + P_3 + \cO_4 \, ,
\end{equation}
where
\begin{equation}\label{Psani}
\begin{aligned}
 & P_0 := P_{0,0} \, , \quad P_1 := \mathcal{P} \big[\cB_1\big], \quad P_2 := \mathcal{P} \big[\cB_2\big] + \mathcal{P} \big[\cB_1,\cB_1 \big] \, ,\\
  &P_3 := \mathcal{P} \big[\cB_3\big] + \mathcal{P} \big[\cB_2,\cB_1 \big] + \mathcal{P} \big[\cB_1,\cB_2 \big]  + \mathcal{P} \big[\cB_1,\cB_1,\cB_1 \big] \, ,
 \end{aligned} 
 \end{equation}
and 
\begin{equation*}
\begin{aligned}
&\mathcal{P} \big[A_1]: = \frac{1}{2\pi\im} \oint_\Gamma (\sL_{0,0}-\lambda)^{-1} \cJ A_1 (\sL_{0,0}-\lambda)^{-1} \de\lambda \, ,\quad \text{and for }k\geq 2 \\
&\mathcal{P} \big[A_1,\dots, A_k \big] := \frac{(-1)^{k+1} }{2\pi \im } \oint_\Gamma (\sL_{0,0}-\lambda)^{-1} \cJ A_1 (\sL_{0,0}-\lambda)^{-1} \dots \cJ A_k (\sL_{0,0}-\lambda)^{-1} \de\lambda \, .
 \end{aligned}
\end{equation*}
Similarly to \eqref{expP}, we can write
\begin{equation*} 
\sL_{\mu,\e,\tts} = \sL_0 + \sL_1 + \sL_2 + \sL_3 + \cO_4 \, ,  \quad
U_{\mu,\e,\tts} = U_0 + U_1 + U_2 + U_3 + \cO_4 \, ,
\end{equation*}
where we define the terms $\sL_j$ and $U_j$ as in \eqref{Taylor}.

We first Taylor-expand the transformation operators $U_{\mu,\e,\tts}$ defined in \eqref{U transformation operators}.
The next lemma follows as \cite[Lemma A.1]{BMV1}
\begin{lemma} \label{U mu epsilon P00}
    The first jets of $U_{\mu,\epsilon} P_{0,0}$ are
\begin{align} 
    &U_{0} P_{0} = P_{0,0}, \quad\quad U_{0,1} P_{0,0} = P_{0,1} P_{0,0}, \quad\quad U_{1,0} P_{0,0} = P_{1,0} P_{0,0},\notag\\ \label{A2}
    &U_{1,1} P_{0,0} = \left( P_{1,1} - \frac{1}{2} P_{0,0} P_{1,1} \right) P_{0,0}, \quad U_{0,2} P_{0,0} = \left(P_{0,2} - \frac{1}{2} P_{0,0} P_{0,2} \right) P_{0,0}, 
\end{align}
where
\begin{align} \label{A3}
    P_{0,1} &= \frac{1}{2\pi \im} \oint_{\Gamma} (\mathscr{L}_{0,0} - \lambda)^{-1} \sL_{0,1} (\mathscr{L}_{0,0} - \lambda)^{-1} \mathrm{d}\lambda\,,  \\ 
    P_{1,0} &= \frac{1}{2\pi \im} \oint_{\Gamma} (\mathscr{L}_{0,0} - \lambda)^{-1} \sL_{1,0} (\mathscr{L}_{0,0} - \lambda)^{-1} \mathrm{d}\lambda\,,\notag 
\end{align}
and
\begin{align} \label{A5a}
    P_{1,1} &= -\frac{1}{2\pi \im} \oint_{\Gamma} (\mathscr{L}_{0,0} - \lambda)^{-1} \sL_{1,0} (\mathscr{L}_{0,0} - \lambda)^{-1} \sL_{0,1} (\mathscr{L}_{0,0} - \lambda)^{-1} \mathrm{d}\lambda \\ \label{A5b}
    &\quad -\frac{1}{2\pi \im} \oint_{\Gamma} (\mathscr{L}_{0,0} - \lambda)^{-1} \sL_{0,1} (\mathscr{L}_{0,0} - \lambda)^{-1} \sL_{1,0} (\mathscr{L}_{0,0} - \lambda)^{-1} \mathrm{d}\lambda \\ \label{A5c}
    &\quad + \frac{1}{2\pi \im} \oint_{\Gamma} (\mathscr{L}_{0,0} - \lambda)^{-1} \sL_{1,1} (\mathscr{L}_{0,0} - \lambda)^{-1} \mathrm{d}\lambda,
\end{align}

\begin{align*}
P_{0,2}
&=
-\frac{1}{2\pi \im}
\oint_{\Gamma}
(\mathscr{L}_{0,0}-\lambda)^{-1}
\sL_{0,1}
(\mathscr{L}_{0,0}-\lambda)^{-1}
\sL_{0,1}
(\mathscr{L}_{0,0}-\lambda)^{-1}
\,\mathrm d\lambda
\notag\\
&\quad
+\frac{1}{2\pi \im}
\oint_{\Gamma}
(\mathscr{L}_{0,0}-\lambda)^{-1}
\sL_{0,2}
(\mathscr{L}_{0,0}-\lambda)^{-1}
\,\mathrm d\lambda .
\end{align*}

The operators $\sL_{0,1}$ and $\sL_{1,0}$ are
\begin{align} \label{A6}
    \sL_{0,1}=\begin{bmatrix}
        \partial_x\circ p_1(x)  &0\\
        -a_1(x)+\kappa \Sigma_1 &p_1(x)\circ \partial_x
    \end{bmatrix},~~\sL_{1,0}=\begin{bmatrix}
        0 &\mathrm{sgn}(D)+\Pi_0\\
        2\im \kappa \partial_{x} & 0
    \end{bmatrix},
\end{align}
where $\mathrm{sgn}(D)$ is defined in \eqref{D+mu}, \eqref{sgn D} 
and $a_1(x)$, $p_1(x)$, $\Sigma_1$ are given in Lemma \ref{lem:pa.exp}.

The operators $\sL_{1,1}$ and $\sL_{0,2}$ are
\begin{align} \label{A8}
    \sL_{1,1}=
    \begin{bmatrix}
        \im p_1(x) & 0\\
        \im \kappa\left(2 d_1(x) \pa_x + e_1(x)\right)        & \im p_1(x)
    \end{bmatrix},~~
    \sL_{0,2}=\frac{1}{2}
    \begin{bmatrix}
        2\partial_x\circ p_2(x)  &0\\
        -2a_2(x)+2\kappa \Sigma_2 &2p_2(x)\circ \partial_x
    \end{bmatrix},
\end{align}
with $d_1(x), e_1(x)$ in \eqref{d1x}. 
\end{lemma}

By the Lemma \ref{U mu epsilon P00}, we have the Taylor expansion
\begin{equation} \label{the expansion of f sigma mu}
    \begin{aligned}
    f_k^\sigma (\mu, \epsilon) &= f_k^\sigma + \epsilon P_{0,1} f_k^\sigma + \mu P_{1,0} f_k^\sigma  \\
    &+ \mu \epsilon \big( P_{1,1} - \frac{1}{2} P_{0,0} P_{1,1} \big) f_k^\sigma +  \epsilon^2 \big( P_{0,2} - \frac{1}{2} P_{0,0} P_{0,2} \big) f_k^\sigma+ \mathcal{O}(\mu^2, \epsilon^3,\mu\e^2).
\end{aligned}
\end{equation}

\paragraph{Action of $P_{\ell,j}$ on the unperturbed vectors: first order jets.}
We now collect how the operators $P_{\ell,j}$ act on the vectors $(f^{\sigma}_k)_{ \substack{ \sigma = \pm , \\ k=0,1,2} }$. We first consider the first order jets $ P_{0,1} $ and $P_{1,0} $ of $ P_1 $. 


We shall  use the following formulas obtained by \eqref{A6} and \eqref{eigenfunc of mathcall L00 2}:
   \begin{equation}\label{A15prime}
\begin{aligned}
\sL_{0,1} f_2^+
&= \frac{\ck}{2}\Bigl(
    -f_1^- + f_{-1}^- + \tts f_2^-
    - f_3^- + 5f_{-3}^- + 8\tts f_{-4}^-
    \Bigr),\\
\sL_{0,1} f_2^-
&= \frac{\ck}{2}\Bigl(
    f_1^+ + f_{-1}^+ - \tts f_2^+
    + f_3^+ + 5f_{-3}^+ + 8\tts f_{-4}^+
    \Bigr),\\
\sL_{0,1} f_1^+
&= \frac{\ck}{4}\Bigl(
    3\tts f_{-1}^- - 2f_2^- + 6f_{-2}^-
    - \tts f_3^- + 11\tts f_{-3}^-
    \Bigr),\\
\sL_{0,1} f_1^-
&= \frac{\ck}{4}\Bigl(
    -2\tts f_1^+ - 3\tts f_{-1}^+
    + 2f_2^+ + 6f_{-2}^+
    + \tts f_3^+ + 11\tts f_{-3}^+
    \Bigr),\\
\sL_{0,1} f_0^+
&= 2\ck^{3/2} f_{-1}^-
    + 4\ck^{3/2}\tts f_{-2}^-,
\qquad
\sL_{0,1} f_0^- = 0,
\end{aligned}
\end{equation}
and
\begin{equation}\label{A15dot}
\begin{aligned}
\sL_{1,0} f_2^+
&= -\frac{\im\ck}{6}\bigl(7f_2^+ - f_{-2}^+\bigr),
\qquad
\sL_{1,0} f_2^-
= -\frac{\im\ck}{6}\bigl(7f_2^- + f_{-2}^-\bigr),\\
\sL_{1,0} f_1^+
&= -\frac{\im\ck}{6}\bigl(5f_1^+ + f_{-1}^+\bigr),
\qquad
\sL_{1,0} f_1^-
= -\frac{\im\ck}{6}\bigl(5f_1^- - f_{-1}^-\bigr),\\
\sL_{1,0} f_0^+
&= 0,
\qquad
\sL_{1,0} f_0^- = f_0^+.
\end{aligned}
\end{equation}
 By Lemma \ref{lem:Dec}, a straightforward calculation proves the following results.
\begin{lemma} \label{C3 lemma}
   One has
    \begin{align*}
& \frac{1}{2\pi\im}\oint_\Gamma
\frac{1}{\lambda}(\sL_{0,0}-\lambda)^{-1} f_3^+ \,\mathrm{d}\lambda
=\frac{1}{6\ck}\bigl(19f_3^-+f_{-3}^-\bigr),
\qquad 
\frac{1}{2\pi\im}\oint_\Gamma
\frac{1}{\lambda}(\sL_{0,0}-\lambda)^{-1} f_3^-
\,\mathrm{d}\lambda =\frac{1}{6\ck}\bigl(-19f_3^++f_{-3}^+\bigr),\\
& \frac{1}{2\pi\im}\oint_\Gamma
\frac{1}{\lambda}(\sL_{0,0}-\lambda)^{-1} f_{-3}^+
\,\mathrm{d}\lambda
=\frac{1}{6\ck}\bigl(f_3^-+f_{-3}^-\bigr), \qquad 
\frac{1}{2\pi\im}\oint_\Gamma
\frac{1}{\lambda}(\sL_{0,0}-\lambda)^{-1} f_{-3}^-
\,\mathrm{d}\lambda
=\frac{1}{6\ck}\bigl(f_3^+-f_{-3}^+\bigr),\\[1ex]
& \frac{1}{2\pi\im}\oint_\Gamma
\frac{1}{\lambda}(\sL_{0,0}-\lambda)^{-1} f_4^+
\,\mathrm{d}\lambda
=\frac{1}{8\ck}\bigl(9f_4^-+f_{-4}^-\bigr), \qquad 
\frac{1}{2\pi\im}\oint_\Gamma
\frac{1}{\lambda}(\sL_{0,0}-\lambda)^{-1} f_4^-
\,\mathrm{d}\lambda
=\frac{1}{8\ck}\bigl(-9f_4^++f_{-4}^+\bigr),\\
& \frac{1}{2\pi\im}\oint_\Gamma
\frac{1}{\lambda}(\sL_{0,0}-\lambda)^{-1} f_{-4}^+
\,\mathrm{d}\lambda
=\frac{1}{8\ck}\bigl(f_4^-+f_{-4}^-\bigr), \qquad 
\frac{1}{2\pi\im}\oint_\Gamma
\frac{1}{\lambda}(\sL_{0,0}-\lambda)^{-1} f_{-4}^-
\,\mathrm{d}\lambda
=\frac{1}{8\ck}\bigl(f_4^+-f_{-4}^+\bigr). 
\end{align*}
Also, one has
\begin{align*}
& \frac{1}{2\pi\im}\oint_\Gamma
\frac{1}{\lambda}(\sL_{0,0}-\lambda)^{-1} f_{-1}^+
\,\mathrm{d}\lambda
=\frac{1}{2\ck}f_{-1}^-, \qquad 
\frac{1}{2\pi\im}\oint_\Gamma
\frac{1}{\lambda}(\sL_{0,0}-\lambda)^{-1} f_{-1}^-
\,\mathrm{d}\lambda
=-\frac{1}{2\ck}f_{-1}^+,\\
& \frac{1}{2\pi\im}\oint_\Gamma
\frac{1}{\lambda}(\sL_{0,0}-\lambda)^{-1} f_{-2}^+
\,\mathrm{d}\lambda
=\frac{1}{4\ck}f_{-2}^-, \qquad 
\frac{1}{2\pi\im}\oint_\Gamma
\frac{1}{\lambda}(\sL_{0,0}-\lambda)^{-1} f_{-2}^-
\,\mathrm{d}\lambda
=-\frac{1}{4\ck}f_{-2}^+.
\end{align*}
\end{lemma}
    
We finally calculate $P_{0,1} f^\sigma_k$ and $P_{1,0} f^\sigma_k$.

\begin{lemma} {\bf (First order jets)} The action of the jets $ P_{0,1} $ and $P_{1,0} $ of $ P_1 $
 in \eqref{Psani} on the basis $(f^{\sigma}_k)_{ \substack{ \sigma = \pm , \\ k=0,1,2} }$ is given by 
    \begin{equation} \label{A16prime}
        \begin{aligned}
         &   P_{0,1}f^\pm_2 = \alpha^{\pm}_{2} f^\pm_{-1} + \beta^{\pm}_{2} f^\pm_{-2} + \gamma^{\pm}_{2} f^{\pm}_{3} + \delta^{\pm}_{2} f^{\pm}_{-3} +  \theta^{\pm}_{2} f^{\pm}_{4} + \phi^{\pm}_{2} f^{\pm}_{-4}  ,  \\
         &   P_{0,1}f^\pm_1 = \alpha^{\pm}_{1} f^\pm_{-1} + \beta^{\pm}_{1} f^\pm_{-2}  + \gamma^{\pm}_{1} f^{\pm}_{3} + \delta^{\pm}_{1} f^{\pm}_{-3} ,  \\
        &    P_{0,1}f^+_0 = \alpha^{+}_{0} f^+_{-1} + \beta^{+}_{0} f^+_{-2} ,  \qquad      P_{0,1}f^-_0 = 0 , 
        \end{aligned}
    \end{equation}
    where the coefficients are given by 
    \begin{equation*} 
        \begin{aligned}
        &    \alpha^{\pm}_{2}  = \pm \frac{1}{4} ,
        \qquad \beta^{\pm}_{2} = 0 , 
        \qquad     \gamma^{\pm}_{2}  = -2 , 
        \qquad     \delta^{\pm}_{2}  = \pm \frac{1}{2} , 
        \qquad    \theta^{\pm}_{2} = -\frac{1}{2} \tts , 
        \qquad    \phi^{\pm}_{2} =  \pm  \frac{1}{2}  \tts, \\
    &   \alpha^{\pm}_{1} = \frac{3}{8} \tts , 
       \qquad   \beta^{\pm}_{1}  = \pm \frac{3}{8} , 
        \qquad  \gamma^{\pm}_{1}  = -\frac{5}{4} \tts , 
      \qquad  \delta^{\pm}_{1}  = \pm \frac{1}{2} \tts , 
    \qquad  \alpha^{+}_{0} =  \ck^{1/2} , 
    \qquad \beta^{+}_{0} =  \ck^{1/2} \tts ,
        \end{aligned}
    \end{equation*}
    and    
    \begin{equation} \label{A16dot}
        \begin{aligned}
            P_{1,0}f^+_0 =0, \qquad P_{1,0}f^-_0=0,  \qquad 
            P_{1,0}f^+_j = \im \rho^+_{j} f^-_{-j} , \qquad  P_{1,0}f^-_j= \im \rho^-_{j} f^+_{-j}  , \quad j=1,2,
        \end{aligned}
    \end{equation}
    where
    \begin{equation*}
        \begin{aligned}
            \rho^\pm_{j} = \frac{1}{2j\ck} \left(  \frac{\ck}{2} - \frac{j\kappa}{\ck} \right) = -\frac{1}{6} + \frac{1}{4j} ,  \quad j=1,2.
        \end{aligned}
    \end{equation*}

\end{lemma}

\begin{proof}

We first calculate $P_{0,1} f^+_2$. By \eqref{A3}, Lemma \ref{lem:Dec}, and \eqref{A15prime} we deduce
    \begin{align*}
        P_{0,1} f^+_2=-\frac{1}{2\pi \im}\oint_{\Gamma} \frac{1}{\lambda} \left(\mathscr{L}_{0,0}-\lambda\right)^{-1} \sL_{0,1} f^+_2 \mathrm{d}\lambda ,
    \end{align*}
where
\begin{align*}
    \sL_{0,1} f^+_2 = \frac{\ck}{2}\Bigl(
    -f_1^- + f_{-1}^- + \tts f_2^-
    - f_3^- + 5f_{-3}^- + 8\tts f_{-4}^-
    \Bigr).
\end{align*}
Therefore by Lemma \ref{C3 lemma}, we obtain
\begin{align*}
         P_{0,1}f^+_2 = \alpha^{+}_{2} f^+_{-1} + \beta^{+}_{2} f^+_{-2} + \gamma^{+}_{2} f^{+}_{3} + \delta^{+}_{2} f^{+}_{-3} +  \theta^{+}_{2} f^{+}_{4} + \phi^{+}_{2} f^{+}_{-4}.
\end{align*}
so that by the residue theorem we get our first identity in \eqref{A16prime}. Similarly, one may compute the other identities in \eqref{A16prime}.

Now we prove the identities in \eqref{A16dot}. In view of \eqref{A15dot}, we have by Lemma \ref{lem:Dec} that

\begin{equation*}
    \begin{aligned}
        P_{1,0}f^+_0=&\frac{1}{2\pi \im}\oint_\Gamma \left(\mathscr{L}_{0,0}-\lambda\right)^{-1}\sL_{1,0}\left(\frac{1}{\lambda^2}f^-_0-\frac{1}{\lambda}f^+_0\right) \mathrm{d}\lambda=0,\\
        P_{1,0}f^-_0=&\frac{1}{2\pi \im}\oint_\Gamma \left(\mathscr{L}_{0,0}-\lambda\right)^{-1}\sL_{1,0}\left(\frac{-1}{\lambda}f^-_0\right) \mathrm{d}\lambda=0.
    \end{aligned}
\end{equation*}
Finally, for $j=1,2$ we have
\begin{align*}
    P_{1,0} f^\pm_j=\frac{1}{2\pi \im}\oint_\Gamma -\frac{1}{\lambda}\left(\mathscr{L}_{0,0}-\lambda\right)^{-1} \sL_{1,0} f^{\pm}_j \mathrm{d}\lambda.
\end{align*}
By \eqref{A15dot} and Lemma~\ref{C3 lemma}, all formulae in \eqref{A16dot} are proven.
\end{proof}

So far we have obtained the linear terms of the expansions \eqref{f+2}, \eqref{f+1}, \eqref{f-2}, \eqref{f-1}, \eqref{f+0}, and \eqref{f-0}. We now provide further information about the expansion of the basis at $\mu = 0$. 

\begin{lemma} 
     The basis $\{ f_k^{\sigma}(0,\e,\tts), \, k = 0,1,2 ,\, \sigma = \pm \}$ is real. For any $\epsilon$ it results $f_0^-(0,\e,\tts) \equiv f_0^-$. The property \eqref{f-0=01} holds.
\end{lemma}

\begin{proof}
    The reality of the basis $f_k^\sigma(0,\e,\tts)$ is a consequence of Lemma \ref{properties of U and P}-(iii). Since, recalling \eqref{mathscr L mu e} and \eqref{mathcal B mu e}, $\mathscr{L}_{0,\epsilon} f_0^- = 0$ for any $\epsilon$, we deduce $(\mathscr{L}_{0,\epsilon} - \lambda)^{-1} f_0^- = -\frac{1}{\lambda} f_0^-$ and then, using also the residue theorem,
\begin{align} \label{pa e f01}
P_{0,\e,\tts} f_0^- = -\frac{1}{2\pi \im} \oint_\Gamma (\mathscr{L}_{0,\e,\tts} - \lambda)^{-1} f_0^- \, \mathrm{d}\lambda = f_0^-.
\end{align}
In particular $P_{0,\e,\tts} f_0^- = P_{0,0} f_0^-$, for any $\epsilon$ and we get, by \eqref{U transformation operators}, $f_0^-(0,\epsilon) = U_{0,\epsilon} f_0^- = f_0^-$, for any $\epsilon$.

Let us prove property \eqref{f-0=01}. In view of \eqref{Parity structure} and since the basis is real, we know that
\[
f_k^+(0,\e,\tts) = \begin{bmatrix}\mathit{even}(x)\\ \mathit{odd}(x)\end{bmatrix}, \quad f_k^-(0,\e,\tts) = \begin{bmatrix}\mathit{odd}(x)\\ \mathit{even}(x)\end{bmatrix},
\]
for any $k = 0, 1, 2$. By Lemma \ref{F is symplectic and reversible} the basis $\{f_k^\sigma(0,\e,\tts)\}$ is symplectic (c.f. \eqref{basis is symplectic}) and, since
\[
\mathcal{J} f_0^-(0,\e,\tts) = \mathcal{J} f_0^- = \begin{bmatrix}1\\0\end{bmatrix},
\]
for any $\epsilon$, we get, for $k=1,2$,
\[
0 = (\mathcal{J} f_0^-(0,\e,\tts),\, f_k^+(0,\e,\tts)) = \left(\begin{bmatrix}1\\0\end{bmatrix}, f_k^+(0,\e,\tts)\right), \qquad 
1 = (\mathcal{J} f_0^-(0,\e,\tts),\, f_0^+(0,\e,\tts)) = \left(\begin{bmatrix}1\\0\end{bmatrix}, f_0^+(0,\e,\tts)\right).
\]
Thus the first component of both $f_k^+(0,\e,\tts)$, $k=1,2$, and $f_0^+(0,\e,\tts) - \begin{bmatrix}1\\0\end{bmatrix}$ has zero average, proving \eqref{f-0=01}.
\end{proof}

In the following we denote by $f_j^{\sigma}(\mu,0)$ the vector $f_j^{\sigma}(\mu,0,\tts)$ ($\sigma=\pm$,$j=0,1,2$), as $f_j^{\sigma}(\mu,0,\tts)$ does not depend on $\tts$. We also denote by $\sL_{\mu,0}$ the operator $\sL_{\mu,\e,\tts}|_{\e=0}$, and similarly for $P_{\mu,0}$, $U_{\mu,0}$ and $\cV_{\mu,0}$.

\begin{lemma} \label{fmu0 eq f0} 
For any small $\mu$, we have $f_0^+(\mu,0) \equiv f_0^+$ and $f_0^-(\mu,0) \equiv f_0^-$. Moreover, for any $j=1,2$ the vectors $f_j^+(\mu,0)$ and $f_j^-(\mu,0)$ have both components with zero space average.
\end{lemma}
\begin{proof}
The operator 
\begin{equation*}
\mathscr{L}_{\mu,0} = 
\begin{bmatrix}
\ck\partial_x & |D+\mu| \\
-1+\kappa(\pa^2_{x}+2\im \mu\pa_x-\mu^2) & \ck\partial_x
\end{bmatrix}
\end{equation*}
leaves invariant the subspace $\mathcal{Z} := \text{span}\{ f_0^+, f_0^- \}$ since 
$\mathscr{L}_{\mu,0} f_0^+ = -(1+\kappa\mu^2)f_0^-$ and $\mathscr{L}_{\mu,0} f_0^- = \mu f_0^+$. 
The operator $\mathscr{L}_{\mu,0}|_{\mathcal{Z}}$ has the two eigenvalues $\pm \im \sqrt{(1+\kappa\mu^2)\mu}$, 
which, for small $\mu$, lie inside the loop $\Gamma$ around $0$ in \eqref{Projection P mu e}. Then, by \eqref{spectrum separated by Gamma}, 
we have $\mathcal{Z} \subseteq \mathcal{V}_{\mu,0} = \text{Rg}(P_{\mu,0})$ and 
\begin{align} \label{pa mu f01}
P_{\mu,0} f_0^\pm = f_0^\pm, \quad 
f_0^\pm (\mu,0) = U_{\mu,0} f_0^\pm = f_0^\pm, 
\quad \text{for any $\mu$ small.}
\end{align}

The basis $\{ f_k^\sigma (\mu,0) \}$ is symplectic (c.f. \eqref{basis is symplectic}). Then, since 
$\mathcal{J} f_0^+ = \begin{bmatrix} 0 \\ -1 \end{bmatrix}$ 
and $\mathcal{J} f_0^- = \begin{bmatrix} 1 \\ 0 \end{bmatrix}$, we have

\[
\begin{aligned}
0 &= \big( \mathcal{J} f_0^+ (\mu,0), f_j^\sigma (\mu,0) \big) 
= \Big( \begin{bmatrix} 0 \\ -1 \end{bmatrix}, f_j^\sigma (\mu,0) \Big), \quad 
0 = \big( \mathcal{J} f_0^- (\mu,0), f_j^\sigma (\mu,0) \big) 
= \Big( \begin{bmatrix} 1 \\ 0 \end{bmatrix}, f_j^\sigma (\mu,0) \Big), \quad j=1,2,
\end{aligned} 
\]
namely both the components of $f_j^\pm (\mu,0)$ ($j=1,2$) have zero average.
\end{proof}

Next, we  consider the $\cO(\mu \epsilon)$ term in the expansion \eqref{the expansion of f sigma mu}.

\begin{lemma} 
The vectors $\partial_\mu \partial_\epsilon f_k^{\sigma}(0,0) = \left( P_{1,1} - \frac{1}{2} P_{0,0} P_{1,1} \right) f_k^{\sigma}$ satisfy
\begin{equation} \label{A17}
\begin{aligned}
(\partial_\mu \partial_\epsilon f_j^+(0,0)) &= \im \begin{bmatrix} \mathit{odd}(x) \\ \mathit{even}(x) \end{bmatrix}, \qquad (\partial_\mu \partial_\epsilon f_j^-(0,0)) = \im \begin{bmatrix} {\mathit{even}_0}(x) \\ \mathit{odd}(x) \end{bmatrix}, \quad j=1,2, \\
(\partial_\mu \partial_\epsilon f_0^+(0,0)) &= \frac{1}{4} \ck^{-3/2} f^+_{-1} + \frac{\tts}{16} \ck^{-3/2} f^+_{-2}+\im \begin{bmatrix} \mathit{odd}(x) \\ \mathit{even}_0(x) \end{bmatrix}, \\
(\partial_\mu \partial_\epsilon f_0^-(0,0)) &= \frac{1}{2} \ck^{-1/2} f^-_{-1} + \frac{1}{4} \ck^{-1/2} \tts f^-_{-2} +\im \begin{bmatrix} \mathit{even}_0(x) \\ \mathit{odd}(x) \end{bmatrix}.
\end{aligned}
\end{equation}
\end{lemma}

\begin{proof}
 Recall \eqref{A6} and decompose the Fourier multiplier operator $\sL_{1,0}$ as
\begin{align*}
    \sL_{1,0}=\sL_{1,0}^{I}+\sL_{1,0}^{II}\,,\quad \sL_{1,0}^{I}:=\begin{bmatrix}
        0 &\mathrm{sgn}(D)\\
        2\im \kappa \partial_{x} &0
    \end{bmatrix}\,,\quad
    \sL_{1,0}^{II}:=\begin{bmatrix}
        0 &\Pi_0\\
       0 &0
    \end{bmatrix}.
        \end{align*}
Moreover, we write 
 \begin{align*}
      P_{1,1} = \eqref{A5a}^{I}+\eqref{A5a}^{II}+\eqref{A5b}^{I}+\eqref{A5b}^{II}+\eqref{A5c}
 \end{align*}
 where
\begin{align*}
   \eqref{A5a}^{I} &:= -\frac{1}{2\pi \im} \oint_{\Gamma} (\mathscr{L}_{0,0} - \lambda)^{-1} \sL_{1,0}^{I} (\mathscr{L}_{0,0} - \lambda)^{-1} \sL_{0,1} (\mathscr{L}_{0,0} - \lambda)^{-1} \mathrm{d}\lambda\,,\\
  \eqref{A5a}^{II} &:= -\frac{1}{2\pi \im} \oint_{\Gamma} (\mathscr{L}_{0,0} - \lambda)^{-1} \sL_{1,0}^{II} (\mathscr{L}_{0,0} - \lambda)^{-1} \sL_{0,1} (\mathscr{L}_{0,0} - \lambda)^{-1} \mathrm{d}\lambda\,,\\
  \eqref{A5b}^{I} &:= -\frac{1}{2\pi \im} \oint_{\Gamma} (\mathscr{L}_{0,0} - \lambda)^{-1} \sL_{0,1} (\mathscr{L}_{0,0} - \lambda)^{-1}\sL_{1,0} ^{I}(\mathscr{L}_{0,0} - \lambda)^{-1}  \mathrm{d}\lambda\,,\\
  \eqref{A5b}^{II} &:= -\frac{1}{2\pi \im} \oint_{\Gamma} (\mathscr{L}_{0,0} - \lambda)^{-1} \sL_{0,1} (\mathscr{L}_{0,0} - \lambda)^{-1}\sL_{1,0} ^{II}(\mathscr{L}_{0,0} - \lambda)^{-1}  \mathrm{d}\lambda\,.
        \end{align*}
 Note that $\eqref{A5a}^{I}$, $\eqref{A5b}^{I}$ and \eqref{A5c} are purely imaginary. Indeed, $\sL_{1,0}^{I}$ is purely imaginary, whereas $\mathscr{L}_{0,1}$ in \eqref{A6} is real, and $\sL_{1,1}$ in \eqref{A8} is purely imaginary. Consequently, when applied to the real vectors $f_k^{\sigma}$, for $k = 0,1,2$ and $\sigma = \pm$, these operators produce purely imaginary vectors.
 
We first compute $(\partial_\mu \partial_\epsilon f_j^{\sigma}(0,0))$ ($\sigma =\pm$, $j=1,2$). Using Lemma \ref{lem:Dec} and \eqref{A15prime}-\eqref{A15dot} we obtain
\begin{align*}
   \eqref{A5a}^{II} f_j^+ = \frac{1}{2\pi \im} \oint_{\Gamma} (\mathscr{L}_{0,0} - \lambda)^{-1} \sL_{1,0}^{II} (\mathscr{L}_{0,0} - \lambda)^{-1} \frac{1}{\lambda} \, \sL_{0,1} f^+_j \mathrm{d}\lambda .
\end{align*}
Notice that the above term vanishes, since $(\mathscr{L}_{0,0} - \lambda)^{-1}  \sL_{0,1} f^+_j$ has zero average, and since $\sL_{1,0}^{II}$ maps terms with zero average to zero.

Moreover, $\eqref{A5b}^{II} f_j^+=0$, since $\sL_{1,0} ^{II}(\mathscr{L}_{0,0} - \lambda)^{-1}f_j^+=0$. Combining the above identities, we conclude that $  P_{1,1}f_j^+$
is a purely imaginary vector. Since $P_{0,0}$ is a real operator, it follows that $ \left( P_{1,1} - \frac{1}{2} P_{0,0} P_{1,1} \right)f_j^+$ is also purely imaginary. 
Finally, \eqref{Parity structure} implies that $\partial_\mu \partial_\epsilon f_j^+(0,0)$ has the structure claimed in \eqref{A17}. The structure for $\partial_\mu \partial_\epsilon f_j^-(0,0)$ follows by a similar argument, using that 
\begin{align*}
   \eqref{A5a}^{II} f_j^- = \frac{1}{2\pi \im} \oint_{\Gamma} (\mathscr{L}_{0,0} - \lambda)^{-1} \sL_{1,0}^{II} (\mathscr{L}_{0,0} - \lambda)^{-1} \frac{1}{\lambda} \, \sL_{0,1} f^-_j \mathrm{d}\lambda = 0 ,
\end{align*}
since by \eqref{A15prime} and by Lemma \ref{lem:Dec} we can deduce that $(\mathscr{L}_{0,0} - \lambda)^{-1}  \sL_{0,1} f^-_j$ has zero average. 

Next, we prove that  $\partial_\mu \partial_\epsilon f_0^+(0,0)$ has zero average. By Lemma \ref{lem:Dec} and \eqref{A15prime}-\eqref{A15dot}, we have
\begin{align*}
   \eqref{A5a}^{I}f_0^+& = -\frac{1}{2\pi \im} \oint_{\Gamma} (\mathscr{L}_{0,0} - \lambda)^{-1} \sL_{1,0}^{I} (\mathscr{L}_{0,0} - \lambda)^{-1} \sL_{0,1} \left(\frac{-1}{\lambda}f_0^++\frac{1}{\lambda^2}f_0^-\right) \mathrm{d}\lambda\\
   & = -\frac{1}{2\pi \im} \oint_{\Gamma} (\mathscr{L}_{0,0} - \lambda)^{-1} \sL_{1,0}^{I} (\mathscr{L}_{0,0} - \lambda)^{-1}  \frac{-1}{\lambda} \left( 2 \ck^{3/2} f^-_{-1} + 4 \ck^{3/2} \tts f^-_{-2} \right) \mathrm{d}\lambda\,. \end{align*}
Since $ (\mathscr{L}_{0,0} - \lambda)^{-1} $ and $\sL_{1,0}^{I}$ are Fourier multiplier operators, they preserve the zero-average property of functions. Therefore, $\eqref{A5a}^{I}f_0^+$ has zero average.
Moreover, 
\begin{align*}
   \eqref{A5a}^{II}f_0^+& = 0 \ \ \text{since}\ \  \sL_{1,0}^{II} (\mathscr{L}_{0,0} - \lambda)^{-1}  \left( 2 \ck^{3/2} f^-_{-1} + 4 \ck^{3/2} \tts f^-_{-2} \right) =0\,.
   \end{align*}
Next, $\eqref{A5b}^{I}f_0^+=0$ because $\sL_{1,0}^{I}f_0^+=\sL_{1,0}^{I}f_0^-=0$.
Using also that $\sL_{1,0}^{II}f_0^+=0$ and $\sL_{1,0}^{II}f_0^-=f_0^+$, we obtain
\begin{align*}
   \eqref{A5b}^{II}f_0^+& =-\frac{1}{2\pi \im} \oint_{\Gamma} (\mathscr{L}_{0,0} - \lambda)^{-1}\sL_{0,1}  (\mathscr{L}_{0,0} - \lambda)^{-1}  \sL_{1,0}^{II}\left(\frac{-1}{\lambda}f_0^++\frac{1}{\lambda^2}f_0^-\right) \mathrm{d}\lambda\\
   & =-\frac{1}{2\pi \im} \oint_{\Gamma} (\mathscr{L}_{0,0} - \lambda)^{-1}\sL_{0,1}  (\mathscr{L}_{0,0} - \lambda)^{-1}  \frac{1}{\lambda^2}f_0^+ \mathrm{d}\lambda\\
   & =-\frac{1}{2\pi \im} \oint_{\Gamma} (\mathscr{L}_{0,0} - \lambda)^{-1}\sL_{0,1} \frac{1}{\lambda^2}\left(\frac{-1}{\lambda}f_0^++\frac{1}{\lambda^2}f_0^-\right) \mathrm{d}\lambda\\
   & = -\frac{1}{2\pi \im} \oint_{\Gamma} (\mathscr{L}_{0,0} - \lambda)^{-1}  \frac{-1}{\lambda^3} \left( 2 \ck^{3/2} f^-_{-1} + 4 \ck^{3/2} \tts f^-_{-2} \right) \mathrm{d}\lambda= \frac{1}{4} \ck^{-3/2} f^+_{-1} + \frac{\tts}{16} \ck^{-3/2} f^+_{-2} \,,
   \end{align*}
where the last equality follows from Lemma \ref{lem:Dec} and from the Residue Theorem. Finally, by Lemma \ref{lem:Dec} and by \eqref{A8}, where $p_1(x),d_1,e_1$ are given by \eqref{pino1fd} and \eqref{d1x}, we obtain
\begin{align*}
   \eqref{A5c}f_0^+& = \frac{1}{2\pi \im} \oint_{\Gamma} (\mathscr{L}_{0,0} - \lambda)^{-1} \sL_{1,1}\left(\frac{-1}{\lambda}f_0^++\frac{1}{\lambda^2}f_0^-\right) \mathrm{d}\lambda\, .
   \end{align*}
Notice that the above vector has zero average. Therefore $P_{1,1} f_0^+$ is a vector with zero average, and the same holds for $(\partial_\mu \partial_\epsilon f_0^+(0,0))$. Finally, by \eqref{Parity structure}, $(\partial_\mu \partial_\epsilon f_0^+(0,0))$ has the structure claimed in \eqref{A17}.
We finally consider $(\partial_\mu \partial_\epsilon f_0^-(0,0))$. By Lemma \ref{lem:Dec} and by the fact that $\mathscr{L}_{0,1} f_0^- = 0$, it follows that for M=I,II,
\begin{equation*}
\eqref{A5a}^{M}f_0^- = -\frac{1}{2\pi \im} \oint_{\Gamma} (\mathscr{L}_{0,0} - \lambda)^{-1} \sL_{1,0}^{M} (\mathscr{L}_{0,0} - \lambda)^{-1} \sL_{0,1} \frac{-1}{\lambda} f_0^-\mathrm{d}\lambda=0\,.
\end{equation*}
Next, by Lemma \ref{lem:Dec} and $\sL_{1,0}^{I} f_0^- = 0$, we obtain $\eqref{A5b}^{I}f_0^- = 0$. Moreover, since $\sL_{1,0}^{II}f_0^-=f_0^+$, we have
\begin{align*}
   \eqref{A5b}^{II}f_0^-& =-\frac{1}{2\pi \im} \oint_{\Gamma} (\mathscr{L}_{0,0} - \lambda)^{-1}\sL_{0,1}  (\mathscr{L}_{0,0} - \lambda)^{-1}  \sL_{1,0}^{II}\frac{-1}{\lambda}f_0^-\mathrm{d}\lambda\\
   & =-\frac{1}{2\pi \im} \oint_{\Gamma} (\mathscr{L}_{0,0} - \lambda)^{-1}\sL_{0,1}  \frac{-1}{\lambda}\left(\frac{-1}{\lambda}f_0^++\frac{1}{\lambda^2}f_0^-\right)\mathrm{d}\lambda\\
   & =-\frac{1}{2\pi \im} \oint_{\Gamma} (\mathscr{L}_{0,0} - \lambda)^{-1}\frac{1}{\lambda^2} \left( 2 \ck^{3/2} f^-_{-1} + 4 \ck^{3/2} \tts f^-_{-2} \right) \mathrm{d}\lambda = \frac{1}{2} \ck^{-1/2} f^-_{-1} + \frac{1}{4} \ck^{-1/2} \tts f^-_{-2}  .
   \end{align*}
Finally, by Lemma \ref{lem:Dec} and \eqref{A8},
\begin{equation*}
 \eqref{A5c} f_0^- = -\frac{1}{2\pi \im} \oint_{\Gamma} (\mathscr{L}_{0,0} - \lambda)^{-1} \frac{1}{\lambda} \begin{bmatrix} 0 \\ \im p_1(x)\end{bmatrix} d \lambda\,.
\end{equation*} 
The above vector has zero average.

It remains to prove that the first component of
$\partial_\mu\partial_\epsilon f_j^-(0,0)$ has zero average for $j=1,2$. Recalling \eqref{A16prime} and \eqref{A16dot}, we have $P_{0,1}f_0^-=P_{1,0}f_0^-=0$.
Taking the coefficient of $\mu\epsilon$ in
$P_{\mu,\epsilon}^2=P_{\mu,\epsilon}$ gives
\begin{align*}
P_{0,0}P_{1,1}+P_{1,1}P_{0,0}
+P_{1,0}P_{0,1}+P_{0,1}P_{1,0}
=P_{1,1}.
\end{align*}
Applying this identity to $f_0^-$ and using
$P_{0,0}f_0^-=f_0^-$, we obtain $P_{0,0}P_{1,1}f_0^-=0$. Therefore, by \eqref{A2}, \eqref{the expansion of f sigma mu},
\begin{equation}
\label{eq:mixed-zero-mode-complement}
P_{0,0}\partial_\mu\partial_\epsilon f_0^-(0,0)=P_{0,0}\left(P_{1,1}-\frac12P_{0,0}P_{1,1}\right)f_0^-=0.
\end{equation}
Since $P_{0,0}f_j^-=f_j^-$ and
$P_{0,0}^*\mathcal J=\mathcal JP_{0,0}$ by
Lemma~\ref{properties of U and P}, we deduce that
\begin{equation}\label{1711:1002}
\begin{aligned}
(\mathcal J (\partial_\mu\partial_\epsilon f_0^-(0,0)),f_j^-)=(P_{0,0}^*\mathcal J(\partial_\mu\partial_\epsilon f_0^-(0,0)),f_j^-)=(\mathcal JP_{0,0}(\partial_\mu\partial_\epsilon f_0^-(0,0)),f_j^-)=0,
\end{aligned}
\end{equation}
for $j=1,2$. On the other hand, Lemma~\ref{F is symplectic and reversible}
gives $\bigl(\mathcal Jf_0^-(\mu,\epsilon),
f_j^-(\mu,\epsilon)\bigr)=0$.
 Differentiating $\bigl(\mathcal Jf_0^-(\mu,\epsilon),
f_j^-(\mu,\epsilon)\bigr)=0$ with respect
to $\mu$ and $\epsilon$, and evaluating at $(0,0)$, thus yields
\begin{align*}
\underbrace{(\mathcal J(\partial_\mu\partial_\epsilon f_0^-(0,0)),f_j^-)}_{ = 0 \mbox{ by } \eqref{1711:1002}}
+
(\mathcal J \underbrace{\partial_\epsilon f_0^-(0,0)}_{ = 0 \mbox{ by } \eqref{pa e f01} },\partial_\mu f_j^-(0,0))
+
(\mathcal J \underbrace{\partial_\mu f_0^-(0,0)}_{ = 0 \mbox{ by } \eqref{pa mu f01} },\partial_\epsilon f_j^-(0,0))
+(\mathcal Jf_0^-(0,0),\partial_\mu\partial_\epsilon f_j^-(0,0))=0.
\end{align*}
As $\mathcal Jf_0^-=f_0^+$, it follows that $(f_0^+,\partial_\mu\partial_\epsilon f_j^-(0,0))=0$. This verifies $\partial_\mu\partial_\epsilon f_j^-(0,0)$ has zero average for $j=1,2$.
\end{proof}

This completes the proof of Lemma \ref{expansion of the basis F}.

\begin{lemma} \label{eq:appendix-mixed-basis-cancellation}
    Fix a branch $\tts=\pm 1$. For each $j=1,2$ and each $\sigma=\pm$, we have 
    \begin{equation} \label{lemma C1}
 \bigl(\partial_\mu\partial_\e f_0^-(0,0,\tts),
       \mathcal B_1f_j^{\sigma}
       +\mathcal B_{0,0}\partial_\e  f_j^{\sigma}(0,0,\tts)\bigr)=0,
\end{equation}
where $\mathcal{B}_{0}$ and $\mathcal{B}_1$ are defined in \eqref{exp B0e}.
\end{lemma}
\begin{proof}
Recall \eqref{eq:mixed-zero-mode-complement}. By Lemma~\ref{properties of U and P},
$P_{0,0}^*\mathcal J=\mathcal JP_{0,0}$. Hence, for every
$v\in\mathcal V_{0,0}$,
\begin{align*}
 \bigl(\partial_\mu\partial_\e  f_0^-(0,0,\tts),\mathcal Jv\bigr)=\bigl(\partial_\mu\partial_\e  f_0^-(0,0,\tts),
                P_{0,0}^*\mathcal Jv\bigr)=\bigl(P_{0,0}\partial_\mu\partial_\e  f_0^-(0,0,\tts),
                \mathcal Jv\bigr)=0.
\end{align*}
For $j=1,2$, differentiating the invariance identity
\[
 P_{0,\e,\tts}\mathscr L_{0,\e,\tts} f_j^{\sigma}(0,\e,\tts)
 =\mathscr L_{0,\e,\tts} f_j^{\sigma}(0,\e,\tts)
\]
at $\e=0$, and using $\mathscr L_{0,0}f_j^{\sigma}=0$, shows that
\[
 \mathscr L_{0,1}f_j^{\sigma}
 +\mathscr L_{0,0}\partial_\e f_j^{\sigma}(0,0,\tts)
 \in\mathcal V_{0,0}.
\]
Using $\mathcal B_1$ as defined in \eqref{exp B0e}, we have
$\mathscr L_{0,1}=\mathcal J\mathcal B_1$ and
$\mathscr L_{0,0}=\mathcal J\mathcal B_{0,0}$. It follows that
\begin{equation*}
 \bigl(\partial_\mu\partial_\e  f_0^-(0,0,\tts),
       \mathcal B_1f_j^{\sigma}
       +\mathcal B_{0,0}\partial_\e  f_j^{\sigma}(0,0,\tts)\bigr)=0,
 \qquad j=1,2.
\end{equation*}
\end{proof}
\begin{lemma}
\label{lem:appendix-pure-cubic-cross}
The entries $F_{11}$
and $L_{11}$ in \eqref{F} and \eqref{L} satisfy
\begin{align}
 F_{11}(0,\e,\tts)= \frac{3\tts \ck^{3/2}}{16}\e^3+\cO(\e^4),\qquad
 L_{11}(0,\e,\tts)=- \frac{\ck^{3/2}}{8} \e^3+\cO(\e^4).
 \label{eq:appendix-L11-cubic}
\end{align}
\end{lemma}

\begin{proof}
Let us define
\begin{equation}\label{eq:appendix-cubic-vertical-vector}
 \tilde{U}_4:=\mathscr T_{\e,\tts}\vet{1}{0}
 =\begin{bmatrix}
 1+ (\mathfrak p_{\e,\tts})_x\\
 -B_{\e,\tts}\circ(\mathrm{Id}+\mathfrak p_{\e,\tts})(x)
 \end{bmatrix},
 \qquad \mathscr T_{\e,\tts}=\mathcal P_{\e,\tts} Z_{\e,\tts}^{-1},
\end{equation}
where $B_{\e,\tts}$ is the function $B$ in \eqref{espB}, and the coordinate
changes are those in \eqref{alinhac} and \eqref{LC}. Fix $c=c_{\e,\tts}$. Denote by $F_c(q)$ the left-hand side of
\eqref{travelingWWstokes} for $q=(\eta,\psi)$.
In infinite depth $G(\eta+a)\psi=G(\eta)\psi$ for any constant $a$,
by vertical translation of the harmonic extension, and $G(\eta)1=0$.
Thus, for constants $a,b$,
\[
 F_c(q+a f_0^+ + b f_0^-)=F_c(q)-a f_0^- .
\]
Differentiating in $a$ and $b$ at zero gives
$D F_c(q)f_0^+=-f_0^-$ and $D F_c(q)f_0^-=0$.
At the wave $q=(\eta_{\e,\tts},\psi_{\e,\tts})$,
\eqref{alinhac} and \eqref{LC} give
$\mathscr L_{0,\e,\tts}=\mathscr T_{\e,\tts}D F_c(q)\mathscr T_{\e,\tts}^{-1}$
and $\mathscr T_{\e,\tts} f_0^-=f_0^-$. Consequently
$\mathscr L_{0,\e,\tts}\tilde{U}_4=-f_0^-$ and
$\mathscr L_{0,\e,\tts}f_0^-=0$. Hence
$\tilde{U}_4\in\ker(\mathscr L_{0,\e,\tts}^2)\subset
\mathcal V_{0,\e,\tts}$ \eqref{Projection P mu e}, and, since
$\mathscr L_{0,\e,\tts}=\mathcal J\mathcal B_{0,\e,\tts}$,
\begin{equation}\label{eq:appendix-cubic-vertical-identity}
 \mathcal B_{0,\e,\tts}\tilde{U}_4=f_0^+.
\end{equation}
The vector $\tilde{U}_4$ is real and has an even first component and an
odd second component. By  \eqref{f-0=01}, it therefore
has an expansion
\begin{equation}\label{eq:appendix-cubic-vertical-coordinates}
 \tilde{U}_4= f_0^+(0,\e,\tts)+q_2(\e,\tts) f_2^+(0,\e,\tts)+q_1(\e,\tts) f_1^+(0,\e,\tts).
\end{equation}
The coefficient of $ f_0^+$ is exactly one: the first component of
$\tilde{U}_4$ has average one, those of $f_1^+(0,\e,\tts),f_2^+(0,\e,\tts)$ have average zero,
and that of $f_0^+(0,\e,\tts)$ has average one.

Symplecticity gives
\begin{equation}\label{eq:appendix-cubic-coordinate-pairing}
 q_j(\e,\tts)=-\bigl(\mathcal J\tilde{U}_4,f_j^-(0,\e,\tts)\bigr),
 \qquad j=1,2.
\end{equation}

We compute these coordinates to second order. Expanding
\eqref{eq:appendix-cubic-vertical-vector} using
\eqref{pexp}--\eqref{fexp} and \eqref{espB1}, we obtain
\begin{equation*}
 \tilde{U}_4=f_0^++\e \tilde{U}_{4_1}+\e^2\tilde{U}_{4_2}+\cO(\e^3),
\end{equation*}
where
\begin{align*}
 \tilde{U}_{4_1}
 =\begin{bmatrix}
 \cos (x)+\tts\cos(2x)\\
 -\ck\bigl(\sin (x)+\tts\sin(2x)\bigr)
 \end{bmatrix}, \qquad \tilde{U}_{4_2}
 =\begin{bmatrix}
 -\frac34 \tts\cos (x)-\cos(2x)-\frac94 \tts\cos(3x)\\
 \ck\bigl(\frac12 \tts\sin (x)+\frac74\sin(2x)
              +\frac{17}{4} \tts\sin(3x)+\sin(4x)\bigr)
 \end{bmatrix}.
\end{align*}
Indeed, the second-order coefficient in the second component is
$-B_2-\mathfrak p_1B_1'$. By \eqref{A2} and
$ P_{0,\e,\tts}^2=P_{0,\e,\tts}$,
\[
  P_{0,0}U_{0,2}P_{0,0}
 =-\frac12P_{0,0}P_{0,1}^2P_{0,0}.
\]
Recall \eqref{basis first derivatives quadratic calculation}. Since $\tilde U_{4_1}=P_{0,1}f_0^+$ and
$\mathcal JP_{0,1}=P_{0,1}^*\mathcal J$, it follows that
\[
  (\mathcal Jf_0^+,U_{0,2}f_j^-)
 =-\frac12(\mathcal J\tilde U_{4_1},\partial_\e  f_j^-(0,0,\tts)).
\]
A direct computation reveals
\begin{align*}
 (\mathcal J \tilde{U}_{4_1},f_j^-)&=0,\qquad j=1,2,\\
 (\mathcal J \tilde{U}_{4_2},f_1^-)&=\frac{\tts}{8}\sqrt{\ck},\quad (\mathcal J \tilde{U}_{4_1},\partial_\e  f_1^-(0,0,\tts))=0,\quad(\mathcal J \tilde{U}_{4_2},f_2^-)=-\frac38\sqrt{\ck},\quad (\mathcal J \tilde{U}_{4_1},\partial_\e  f_2^-(0,0,\tts))=\frac14\sqrt{\ck}.
\end{align*}
Consequently, expanding \eqref{eq:appendix-cubic-coordinate-pairing} and using
$ P_{0,0}P_{0,1}P_{0,0}=0$, we obtain
\[
  q_j(\e,\tts)=-\e^2\left[
 (\mathcal J\tilde U_{4_2},f_j^-)
 +\frac12(\mathcal J\tilde U_{4_1},\partial_\e  f_j^-(0,0,\tts))
 \right]+\cO(\e^3),\qquad j=1,2.
\]
Hence
\begin{align}
 q_1(\e,\tts)=-\tts\frac{\sqrt{\ck}}{8} \e^2+\cO(\e^3),\qquad q_2(\e,\tts)= \frac{\sqrt{\ck}}{4} \e^2+\cO(\e^3).
 \label{eq:appendix-cubic-q-expansion}
\end{align}
Finally, pair \eqref{eq:appendix-cubic-vertical-identity} with
$ f_2^+(0,\e,\tts)$ and $ f_1^+(0,\e,\tts)$. Their first components have zero average,
and \eqref{eq:appendix-cubic-vertical-coordinates} implies the
exact identities
\begin{align*}
 F_{11}(0,\e,\tts)
 =-E_{11}(0,\e,\tts)q_2(\e,\tts)-I_{11}(0,\e,\tts)q_1(\e,\tts),\quad
 L_{11}(0,\e,\tts)
 =-I_{11}(0,\e,\tts)q_2(\e,\tts)-H_{11}(0,\e,\tts)q_1(\e,\tts).
\end{align*}
Using that
\[
 E_{11}(0,\e,\tts)=-\tts\frac{\ck}2\e+\cO(\e^2),\qquad
 I_{11}(0,\e,\tts)=\frac{\ck}2\e+\cO(\e^2),\qquad
 H_{11}(0,\e,\tts)=\cO(\e^2),
\]
together with \eqref{eq:appendix-cubic-q-expansion}, we obtain \eqref{eq:appendix-L11-cubic}.
\end{proof}

\section{Analysis of the branch $\tts = -1$} \label{app:s-1}

Here we prove that the six eigenvalues of $\mathcal L_{\mu,\e,-1}$ near the origin are purely imaginary for sufficiently small $\e>0$ and $\e^2\leq\mu\leq2\e^2$.

\begin{proposition}[Non perturbative block-decoupling, $\tts =-1$]\label{prop:first-block-decoupling 2}
There exists $\e_0>0$ such that for any 
$(\tau,\ve)\in[1,2^{1/4}]\times(0,\e_0^{1/2})$, the following holds true. 
There exists a $4\times2$ reversibility-preserving matrix $X_-$
whose upper and lower $2\times2$ blocks, denoted by $X_{2,-1}$
and $X_{3,-1}$, respectively, are real analytic in $(\tau,\ve)$ and expand as
\begin{equation}\label{eq:Xm1-sylvester}
\begin{aligned}
X_{2,-1} =
\begin{pmatrix}
-\fa_2\ve^2-\fc\ve^3+r(\ve^4)
&
-\im\fb_2\ve^2+\im r(\ve^4)
\\[1ex]
\im\fa_2\ve^2-\im\fc\ve^3+\im r(\ve^4)
&
-\fb_2\ve^2+r(\ve^4)
\end{pmatrix},\, X_{3,-1} =
\begin{pmatrix}
-\fa_3\ve^2+\fc\ve^3+r(\ve^4)
&
\im\fb_3\ve^2+\im r(\ve^4)
\\[1ex]
-\im\fa_3\ve^2-\im\fc\ve^3+\im r(\ve^4)
&
-\fb_3\ve^2+r(\ve^4)
\end{pmatrix},
\end{aligned}
\end{equation}
with $\fa_2,\fb_2,\fa_3,\fb_3$ and $\fc$
as in \eqref{X.coeff}. By conjugating the Hamiltonian and reversible
matrix $\mathsf L^{(1)}_{\tau,\ve,-1}$, defined in \eqref{eq:L1-first-change}, with the symplectic and reversibility-preserving $6\times 6$
matrix $\exp(S_{-1}^{(1)})$, where
\begin{equation}\label{Sminus}
S_{-1}^{(1)}=\mathsf{J}_6\begin{pmatrix}0_{4}&\mathsf{J}_4X_-\\(\mathsf{J}_4X_-)^*&0_{2}\end{pmatrix}
,\quad 
X_-:=\begin{pmatrix}X_{2,-1}\\X_{3,-1}\end{pmatrix} , \quad \mathsf{J}_4 := \begin{pmatrix}
    \mathsf{J}_2 & 0 \\
    0 & \mathsf{J}_2
\end{pmatrix}
,
\end{equation}
we get the Hamiltonian and reversible matrix 
\begin{equation}\label{minusresult}
\begin{aligned}
 \mathsf L^{(2)}_{\tau,\ve,-1}
 &:=\exp(S_{-1}^{(1)})\mathsf L^{(1)}_{\tau,\ve,-1}\exp(-S_{-1}^{(1)})
 =\mathsf J_6\mathsf B^{(2)}_{\tau,\ve,-1} , \\
\mathsf B^{(2)}_{\tau,\ve,-1}
 &:= \begin{pmatrix}A^{(2)}& C^{(2)}\\
[ C^{(2)} ]^*&G^{(2)}\end{pmatrix}, \quad
A^{(2)}:=\begin{pmatrix}E^{(2)} &I^{(2)} \\ [ I^{(2)} ]^*&H^{(2)} \end{pmatrix} , 
\end{aligned}
\end{equation}
where the blocks $E^{(2)}$, $H^{(2)}$ and $G^{(2)}$ are self-adjoint, and 
\begin{align} \label{eq:E2-final-minus}
    E^{(2)}=\ve^2\begin{pmatrix}
\beta+r(\ve^2)
&
  -\im\left(\alpha+r(\ve^2)\right)
\\
 \im \left(\alpha+ r(\ve^2)\right)
&
\beta+ r(\ve^2)
\end{pmatrix}, \qquad H^{(2)}=\ve^2\begin{pmatrix}
\beta+r(\ve^2)
&
  \im\left(\alpha+r(\ve^2)\right)
\\
 -\im \left(\alpha+ r(\ve^2)\right)
&
\beta+ r(\ve^2)
\end{pmatrix},
\end{align}
and
\begin{align} \label{eq:G2-final-minus}
    G^{(2)}=\begin{pmatrix}
G^{(2)}_{11}
&
  \im G^{(2)}_{12}
\\
-\im G^{(2)}_{12}
&
G^{(2)}_{22}
\end{pmatrix}=\begin{pmatrix}
-\frac{\ck}{4}\ve^2+r(\ve^4)
&
  \im r(\ve^4)
\\
-\im r(\ve^4)
&
-\frac{4\ck}{81}\tau^8\ve^6+r(\ve^8)
\end{pmatrix}, \quad I^{(2)}=\cO(\ve^4), \quad C^{(2)} =
\begin{pmatrix}
F^{(2)}\\
L^{(2)}
\end{pmatrix}=\cO(\ve^8),
\end{align}
where $\alpha$, $\beta$ are as in \eqref{final.coeff}.
\end{proposition}

\begin{proof}
We split the matrix
\eqref{eq:L1-first-change} into a $(4\times4) + (2\times2)$ block-diagonal part, and the relative
off-diagonal parts:

\begin{equation*} 
\mathsf{L}_{\tau,\ve,-1}^{(1)}
=D_{-1}^{(1)}+R_{-1}^{(1)},
\quad
D_{-1}^{(1)}:=
\begin{pmatrix}
\mathsf{J}_4 A^{(1)} & 0\\
0 & \mathsf{J}_2 G^{(1)}
\end{pmatrix},
\qquad
R_{-1}^{(1)}:=
\begin{pmatrix}
0 & \mathsf{J}_4 C^{(1)}\\
\mathsf{J}_2 [C^{(1)}]^* & 0
\end{pmatrix},
\end{equation*}
where
\[
A^{(1)} =
\begin{pmatrix}
E^{(1)} & I^{(1)}\\
[I^{(1)}]^* & H^{(1)}
\end{pmatrix},
\qquad
C^{(1)} =
\begin{pmatrix}
F^{(1)}\\
L^{(1)}
\end{pmatrix},
\qquad 
\mathsf{L}_{\tau,\ve,-1}^{(1)}
=D_{-1}^{(1)}+R_{-1}^{(1)}.
\]

By \eqref{eq:E1-first-change}-\eqref{eq:Lblock1-first-change}, we have that $D_{-1}^{(1)} = \cO(\ve^2)$ and $R_{-1}^{(1)} = \cO(\ve^4)$.
We seek a Hamiltonian and reversibility-preserving matrix
$S_{-1}^{(1)}$ of the form \eqref{Sminus}, satisfying the homological equation
\[
R_{-1}^{(1)}+[S_{-1}^{(1)},D_{-1}^{(1)}]=0.
\]
Writing $X_-$ as in \eqref{Sminus}, the homological equation is equivalent to
\begin{equation} \label{eq:Syslvester4x2}
(\mathsf{J}_4 A^{(1)})X_- - X_-(\mathsf{J}_2 G^{(1)}) = -\mathsf{J}_4 C^{(1)},
\end{equation}
or, in terms of the corresponding $2\times2$ blocks,
\begin{equation*}
\begin{aligned}
(\mathsf{J}_2 E^{(1)})X_{2,-1} +\underbrace{\mathsf{J}_2 I^{(1)}}_{\cO(\ve^4)}X_{3,-1}
- X_{2,-1}(\mathsf{J}_2 G^{(1)}) &=-\mathsf{J}_2 F^{(1)},\\
\underbrace{\mathsf{J}_2 [I^{(1)}]^*}_{\cO(\ve^4)}X_{2,-1} +(\mathsf{J}_2 H^{(1)})X_{3,-1}
- X_{3,-1} (\mathsf{J}_2 G^{(1)}) &= -\mathsf{J}_2 L^{(1)}.
\end{aligned}
\end{equation*}
Arguing as in \eqref{eq:three-real-systems}, we can rewrite \eqref{eq:Syslvester4x2} in terms of an $8\times8$ system, whose coefficient matrix has the form
\[
A_-=\operatorname{diag}(A_2,A_3)+A_{\mathrm{r}},
\qquad
\|A_{\mathrm{r}}\|\leq C\|I^{(1)}\|\leq C\ve^4,
\]
where $A_2,A_3$ are defined by \eqref{eq:Ak-three-blocks} (and where their entries are evaluated at $\tts=-1$). The matrix $A_{\mathrm r}$ has zero diagonal
$4\times4$ blocks, while its off-diagonal blocks
arise from the coupling terms involving
$I^{(1)}$ and $[I^{(1)}]^*$.

Using \eqref{eq:E1-first-change}-\eqref{eq:G1-first-change} and \eqref{final.coeff}, \eqref{eq:coeff-first-change}, we have that the matrices $A_2,A_3$ satisfy
\begin{equation*}
 \det A_k=\frac{\ck^4}{2}\tau^4\ve^8\bigl(1+r(\ve^2)\bigr) , \quad
 \|A_k^{-1}\|\leq C\ve^{-2}, k=2,3,
\end{equation*}
uniformly on \eqref{eq:sylvester-band-new}.
We factor $A_-=\operatorname{diag}(A_2,A_3)
\left[
I_8+\operatorname{diag}(A_2^{-1},A_3^{-1})
A_{\mathrm r}
\right]$. Since $\| \operatorname{diag}(A_2^{-1},A_3^{-1}) A_{\mathrm{r}} \| \leq C \ve^2$, the matrix $A_-$ is invertible by a Neumann series, and  $\| A_-^{-1} \| \leq C \ve^{-2}$.

Therefore, the system \eqref{eq:Syslvester4x2} has a unique solution which is real analytic in $(\tau,\ve)$, and such that
\begin{equation*}
    X_{2,-1} , X_{3,-1} = \cO(\ve^2) , \quad S_{-1}^{(1)} = \cO(\ve^2) .
\end{equation*}
Since $I^{(1)} = \cO(\ve^4)$, the solutions of the coupled system \eqref{eq:Syslvester4x2} differ from the solution of the decoupled system by $\cO(\ve^4)$.
Substituting \eqref{eq:Lblock1-first-change} into the decoupled systems gives
$X_{k,-1}=-X_{k,1}+\cO(\ve^4)$, $k=2,3$, which proves \eqref{eq:Xm1-sylvester}.

The Lie expansion and the homological equation give the exact identity
\begin{equation}\label{d:lie}
 \mathsf{L}_{\tau,\ve,-1}^{(2)}
 =D^{(1)}_{-1}+\frac12[S_{-1}^{(1)},R_{-1}^{(1)}]+\mathcal T,\qquad
 \mathcal T:=\frac12\int_0^1(1-t^2)e^{t S_{-1}^{(1)}}
 \ad_{S_{-1}^{(1)}}^2(R_{-1}^{(1)})e^{-t S_{-1}^{(1)}}\,\mathrm{d}t.
\end{equation}
Let $\Pi_D,\Pi_{\emptyset}$ be the diagonal and off-diagonal projections
for the $4+2$ decomposition. Since $S_{-1}^{(1)}$, $R_{-1}^{(1)}$ are block off-diagonal,
$[S_{-1}^{(1)},R_{-1}^{(1)}]$ is block-diagonal and
$N:=\ad_{S_{-1}^{(1)}}^2(R_{-1}^{(1)})$ is block off-diagonal, with $N=\cO(\ve^8)$.
Thus $\Pi_{\emptyset}\mathcal T=\cO(\ve^8)$. Also
\[
 \Pi_D(e^{t S_{-1}^{(1)}} N e^{-t S_{-1}^{(1)}})
 =\int_0^t\Pi_D\left(e^{u S_{-1}^{(1)}}[S_{-1}^{(1)},N]e^{-u S_{-1}^{(1)}}\right) \mathrm{d}u =O(\ve^{10}),
\]
so $\Pi_D\mathcal T=\cO(\ve^{10})$. Computing the commutator gives
\begin{equation*}
\begin{aligned}
A^{(2)}=A^{(1)}+\mathrm{Sym}(\mathsf{J}_4 X_-\mathsf{J}_2 ( C^{(1)} )^*)+\cO(\ve^{10}), \quad
G^{(2)}=G^{(1)}+\mathrm{Sym}(X_-^* C^{(1)})+\cO(\ve^{10}) .
\end{aligned}
\end{equation*}
In particular, $A^{(2)}-A^{(1)}=\cO(\ve^6)$
whereas by \eqref{d:lie} the remaining $4\times2$ off-diagonal block $C^{(2)}=\cO(\ve^8)$, since $\Pi_{\emptyset}\mathcal T=\cO(\ve^8)$. Also, a straightforward calculation reveals that $G^{(2)}_{22}=-\tfrac{4\ck}{81}\tau^8\ve^6+r(\ve^8)$. This verifies \eqref{eq:G2-final-minus}.

Finally, $e^{S_{-1}^{(1)}}$ is symplectic and reversibility preserving, so the
transformed matrix retains its Hamiltonian and reversible structure.
\end{proof}
\begin{proposition}[Spectral stability of the negative branch] \label{final proposition}
On the regime \eqref{eq:sylvester-band-new}, the six eigenvalues of
$\mathcal L_{\mu,\e,-1}$ near the origin are purely imaginary for
sufficiently small $\ve>0$.
\end{proposition}

\begin{proof}
By \eqref{mathcal L= ichmu+mathscr L} and \eqref{minusresult}, it
suffices to prove that the eigenvalues of $\mathsf L^{(2)}_{\tau,\ve,-1}$
are purely imaginary. By \eqref{eq:E2-final-minus} and
\eqref{eq:G2-final-minus}, $A^{(2)}\geq m\ve^2 I_4$ for some $m>0$
uniformly on the band. Fix $0<\rho<m/2$. By \eqref{minusresult} and
Remark~\ref{rem:-s}, the eigenvalues
with $|\lambda|\geq\rho\ve^2$ are purely imaginary.

Let $\lambda$ be an eigenvalue with $|\lambda|<\rho\ve^2$ and
$w=(u,v)$ an associated eigenvector. The upper block equation is
\[
(\mathsf J_4A^{(2)}-\lambda I_4)u=-\mathsf J_4C^{(2)}v.
\]
Since $\|(\mathsf J_4A^{(2)})^{-1}\|\leq m^{-1}\ve^{-2}$, a
Neumann series bounds the inverse on the left by $2m^{-1}\ve^{-2}$.
As $C^{(2)}=\cO(\ve^8)$, we obtain $\|u\|\lesssim\ve^6\|v\|$.

Let $g=G^{(2)}_{12}$ be the real scalar in the off-diagonal entry
$\im G^{(2)}_{12}$ of \eqref{eq:G2-final-minus}. The shifted matrix
$\mathsf L^{(2)}_{\tau,\ve,-1}+\im g I_6=\mathsf J_6\widetilde B$,
where $\widetilde B=\mathsf B^{(2)}_{\tau,\ve,-1}-\im g\mathsf J_6$,
has the same eigenvector $w$. The lower diagonal block of
$\widetilde B$ is $\operatorname{diag}(G^{(2)}_{11},G^{(2)}_{22})$.
By \eqref{eq:G2-final-minus}, it is bounded above by
$-c\ve^6 I_2$ for some $c>0$. The upper diagonal block is
$\cO(\ve^2)$, and the off-diagonal block is still $C^{(2)}$.
Consequently,
\[
w^*\widetilde B w
\leq -c\ve^6\|v\|^2
+\cO(\ve^2)\|u\|^2+\cO(\ve^8)\|u\|\|v\|
\leq \bigl(-c\ve^6+\cO(\ve^{14})\bigr)\|v\|^2<0.
\]
The Hamiltonian identity for the shifted matrix implies
$2\operatorname{Re}(\lambda+\im g)\,\langle w, \widetilde B w\rangle=0$.
Thus $\operatorname{Re}\lambda=0$.
\end{proof}
\footnotesize

\vspace{1em}
\noindent{\bf Acknowledgments:}
T.-Y. Hsiao, A. Maspero and S. Pasquali are  supported by   the European Union  ERC CONSOLIDATOR GRANT 2023 GUnDHam, Project Number: 101124921.
 Views and opinions expressed are however those of the authors only and do not necessarily reflect those of the European Union or the European Research Council. Neither the European Union nor the granting authority can be held responsible for them. The work of C. Sun is partially supported by the ANR project ANR-24-CE40-3260.

\vspace{1em}
\noindent{\bf Declaration on the use of AI:}
The mathematical ideas and proof strategies were developed by the authors. Astra 6.1 was used to support algebraic and symbolic computations and to translate their formulas into MATLAB code for the schematic plots. The authors take full responsibility for the mathematical correctness and scientific content of the paper.


\begin{thebibliography}{99}
\begingroup\hypersetup{urlcolor=black}


\bibitem{akers2021wilton}
B. Akers and D. P. Nicholls. \href{https://doi.org/10.1007/s42286-020-00034-w}{\emph{Wilton ripples in weakly nonlinear dispersive models of water waves: existence and analyticity of solution branches}}. Water Waves \textbf{3} (2021), no.~1, 25--47.


\bibitem{BBMM}
T. Barbieri, M. Berti, A. Maspero and M. Mazzucchelli. \href{https://doi.org/10.1016/j.jde.2025.113753}{\emph{Bifurcation of gravity-capillary Stokes waves with constant vorticity}}. J. Differential Equations \textbf{451} (2026), 113753.

\bibitem{BarbieriBertiMazzucchelli2026_3D}
T. Barbieri, M. Berti and M. Mazzucchelli.
\href{https://arxiv.org/abs/2603.27847}
{\emph{Equivariant critical point theory and bifurcation
of 3d gravity-capillary Stokes waves}}.
Preprint, arXiv:2603.27847 (2026).

\bibitem{BF}
T. B. Benjamin and J. E. Feir. \href{https://doi.org/10.1017/S002211206700045X}{\emph{The disintegration of wave trains on deep water. Part~1. Theory}}. J. Fluid Mech. \textbf{27} (1967), no.~3, 417--430.

\bibitem{BCMV}
M. Berti, L. Corsi, A. Maspero and P. Ventura. \emph{Infinitely many isolas of modulational instability for Stokes waves}. Preprint, \href{https://arxiv.org/abs/2405.05854v3}{arXiv:2405.05854v3} (2024).



\bibitem{BMR}
M. Berti, A. Maspero and A. M. Radakovic. \href{https://arxiv.org/abs/2607.01939}{\emph{McLean resonances and 3d spectral instability of Stokes waves}}. \texttt{arXiv:2607.01939}, 2026.


\bibitem{BMV1}
M. Berti, A. Maspero and P. Ventura. \href{https://doi.org/10.1007/s00222-022-01130-z}{\emph{Full description of Benjamin--Feir instability of Stokes waves in deep water}}. Invent. Math. \textbf{230} (2022), no.~2, 651--711.

\bibitem{BMV3}
M. Berti, A. Maspero and P. Ventura. \href{https://doi.org/10.1007/s00205-023-01916-2}{\emph{Benjamin--Feir instability of Stokes waves in finite depth}}. Arch. Ration. Mech. Anal. \textbf{247} (2023), article~91.

\bibitem{BMV_ed}
M. Berti, A. Maspero and P. Ventura. \href{https://doi.org/10.1007/s00220-023-04928-x}{\emph{Stokes waves at the critical depth are modulationally unstable}}. Comm. Math. Phys. \textbf{405} (2024), article~56.

\bibitem{BMV4}
M. Berti, A. Maspero and P. Ventura. \href{https://doi.org/10.4171/EMSS/91}{\emph{First isola of modulational instability of Stokes waves in deep water}}. EMS Surv. Math. Sci. \textbf{12} (2025), no.~1, 79--122.


\bibitem{BBHZ}
B. Barker, J. C. Bronski, V. M. Hur and Z. Yang. \href{https://doi.org/10.1016/j.jde.2025.113550}{\emph{Asymptotic stability of sharp fronts: Analysis and rigorous computation}}. J. Differential Equations \textbf{444} (2025), 113550.

\bibitem{BMP}
R. Bianchini, A. Maspero and S. Pasquali. \emph{Instabilities of internal gravity waves in the two-dimensional Boussinesq system}. Accepted for publication in SIAM J. Math. Anal. (2026), \href{https://arxiv.org/abs/2507.10390}{arXiv:2507.10390}.

\bibitem{BrM}
T. J. Bridges and A. Mielke. \href{https://doi.org/10.1007/BF00376815}{\emph{A proof of the Benjamin--Feir instability}}. Arch. Ration. Mech. Anal. \textbf{133} (1995), 145--198.

\bibitem{BHJ}
J. C. Bronski, V. M. Hur and M. A. Johnson. \href{https://doi.org/10.1007/978-3-319-20690-5_4}{\emph{Modulational instability in equations of KdV type}}. In \emph{New Approaches to Nonlinear Waves}, E. Tobisch (ed.), Lecture Notes in Physics 908, Springer, Cham (2016), 83--133.

\bibitem{BHR}
J. C. Bronski, V. M. Hur and R. Marangell. \href{https://doi.org/10.1088/1361-6544/ad8698}{\emph{Floquet theory and stability analysis for Hamiltonian PDEs}}. Nonlinearity \textbf{37} (2024), 125010.

\bibitem{BJ}
J. C. Bronski and M. A. Johnson. \href{https://doi.org/10.1007/s00205-009-0270-5}{\emph{The modulational instability for a generalized Korteweg--de Vries equation}}. Arch. Ration. Mech. Anal. \textbf{197} (2010), 357--400.

\bibitem{BDS}
E. Byrnes, B. Deconinck and A. Semenova. \href{https://doi.org/10.1017/jnw.2026.10031}{\emph{Stokes waves in water of finite depth}}. J. Nonlinear Waves \textbf{2} (2026), e7.

\bibitem{CN1}
B. Chen and H. Q. Nguyen.
\emph{Suppression of the Benjamin--Feir instability by constant vorticity}.
\href{https://arxiv.org/abs/2610.39000}{arXiv:2610.39000} (2026).

\bibitem{CN2}
B. Chen and H. Q. Nguyen.
\emph{The Benjamin--Feir stability diagram for Stokes waves with constant vorticity in finite depth}.
\href{https://arxiv.org/abs/2611.18203}{arXiv:2611.18203} (2026).

\bibitem{ChenSaffman1979}
B. Chen and P. G. Saffman.
\href{https://doi.org/10.1002/sapm1979603183}
{\emph{Steady gravity-capillary waves on deep water.
I. Weakly nonlinear waves}}.
Stud. Appl. Math. \textbf{60} (1979), no.~3, 183--210.

\bibitem{ChenSu}
G. Chen and Q. Su. \href{https://doi.org/10.1007/s00220-023-04747-0}{\emph{Nonlinear modulational instabililty of the Stokes waves in 2D full water waves}}. Comm. Math. Phys. \textbf{402} (2023), 1345--1452.


\bibitem{CraigNicholls2000}
W. Craig and D. P. Nicholls.
\href{https://doi.org/10.1137/S0036141099354181}
{\emph{Traveling two and three dimensional capillary gravity water waves}}.
SIAM J. Math. Anal. \textbf{32} (2000), no.~2, 323--359.

\bibitem{CS}
W. Craig and C. Sulem. \href{https://doi.org/10.1006/jcph.1993.1164}{\emph{Numerical simulation of gravity waves}}. J. Comput. Phys. \textbf{108} (1993), 73--83.

\bibitem{CD}
R. P. Creedon and B. Deconinck. \href{https://doi.org/10.1017/jfm.2022.1031}{\emph{A high-order asymptotic analysis of the Benjamin--Feir instability spectrum in arbitrary depth}}. J. Fluid Mech. \textbf{956} (2023), A29.

\bibitem{CDT}
R. P. Creedon, B. Deconinck and O. Trichtchenko. \href{https://doi.org/10.1017/jfm.2021.1119}{\emph{High-frequency instabilities of Stokes waves}}. J. Fluid Mech. \textbf{937} (2022), A24.

\bibitem{CNS}
R. P. Creedon, H. Q. Nguyen and W. A. Strauss. \href{https://doi.org/10.1007/s40818-024-00188-7}{\emph{Proof of the transverse instability of Stokes waves}}. Ann. PDE \textbf{11} (2025), 4.

\bibitem{CNS2}
R. P. Creedon, H. Q. Nguyen and W. A. Strauss. \href{https://doi.org/10.1137/24M1684372}{\emph{Transverse instability of Stokes waves at finite depth}}. SIAM J. Math. Anal. \textbf{58} (2026), no.~1, 661--697.



\bibitem{DT}
B. Deconinck and O. Trichtchenko. \href{https://doi.org/10.1016/j.euromechflu.2014.02.010}{\emph{Stability of periodic gravity waves in the presence of surface tension}}. Eur. J. Mech. B Fluids \textbf{46} (2014), 97--108.


\bibitem{FS}
L. Fan and C. Sun.
\emph{Modulational stability of finite-depth Stokes waves with constant vorticity}.
\href{https://arxiv.org/abs/2611.05237}{arXiv:2611.05237} (2026).

\bibitem{GH}
T. Gallay and M. Haragus. \href{https://doi.org/10.1016/j.jde.2006.12.007}{\emph{Stability of small periodic waves for the nonlinear Schr\"odinger equation}}. J. Differential Equations \textbf{234} (2007), no.~2, 544--581.

\bibitem{HK}
M. Haragus and T. Kapitula. \href{https://doi.org/10.1016/j.physd.2008.03.050}{\emph{On the spectra of periodic waves for infinite-dimensional Hamiltonian systems}}. Phys. D \textbf{237} (2008), no.~20, 2649--2671.

\bibitem{HTW}
M. Haragus, T. Truong, and E. Wahl\'en. \href{https://doi.org/10.1007/s42286-023-00074-y}{\emph{Transverse dynamics of two-dimensional traveling periodic gravity-capillary water waves}}. Water Waves \textbf{5} (2023), no.~1, 65--99.

\bibitem{ebb}
S. V. Haziot, V. M. Hur, W. A. Strauss, J. F. Toland, E. Wahl\'en, S. Walsh, and M. H. Wheeler. \href{https://doi.org/10.1090/qam/1614}{\emph{Traveling water waves---the ebb and flow of two centuries}}. Quart. Appl. Math. \textbf{80} (2022), no.~2, 317--401.

\bibitem{HJ}
V. M. Hur and M. A. Johnson. \href{https://doi.org/10.1111/sapm.12061}{\emph{Modulational instability in the Whitham equation for water waves}}. Stud. Appl. Math. \textbf{134} (2015), no.~1, 120--143.

\bibitem{HsiaoThesis}
T.-Y. Hsiao.
\href{https://hdl.handle.net/2142/127490}{\emph{Unstable Stokes waves in constant vorticity flows}}.
Ph.D. thesis, University of Illinois at Urbana-Champaign (2024).

\bibitem{HHY1}
T.-Y. Hsiao, V. M. Hur and Z. Yang.
\emph{Unstable Stokes waves in constant vorticity flows of finite depth}.
\href{https://arxiv.org/abs/2611.04509}{arXiv:2611.04509} (2026).

\bibitem{HHY2}
T.-Y. Hsiao, V. M. Hur and Z. Yang.
\emph{Unstable Stokes waves in constant vorticity flows of infinite depth}.
\href{https://arxiv.org/abs/2611.04528}{arXiv:2611.04528} (2026).

\bibitem{HM}
T.-Y. Hsiao and A. Maspero.
\href{https://doi.org/10.1007/s00220-026-05776-1}{\emph{Full Benjamin--Feir instability of capillary-gravity Stokes waves in finite depth}}.
Comm. Math. Phys. \textbf{408} (2027), article 16.

\bibitem{HW}
T.-Y. Hsiao and X. Wang. \emph{Spectral structure of the Benjamin--Feir instability in deep-water gravity-capillary Stokes waves}. Accepted for publication in SIAM J. Math. Anal. (2026), \href{https://arxiv.org/abs/2604.16724v3}{arXiv:2604.16724v3}.

\bibitem{HP}
V. M. Hur and A. K. Pandey. \href{https://doi.org/10.1016/j.physd.2016.03.005}{\emph{Modulational instability in nonlinear nonlocal equations of regularized long wave type}}. Phys. D \textbf{325} (2016), 98--112.

\bibitem{HY}
V. M. Hur and Z. Yang. \href{https://doi.org/10.1007/s00205-023-01889-2}{\emph{Unstable Stokes waves}}. Arch. Ration. Mech. Anal. \textbf{247} (2023), no.~4, article~62.

\bibitem{HY2}
V. M. Hur and Z. Yang. \emph{Stable and unstable capillary-gravity waves}. \href{https://arxiv.org/abs/2311.01368v4}{arXiv:2311.01368v4} (2026).

\bibitem{JRSY}
Z. Jiao, L. M. Rodrigues, C. Sun, and Z. Yang. \href{https://doi.org/10.1007/s00220-025-05428-w}{\emph{Small-amplitude finite-depth Stokes waves are transversally unstable}}. Comm. Math. Phys. \textbf{406} (2025), no.~10, article~255.

\bibitem{Jones1992}
M. C. W. Jones. \href{https://doi.org/10.1016/0165-2125(92)90011-P}{\emph{Nonlinear stability of resonant capillary-gravity waves}}. Wave Motion \textbf{15} (1992), no.~3, 267--283.

\bibitem{langer2024wilton}
R. Langer, O. Trichtchenko, and B. Akers. \href{https://doi.org/10.1007/s42286-024-00086-2}{\emph{Wilton ripples with high-order resonances in weakly nonlinear models}}. Water Waves \textbf{6} (2024), no.~1, 97--126.



\bibitem{LD book}
D. Lannes. \href{https://doi.org/10.1090/surv/188}{\emph{The Water Waves Problem: Mathematical Analysis and Asymptotics}}. Mathematical Surveys and Monographs, vol.~188, American Mathematical Society, Providence, RI, 2013.



\bibitem{Kato1966}
T. Kato. \href{https://doi.org/10.1007/978-3-662-12678-3}{\emph{Perturbation Theory for Linear Operators}}. Grundlehren der mathematischen Wissenschaften, vol.~132, Springer-Verlag, 1966.

\bibitem{Li}
M. J. Lighthill. \href{https://doi.org/10.1093/imamat/1.3.269}{\emph{Contributions to the theory of waves in non-linear dispersive systems}}. IMA J. Appl. Math. \textbf{1} (1965), no.~3, 269--306.

\bibitem{MaRa}
A. Maspero and A. M. Radakovic.
\href{https://doi.org/10.1137/24M1652076}
{\emph{Full description of Benjamin--Feir instability
for generalized Korteweg--de Vries equations}}.
SIAM J. Math. Anal. \textbf{57} (2025), no.~3,
3030--3070.

\bibitem{Mc1}
J. W. McLean. \href{https://doi.org/10.1017/S0022112082000172}{\emph{Instabilities of finite-amplitude water waves}}. J. Fluid Mech. \textbf{114} (1982), 315--330.

\bibitem{Mc2}
J. W. McLean. \href{https://doi.org/10.1017/S0022112082000184}{\emph{Instabilities of finite-amplitude gravity waves on water of finite depth}}. J. Fluid Mech. \textbf{114} (1982), 331--341.



\bibitem{NS}
H. Q. Nguyen and W. A. Strauss. \href{https://doi.org/10.1002/cpa.22073}{\emph{Proof of modulational instability of Stokes waves in deep water}}. Comm. Pure Appl. Math. \textbf{76} (2023), no.~5, 1035--1084.



\bibitem{PiersonFife1961}
W. J. Pierson, Jr. and P. Fife.
\href{https://doi.org/10.1029/JZ066i001p00163}
{\emph{Some nonlinear properties of long-crested periodic waves
with lengths near 2.44 centimeters}}.
J. Geophys. Res. \textbf{66} (1961), no.~1, 163--179.

\bibitem{reeder1981wiltonI}
J. Reeder and M. Shinbrot. \href{https://doi.org/10.1016/0165-2125(81)90001-9}{\emph{On Wilton ripples, I: Formal derivation of the phenomenon}}. Wave Motion \textbf{3} (1981), no.~2, 115--135.

\bibitem{reeder1981wiltonII}
J. Reeder and M. Shinbrot. \href{https://doi.org/10.1007/BF00280641}{\emph{On Wilton ripples, II: Rigorous results}}. Arch. Rational Mech. Anal. \textbf{77} (1981), no.~4, 321--347.

\bibitem{ShavitWangBuckmaster2026}
M. Shavit, Y. Wang and T. Buckmaster.
\href{https://arxiv.org/abs/2611.07395}
{\emph{Rigorously certified neural-network computation
of nonlinear water waves at near-machine precision}}.
Preprint, arXiv:2611.07395 (2026).

\bibitem{stokes}
G. G. Stokes. \emph{On the theory of oscillatory waves}. Trans. Cambridge Philos. Soc. \textbf{8} (1847), 441--455.

\bibitem{SW}
C. Sun and E. Wahl\'en. \emph{On the spectral stability of periodic capillary-gravity waves}. \href{https://arxiv.org/abs/2509.17534v3}{to appear in Comm. Math. Phys.} 

\bibitem{TolandJones1985}
J. F. Toland and M. C. W. Jones.
\href{https://doi.org/10.1098/rspa.1985.0063}
{\emph{The bifurcation and secondary bifurcation
of capillary-gravity waves}}.
Proc. Roy. Soc. London Ser. A \textbf{399} (1985),
no.~1817, 391--417.


\bibitem{trichtchenko2016instability}
O. Trichtchenko, B. Deconinck, and J. Wilkening. \href{https://doi.org/10.1016/j.wavemoti.2016.06.004}{\emph{The instability of Wilton ripples}}. Wave Motion \textbf{66} (2016), 147--155.



\bibitem{Wilton}
J. R. Wilton. \href{https://doi.org/10.1080/14786440508635350}{\emph{LXXII. On ripples}}. Philos. Mag. (6) \textbf{29} (1915), no.~173, 688--700.


\bibitem{Zak68}
V. E. Zakharov. \href{https://doi.org/10.1007/BF00913182}{\emph{Stability of periodic waves of finite amplitude on the surface of a deep fluid}}. J. Appl. Mech. Tech. Phys. \textbf{9} (1968), no.~2, 190--194.

\endgroup
\end{thebibliography}
\end{document}